**Université catholique de Louvain**

Institute of Information and Communication
Technologies, Electronics and Applied Mathematics
Louvain-la-Neuve, Belgium.

# Transceiver design for single-cell and multi-cell downlink multiuser MIMO systems

Tadilo Endeshaw Bogale

Thesis presented in partial fulfillment of the requirements
for the PhD degree in Engineering Sciences

PhD Committee:

| | |
|---|---|
| Luc Vandendorpe, *Supervisor* | Université catholique de Louvain, Belgium |
| Michel Verleysen, *President* | Université catholique de Louvain, Belgium |
| Jérôme Louveaux | Université catholique de Louvain, Belgium |
| Francois Horlin | Université Libre de Bruxelles, Belgium |
| Bjorn Ottersten | Université du Luxembourg, Luxembourg |
| Marc Moonen | Katholieke Universiteit Leuven, Belgium |

December 2013

# Abstract

This thesis designs linear transmitters and receivers (i.e., transceivers) for the downlink multiuser multiple input multiple output (MIMO) uncoordinated and coordinated base station (BS) systems. The transmitters and receivers are designed by assuming perfect and imperfect channel state information (CSI) at the BSs and mobile stations (MSs). Different signal to interference plus noise ratio (SINR), mean square error (MSE) and rate-based design criteria are considered. These design criteria are formulated by considering total BS, per BS antenna, per user, per symbol or a combination of per BS antenna and per user (symbol) power constraints. Centralized and distributed algorithms are proposed to examine the design criteria.

For the total BS power constrained robust transceiver design problems, we propose uplink-downlink duality based solutions. These duality are established just by transforming the power allocation matrices (which are diagonal) from uplink to downlink channel and vice-versa. And for the more generalized power constrained robust transceiver design problems, we propose new MSE downlink-interference duality based solutions. The new MSE downlink-interference duality are established by formulating the noise covariance matrices of the interference channels as fixed point and marginally stable (convergent) discrete time switched systems for weighted sum rate/MSE and rate/SINR/MSE-based problems, respectively.

We have shown that the weighted sum rate maximization problem can be equivalently formulated as weighted sum MSE minimization problem with additional optimization variables and constraints. We also develop distributed transceiver design algorithms to solve weighted sum rate and MSE optimization problems for coordinated BS systems. The distributed transceiver design algorithms employ modified matrix fractional minimization and Lagrangian dual decomposition methods.

# Acknowledgments

First of all I would like to express my deepest gratitude to my supervisor Professor Luc Vandendorpe for giving me the opportunity to work in such ambitious research project. I am especially grateful for his confidence and the freedom he gave me to do this research. I also appreciate his stimulating suggestions and guidance during this research work. I would also like to thank Dr. Batu Krishna Chalise who was working with me closely when he was at the University Catholique de Louvain (UCL). He is now with the Center for Advanced Communications, Villanova University, Villanova, PA 19085 USA. I would like to thank all the jury members for agreeing to read this thesis and participate in both the private and public defence.

I would like to acknowledge the financial support of the Region Wallonne for the project MIMOCOM in the framework of which this work has been achieved. I also thank SES and BELSPO for the financial support of the COGRADIO project and IAP project BESTCOM.

I thank the UCL Institute of Information and Communication Technologies, Electronics and Applied Mathematics (ICTEAM) research members, Jin Zhiwen, Ivan Stupia, Achraf Mallat, Adrià Gusi Amigó, Tao Wang, Zohaib Hassan Awan, Mohieddine El Soussi, Achraf Mallat and Felix Brah, whom I had the opportunity to discuss with in a series of seminars and conferences.

I thank the secretaries, Jean Deschuyter and Marie Helene, who arrange me several travels to present research papers at different conferences. I also thank Francois Hubin for his technical support related to computer and networks.

Lastly, I would like to thank my family for giving me unconditional support during the period of my PhD. They also encouraging me to be patient, strong and consistence.

# Notations

| | |
|---|---|
| $(.)_{(n,n)}$ | denotes the $(n, n)$ element of a matrix. |
| $(.)_{(n,:)}$ | denotes the $(n, :)$th row of a matrix. |
| $(.)^T$ | denotes transpose. |
| $(.)^H$ | denotes conjugate transpose. |
| $(.)^*$ | denotes conjugate. |
| $(.)^\star$ | denotes optimal solution. |
| $(.)^{DL}/(.)^{UL}/(.)^I$ | denotes downlink/uplink/interference. |
| $\|(.)\|_n / \|\|\|(.)\|\|\|_n$ | denotes the $n$th norm of a vector (matrix). |
| $\mathbb{C}^{N\times M}$ ($\Re^{N\times M}$) | denotes an $N \times M$ matrix with complex (real) entries. |
| diag(.) /blkdiag(.) | denotes diagonal/block diagonal. |
| $\mathrm{E}\{(.)\}$ | denotes expected value. |
| $\mathbf{I}_n(\mathbf{I})$ | denotes an identity matrix of size n (appropriate size). |
| Lower (upper) case letters | denotes column vectors (matrices). |
| $\Re\{.\}(\Im\{.\})$ | denotes real (imaginary). |
| s.t | denotes subject to. |
| $\mathrm{tr}\{(.)\}$ | denotes trace. |
| vec(.) | denotes vectorization of a matrix. |

# Acronyms

| | |
|---|---|
| AMSE | Average mean square error. |
| ASER | Average symbol error rate. |
| AWGN | Additive white Gaussian noise. |
| BER | Bit error rate. |
| BC | Broadcast channel. |
| BS | Base station. |
| CSI | Channel state information. |
| DPC | Dirty paper coding. |
| FDD | Frequency division duplex. |
| GM | Global minimum. |
| GP | Geometric program. |
| GSM | Global system for mobile communications. |
| i.i.d | Independent and identically distributed. |
| ISI | Inter symbol interference. |
| MAC | Multiple access channel. |
| MF | Matched filtering. |
| IZF | Improved zero forcing. |
| KKT | Karush Kuhn Tucker. |
| MAMSE | Minimum average mean square error. |
| MGO | Monotonic global optimization. |
| MIMO | Multiple input multiple output. |
| MISO | Multiple input single output. |
| ML | Maximum likelihood. |
| MMSE | Minimum mean square error. |
| MS | Mobile station. |
| MSE | Mean square error. |
| QoS | Quality of service. |
| SIC | Successive interference cancelation. |
| SINR | Signal to interference plus noise ratio. |
| SISO | Single input single output. |
| SIMO | Single input multiple output. |
| SVD | Singular value decomposition. |
| SNR | Signal to noise ratio. |
| SOC | Second order cone. |
| SOCP | Second order cone programming. |
| SQP | Sequential quadratic programming. |
| TDD | Time division duplex. |
| WMSE | Weighted mean square error. |
| WSMSE | Weighted sum mean square error. |
| ZMCSCG | Zero mean circularly symmetric complex Gaussian. |

# Contents

# Introduction

# 1

## 1.1 Introduction to digital communication

Now a days, communities in almost all countries of the world are made up of a people from varied backgrounds. Our colleagues may be from different races, faiths, cultures and experiences. We have different preferences in living style, music, food etc. Given these huge differences, human relation activities help communities become more harmonious, respectful, and cohesive. And, mutual understanding is the core of human relations. It is apparent that there can be no mutual understanding without communication.

Communication enables us to do plenty of important things: to grow, to be aware of ourselves, to adjust to our environment, to educate people and to express a myriad of emotions and needs. Excellent communication skill is of course central to a civilized society.

In engineering, the various communication disciplines have the purpose of providing technological aids to human communication. One could view the drum rolls of primitive societies as being technological aids to communication, but communication technology as we view it today became important with telegraphy, then telephony, then video, then computer communication, then internet and so on. Surprisingly, all of these services are currently available in almost all personal mobile telephones.

Initially all of these technologies were developed as separate networks and were viewed as having little in common. As these networks grew, however, the fact that all parts of a given network had to work together, coupled with the fact that different components were developed at different times using different design methodologies, caused an increased focus on the underlying principles and architectural understanding required for continued system evolution.

Perhaps the greatest contribution for the principles of digital communication was the creation of Information Theory by Claude Shannon [Sha48]. After this theory, both the device technology and the engineering understanding of the theory were sufficient to enable system developments follow information theoretic principles. One key benefit of this theory is it help the engineers view all communication sources, e.g., text, audio, video, and image signals, as being representable in binary sequences. The other key benefit is it helps the engineers design communication systems that first convert the source signal into a binary sequence and then convert that binary sequence into a form suitable for transmission over a particular physical media such as cable, twisted wire pair, optical fiber or electromagnetic radiation through space. The succuss of this theory (and many other communication theories and principles) allows end user community enjoy compact and low cost devices with amazing mixtures of services, like data, audio and video (for instance global system for mobile communications (GSM) phones, smart-phones, I-phones, Laptop computers, etc).

## 1.2 Baseband digital communication system

A digital communication requires mapping information bits into symbols. As an example, to transmit digital information bits 1 and 0, one can map these two information bits as one complex symbol $d = 1 - j1$. Such a mapping is termed as "modulation", Quadrature phase shift keying (QPSK) modulation [BL04, Pro01]. For a given modulation scheme, there is strictly a one to one relation between a symbol and a sequence of information bits. Thus, without loss of generality, a communication system can be designed to transmit a symbol and reliably recover the transmitted symbol from the received signal. If the recovered symbol is error free so is the recovered information bits.

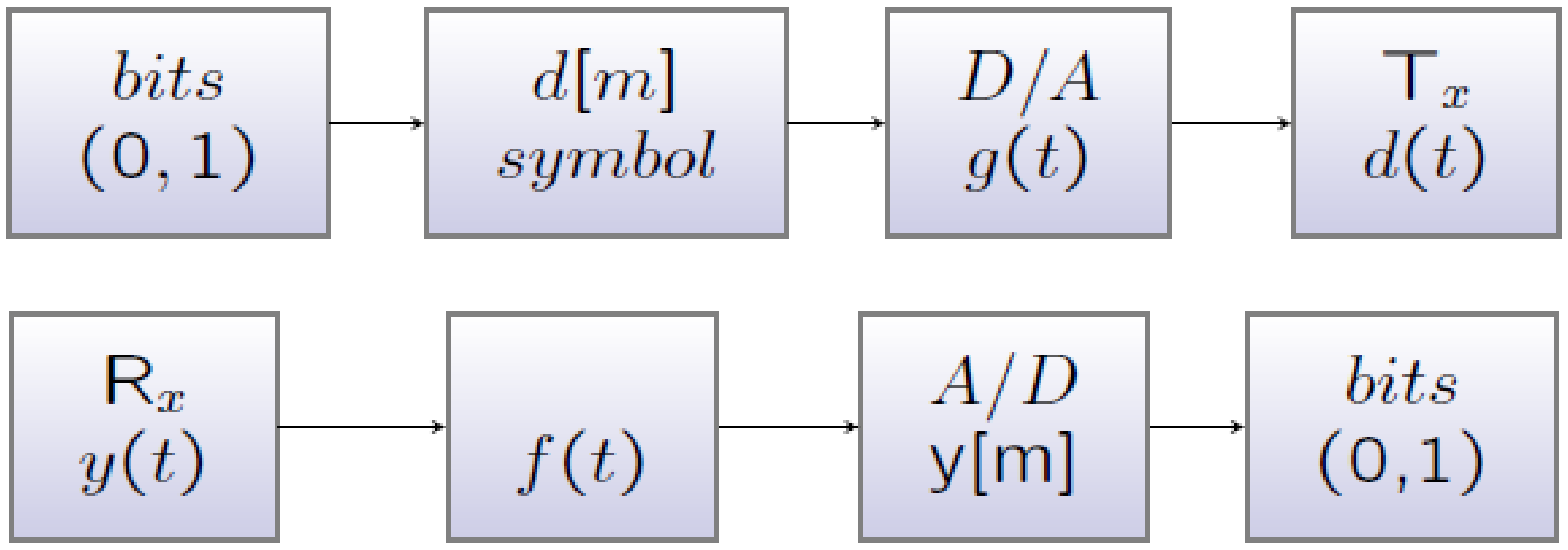

**Figure 1.1** Basic building blocks of baseband digital communication system: [upper] Transmitter. [lower] Receiver.

The basic building blocks of a baseband digital communication system is shown in Fig. 1.1. At the transmitter, digital information bits are mapped into a discrete sequence of symbols (i.e., modulation). Next, these symbols are interpolated to produce filtered analog waveform which is then transmitted over a physical medium: this is the D/A (digital to analog conversion) step.

At the receiver, the received analog signal is first filtered. Then, the filtered signal is sampled to produce a discrete sequence of symbols; this is the A/D (analog to digital conversion) step. Finally, the discrete sequence of symbols are de-mapped to the information bits; this is the demodulation step.

For better exposition, we assume that the transmission medium is noise free. By applying the principles of signals and systems, the baseband transmitted signal $d(t)$ is expressed as

$$d(t) = \sum_{k=-\infty}^{\infty} d[k]g(t - kT_s) \tag{1.1}$$

where $T_s$ is the symbol period (the bandwidth of the transmitted signal is $B_s = \frac{1}{T_s}$Hz). And the received signal after filtered by $f(t)$ is given as[1]

$$\begin{aligned}\tilde{y}(t) &= \int_{-\infty}^{\infty} f^*(\tau)y(t-\tau)d\tau \\ &= \int_{-\infty}^{\infty} f(\tau)d(t-\tau)d\tau\end{aligned}$$

[1]It is assumed that the bandwidth of $f(t)$ and $g(t)$ are the same.

$$
\begin{aligned}
&= \sum_{k=-\infty}^{\infty} d[k] \int_{-\infty}^{\infty} f(\tau) g(t-\tau-kT_s) d\tau \\
&= \sum_{k=-\infty}^{\infty} d[k] \tilde{g}(t-kT_s)
\end{aligned}
\tag{1.2}
$$

where $\tilde{g}(t) = \int_{-\infty}^{\infty} f(\tau)g(t-\tau)d\tau$. As the transmission is noise free, $\tilde{g}(t)$ should be selected such that $\tilde{y}(mT_s) = d[m], \forall m$ holds true (i.e., after the analog to digital conversion). The sampled signal at $t = mT_s$ is expressed as

$$
\begin{aligned}
y[m] =& \tilde{y}(mT_s) = \sum_{k=-\infty}^{\infty} d[k]\tilde{g}((m-k)T_s) \\
=& d[m]\tilde{g}(0) + \sum_{k=-\infty, k\neq m}^{\infty} d[k]\tilde{g}((m-k)T_s).
\end{aligned}
$$

As can be seen from this expression, to achieve $y[m] = d[m]$, we should select $\tilde{g}(t)$ such that $\tilde{g}(0) = 1, \tilde{g}(kT_s) = 0, k \geq 1$ with accurate initial sampling time (i.e., this holds true when the transmitter and receiver are symbol synchronous). A filter satisfying these criteria is called inter symbol interference (ISI) free filter. Example of such a filter satisfying the bandwidth requirement $B_s = \frac{1}{T_s}$ is $\tilde{g}(t) = \text{sinc}(\frac{t}{T_s})$ (i.e., a sinc filter). However, since this filter requires infinite time duration, it is not suitable for practical realization [BL04, Gol05]. For this reason, a modified version of a sinc function, raise cosine filter, is commonly used in practice which is mathematically expressed as

$$
\tilde{g}(t) = \text{sinc}\left(\frac{t}{T_s}\right) \frac{\cos\left(\frac{\pi\beta t}{T_s}\right)}{1 - \frac{4\pi^2\beta^2 t^2}{T_s^2}}. \tag{1.3}
$$

The bandwidth of this filter is given by $B = \frac{1}{T_s}(1+\beta)$, where $\beta$ is termed as "excess bandwidth". By appropriately choosing $\beta$, one can control both the bandwidth and the number of interpolated symbols; increasing $\beta$ will increase the bandwidth while reducing the number of interpolated symbols, whereas decreasing $\beta$ will decrease the bandwidth while increasing the number of interpolated symbols. This scenario is a byproduct of time-bandwidth product analysis of a signal [SM05] (i.e., a signal can not be both time limited and band limited simultaneously).

In a practical communication (i.e., under noise and wireline (wireless) channel environment), for a given transmitter filter $g(t)$, choosing $f^*(t) = g(t)$ (called matched filter) achieves the best performance (in terms of signal to noise ratio (SNR)). Thus, when $\tilde{g}(t)$ is raised cosine filter, $f^*(t) = g(t)$ becomes a square root raised cosine filter [Gol05].

We would like to mention here that there are also several other (non) ISI free filters that are applicable for practical purpose, each of these filters has both advantages and disadvantages [BL04, Gol05].

From these explanations, we can understand that a transmitted symbol can be recovered without any error (in noise free transmission) by appropriately designing each sub-block of Fig. 1.1. Furthermore, for practical channels, we can study the transmit and receive signals in terms of discrete time transmitted and received symbols. In the above example, since we transmit and receive one symbol at a time, we call such a communication system as single input single output (SISO) [BL04, Gol05]. In the following section, we examine the performance of SISO systems for practical channels.

## 1.3 SISO

Let us assume a simple communication system where the transmitted symbols $d[m], \forall m$ are corrupted by noise. By incorporating the effect of noise, the discrete time received signal can be expressed as (ignoring the time index)

$$y = d + n. \tag{1.4}$$

When the transmission medium is corrupted only by noise, the transmission channel is termed as additive white Gaussian noise (AWGN) channel.

For better exposition, let us assume that $n \sim \mathcal{CN}(0, N_0)$ (i.e., a zero mean Circularly symmetric complex Gaussian noise with variance $N_0$) and $d$ is a binary phase shift keying (BPSK) symbol (i.e., $d = \pm a$) [BL04, Gol05]. For BPSK symbol transmission, the received signal is given by

$$y = \begin{cases} a + n, & d = a \\ -a + n, & d = -a. \end{cases} \tag{1.5}$$

Now let us detect the transmitted symbol ($d$) from the received signal $y$. This problem can be considered as a simple hypothesis testing problem. The max-imum likelihood (ML) detection rule then becomes [BL04, Kay98]

$$d = \begin{cases} a, & \text{if } \Re\{y\} \geq 0 \\ -a, & \text{if } \Re\{y\} < 0. \end{cases} \tag{1.6}$$

Now if we use the most common definition of SNR [BL04, Gol05]

$$SNR = \frac{\text{Average received signal power per complex symbol time}}{\text{Noise power per complex symbol time}} \tag{1.7}$$

the SNR of $y$ can be expressed as

$$SNR = \frac{a^2}{N_0}. \tag{1.8}$$

Thus, the probability of error $(P_e)^2$ is given as [BL04, Pro01]

$$P_e = Q(\sqrt{2SNR}) \tag{1.9}$$

where

$$Q(x) = \frac{1}{2\pi}\int_x^{\infty} \exp(-\frac{u^2}{2})du.$$

By applying fundamental calculus, one can get the following bounds

$$\begin{aligned} Q(x) &< \exp^{-x^2/2}, & x > 0, \\ Q(x) &> \frac{1}{\sqrt{2\pi x^2}}\left(1 - \frac{1}{x^2}\right)\exp^{-x^2/2}, & x > 1. \end{aligned}$$

In most practical communication systems, $SNR >> 1$. Thus, in an AWGN channel, $P_e$ decays exponentially, i.e.,

$$P_{e_{AWGN}} \sim \exp^{-SNR}. \tag{1.10}$$

For both wireline and wireless communication systems, the transmitted signal is not only corrupted by noise but it may also be faded by the propagation environment [BL04, Pro01]. Thus, it is reasonable to examine the error probability of (1.4) for a fading channel $h$. For this channel, the received signal $y$ can be expressed as

$$y = hd + n. \tag{1.11}$$

For better exposition, let us assume that $h \sim \mathcal{CN}(0, 1)$ (i.e., a unit variance Rayleigh fading channel) [Pro01, BL04, Gol05].

[2] For this example, $P_e$ is the probability of getting $\Re\{y\} > 0$ (i.e., deciding $d$ as $d = a$) when $d = -a$ is transmitted and vice versa.

### 1.3.1 Unknown $h$

In this subsection, we examine the error probability of (1.11) for the scenario where the fading coefficient $h$ is unknown for both the transmitter and receiver.

As we can see from the above signal model, the phase of the transmitted signal $d$ is affected by the channel coefficient. Thus, if we do not know $h$, BPSK symbols can not be reliably detected at the receiver. Due to this reason, when $h$ is not known, one common approach is to incorporate mapping at the transmitter. Now, let us consider the following orthogonal mapping: rather than transmitting $\pm a$ and detect individually, we will map $\pm a$ as $[a, 0]$ and transmit these mapped symbols in two symbol periods, i.e.,

$$\mathbf{d}_0 = [a,\ 0],\ \ \mathbf{d}_1 = [0,\ a]. \tag{1.12}$$

Now in two consecutive symbol periods (indexed as 0 and 1), we will have

$$\begin{aligned} y[0] =& h[0]a + n[0], \quad y[1] = n[1], \qquad\qquad \text{If } \mathbf{d}_0 \text{ is Transmitted} \\ y[0] =& n[0], \qquad\qquad\quad y[1] = h[1]a + n[1], \ \ \text{If } \mathbf{d}_1 \text{ is Transmitted.} \end{aligned}$$

For the received signals $y[0]$ and $y[1]$, the ML detection rule turns out to be a simple comparison of the powers of the received signals $y[0]$ and $y[1]$. And the decision rule then becomes [Pro01, BL04, Gol05]

$$\begin{aligned} [y[0],\ y[1]] &= \mathbf{d}_0,\ \ \text{If } |y[0]|^2 \geq |y[1]|^2 \\ [y[0],\ y[1]] &= \mathbf{d}_1,\ \ \text{If } |y[0]|^2 < |y[1]|^2. \end{aligned} \tag{1.13}$$

As this detector does not exploit the channel state information, it is called "Non coherent" detector [TV05]. For this detector, the $P_e$ is given as [TV05, Pro01, BL04, Gol05]

$$P_e = \frac{1}{2(1 + SNR)}. \tag{1.14}$$

We would like to mention here that the $P_e$ of this detector is analyzed by applying geometric approach of [TV05].

### 1.3.2 Known $h$

As we can see from the above subsection, we employ symbol mapping to counteract the effect of phase change due to the channel. From fundamental digital communication, it is well known that channel gain $h$ can be estimated at the receiver by employing a training (pilot) symbols [Pro01, BL04, Gol05]. In this section, we will see the $P_e$ of (1.11) by assuming that $h$ is perfectly estimated at the receiver.

If $h$ is known at the receiver, one straightforward approach of eliminating the effect of phase change of the transmitted symbol $d$ is just to multiply the received signal $y$ by the conjugate of the channel $h$ (it is called matched filtering (MF)). By doing so, we will get

$$h^* y = |h|^2 d + h^* n. \tag{1.15}$$

As we can see, such multiplication will not affect the phase of transmitted signal $d$. For a given $h$, the decision statistics can be obtained like in the AWGN channel and is given as

$$d = \begin{cases} a, & \text{if } \Re\{h^\star y\} \geq 0 \\ -a, & \text{if } \Re\{h^\star y\} < 0. \end{cases} \tag{1.16}$$

For the given $h$, the $P_e$ is then becomes

$$P_{e|h} = Q\left(\frac{a|h|}{\sqrt{N_0/2}}\right) = Q\left(\sqrt{2|h|^2 SNR}\right). \tag{1.17}$$

As this detector does exploits the channel state information, it is called "coherent" detector [TV05]. For a unit variance Rayleigh fading channel $h$, the average $P_e$ is given by [Pro01, BL04]

$$\begin{aligned} P_e = \mathrm{E}_h\{P_{e|h}\} &= \mathrm{E}_h\left[Q\left(\sqrt{2|h|^2 SNR}\right)\right] \\ &= 0.5\left(1 - \sqrt{\frac{SNR}{1+SNR}}\right). \end{aligned}$$

Fig. 1.2 shows the $P_e$ for different SNR for AWGN and Rayleigh fading channels. From this figure, we can observe that the performance of AWGN channel is significantly better than that of the Rayleigh fading channel. This is evidently seen from the $P_e$ expressions that for AWGN channel $P_e$ decays exponentially, whereas for the Rayleigh fading channel, $P_e$ decay only linearly. From this figure, we can also observe that for a Rayleigh fading channel, to achieve $P_e = 10^{-4}$, the AWGN channel requires 8dB, whereas the Rayleigh fading channel with and without the channel knowledge $h$ require 34dB and 37dB, respectively. This fact shows that for a Rayleigh fading channel, the gain due to the channel knowledge at the receiver is not significant (3dB). However, the gain between the AWGN and Rayleigh fading channels is around 26dB (which is huge).

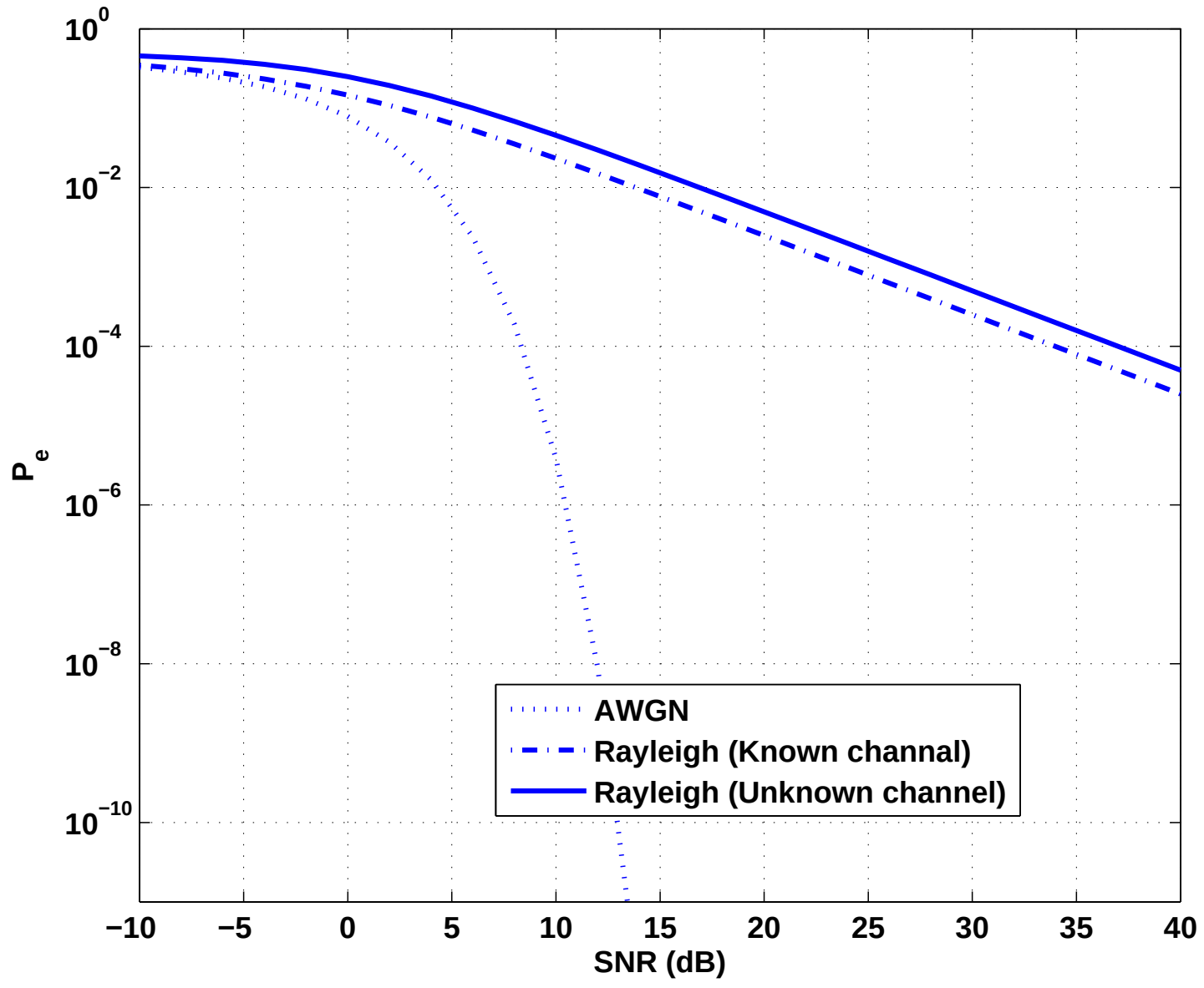

**Figure 1.2** Probability of error for SISO system under BPSK transmission.

Now let us examine closely the reason why $P_e$ is very high in a Rayleigh fading channel even under the assumption of known channel knowledge at the receiver. By considering (1.15), the received signal can be expressed as

$$h^* y = |h|^2 d + h^* n. \tag{1.18}$$

When $h$ is a Rayleigh distributed random variable, $|h|^2$ will have chi-square distributed random variable with 2 degrees of freedom i.e., $|h|^2 \sim \chi_2^2$. The probability density function of $\chi_2^2$ distributed random variable is plotted in Fig. 1.3. As can be seen from this figure, there is a significant probability that $|h|^2$ is closer to zero (i.e., the channel experience "deep fading"). From this explanation, we can understand that signal detection in a fading channel has poor performance. And the reason why detection in the fading channel has poor performance is not because of the lack of knowledge of the channel at the receiver. It is due to the fact that the channel gain is random and there is a significant probability that the channel is in a "deep fade".

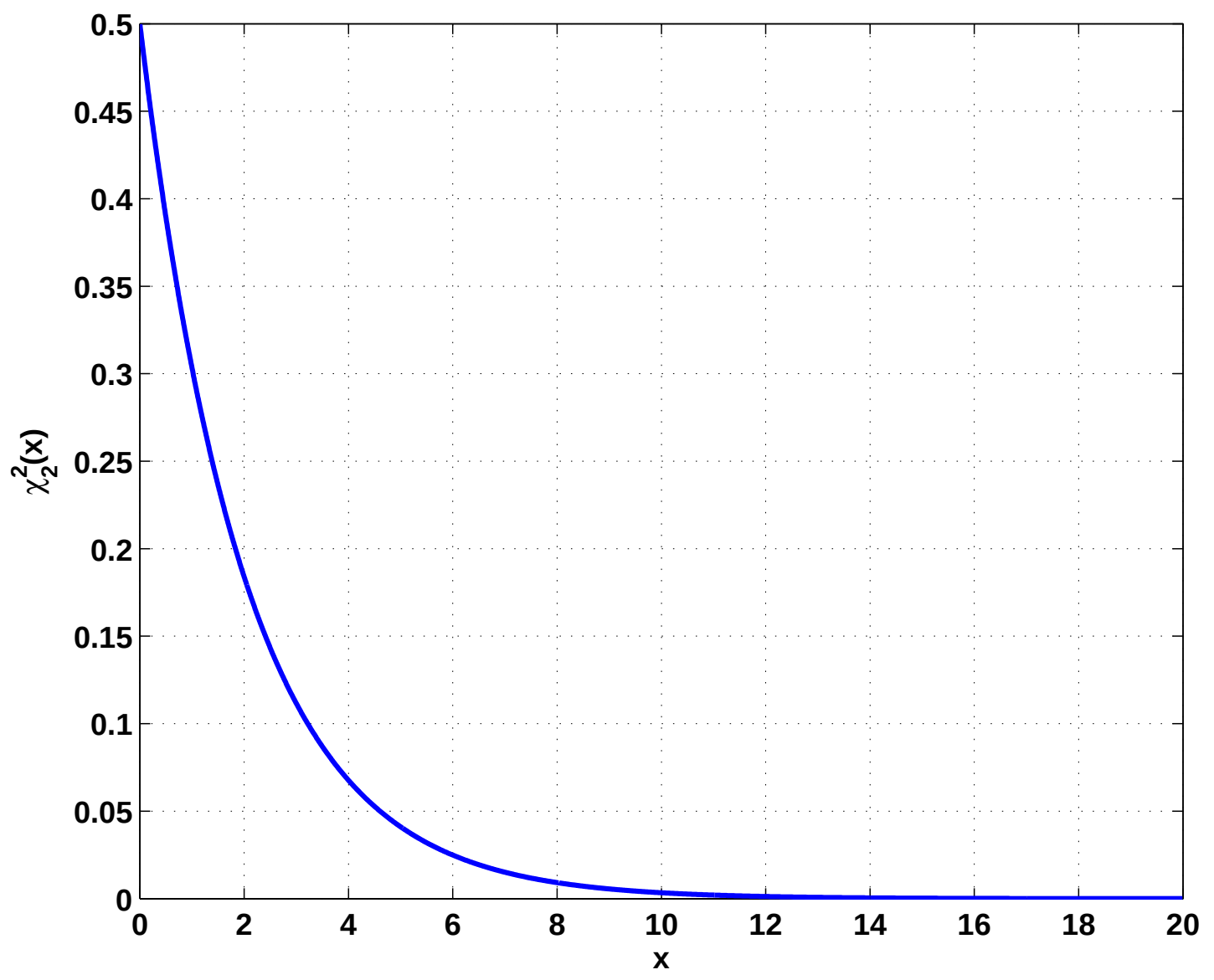

**Figure 1.3** The probability density function of chi-square distribution with 2 degrees of freedom.

It can be easily seen that the root cause of this poor performance is that reliable communication depends on the strength of only one signal path. There is a significant probability that this path will be in a deep fade. When the path is in a deep fade, any communication scheme will likely suffer from errors. A natural solution to improve the performance is to ensure that the information symbols pass through multiple signal paths, each of which fades independently, making sure that reliable communication is possible as long as one of the paths is strong. This technique is called "diversity", and as will be clear in the sequel, diversity can dramatically improve the performance of a wireless system in fading channels [Pro01, BL04].

## 1.4 Diversity

There are many approaches to achieve diversity. The most widely used approaches are diversity over time, diversity over frequency and diversity over space [Gol05]. Diversity over time can be achieved by employing coding and interleaving of the original message symbols. Frequency diversity can also be achieved when the channel is frequency selective. And space diversity is achieved by applying multiple transmit and (or) receive antennas. In this thesis, we will provide a brief summary of time and space diversity techniques as these techniques can be applied for flat fading channels (which is the focus of this thesis). To keep the discussion simple, we assume that the receiver has perfect knowledge of the channel gains and can coherently combine the received signals in the diversity paths.

### 1.4.1 Time diversity

Time diversity can be exploited by several ways. One simple time diversity scheme is based on repetition coding: the same information symbol is transmitted over several signal paths. Time diversity is achieved by averaging the fading of the channel over time. In a typical communication system, the channel coherence time is of the order of tens to hundreds of symbols, and therefore the channel is highly correlated across consecutive symbols. To ensure that the symbols are transmitted through independent (nearly independent) fading gains, interleaving of the message symbols is employed [BL04, Gol05].

Assuming ideal interleaving so that consecutive symbols $\{d_l\}_{l=1}^{L}$ are transmitted sufficiently far apart in time[3], the $l$th received signal can be expressed as

$$y_l = h_l d_l + n_l, \qquad l = 1, \cdots, L. \tag{1.19}$$

[3]In practice, usually $\{d_l\}_{l=1}^{L}$ are codewords. But here we consider it as symbols for better understanding about time diversity.

In a repetition coding, we will transmit $\{x_l = x_1\}_{l=1}^{L}$. Under such transmission, the overall received signal can be expressed as

$$\mathbf{y} = \mathbf{h}d_1 + \mathbf{n} \tag{1.20}$$

where $\mathbf{y} = [y_1, \cdots, y_L]^T$, $\mathbf{h} = [h_1, \cdots, h_L]^T$ and $\mathbf{n} = [n_1, \cdots, n_L]^T$. Applying matched filtering at the receiver gives us

$$\mathbf{h}^H\mathbf{y} = |\mathbf{h}|^2 d_1 + \mathbf{h}^H\mathbf{n}. \tag{1.21}$$

Now if $d_1$ is a BPSK signal, the probability of error can be computed as

$$P_{e|\mathbf{h}} = Q\left(\sqrt{2|\mathbf{h}|^2 SNR}\right). \tag{1.22}$$

Under the assumption of a unit variance Rayleigh fading channel (i.e., $\{h_l \sim \mathcal{CN}(0,1)\}_{l=1}^{L}$), $|\mathbf{h}|^2$ will have a Chi-square random variable with $2L$ degrees of freedom. The probability density function of a Chi-square random variable with $2L$ degrees of freedom is given by

$$f(x) = \frac{1}{(L-1)!}x^{L-1}\exp^{-x}, \qquad x \geq 0.$$

The average probability of error can be computed as

$$\begin{aligned} P_e =& \int_0^\infty Q(\sqrt{2xSNR})f(x)dx \\ =& \left(\frac{1-\mu}{2}\right)^L \sum_{l=0}^{L-1} \left(\begin{array}{c} L-1+l \\ l \end{array}\right) \left(\frac{1+\mu}{2}\right)^l \\ \approx& O\left(\frac{1}{SNR^L}\right) \end{aligned} \tag{1.23}$$

where $\mu = \sqrt{\frac{SNR}{1+SNR}}$. The $P_e$ is plotted for different $L$ in Fig. 1.4. As can be seen from this figure, increasing $L$ dramatically decreases the error probability ($P_e$). This is due to the fact that the tail of $|\mathbf{h}|^2$ near zero decreases as $L$ increases.

As we can see, this time diversity achieves a diversity gain. However, such diversity gain is at the expense of a more transmission power and reduced data transmission rate (i.e., to transmit one symbol we need to spend L times the power of each symbol and we also require L symbol periods). Thus, although such a diversity has some gain in $P_e$, it is not power and spectrally efficient. There are also other complicated time diversity schemed. However, those schemes will not be presented in this thesis (The reader can see [BL04, Gol05] for many other time diversity schemes).

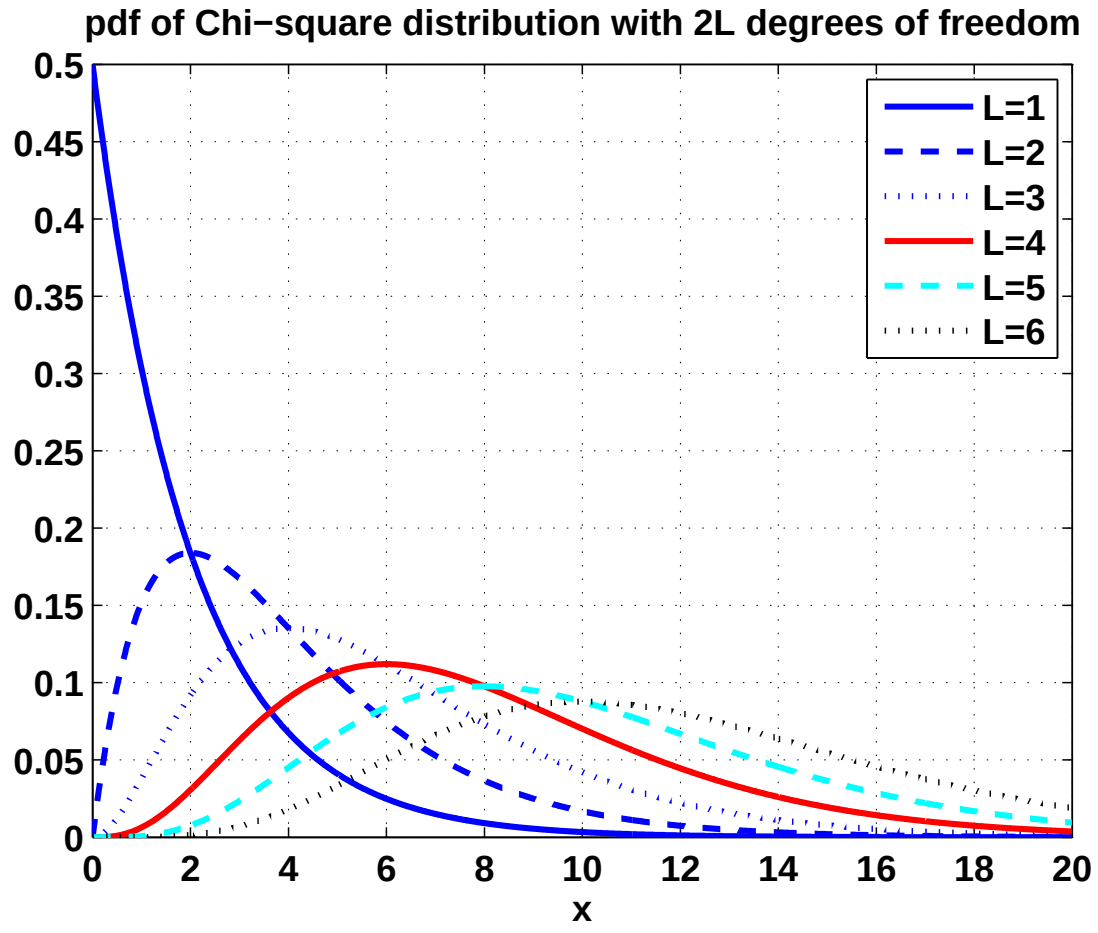

**Figure 1.4** The probability density function of Chi-square distribution with different degrees of freedom.

### 1.4.2 Space diversity

As can be understood from the above discussion, to exploit time diversity, interleaving and coding over several channel coherence time periods is necessary. Furthermore, time diversity is not power and spectrally efficient. When the coherence time of the channel is large (and when we need more power and spectrally efficient communication), time diversity can not be applied. In such a case spatial diversity can be obtained by placing multiple antennas at the transmitter and (or) the receiver. If the antennas are placed sufficiently far apart, the channel gains between different antenna pairs fade almost independently, which can consequently create independent signal paths. The required antenna separation depends on the local scattering environment as well as on the carrier frequency. For a mobile which is near the ground with many scatterers around, the channel becomes uncorrelated in a shorter spatial distances, and typical antenna separation of half to one carrier wavelength is sufficient. For base stations on high towers, larger antenna separation of several to tens of wavelengths may be required [BL04, Gol05].

#### 1.4.2.1 Receive diversity

In a system with one transmit antenna and $L$ receive antennas (such a system is termed as single input multiple output (SIMO)), the relation between the transmitted signal $d$ and received signal $y_l$ is given as

$$\begin{aligned} y_l =& h_l d + n_l, \qquad l = 1, \cdots, L \\ \mathbf{y} =& \mathbf{h} d + \mathbf{n} \end{aligned} \tag{1.24}$$

where $h_l$, $y_l$ and $n_l$ are the channel gain, received signal and noise at the $l$th receive antenna, respectively, $d$ is the transmitted symbol, $\mathbf{y} = [y_1, \cdots, y_L]^T$, $\mathbf{h} = [h_1, \cdots, h_L]^T$ and $\mathbf{n} = [n_1, \cdots, n_L]^T$.

This is exactly the same detection problem as in the use of a repetition code and interleaving over time, with $L$ diversity branches now over space instead of over time. Thus, the detection of $d$ from the received signal $\mathbf{y}$ can be performed like in the above section. By doing so, we can achieve the following error probability:

$$P_e = Q(\sqrt{2|\mathbf{h}|^2 SNR}).$$

The gain achieved by employing multiple receive antenna can be decomposed as

$$|\mathbf{h}|^2 SNR = L.SNR \times \frac{1}{L}|\mathbf{h}[m]|^2.$$

One can interpret the above gain as two gains one is due to $L.SNR$ and the other is due to $\frac{1}{L}|\mathbf{h}[m]|^2$. The former gain is termed as an "array gain" (it is always achieved when $L > 1$, even in an AWGN channel)[4], the latter gain is diversity gain (full diversity gain is achieved when all of the $L$ paths are independent) [BL04, Gol05]. By examining $\frac{1}{L}|\mathbf{h}[m]|^2$ for large $L$, one can obtain that $\frac{1}{L}|\mathbf{h}[m]|^2$ tends to one. Consequently, for large $L$, applying multiple receive antenna will only achieve array gain.

[4]We would like to mention here that although the $P_e$ of the time diversity repetition coding is decomposed as the above expression, one can not call the term $L \times SNR$ of repetition coding as "power gain". This is due to the fact that in the repetition coding case, the increase in received SNR comes from increasing the total transmit power required to send a single symbol.

#### 1.4.2.2 Transmit diversity

Now consider the case when there are $L$ transmit antennas and 1 receive antenna (such a system is literary termed as multiple input single output (MISO)). For such a system, it is easy to get a diversity gain of $L$: simply transmit the same symbol over the L different antennas during $L$ symbol times. At any one symbol time, only one antenna is turned on and the rest are silent. This kind of transmission is similar to a repetition code, and, as we have seen in the previous section, repetition codes are not power and spectrally efficient.

There have been a lot of research activities on the diversity achieving transmission strategies for MISO systems. One of the most simple and well known diversity achieving strategy is the so-called Alamouti space-time coding scheme [BL04, Ala98]. This scheme is designed for 2 transmit and $L$ receive antenna scenarios. The generalization of Alamouti scheme for any transmit and receive antennas is still an ongoing research.

#### 1.4.2.3 Transmit and receive diversity

Now suppose that there are $N$ transmit and $L$ receive antennas (such a system is termed as multiple input multiple output (MIMO)). In such a case, the same repetition scheme described in the aforementioned section can achieve diversity: transmit the same symbol over the $N$ antennas in $N$ consecutive symbol times (at each time, nothing is sent over the other antennas). By doing so, we will have a diversity order of $NL$. However, this kind of transmission is similar to a repetition code which is not power and spectrally efficient. Thus, one can improve the performance by employing Alamouti scheme [Ala98]: transmit 2 symbols over the $N$ antennas in $N/2$ consecutive symbol times (i.e., at each time, only 2 transmit antennas are active). By doing so, one may get better power and (or) diversity gain, and spectral efficiency compared to SISO, SIMO and MISO scenarios.

For better explanation, let us consider a MIMO system with 2 transmit and 2 receive antennas. For this system, the Alamouti scheme achieve a diversity order of 4 [TV05, BL04]. However, this scheme effectively transmits 1 symbol in 1 symbol period[5]. Now, obviously one may ask the following question: For a 2 by 2 channel, is there any transmission scheme that can transmit 2 symbols in 1 symbol period while ensuring a diversity order of 4? The following discussion addresses this issue.

[5] This is due to the fact that this scheme transmits 2 symbols in 2 symbol periods.

## 1.5 2 by 2 MIMO: Multiplexing versus diversity

For a 2 by 2 MIMO system, we will get the following input output relation

$$\mathbf{y} = \mathbf{H}\mathbf{d} + \mathbf{n} \tag{1.25}$$

where $\mathbf{d}$ ($\mathbf{y}$) is the transmitted (received) signal vector, $\mathbf{n}$ is the additive noise vector and $\mathbf{H}$ is the channel matrix between the transmitter and receiver, i.e.,

$$\mathbf{y} = \begin{bmatrix} y_1 \\ y_2 \end{bmatrix}, \ \mathbf{H} = \begin{bmatrix} h_{11} & h_{21} \\ h_{12} & h_{22} \end{bmatrix}, \mathbf{d} = \begin{bmatrix} d_1 \\ d_2 \end{bmatrix}, \mathbf{n} = \begin{bmatrix} n_1 \\ n_2 \end{bmatrix}.$$

It is assumed that $d_1$ and $d_2$ are i.i.d symbols. This is an example of a spatial multiplexing scheme: independent data streams are multiplexed in space. Thus, the aim is therefore to design a transmission and reception scheme achieving a diversity gain of 4 and multiplexing gain of 2.

### 1.5.1 H is known only at the receiver

If $\mathbf{H}$ is known only at the receiver, one simple approach of recovering the transmitted symbols $d_1$ and $d_2$ is to apply channel inversion, i.e.,

$$\tilde{\mathbf{y}} = \mathbf{H}^{-1}\mathbf{y} = \mathbf{d} + \mathbf{H}^{-1}\mathbf{n}. \tag{1.26}$$

Now let us examine the performance of this receiver approach. For simplicity, we examine the first term of $\tilde{\mathbf{y}}$ (i.e., $\tilde{y}_1$)

$$\begin{aligned} \tilde{y}_1 =& d_1 + \mathbf{H}^{-1}\mathbf{n} \\ =& d_1 + \frac{h_{22}n_1 - h_{12}n_2}{h_{11}h_{22} - h_{12}h_{21}}. \end{aligned} \tag{1.27}$$

The scaled version of $\tilde{y}_1$ is

$$\begin{aligned} \tilde{y}_1' =& \frac{h_{11}h_{22} - h_{12}h_{21}}{\sqrt{|h_{22}|^2 + |h_{12}|^2}}\tilde{y}_1 \\ =& \frac{h_{11}h_{22} - h_{12}h_{21}}{\sqrt{|h_{22}|^2 + |h_{12}|^2}}d_1 + \frac{h_{22}n_1 - h_{12}n_2}{\sqrt{|h_{22}|^2 + |h_{12}|^2}}. \end{aligned}$$

This equation has the same mathematical expression as a Rayleigh faded SISO channel. Therefore, with the channel inversion approach, it is possible to transmit 2 symbols in one symbol period (i.e., we achieve a multiplexing gain of 2) but with a unit diversity gain[6].

The other most common approaches of recovering the transmitted symbols $d_1$ and $d_2$ are minimum mean square error (MMSE), MF, combination of MMSE and successive interference cancelation (SIC), and the ML detection approach [BL04, Fos96]. Each of these approaches has its own advantage and disadvantages. Out of these approaches, the ML detector achieves the best performance: it achieves multiplexing and diversity gain of 2 [BL04, Fos96]. The main disadvantage of this approach is that its complexity grows exponentially with the number of antennas for general MIMO channel [TV05, BL04]. Thus, the aim should be the development of practically realizable algorithm achieving full diversity and multiplexing gain.

From fundamental digital communication theory [Pro01, BL04], we know that the transmitter can get the CSI knowledge by exploiting the channel reciprocity in the case of time division duplex (TDD) system and through feedback in the case of frequency division duplex (FDD) system. In the following, we examine MIMO transmission schemes by assuming perfect CSI both at the transmitter and receiver.

### 1.5.2 H is known both at the transmitter and receiver

In this subsection, we examine the achievable performance of (1.25) by assuming perfect CSI both at the transmitter and receiver. Taking the singular value decomposition (SVD) of $\mathbf{H} = \mathbf{U}_h\mathbf{D}_h\mathbf{V}_h^H$, where $\mathbf{U}_h$ and $\mathbf{V}_h$ are the unitary matrices and $\mathbf{D}_h$ is the diagonal matrix containing the singular values of $\mathbf{H}$, we can write (1.25) as

$$\mathbf{y} = \mathbf{H}\mathbf{d} + \mathbf{n} = \mathbf{U}_h\mathbf{D}_h\mathbf{V}_h^H\mathbf{d} + \mathbf{n}. \tag{1.28}$$

[6]We would like to mention here that as $d_1$ and $d_2$ are transmitted without coding, these two symbols can be transmitted from 2 geographically separated users. In such a multiuser context, such channel inversion receiver design approach is termed as "zero forcing" or "interference nulling".

Now if the transmitter preprocess the symbols $\mathbf{d}$ as $\tilde{\mathbf{d}} = \mathbf{V}_h\mathbf{d}$, then $|\mathbf{d}|^2 = |\tilde{\mathbf{d}}|^2$ and received signal becomes

$$\begin{aligned}\mathbf{y} =& \mathbf{H}\mathbf{d} + \mathbf{n} \\ =& \mathbf{U}_h\mathbf{D}_h\mathbf{d} + \mathbf{n}.\end{aligned} \tag{1.29}$$

At the receiver, we can apply the following transformation

$$\tilde{\mathbf{y}} = \mathbf{U}^H\mathbf{y} = \mathbf{D}_h\mathbf{d} + \mathbf{U}^H\mathbf{n}. \tag{1.30}$$

As can be seen from this equation, if $n_1$ and $n_2$ are i.i.d Gaussian noise, each entry of $\mathbf{U}^H\mathbf{n}$ is also i.i.d Gaussian. Thus, this expression can be interpreted as a 2 Gaussian parallel channels. From matrix theory, we have

$$\mathbf{D}_{h_{11}}^2 + \mathbf{D}_{h_{22}}^2 = \mathrm{tr}\{\mathbf{H}\mathbf{H}^H\} = |h_{11}|^2 + |h_{12}|^2 + |h_{21}|^2 + |h_{22}|^2.$$

The condition number of $\mathbf{H}$ is expressed as [HJ85]

$$C(\mathbf{H}) = \frac{\mathbf{D}_{h_{11}}}{\mathbf{D}_{h_{22}}}.$$

We say that $\mathbf{H}$ is well conditioned when $C(\mathbf{H}) \approx 1$. Thus, for a well conditioned channel, this transmission and reception scheme achieve a multiplexing and diversity gain of 2. This scenario exhibits with high probability when the channel gain between each transmit antenna and each receiver antenna fade independently. One key advantage of this approach is it is simple to implement and it can be extended to an arbitrary number of transmit and receive antennas [TV05, BL04].

From this discussion, one may come up with the following question: Given a general MIMO channel $\mathbf{H}$, how much multiplexing and diversity gain can be achieved simultaneously. The detailed analysis of this question has been discussed from the information theoretic point of view in [ZT03]. The reader can refer this paper for more details about the tradeoff between multiplexing gain and diversity gain of a MIMO channel.

Hence, by applying multiple antenna both at the transmitter and receiver, one can achieve diversity, multiplexing and (or) power gains. Due to these key advantages, MIMO systems are proposed in several wireless network standards [BL04, Gol05].

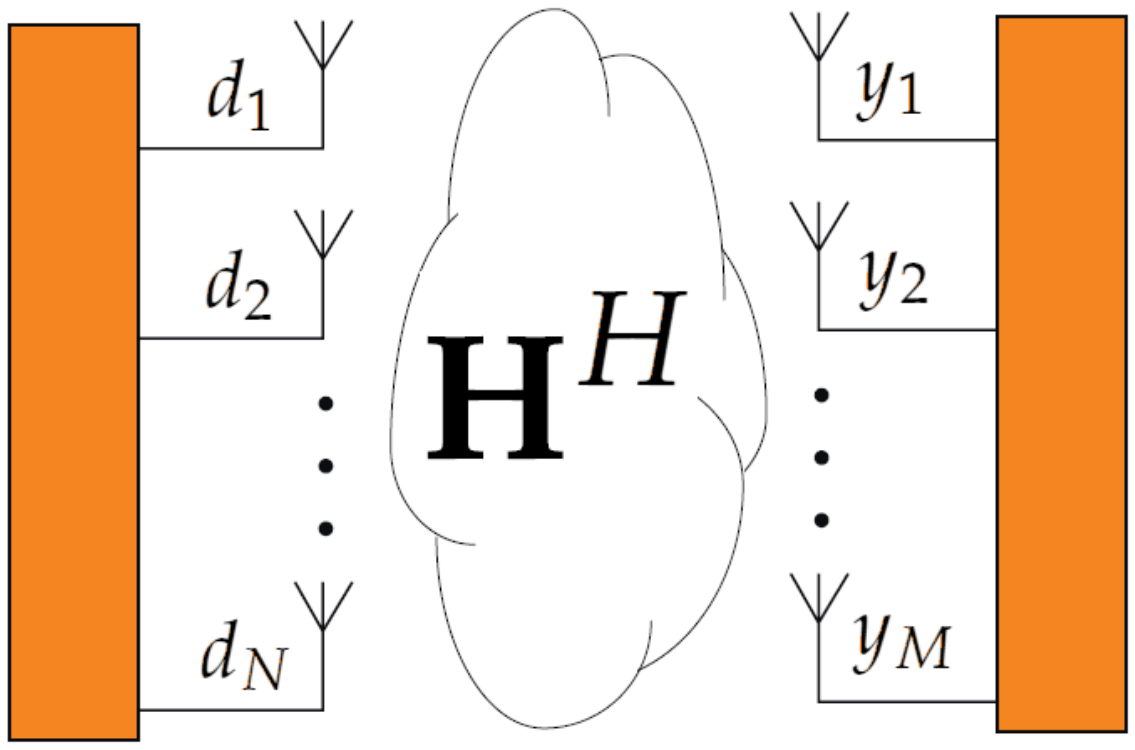

**Figure 1.5** A MIMO system with N transmit and M receive antennas.

## 1.6 MIMO

A MIMO system model with $N$ transmit and $M$ receive antennas is shown in Fig. 1.5. The general mathematical representation of this system is given as [PCL03]

$$\mathbf{y} = \mathbf{H}^H \mathbf{d} + \mathbf{n} \tag{1.31}$$

where $\mathbf{d}$ ($\mathbf{y}$) is the transmitted (received) signal vector, $\mathbf{n}$ is the additive noise vector and $\mathbf{H}^H$ is the channel matrix between the transmitter and receiver[7].

For better flow of the text, we consider the case where the transmitter and receiver have CSI information. In this case the transmitter can preprocess (precod) the data ($\mathbf{d}$) and the receiver can perform postprocessing (decoding) of the received signal ($\mathbf{y}$) to recover $\mathbf{d}$. When the precoding and decoding operations are incorporated, (1.31) can be reexpressed as[8]

$$\widehat{\mathbf{d}} = D(\mathbf{y}) = D(\mathbf{H}^H P(\mathbf{d}) + \mathbf{n}) \tag{1.32}$$

where $\widehat{\mathbf{d}}$ is the estimate of $\mathbf{d}$, and $P(.)$ and $D(.)$ are the precoding and decoding operations, respectively.

[7]We would like to mention here that many other communication systems, like orthogonal frequency division multiplexing (OFDM), digital subscriber line (DSL) and code division multiple access (CDMA) systems can be well represented by (1.31). Thus, this input output relation can describe many classes of practically relevant communication systems. For more detailed explanation, the reader can refer [PCL03]. However, since the aim of this thesis is about transceiver design algorithms for multiuser MIMO system, transceiver design strategies for OFDM, DSL and CDMA systems will not be discussed in this thesis.

When $P(\mathbf{d}) = \mathbf{V}_h\mathbf{d}$ and $D(\mathbf{y}) = \mathbf{U}^H\mathbf{y}$, this expressions turns to that of (1.30). When $P(\mathbf{d}) = \mathbf{d}$, $\mathbf{d}$ is i.i.d and $D(\mathbf{y})$ is a ML estimator of $\mathbf{d}$, the above precoding decoding operation turns out to be that of the approach of [Fos96][9]. Therefore, to exploit the benefits of a MIMO channel, one may need to choose appropriate $P(.)$ and $D(.)$. Now the critical question is how to design $P(.)$ and $D(.)$ for the given channel matrix $\mathbf{H}$ (and the statistical properties of the noise, like covariance matrix)? And are there any general $P(.)$ and $D(.)$ which are optimal for all criteria. To address this, let us consider the following practical design problem:

Assume that we would like to transmit different types of information (for example, text, audio and video information) over a MIMO channel with 3 antennas. We map the text, audio and video information into symbols $d_1$, $d_2$ and $d_3$. It is evident that the bit error probability ($P_e$) requirements of all of these information will not be the same (the $P_e$ corresponding to the video information should be lower than that of the text and audio information). Moreover, as pointed out in [YL07], in a multi-antenna system, each antenna has its own power amplifier and the maximum power of each antenna is limited by some value. Taking into account these quality of service (QoS) requirements and power constraints, one may design the precoders and decoders to minimize the total transmission power while taking into account each symbol QoS requirement and each antenna power constraint. This problem can be formulated as

$$\min_{P(.),\, D(.)} P_1 + P_2 + P_3,$$
$$\text{s.t } P_i \leq P_{am}, \;\; P_{ed_i} \leq P_{emi}, \;\; i = 1, 2, 3$$

[8] Note that the precoder design and decoder design are also termed as transmitter design and receiver design, respectively.

[9] In this paper, such transmission and reception approach is termed as space time coder and decoder

where $\{P_i$ is the power utilized at the $i$th$\}_{i=1}^{3}$ antenna, $P_{am}$ is the maximum power available at each antenna and $\{P_{emi}\}_{i=1}^{3}$ are the maximum allowed $P_e$ for the symbols $\{d_i\}_{i=1}^{3}$.

From these explanations we can understand that (1.32) can act as a general precoding and decoding operation for any MIMO channel. And, all the aforementioned benefits of a MIMO channel (i.e., multiplexing, diversity and power gains) can be appropriately formulated as a precoder/decoder optimization problem. Furthermore, the precoder/decoder design problems and solutions can vary from one design criterion to another.

The aim of this thesis is to examine several practically relevant precoder/decoder design strategies and algorithms for downlink multiuser MIMO systems briefly (i.e., a particular class of MIMO system). We would like to stress here that our objective is not to exploit all the available multiplexing, diversity and power gains of a MIMO channel, rather, we aim to formulate practically relevant precoder/decoder design criteria as optimization problems and try to get the optimal (suboptimal) solution for each problem. In general, the optimal (suboptimal) solution of each problem may not exploit all the available multiplexing, diversity and power gains of the MIMO channel; in other words, the optimal (suboptimal) solution of one problem may not be optimal (even suboptimal) for the other problem.

The precoder and decoder operations can be linear or nonlinear[10]. In general, the best system performance can be achieved by utilizing a nonlinear precoder and decoder. However, since the complexity of such precoder and decoder grows exponentially with the number of antennas (transmitted symbols), nonlinear precoder and decoders are not suitable for practical realization. Due to this fact, linear precoder and decoders are motivated as they are simple to implement at the expense of performance loss [BL04, Fos96, DB09, PLC04]. Thus, the tradeoff between linear and nonlinear precoder(decoder) is complexity versus performance.

[10] An operation is linear when its complexity scales linearly with the number of variables, otherwise it is called a non linear operation.

Under linear precoder and decoder operations, (1.32) can be elegantly expressed as

$$\widehat{\mathbf{d}} = \mathbf{W}^H(\mathbf{H}^H\mathbf{B}\mathbf{d} + \mathbf{n}) \tag{1.33}$$

where $\mathbf{B}$ and $\mathbf{W}^H$ are the precoder and decoder matrices with appropriate dimensions, respectively.

## 1.7 Linear transceiver design for MIMO systems

In this section, we summarize linear precoder and decoder design algorithms for MIMO systems. For better insight of these algorithms, we will provide a concise mathematical description for a particular type of mean square error (MSE)-based design criteria which is presented as follows:

For the given channel $\mathbf{H}$, the MSE between $\mathbf{d}$ and $\widehat{\mathbf{d}}$ of (1.33) is given by

$$\begin{aligned}\xi =& \mathrm{E}\{(\mathbf{d} - \widehat{\mathbf{d}})(\mathbf{d} - \widehat{\mathbf{d}})^H\} \\ =& \mathbf{R}_d + \mathbf{W}^H(\mathbf{H}^H\mathbf{B}\mathbf{R}_d\mathbf{B}^H\mathbf{H} + \mathbf{R}_n)\mathbf{W} - \mathbf{W}^H\mathbf{H}^H\mathbf{B}\mathbf{R}_d - \mathbf{R}_d^H\mathbf{B}^H\mathbf{H}\mathbf{W}\end{aligned} \tag{1.34}$$

where $\mathbf{R}_d(\mathbf{R}_n)$ is the covariance matrix of $\mathbf{d}(\mathbf{n})$ and it is assumed that $\mathbf{d}$ and $\mathbf{n}$ are uncorrelated, i.e., $\mathbf{R}_d = \mathrm{E}\{\mathbf{d}\mathbf{d}^H\}$, $\mathbf{R}_n = \mathrm{E}\{\mathbf{n}\mathbf{n}^H\}$ and $\mathrm{E}\{\mathbf{n}\mathbf{d}^H\} = 0$. One MSE-based precoder(decoder) design criteria could thus be to minimize the overall MSE under a total transmit power constraint. Mathematically, this problem can be formulated as

$$\begin{aligned}&\min_{\mathbf{B},\mathbf{W}} \mathrm{tr}\{\mathbf{R}_d + \mathbf{W}^H(\mathbf{H}^H\mathbf{B}\mathbf{R}_d\mathbf{B}^H\mathbf{H} + \mathbf{R}_n)\mathbf{W} - \mathbf{W}^H\mathbf{H}^H\mathbf{B}\mathbf{R}_d - \mathbf{R}_d^H\mathbf{B}^H\mathbf{H}\mathbf{W}\} \\ &\quad \text{s.t } \mathrm{tr}\{\mathbf{B}\mathbf{R}_d\mathbf{B}^H\} \leq P_{max}\end{aligned} \tag{1.35}$$

where $P_{max}$ is the total available power at the transmitter.

One may think of solving this matrix valued problem by employing the Lagrangian multiplier and gradient methods which is known from postgraduate numerical mathematics course. However, by employing modern convex optimization theory (see Appendix A for the basics on convex optimization), one can show that the global optimality of this problem can not be ensured by applying the Lagrangian multiplier and gradient methods [BV04]. This is simply because the convexity of this problem is not exploited[11].

[11] To the best of our knowledge, we are not aware of any work showing the convexity of the above problem.

This problem and several other data transmission rate, signal to interference plus noise ratio (SINR) and MSE-based design problems have been extensively considered in [PLC04, PCL03] (see also the references in [PLC04, PCL03]). And it is shown that most of the practically relevant transceiver design problems for MIMO systems are not convex. Hence convex optimization approaches can not be applied to solve them. Despite this challenge, [PLC04] comes up with a novel and unified Majorization theory to solve the aforementioned classes of transceiver design problems for MIMO systems.

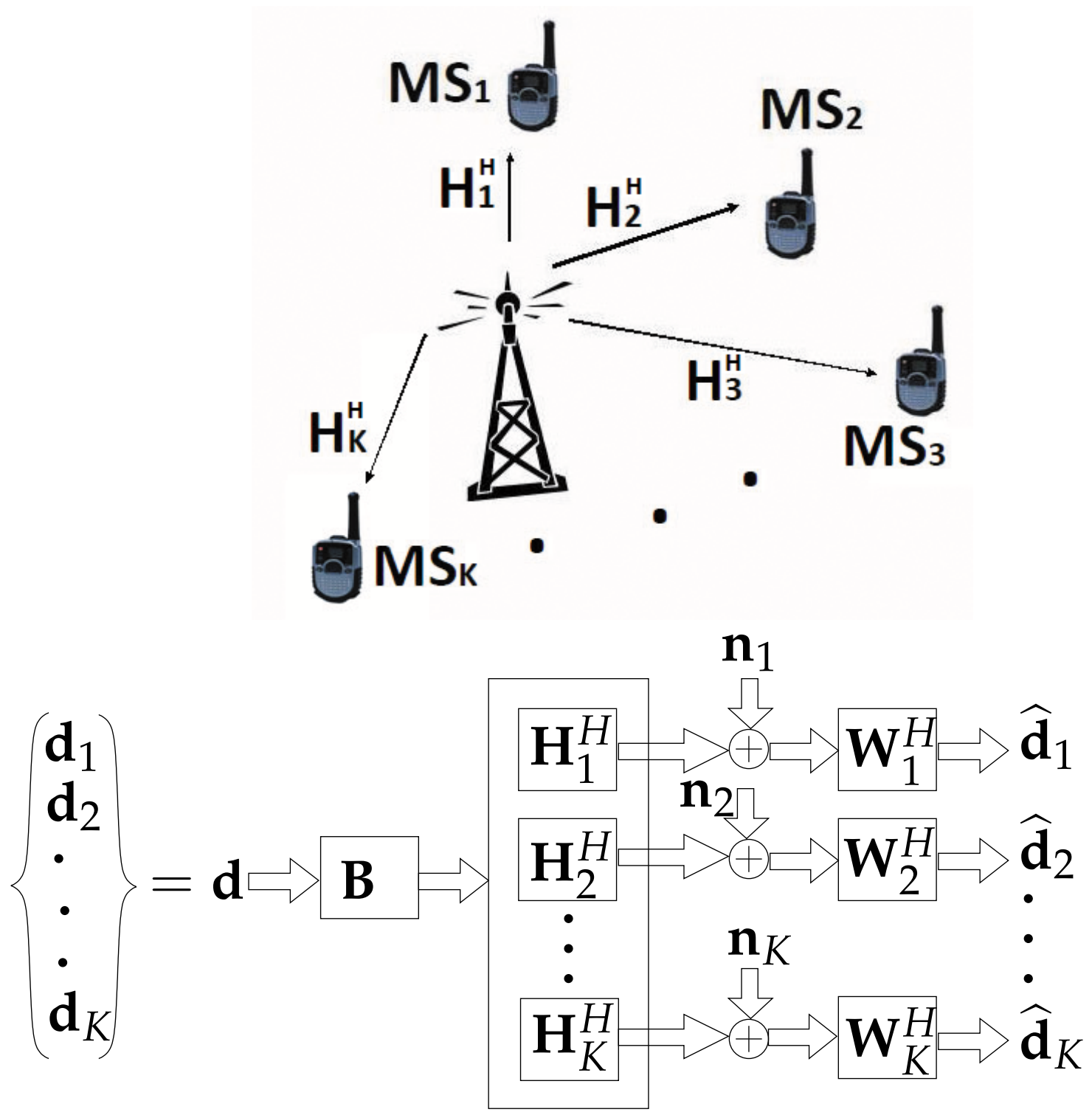

**Figure 1.6** [upper]: Conventional downlink communication, [lower]: The equivalent system model for this downlink communication.

## 1.8 Multiuser MIMO

In the previous section, a brief summary of MIMO system is presented for a point to point scenario (i.e., one transmitter and one receiver). Now let us extend the MIMO system model of Fig. 1.5 to the scenario where there is one transmitter (base station (BS)) serving $K$ decentralized receivers (mobile stations (MS)) which is shown in Fig. 1.6 (This system is termed as downlink multiuser MIMO system). For this system, the relationship between the transmitted signal of the $i$th user $\mathbf{d}_i$ and its estimated version $\widehat{\mathbf{d}}_i$ with linear precoders and decoders can be expressed as

$$\widehat{\mathbf{d}}_i = \mathbf{W}_i(\mathbf{H}_i^H \sum_{m=1}^{K} \mathbf{B}_m \mathbf{d}_m + \mathbf{n}_i) \tag{1.36}$$

where $\mathbf{H}_i^H$ is the channel between the BS and $i$th MS and $\mathbf{B}_i(\mathbf{W}_i)$ is the $i$th user precoder(decoder) matrix. The extension of problem (1.35) for the downlink multiuser MIMO system can be formulated as

$$\begin{aligned} \min_{\{\mathbf{B}_i,\mathbf{W}_i\}_{i=1}^K} & \sum_{i=1}^{K} \mathrm{tr}\{\mathbf{R}_{di} + \mathbf{W}_i^H(\mathbf{H}_i^H \sum_{m=1}^{K} \mathbf{B}_m \mathbf{R}_{dm} \mathbf{B}_m^H \mathbf{H}_i + \mathbf{R}_{ni})\mathbf{W}_i - \mathbf{W}_i^H \mathbf{H}_i^H \mathbf{B}_i \mathbf{R}_{di} - \\ & \mathbf{R}_{di}^H \mathbf{B}_i^H \mathbf{H}_i \mathbf{W}_i\} \\ & \text{s.t} \sum_{m=1}^{K} \mathrm{tr}\{\mathbf{B}_m \mathbf{R}_{dm} \mathbf{B}_m^H\} \leq P_{max} \end{aligned} \tag{1.37}$$

where $P_{max}$ is the maximum available power at the BS and $\mathbf{R}_{di}(\mathbf{R}_{ni})$ is the covariance matrix of $\mathbf{d}_i(\mathbf{n}_i)$. As we can see from this problem, the precoders of all users are jointly coupled. Due to this fact, getting the optimal precoder(decoder) pairs of this problem is not trivial. Also for this problem and other multiuser MIMO downlink, rate, SINR and MSE-based problems, the Majorization theory of [PLC04] can not be applied. For this reason, several researchers propose different linear precoder/decoder (i.e., transceiver) design algorithms for solving different objective functions for the downlink multiuser MIMO systems which is summarized in the next section[12].

[12] We would like to mention here that when there are more than one decentralized transmitters and a single receiver. The system is termed as a multiuser MIMO uplink. For this system, since linear precoder/decoder designs are well understood, the precoder/decoder design strategies and algorithms is not discussed for this system.

## 1.9 Existing linear transceiver design algorithms for downlink multiuser MIMO systems

This section discusses existing linear transceiver design algorithms for downlink multiuser MIMO systems. In [SSH04], channel block-diagonalization transceiver design approach is suggested where water-filling algorithm is used to solve the power allocation part of the optimization problem. This method suffers from noise enhancement and has a restriction on the number of transmit and receive antennas. In [SSB07, SSB08a, HJU09, SSJB05], several MSE-based problems are considered in the downlink channel. These papers solve the MSE-based problems by applying uplink-downlink duality (see Appendix B for a brief summary of uplink-downlink duality) solution approach [FLT98, VM99, SB04]. The authors of [SSB07, SSB08a, HJU09, SSJB05] also exploit the fact that solving the downlink MSE-based precoder/decoder design problems by the uplink-downlink duality approach has easier mathematical structure than the direct approach where the downlink MSE-based precoder/decoder design problems are examined directly in the downlink channel. Moreover, for some MSE-based problems, the uplink-downlink duality solution approach can also exploit the hidden convexity of the downlink channel MSE-based problems (see for example [SSB07]).

In [SSB08c] weighted sum rate maximization problem is formulated as the problem of minimizing the geometric product of the minimum mean square errors. This problem is solved by applying MSE uplink-downlink duality approach. Minimizing the product of all users minimum mean square error (MMSE) matrix determinants is proposed as an equivalent formulation for the un-weighted sum rate maximization problem [TA08]. This problem is solved by employing sequential quadratic programming.

## 1.10 Motivation of the thesis

For the downlink multiuser MIMO systems, the above precoder/decoder design problems are examined by assuming perfect CSI at the transmitter and receiver. However, in practice, obtaining perfect CSI at the transmitter and receiver is difficult and often only an imperfect channel is available. Consequently, designing the precoders and decoders without taking into account the CSI imperfections leads to performance degradation [UC08]. Thus, incorporating channel uncertainty in the design problem has significant advantage. This motivates us to examine linear transceiver (i.e., precoder/decoder) design problems for practically relevant capacity, rate, SINR and MSE-based objective functions for the downlink multiuser MIMO systems by introducing the CSI imperfections.

In addition, all of the aforementioned papers examine their problems for conventional downlink networks. In these networks, BSs from different cells communicate with their respective MSs independently. Hence, in the latter network, each BS is obliged to treat its inter-cell interference as a background noise. To show this fact, let us consider a system with two BSs where each BS is serving 2 MSs as shown in Fig. 1.7.

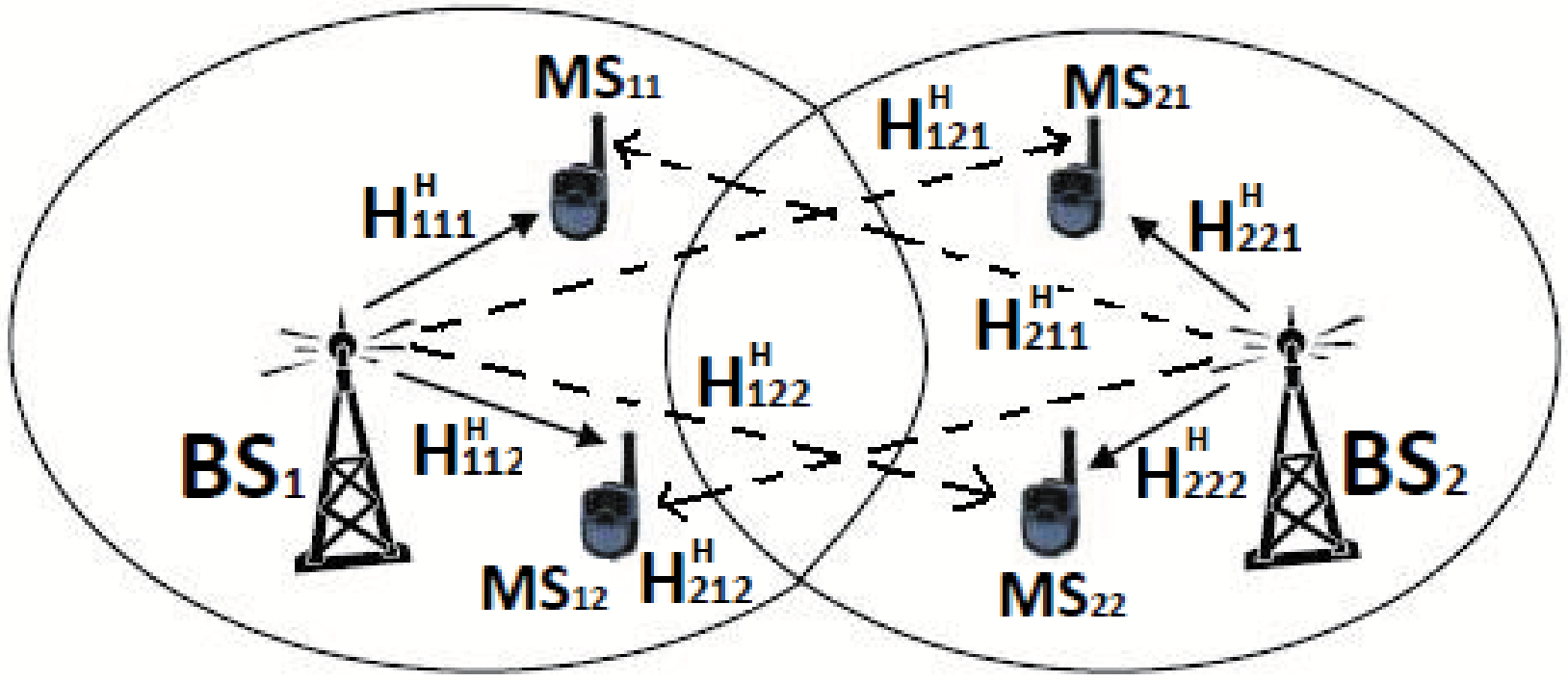

**Figure 1.7** A downlink communication system with 2 uncoordinated base stations. Dashed line denotes the channel due to inter-cell interference.

By employing the same precoder, decoder and channel matrix notations as in Fig. 1.6, the estimated signal at $\text{MS}_{11}$ of Fig. 1.7 becomes

$$\begin{aligned}\hat{\mathbf{d}}_{11} =& \mathbf{W}_{11}^H(\mathbf{H}_{111}^H(\mathbf{B}_{11}\mathbf{d}_{11} + \mathbf{B}_{12}\mathbf{d}_{12}) + \mathbf{H}_{211}^H(\mathbf{B}_{21}\mathbf{d}_{21} + \mathbf{B}_{22}\mathbf{d}_{22}) + \mathbf{n}_{11}) \\ =& \mathbf{W}_{11}^H(\mathbf{H}_{111}^H(\mathbf{B}_{11}\mathbf{d}_{11} + \mathbf{B}_{12}\mathbf{d}_{12}) + \tilde{\mathbf{n}}_{11} + \mathbf{n}_{11}) \end{aligned} \quad (1.38)$$

where $\mathbf{B}_{ij}$ and $\mathbf{W}_{ij}$ are the precoder and decoder matrices corresponding to $\text{MS}_{ij}$[13], and $\tilde{\mathbf{n}}_{11} = \mathbf{H}_{211}^H(\mathbf{B}_{21}\mathbf{d}_{21} + \mathbf{B}_{22}\mathbf{d}_{22})$. As we can see from this equation,

[13] i refers to the BS or cell number and j refers to the MS number.

since $BS_1$ and $MS_{11}$ do not have any knowledge about $\mathbf{d}_{21}$, $\mathbf{d}_{22}$, $\mathbf{H}_{211}$, $\mathbf{B}_{21}$ and $\mathbf{B}_{22}$, $BS_1$ and $MS_{11}$ are obliged to treat $\tilde{\mathbf{n}}_{11}$ as an additional background noise. Consequently, $BS_1$ designs its precoder ($\mathbf{B}_{11}$, $\mathbf{B}_{12}$) and $MS_{11}$ designs its decoder ($\mathbf{W}_{11}$) by considering $\tilde{\mathbf{n}}_{11}$ (inter-cell interference) as an additional background noise. Research results show that such design significantly reduces the performance of precoders and decoders [KFV06].

Recently, it has been shown that BS coordination is a promising technique to significantly improve the spectral efficiency of wireless channels by mitigating (or possibly canceling) inter-cell interference [KFV06, DY10, BZGO10]. As will be detailed later, for such coordinated BS system, the existing papers propose centralized algorithms to solve several classes of transceiver design problems. However, it is evidently seen that such centralized transceiver design algorithm is not feasible for large scale coordinated networks. This fact motivate us to develop distributed transceiver design algorithm for several practically relevant rate, SINR and MSE-based problems in the downlink multiuser MIMO coordinated BS systems.

## 1.11 Outline and history of the thesis

In Chapter 2, three kinds of MSE uplink-downlink duality are established for the multiuser MIMO systems by considering that the BS and MS antennas exhibit spatial correlations and the CSI at both the transmitter and receiver ends are imperfect. These duality are established by extending the three level MSE duality of [HJU09] to imperfect CSI. As application examples of the duality, the joint optimization of transceivers for different MSE-based robust design problems in the downlink channel have been considered. The robustness against imperfect CSI is incorporated into the designs using stochastic approach [DB08].

The work of Chapter 2 has been published in [BCV11] and [ECV09].

In Chapter 2, we show that the MSE uplink-downlink duality under imperfect (also perfect) CSI scenarios can be exploited just by transforming the diagonal power allocation matrices from uplink to downlink channel and vice versa. This simple transformation helps us to get less complexity uplink-downlink duality based iterative solution for several classes of MSE based problems. The main drawback of this chapter is that it solves only total BS power based precoder decoder design problems.

After we complete the works of [BCV11] and [ECV09] (i.e., Chapter 2), we plan to generalize the duality to handle arbitrary power constrained MSE based problems. However, as such generalization is not trivial, we were unable to move forward in this matter. For this reason, we switched to a distributed transceiver design algorithm for coordinated BS system which is motivated as follows:

Chapter 2 examines the problems for conventional downlink networks. Research studies exploit the fact that BS coordination is a promising approach to improve the spectral efficiency of wireless channels [KFV06, DY10, BZGO10]. The BS coordination can be performed by two approaches. In the first approach, BSs are coordinated at the beamforming (precoder) level [DY10] (multi-cell or partially coordinated systems), whereas in the second approach, coordination takes place both at the signal and beamforming (precoder) levels [KFV06], [BZGO10] (network MIMO or fully coordinated systems). It is well know that the latter coordination approach has better performance compared to that of the former one [BZGO10], [BZGO09]. This performance improvement, however, requires additional signal coordination. In [SSVB08], four MSE-based linear transceiver optimization problems have been considered for multiuser MIMO systems with fully coordinated BSs. In [SSVB08], the receiver of each user are optimized independently and distributively. However, the joint optimization of the precoders has been carried out by a centralized controller. When the number of MSs and/or BSs increase, the computational cost of the joint precoder design also increases [TSC07]. Consequently, solving the precoder optimization problem in a centralized manner, especially for large-scale coordinated networks, is not a computationally efficient approach. This limitation has been addressed in Chapter 3 of this thesis.

The work of Chapter 3 has been published in [BV11d], [BVC12] and [BVC10].

From Chapter 2, we realize that the MSE duality solution approach of solving downlink transceiver design problems require less computational cost than that of the direct solution approach. Moreover, the transmit and receive filters of the MSE duality based algorithm can also be implemented distribu-tively which naturally leads to a distributive algorithm (recall **Table 2.1** of Chapter 2).

From the work of Chapter 3, we learn that a fully coordinated BS system can be interpreted as a giant MIMO system. In consequence, a downlink multiuser MIMO fully coordinated BS system and conventional downlink multiuser MIMO system will have similar mathematical structure. In this chapter we also show a clear relationship between weighted sum rate and weighted sum MSE-based optimization problems. On the other hand, according to [YL07], in a practical multi-antenna BS system (either downlink multiuser MIMO coordinated BS system or conventional downlink multiuser MIMO system), the maximum power of each BS antenna is limited. Moreover, MSs are spaced far apart from each other and the noise vector of each MS may include other interference signals [Pal03].

These practically relevant design requirements and the aforementioned advantages of MSE duality solution approach motivate us to put more effort for exploiting the MSE duality for generalized power constraints and noise covariance matrices (i.e., a generalized version of the duality of Chapter 2). And after a lot of effort, we come up with novel downlink-interference MSE duality and the benefits of the new downlink-interference duality have been exploited by showing their importance to solve several practically relevant rate, SINR and MSE-based problems for both conventional downlink multiuser MIMO system and downlink multiuser MIMO coordinated BS systems. Chapter 4 of this thesis discusses the new MSE downlink-interference duality and their applications.

The work of Chapter 4 has been published in [BV13], [BV11b], [BV12] and [BV11c].

Finally in Chapter 5, conclusions and future works has been discussed briefly.

# MSE Uplink-Downlink Duality under Imperfect CSI

# 2

In this chapter, three kinds of MSE uplink-downlink duality are established by considering that the BS and MS antennas exhibit spatial correlations and the CSI at both the transmitter and receiver ends are imperfect. These duality are established by extending the three level MSE duality of [HJU09] to imperfect CSI. As application examples of the duality, we solve two robust design problems: minimization of the weighted sum MSE and minimization of the maximum WMSE. The robustness against imperfect CSI is incorporated into the designs using stochastic approach [DB08]. These problems are solved as follows: first, we establish three kinds of MSE uplink-downlink duality by transforming only the power allocation matrices from uplink channel to downlink channel and vice versa. Second, in the uplink channel, we formulate the power allocation part of each problem ensuring global optimality. Finally, based on the solution of the uplink power allocation and the MSE duality results, for each problem, we propose an iterative algorithm that performs optimization alternatively between the uplink and downlink channels. Computer simulations verify the robustness of the proposed design compared to the non-robust/naive design.

## 2.1 Introduction

In a multiuser network the uplink-downlink duality approach for solving the downlink optimization problems has received a lot of attention. In [SSB07] and [HJU09], MSE based uplink-downlink duality have been exploited. These two papers utilize their duality results to solve MSE-based design problems. These duality are established by assuming that perfect CSI is available at the BS and MSs. However, due to the inevitability of channel estimation error, CSI can never be perfect. This motivates [SD08] to establish the MSE duality under imperfect CSI for MISO systems. The latter work is extended in [UC09] for MIMO case. None of [SD08] and [UC09] incorporates antenna correlation in their channel model and neither of these duality can be applied to symbol wise MSE-based problems for MIMO systems. For instance, the duality of [SD08] and [UC09] can not be used for the robust symbol wise weighted sum MSE problem. Moreover, while solving the robust sum MSE minimization problem, the authors of [SD08] and [UC09] compute $K$ (total number of MSs) scaling factors (see (16) in [SD08] and [UC09]) to transfer the total sum AMSE from uplink to downlink channel which is not computationally efficient. As will be seen later in Section 2.4, we compute only one scaling factor to transfer the sum AMSE from uplink to downlink channel and vice versa. In [DB08], the MSE uplink-downlink duality has been established by considering imperfect CSI both at the BS and MSs, and with antenna correlation only at the BS. The duality is examined by analyzing the KKT conditions for the uplink and downlink channel problems. The latter duality is limited to sum MSE minimization problem.

In [ECV09], three kinds of MSE duality are established by considering that imperfect CSI is available both at the BS and MSs, and with antenna correlation only at the BS. These duality are established by extending the three level MSE duality of [HJU09] to imperfect CSI. Thus, from the MSE duality perspective, the duality of [ECV09] is more general than that of [SD08, UC09, DB08]. In order to solve general MSE-based robust design problems (see for example (14) in **Case 2** of [ECV09]), the approach of [SSB07] and [SSB08c] has been employed where the precoder of each MS is decomposed into a product of unity norm filter and diagonal power allocation matrices, and the decoder of each MS is decomposed into a product of unity norm filter, diagonal scaling factor and the inverse of power allocation matrices (see (15) of [ECV09]). Upon doing so, [ECV09] show that any MSE-based robust design problem can be solved using alternating optimization framework. From (22) of [ECV09], one can also realize that by employing the same filters and scaling factors in both the up-link and downlink channels, three kinds of AMSE uplink-downlink duality can be established just by transforming the power allocation matrices from uplink channel to downlink channel and vice versa. Due to this reason, this chapter employs the system model shown in Fig. 2.1. Note that although this system model is known from [SSB07] and [SSB08c], the authors of these two papers employ another approach to establish the MSE uplink-downlink duality which is computationally costly.

This chapter considers that the BS and MS antennas exhibit spatial correlations and the CSI at both ends is imperfect. The robustness against imperfect CSI is incorporated into our designs using stochastic approach [DB08]. In this regard, first the three kinds of AMSE duality has been established. Then, as application examples, the joint optimization of transceivers for the following MSE-based robust design problems have been examined.

1. The robust minimization of the weighted sum MSE constrained with a total BS power ($\mathcal{P}2.1$).
2. The robust minimization of the maximum weighted MSE (min-max) constrained with a total BS power ($\mathcal{P}2.2$).

As $\mathcal{P}2.1$ and $\mathcal{P}2.2$ are non-convex, convex optimization tools can not be applied to solve them. Due to this, this chapter solves these problems iteratively as follows: First, the power allocation part of each problem has been solved ensuring global optimality. Then, with this solution and the AMSE duality results, like in [ECV09], iterative algorithms are applied for $\mathcal{P}2.1$ and $\mathcal{P}2.2$. The key contributions of this chapter is summarized as follows:

1. By using the system model shown in Fig. 2.1, the three kinds of AMSE duality known from [HJU09][1] for the aforementioned CSI has been established just by transforming the power allocation matrices (which are diagonal) from uplink to downlink channel and vice versa. In contrast to the AMSE duality in [HJU09], [SD08], [UC09] and [DB08], our duality can be used to solve all MSE-based problems by using alternating optimization like in [ECV09]. It is worthwhile to mention that one can also extend the duality approach of [SSB07] to imperfect CSI case as the latter duality also requires only the transformation of powers from uplink to downlink channel and vice versa. However, by utilizing our duality, the computational complexity of the latter power transformation can be reduced (this will be clear later in Section 2.5.1). As a consequence, the overall computational cost of alternating optimization algorithm of [SSB07] reduces. Moreover, this work generalizes the hitherto MSE uplink-downlink duality[2].

[1]Note: The authors of [HJU09] establish the three kinds of duality by transferring the precoder/decoder pairs from uplink to downlink channel and vice versa.

2. It is shown that the uplink power allocation part of each problem can be solved ensuring global optimality.

3. Using the uplink power allocation and AMSE duality results, iterative uplink-downlink duality based algorithms have been proposed for solving $\mathcal{P}2.1$ and $\mathcal{P}2.2$.

4. The effects of channel estimation errors and antenna correlations on the system performance have been examined.

This chapter is organized as follows. In Section 2.2, multiuser MIMO downlink and virtual uplink channel system models are presented. In Section 2.3, brief description of imperfect channel model is given. Section 2.4 presents the proposed AMSE uplink-downlink duality. The applications of the proposed AMSE duality has been discussed in Section 2.5. In Section 2.6, computer simulations are used to compare the performance of the proposed duality algorithms with that of the existing algorithms. Conclusions are presented in Section 2.7.

## 2.2 System model

In this section the MIMO downlink and uplink system models are considered. The BS equipped with $N$ transmit antennas is serving $K$ decentralized

[2]Note that for the considered CSI model, the MSE uplink-downlink duality can be established using the system model like in [DB08] and [ECV09]. However, this system model is not convenient to solve general MSE-based robust design problems (for example $\mathcal{P}2.1$ (**Case 2**) and $\mathcal{P}2.1$).

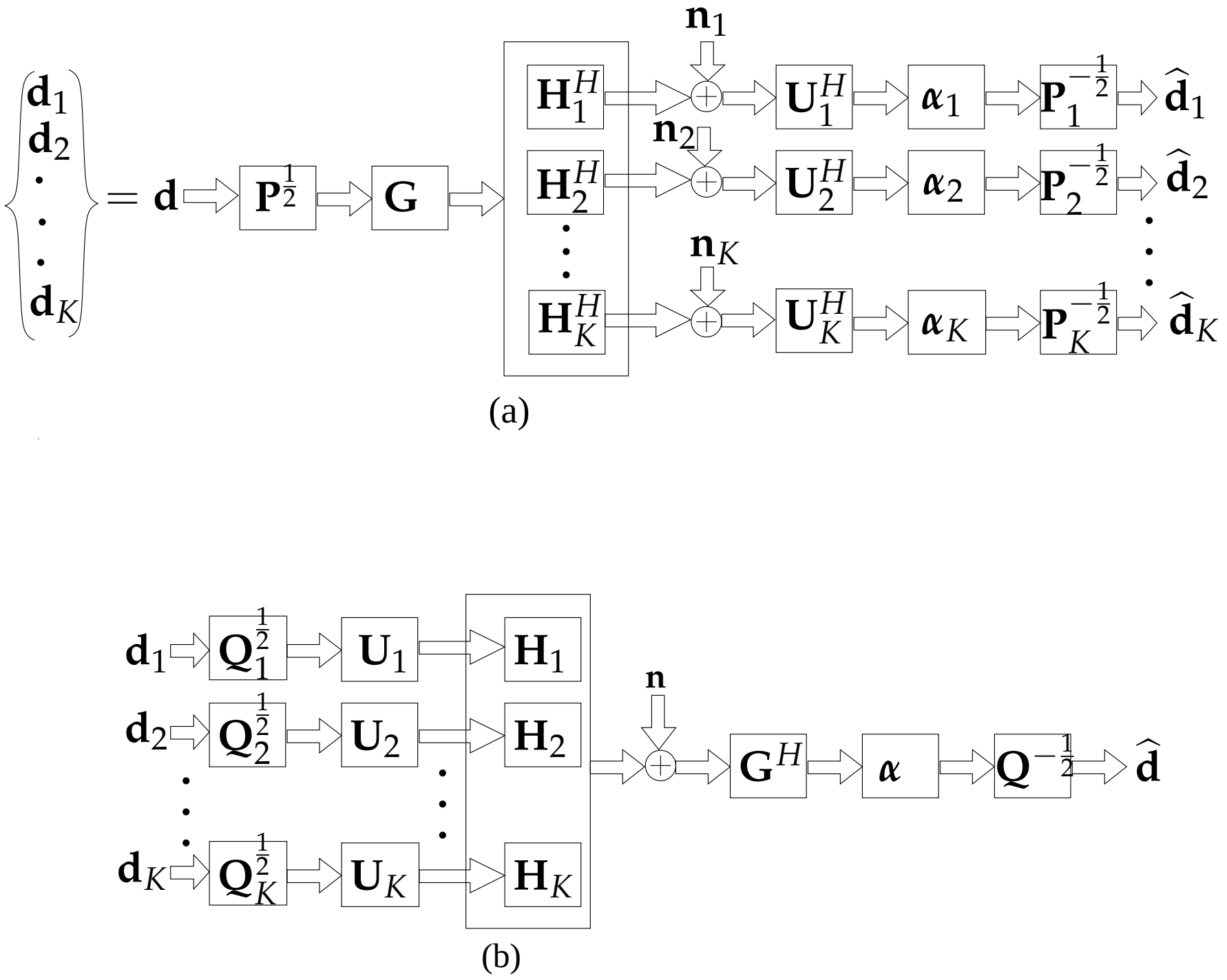

**Figure 2.1** Multiuser MIMO system. (a) downlink channel. (b) uplink channel.

MSs each having $\{M_k\}_{k=1}^K$ antennas to multiplex $S_k$ symbols. The total number of MS antennas and symbols are $M = \sum_{k=1}^K M_k$ and $S = \sum_{k=1}^K S_k$, respectively. All symbols can be stacked in a data vector $\mathbf{d} = [\mathbf{d}_1^T, \cdots, \mathbf{d}_K^T]^T$, where $\mathbf{d}_k \in C^{S_k \times 1}$ is the symbol vector for the $k$th MS. The MAC channel can be expressed as $\mathbf{H} = [\mathbf{H}_1, \cdots, \mathbf{H}_K]$, where $\mathbf{H}_k^H \in C^{M_k \times N}$ is the channel between the BS and $k$th MS.

Using the system model similar to [SSB08c] and as shown in Fig. 2.1, we collect the transmit powers of all users as $\mathbf{P} = \text{blkdiag}(\mathbf{P}_1, \mathbf{P}_2, \cdots, \mathbf{P}_K)$ and $\mathbf{Q} = \text{blkdiag}(\mathbf{Q}_1, \cdots, \mathbf{Q}_K)$, where $\mathbf{P}_k = \text{diag}(p_{k1}, \cdots, p_{kS_k})$, $\mathbf{Q}_k = \text{diag}(q_{k1}, \cdots, q_{kS_k})$ and $p_{ki}$ $(q_{ki})$ is the downlink (uplink) power allocation for the $i$th symbol of the $k$th user. The overall filter matrix at the BS is $\mathbf{G} = [\mathbf{G}_1, \cdots, \mathbf{G}_K]$, where $\mathbf{G}_k = [\mathbf{g}_{k1} \ \cdots \ \mathbf{g}_{kS_k}] \in C^{N \times S_k}$ is the filter matrix for the $k$th user with $\{\mathbf{g}_{ki}^H \mathbf{g}_{ki} = 1\}_{i=1}^{S_k}, k = \{1, \cdots, K\}$. The filters of all users are stacked in a block diagonal matrix $\mathbf{U} = \text{blkdiag}(\mathbf{U}_1, \cdots, \mathbf{U}_K)$, where $\mathbf{U}_k = [\mathbf{u}_{k1} \ \cdots \ \mathbf{u}_{kS_k}] \in C^{M_k \times S_k}$ is the filter matrix for the $k$th user with $\{\mathbf{u}_{ki}^H \mathbf{u}_{ki} = 1\}_{i=1}^{S_k}, \forall k$. The scaling factors are accumulated as $\boldsymbol{\alpha} = \text{blkdiag}(\boldsymbol{\alpha}_1, \cdots, \boldsymbol{\alpha}_K)$, where $\boldsymbol{\alpha}_k = \text{diag}(\alpha_{k1}, \cdots, \alpha_{kS_k})$. The entries of $\mathbf{n} = [\mathbf{n}_1^T, \mathbf{n}_2^T, \cdots, \mathbf{n}_K^T]^T$ are assumed to be i.i.d ZMCSCG random variables all with variance $\sigma^2$. We also assume that $\text{E}\{\mathbf{d}_k \mathbf{d}_k^H\} = \mathbf{I}_{S_k}$, $\text{E}\{\mathbf{d}_k \mathbf{d}_i^H\} = \mathbf{0}, \forall i \neq k$, and $\text{E}\{\mathbf{d}_k \mathbf{n}_k^H\} = \mathbf{0}$.

## 2.3 Channel model

Considering antenna correlation at the BS and MSs, we model the Rayleigh fading MIMO channels as $\mathbf{H}_k^H = \widetilde{\mathbf{R}}_{mk}^{1/2} \mathbf{H}_{wk}^H \mathbf{R}_{bk}^{1/2}, \forall k$, where the elements of $\{\mathbf{H}_{wk}^H\}_{k=1}^K$ are i.i.d ZMCSCG random variables all with unit variance and $\mathbf{R}_{bk} \in C^{N \times N}$, $\widetilde{\mathbf{R}}_{mk} \in C^{M_k \times M_k}$ are antenna correlation matrices at the BS and MSs, respectively [Din08], [YYG04]. The channel estimation is performed on $\{\mathbf{H}_{wk}^H\}_{k=1}^K$ using an orthogonal training method [DB09].

In the training stage, the BS transmit the following training sequence $\widetilde{\mathbf{T}}_k \in \mathcal{C}^{N \times N} = \mathbf{R}_{bk}^{-\frac{1}{2}} \widetilde{\mathbf{T}}_0$ only for the $k$th MS (i.e., all the other MSs are silent), where $\widetilde{\mathbf{T}}_0$ is an $N \times N$ unitary matrix scaled by $\sqrt{\text{tr}\{\widetilde{\mathbf{T}}_k \widetilde{\mathbf{T}}_k^H\} / \text{tr}\{\mathbf{R}_{bk}^{-1}\}}$ . Upon doing so, the received signal of the $k$th user can be expressed as

$$\mathbf{Y}_k = \mathbf{H}_k^H \widetilde{\mathbf{T}}_k + \mathbf{N}_{kt} = \widetilde{\mathbf{R}}_{mk}^{1/2} \mathbf{H}_{wk}^H \widetilde{\mathbf{T}}_0 + \mathbf{N}_{kt}$$

where $\mathbf{N}_{kt}$ is the noise matrix of the $k$th MS during the training phase. Each entry of $\mathbf{N}_{kt}$ is i.i.d with $\mathcal{NC}(0, \sigma^2_{nkt})$. After pre and post-multiplying $\mathbf{Y}_k$ by $\widetilde{\mathbf{R}}_{mk}^{-1/2}$ and $\widetilde{\mathbf{T}}_0^H$, respectively, we get

$$\widetilde{\mathbf{Y}}_k = \widetilde{\mathbf{R}}_{mk}^{-1/2}\mathbf{Y}_k\widetilde{\mathbf{T}}_0^H = \mathbf{H}_{wk}^H + \widetilde{\mathbf{R}}_{mk}^{-1/2}\mathbf{N}_{kt}\widetilde{\mathbf{T}}_0^H = \mathbf{H}_{wk}^H + \widetilde{\mathbf{R}}_{mk}^{-1/2}\widetilde{\mathbf{N}}_{kt}$$

where $\widetilde{\mathbf{N}}_{kt} = \mathbf{N}_{kt}\widetilde{\mathbf{T}}_0^H$. It is easily seen that each element of $\widetilde{\mathbf{N}}_{kt}$ is i.i.d with $\widetilde{\mathbf{N}}_{kt} \sim \mathcal{CN}(0, \sigma^2_{ek})$, where $\sigma^2_{ek} = \sigma^2_{nkt}\text{tr}\{\mathbf{R}_{bk}^{-1}\}/\text{tr}\{\widetilde{\mathbf{T}}_k\widetilde{\mathbf{T}}_k^H\}$.

Now let us estimate $\mathbf{H}_{wk}^H$ using MMSE channel estimation method. By incorporating a linear estimator $\widetilde{\mathbf{W}}_k$, the estimate of $\mathbf{H}_{wk}^H$ can be expressed as

$$\widehat{\mathbf{H}}_{wk}^H = \widetilde{\mathbf{W}}_k^H\widetilde{\mathbf{Y}}_k = \mathbf{W}_k^H(\mathbf{H}_{wk}^H + \widetilde{\mathbf{R}}_{mk}^{-1/2}\widetilde{\mathbf{N}}_{kt}).$$

The error between $\widehat{\mathbf{H}}_{wk}$ and $\mathbf{H}_{wk}$ is expressed as

$$\widetilde{\mathbf{E}}_{wk}^H = \mathbf{H}_{wk}^H - \widehat{\mathbf{H}}_{wk}^H$$

and the MSE matrix becomes

$$\begin{aligned}
\mathbf{MSE}_{chk} =& \text{E}\{\widetilde{\mathbf{E}}_{wk}^H\widetilde{\mathbf{E}}_{wk}\} \\
=& \text{E}\{(\mathbf{H}_{wk}^H - \widehat{\mathbf{H}}_{wk}^H)(\mathbf{H}_{wk}^H - \widehat{\mathbf{H}}_{wk}^H)^H\} \\
=& \text{E}\{(\mathbf{I} - \widetilde{\mathbf{W}}_k^H)\mathbf{H}_{wk}^H\mathbf{H}_{wk}(\mathbf{I} - \widetilde{\mathbf{W}}_k^H)^H + \widetilde{\mathbf{W}}_k^H\widetilde{\mathbf{R}}_{mk}^{-1/2}\widetilde{\mathbf{N}}_{kt}\widetilde{\mathbf{N}}_{kt}^H\widetilde{\mathbf{R}}_{mk}^{-1/2}\widetilde{\mathbf{W}}_k\} \\
=& N((\mathbf{I} - \widetilde{\mathbf{W}}_k^H)(\mathbf{I} - \widetilde{\mathbf{W}}_k^H)^H + \sigma^2_{ek}\widetilde{\mathbf{W}}_k^H\widetilde{\mathbf{R}}_{mk}^{-1}\widetilde{\mathbf{W}}_k).
\end{aligned}$$

In the last equality we have used the fact that $\text{E}\{\mathbf{XAX}^H\} = \sigma_x^2\text{tr}\{\mathbf{A}\}\mathbf{I}$, where each entry of $\mathbf{X}$ is i.i.d with $\mathcal{CN}(0, \sigma_x^2)$. The linear estimator ($\widetilde{\mathbf{W}}_k$) that minimizes $\mathbf{MSE}_{chk}$ is obtained by

$$\frac{\partial \mathbf{MSE}_{chk}}{\partial \widetilde{\mathbf{W}}_k^H} = (\mathbf{I} + \sigma^2_{ek}\widetilde{\mathbf{R}}_{mk}^{-1})\mathbf{W}_k - \mathbf{I} = \mathbf{0} \Rightarrow \widetilde{\mathbf{W}}_k = (\mathbf{I} + \sigma^2_{ek}\widetilde{\mathbf{R}}_{mk}^{-1})^{-1}.$$

It follows

$$\begin{aligned}
\widehat{\mathbf{H}}_{wk}^H =& \widetilde{\mathbf{W}}_k^H\widetilde{\mathbf{Y}}_k \\
=& (\mathbf{I} + \sigma^2_{ek}\widetilde{\mathbf{R}}_{mk}^{-1})^{-1}\widetilde{\mathbf{Y}}_k \\
\mathbf{MSE}_{chk} =& \text{E}\{\widetilde{\mathbf{E}}_{wk}^H\widetilde{\mathbf{E}}_{wk}\} \\
=& N((\mathbf{I} - \widetilde{\mathbf{W}}_k^H)(\mathbf{I} - \widetilde{\mathbf{W}}_k^H)^H + \sigma^2_{ek}\mathbf{W}_k^H\widetilde{\mathbf{R}}_{mk}^{-1}\widetilde{\mathbf{W}}_k) \\
=& N(\mathbf{I} - (\mathbf{I} + \sigma^2_{ek}\widetilde{\mathbf{R}}_{mk}^{-1})^{-1}).
\end{aligned}$$

And $\mathbf{H}_{wk}^H$ can be expressed as

$$\mathbf{H}_{wk}^H = \widehat{\mathbf{H}}_{wk}^H + \widetilde{\mathbf{E}}_{wk}^H$$

where $\widetilde{\mathbf{E}}_{wk}^H \sim \mathcal{CN}(\mathbf{0}, \mathbf{MSE}_{chk} = N(\mathbf{I} - (\mathbf{I} + \sigma_{ek}^2 \widetilde{\mathbf{R}}_{mk}^{-1})^{-1}))$.

Thus, the true channel $\mathbf{H}_k^H$ can be expressed as

$$\mathbf{H}_k^H = \widetilde{\mathbf{R}}_{mk}^{1/2} \mathbf{H}_{wk}^H \mathbf{R}_{bk}^{1/2} = \widetilde{\mathbf{R}}_{mk}^{1/2} \widehat{\mathbf{H}}_{wk}^H \mathbf{R}_{bk}^{1/2} + \widetilde{\mathbf{R}}_{mk}^{1/2} \widetilde{\mathbf{E}}_{wk}^H \mathbf{R}_{bk}^{1/2}.$$

Now by applying matrix inversion lemma, we can equivalently express $\widetilde{\mathbf{R}}_{mk}^{1/2} \widetilde{\mathbf{E}}_{wk}^H$ as

$$\widetilde{\mathbf{R}}_{mk}^{1/2} \widetilde{\mathbf{E}}_{wk}^H = \mathbf{R}_{mk}^{1/2} \mathbf{E}_{wk}^H$$

where $\mathbf{R}_{mk} = (\mathbf{I}_{M_k} + \sigma_{ek}^2 \widetilde{\mathbf{R}}_{mk}^{-1})^{-1}$, $\mathbf{E}_k^H$ is the estimation error and the entries of $\mathbf{E}_{wk}^H$ are i.i.d with $\mathcal{CN}(0, \sigma_{ek}^2)$.

Therefore, the true and estimated channel of the $k$th MS can be related by

$$\mathbf{H}_k^H = \widetilde{\mathbf{R}}_{mk}^{1/2} \widehat{\mathbf{H}}_{wk}^H \mathbf{R}_{bk}^{1/2} + \mathbf{R}_{mk}^{1/2} \mathbf{E}_{wk}^H \mathbf{R}_{bk}^{1/2} = \widehat{\mathbf{H}}_k^H + \mathbf{E}_k^H, \ \forall k. \tag{2.1}$$

We would like to mention here that the above channel estimation approach is an extension of [DB09] to the multiuser scenario (see section II.B of [DB09] for more details). However, the analysis of this chapter can still be applied for other stochastic channel estimation error models.

In this chapter, we consider that $\{\mathbf{E}_{wk}^H\}_{k=1}^K$ are unknown but $\{\widehat{\mathbf{H}}_k^H$, $\mathbf{R}_{bk}$, $\widetilde{\mathbf{R}}_{mk}$ and $\sigma_{ek}^2\}_{k=1}^K$ are available. We assume that each MS estimates its channel and feeds the estimated channel back to the BS without any error. Thus, both the BS and MSs have the same channel imperfections. The $k$th user estimated signal $\widehat{\mathbf{d}}_k^{DL}$ can be expressed as

$$\begin{aligned}\widehat{\mathbf{d}}_k^{DL} =& \mathbf{P}_k^{-1/2} \boldsymbol{\alpha}_k \mathbf{U}_k^H (\mathbf{H}_k^H \mathbf{G} \mathbf{P}^{1/2} \mathbf{d} + \mathbf{n}_k) \\ =& \mathbf{P}_k^{-1/2} \boldsymbol{\alpha}_k \mathbf{U}_k^H (\mathbf{H}_k^H \sum_{k=1}^K \mathbf{G}_k \mathbf{P}_k^{1/2} \mathbf{d}_k + \mathbf{n}_k)\end{aligned} \tag{2.2}$$

where $\mathbf{H}_k^H$ is the channel between the BS and the $k$th user, and $\mathbf{n}_k$ is the additive noise at the $k$th MS. The downlink instantaneous MSE matrix of the $k$th

user is given by

$$\boldsymbol{\xi}_k^{DL} = \mathrm{E}_{\mathbf{d}}\{(\mathbf{d}_k - \widehat{\mathbf{d}}_k^{DL})(\mathbf{d}_k - \widehat{\mathbf{d}}_k^{DL})^H\}$$

$$\boldsymbol{\xi}_k^{DL} = \mathbf{I}_{S_k} + \mathbf{P}_k^{-1/2}\boldsymbol{\alpha}_k\mathbf{U}_k^H(\mathbf{H}_k^H(\sum_{i=1}^{K}\mathbf{G}_i\mathbf{P}_i\mathbf{G}_i^H)\mathbf{H}_k + \sigma^2\mathbf{I}_{M_k})\mathbf{U}_k\boldsymbol{\alpha}_k\mathbf{P}_k^{-1/2} - \mathbf{P}_k^{1/2}\mathbf{G}_k^H\mathbf{H}_k\mathbf{U}_k\boldsymbol{\alpha}_k\mathbf{P}_k^{-1/2} - \mathbf{P}_k^{-1/2}\boldsymbol{\alpha}_k\mathbf{U}_k^H\mathbf{H}_k^H\mathbf{G}_k\mathbf{P}_k^{1/2}.$$

Substituting (2.1) in $\boldsymbol{\xi}_k^{DL}$ and taking the expected value of $\boldsymbol{\xi}_k^{DL}$ over $\mathbf{E}_{wk}^H$, the downlink AMSEs can be expressed as

$$\begin{aligned}
\overline{\boldsymbol{\xi}}_k^{DL} &= \mathrm{E}_{\mathbf{E}_{wk}^H}\{\boldsymbol{\xi}_k^{DL}\} \\
&= \mathbf{I} + \mathbf{P}_k^{-1/2}\boldsymbol{\alpha}_k\mathbf{U}_k^H\boldsymbol{\Gamma}_k^{DL}\mathbf{U}_k\boldsymbol{\alpha}_k\mathbf{P}_k^{-1/2} \\
&\quad - \mathbf{P}_k^{1/2}\mathbf{G}_k^H\widehat{\mathbf{H}}_k\mathbf{U}_k\boldsymbol{\alpha}_k\mathbf{P}_k^{-1/2} - \mathbf{P}_k^{-1/2}\boldsymbol{\alpha}_k\mathbf{U}_k^H\widehat{\mathbf{H}}_k^H\mathbf{G}_k\mathbf{P}_k^{1/2} \\
\overline{\xi}_k^{DL} &= \mathrm{tr}\{\overline{\boldsymbol{\xi}}_k^{DL}\} \\
&= S_k + \mathrm{tr}\{\mathbf{P}_k^{-1}\boldsymbol{\alpha}_k^2\mathbf{U}_k^H\boldsymbol{\Gamma}_k^{DL}\mathbf{U}_k - 2\Re\{\mathbf{G}_k^H\widehat{\mathbf{H}}_k\mathbf{U}_k\boldsymbol{\alpha}_k\}\} \qquad (2.3) \\
\overline{\xi}^{DL} &= \sum_{k=1}^{K}\overline{\xi}_k^{DL} \\
&= S + \sum_{k=1}^{K}\mathrm{tr}\{\mathbf{P}_k^{-1}\boldsymbol{\alpha}_k^2\mathbf{U}_k^H\boldsymbol{\Gamma}_k^{DL}\mathbf{U}_k - 2\Re\{\mathbf{G}_k^H\widehat{\mathbf{H}}_k\mathbf{U}_k\boldsymbol{\alpha}_k\}\} \qquad (2.4)
\end{aligned}$$

where $\boldsymbol{\Gamma}_k^{DL} = (\widehat{\mathbf{H}}_k^H\mathbf{G}\mathbf{P}\mathbf{G}^H\widehat{\mathbf{H}}_k + \sigma_{ek}^2\mathrm{tr}\{\mathbf{R}_{bk}\mathbf{G}\mathbf{P}\mathbf{G}^H\}\mathbf{R}_{mk} + \sigma^2\mathbf{I}_{M_k})$ and we use the fact $\mathrm{E}_{\mathbf{E}}\{\mathbf{E}\mathbf{A}\mathbf{E}^H\} = \sigma_e^2\mathrm{tr}\{\mathbf{A}\}\mathbf{I}$, if the entries of $\mathbf{E}$ are i.i.d with $\mathcal{CN}(0,\sigma_e^2)$ and $\mathbf{A}$ is a given matrix. Like in the downlink channel, by defining $\boldsymbol{\Gamma}_c \triangleq [\sum_{i=1}^{K}(\widehat{\mathbf{H}}_i\mathbf{U}_i\mathbf{Q}_i\mathbf{U}_i^H\widehat{\mathbf{H}}_i^H + \sigma_{ei}^2\mathrm{tr}\{\mathbf{R}_{mi}\mathbf{U}_i\mathbf{Q}_i\mathbf{U}_i^H\}\mathbf{R}_{bi}) + \sigma^2\mathbf{I}_N]$, the uplink channel AMSEs are given by

$$\begin{aligned}
\overline{\boldsymbol{\xi}}_k^{UL} &= \mathbf{I}_{S_k} + \mathbf{Q}_k^{-1/2}\boldsymbol{\alpha}_k\mathbf{G}_k^H\boldsymbol{\Gamma}_c\mathbf{G}_k\boldsymbol{\alpha}_k\mathbf{Q}_k^{-1/2} - \mathbf{Q}_k^{-1/2}\boldsymbol{\alpha}_k\mathbf{G}_k^H\widehat{\mathbf{H}}_k\mathbf{U}_k\mathbf{Q}_k^{1/2} - \\
&\quad \mathbf{Q}_k^{1/2}\mathbf{U}_k^H\widehat{\mathbf{H}}_k^H\mathbf{G}_k\boldsymbol{\alpha}_k\mathbf{Q}_k^{-1/2} \qquad (2.5) \\
\overline{\xi}_k^{UL} &= \mathrm{tr}\{\overline{\boldsymbol{\xi}}_k^{UL}\} \\
&= S_k + \mathrm{tr}\{\mathbf{Q}_k^{-1}\boldsymbol{\alpha}_k^2\mathbf{G}_k^H\boldsymbol{\Gamma}_c\mathbf{G}_k - 2\Re\{\boldsymbol{\alpha}_k\mathbf{G}_k^H\widehat{\mathbf{H}}_k\mathbf{U}_k\}\} \qquad (2.6) \\
\overline{\xi}^{UL} &= \sum_{k=1}^{K}\overline{\xi}_k^{UL} \\
&= S + \sum_{k=1}^{K}\mathrm{tr}\{\mathbf{Q}_k^{-1}\boldsymbol{\alpha}_k^2\mathbf{G}_k^H\boldsymbol{\Gamma}_c\mathbf{G}_k - 2\Re\{\boldsymbol{\alpha}_k\mathbf{G}_k^H\widehat{\mathbf{H}}_k\mathbf{U}_k\}\}. \qquad (2.7)
\end{aligned}$$

## 2.4 Average mean square error uplink-downlink duality

As we mentioned in Section 2.1, our AMSE duality generalizes the work of [ECV09] to the case where the BS and MS antennas are spatially correlated, and both the BS and MSs have imperfect CSI. Thus, in this section, we transfer the sum AMSE, user wise AMSE and symbol wise AMSEs from the uplink to downlink channel and vice versa.

### 2.4.1 AMSE transfer from uplink to downlink channel

#### 2.4.1.1 Total sum AMSE transfer

For a given uplink sum AMSE (with a transmit power $\mathbf{Q}$), we can achieve the same sum AMSE in the downlink channel by using a positive $\beta$ which satisfies $\mathbf{P} = \beta\boldsymbol{\alpha}^2\mathbf{Q}^{-1}$.

Substituting $\mathbf{P}$ in (2.4), equating $\bar{\xi}^{DL} = \bar{\xi}^{UL}$ and after some simple derivations, $\beta$ can be determined as

$$\beta = \text{tr}\{\mathbf{Q}\}/\text{tr}\{\mathbf{Q}^{-1}\boldsymbol{\alpha}^2\}. \tag{2.8}$$

As can be seen from (2.8), the scaling factor $\beta$ does not depend on $\{\sigma_{ek}^2\}_{k=1}^K$. This can be seen from (2.4) and (2.7), after substituting $\{\boldsymbol{\Gamma}_k^{DL}\}_{k=1}^K$ and $\boldsymbol{\Gamma}_c$, where $\{\sigma_{ek}^2\}_{k=1}^K$ are amplified by the same factor. The downlink power is given by $P_{sum}^{DL} = \text{tr}\{\mathbf{P}\} = \text{tr}\{\beta\boldsymbol{\alpha}^2\mathbf{Q}^{-1}\} = \text{tr}\{\mathbf{Q}\} = P_{sum}^{UL}$. Thus, the same sum power is allocated in both channels.

#### 2.4.1.2 User wise AMSE transfer

Given the $k$th user AMSE in the uplink channel with $\{\mathbf{Q}_k\}_{k=1}^K \neq \mathbf{0}$, this user can achieve the same AMSE in the downlink channel if $\mathbf{P}_k$ is computed by

$$\mathbf{P}_k = \beta_k\boldsymbol{\alpha}_k^2\mathbf{Q}_k^{-1}. \tag{2.9}$$

Substituting (2.9) in (2.3), then equating $\bar{\xi}_k^{UL} = \bar{\xi}_k^{DL}$ and after some mathematical manipulations we obtain

$$\sum_{i=1,i\neq k}^{K} \beta_k(\|\mathbf{Q}_k^{-1/2}\boldsymbol{\alpha}_k\mathbf{G}_k^H\widehat{\mathbf{H}}_i\mathbf{U}_i\mathbf{Q}_i^{1/2}\|_F^2 + \sigma_{ei}^2\mathrm{tr}\{\mathbf{R}_{mi}\mathbf{U}_i\mathbf{Q}_i\mathbf{U}_i^H\}\|\mathbf{R}_{bi}^{1/2}\mathbf{G}_k\boldsymbol{\alpha}_k\mathbf{Q}_k^{-1/2}\|_F^2)$$
$$+\beta_k\sigma^2\mathrm{tr}\{\mathbf{Q}_k^{-1}\boldsymbol{\alpha}_k^2\} = \sigma^2\mathrm{tr}\{\mathbf{Q}_k\}+ \quad (2.10)$$
$$\sum_{i=1,i\neq k}^{K} \beta_i(\|\boldsymbol{\alpha}_i\mathbf{Q}_i^{-1/2}\mathbf{G}_i^H\widehat{\mathbf{H}}_k\mathbf{U}_k\mathbf{Q}_k^{1/2}\|_F^2 + \sigma_{ek}^2\mathrm{tr}\{\mathbf{R}_{mk}\mathbf{U}_k\mathbf{Q}_k\mathbf{U}_k^H\}\|\mathbf{R}_{bk}^{1/2}\mathbf{G}_i\boldsymbol{\alpha}_i\mathbf{Q}_i^{-1/2}\|_F^2).$$

After applying (2.10) for all users, we can form the following system of linear equations

$$\mathbf{X}\cdot[\beta_1,\ldots,\beta_K]^T = \sigma^2\,[\mathrm{tr}\{\mathbf{Q}_1\},\ldots,\mathrm{tr}\{\mathbf{Q}_K\}]^T \quad (2.11)$$

$$\text{where } [\mathbf{X}]_{k,l} = \quad (2.12)$$
$$\begin{cases} \sigma^2\mathrm{tr}\{\mathbf{Q}_k^{-1}\boldsymbol{\alpha}_k^2\} + \sum_{i=1,i\neq k}^{K}(\|\mathbf{Q}_k^{-1/2}\boldsymbol{\alpha}_k\mathbf{G}_k^H\widehat{\mathbf{H}}_i\mathbf{U}_i\mathbf{Q}_i^{1/2}\|_F^2 + \\ \sigma_{ei}^2\|\mathbf{R}_{mi}^{1/2}\mathbf{U}_i\mathbf{Q}_i^{1/2}\|_F^2\|\mathbf{R}_{bi}^{1/2}\mathbf{G}_k\boldsymbol{\alpha}_k\mathbf{Q}_k^{-1/2}\|_F^2) & k=l \\ -(\|\boldsymbol{\alpha}_l\mathbf{Q}_l^{-1/2}\mathbf{G}_l^H\widehat{\mathbf{H}}_k\mathbf{U}_k\mathbf{Q}_k^{1/2}\|_F^2 + \\ \sigma_{ek}^2\|\mathbf{R}_{mk}^{1/2}\mathbf{U}_k\mathbf{Q}_k^{1/2}\|_F^2\|\mathbf{R}_{bk}^{1/2}\mathbf{G}_l\boldsymbol{\alpha}_l\mathbf{Q}_l^{-1/2}\|_F^2) & k\neq l. \end{cases}$$

It can be shown that if $\sigma^2 > 0$ then $\{\beta_k\}_{k=1}^K$ of (2.11) are strictly positive [HJU09], [ECV09]. Thus, the $k$th user AMSE can be transfered from uplink to downlink channel. Summing up the left-hand and right-hand sides of (2.11) and cancelling $\sigma^2$ in both sides yields $P_{sum}^{DL} = \sum_{i=1}^K \beta_i\mathrm{tr}\{\mathbf{Q}_i^{-1}\boldsymbol{\alpha}_i^2\} = \sum_{i=1}^K \mathrm{tr}\{\mathbf{P}_i\} = \sum_{i=1}^K \mathrm{tr}\{\mathbf{Q}_i\} = P_{sum}^{UL}$. Thus, the same sum power is allocated in both the uplink and downlink channels.

## 2.4.2 AMSE transfer from downlink to uplink channel

To complete the duality, in this section we examine the AMSE transfer from the downlink to uplink channel.

### 2.4.2.1 Total sum AMSE transfer

Similar to Subsection 2.4.1.1, the sum AMSE can be transferred from the downlink to uplink channel by using a nonzero scaling factor $\widetilde{\beta}$ which satisfies $\mathbf{Q} = \widetilde{\beta}\boldsymbol{\alpha}^2\mathbf{P}^{-1}$.

Substituting $\mathbf{Q}$ in (2.7) and then equating $\bar{\xi}^{UL} = \bar{\xi}^{DL}$, $\widetilde{\beta}$ is determined as

$$\widetilde{\beta} = \text{tr}\{\mathbf{P}\}/\text{tr}\{\mathbf{P}^{-1}\boldsymbol{\alpha}^2\}. \tag{2.13}$$

#### 2.4.2.2 User wise AMSE transfer

Given the $k$th user downlink AMSE with $\{\mathbf{P}_k\}_{k=1}^K \neq \mathbf{0}$, this user can achieve the same AMSE in the uplink channel if $\mathbf{Q}_k$ is computed by $\mathbf{Q}_k = \widetilde{\beta}_k\boldsymbol{\alpha}_k^2\mathbf{P}_k^{-1}$. Like in Subsection 2.4.1.2, by substituting $\mathbf{Q}_k$ in (2.6), equating $\bar{\xi}_k^{UL} = \bar{\xi}_k^{DL}$ and after some mathematical manipulations, the scaling factors $\{\widetilde{\beta}_k\}_{k=1}^K$ are determined by solving the following system of linear equations.

$$\mathbf{T}\cdot[\widetilde{\beta}_1,\ldots,\widetilde{\beta}_K]^T = \sigma^2\,[\text{tr}\{\mathbf{P}_1\},\ldots,\text{tr}\{\mathbf{P}_K\}]^T \tag{2.14}$$

where $[\mathbf{T}]_{k,l} =$ (2.15)

$$\begin{cases} \sigma^2\text{tr}\{\mathbf{P}_k^{-1}\boldsymbol{\alpha}_k^2\} + \sum_{i=1,i\neq k}^K(\|\mathbf{P}_k^{-1/2}\boldsymbol{\alpha}_k\mathbf{U}_k^H\widehat{\mathbf{H}}_k^H\mathbf{G}_i\mathbf{P}_i^{1/2}\|_F^2 + \\ \sigma_{ek}^2\|\mathbf{R}_{bk}^{1/2}\mathbf{G}_i\mathbf{P}_i^{1/2}\|_F^2\|\mathbf{R}_{mk}^{1/2}\mathbf{U}_k\mathbf{P}_k^{-1/2}\boldsymbol{\alpha}_k\|_F^2), & k = l \\ -(\|\mathbf{P}_l^{-1/2}\boldsymbol{\alpha}_l\mathbf{U}_l^H\widehat{\mathbf{H}}_l^H\mathbf{G}_k\mathbf{P}_k^{1/2}\|_F^2 + \sigma_{el}^2\|\mathbf{R}_{ml}^{1/2}\mathbf{U}_l\boldsymbol{\alpha}_l\mathbf{P}_l^{-1/2}\|_F^2\|\mathbf{R}_{bl}\mathbf{G}_k\mathbf{P}_k^{1/2}\|_F^2), & k \neq l. \end{cases}$$

The symbol wise AMSE transfer (from the uplink channel to downlink channel and vice versa) can be examined similar to Subsections 2.4.1.2 and 2.4.2.2. The details are omitted for conciseness.

## 2.5 Applications of AMSE duality

To show the applications of our AMSE duality, in this section, we examine the problem of jointly designing the precoders and decoders for the downlink multiuser MIMO systems to minimize: (i) the weighted sum MSE under a total BS power constraint ($\mathcal{P}$2.1) and (ii) the maximum weighted user AMSE constrained with a total BS power ($\mathcal{P}$2.2). Both design problems provide robustness against the channel uncertainties.

### 2.5.1 The robust weighted sum MSE minimization problem ($\mathcal{P}$2.1)

In the downlink channel, the robust weighted sum MSE minimization problem ($\mathcal{P}$2.1) can be formulated by

$$\min_{\mathbf{G}_k,\mathbf{U}_k,\boldsymbol{\alpha}_k,\mathbf{P}_k} \mathrm{E}_{\mathbf{E}_{wk}^H} \sum_{k=1}^{K} \tau_k \mathrm{tr}\{\boldsymbol{\xi}_k^{DL}\} = \sum_{k=1}^{K} \tau_k \overline{\xi}_k^{DL}$$
$$\text{s.t} \sum_{k=1}^{K} \mathrm{tr}\{\mathbf{P}_k\} \leq P_{max}, \{\mathbf{g}_{ki}^H \mathbf{g}_{ki} = \mathbf{u}_{ki}^H \mathbf{u}_{ki} = 1, p_{ki} > 0\}_{i=1}^{S_k}, \forall k \tag{2.16}$$

where $\tau_k$ is the AMSE weighting factor of the $k$th user (when $\{\tau_k = 1\}_{k=1}^K$, (2.16) simplifies to sum AMSE problem). In $\sum_{k=1}^K \tau_k \overline{\xi}_k^{DL}$, since the precoders of all users are coupled, $\mathcal{P}$2.1 has more complicated mathematical structure than its dual uplink problem [SSB07], [HJU09]. Due to this, we examine the dual uplink problem of (2.16) which is expressed as

$$\min_{\{\mathbf{G}_k,\mathbf{U}_k,\boldsymbol{\alpha}_k,\mathbf{Q}_k\}_{k=1}^K} \sum_{k=1}^{K} \tau_k \overline{\xi}_k^{UL}$$
$$\text{s.t} \sum_{k=1}^{K} \mathrm{tr}\{\mathbf{Q}_k\} \leq P_{\max}, \ \{\mathbf{g}_{ki}^H \mathbf{g}_{ki} = \mathbf{u}_{ki}^H \mathbf{u}_{ki} = 1, q_{ki} > 0\}_{i=1}^{S_k}, \ \forall k. \tag{2.17}$$

For convenience, we consider (2.17) for the following two cases.

#### Case 1: When $\{\tau_k = 1, \widetilde{\mathbf{R}}_{mk} = \mathbf{I}_{M_k}, \mathbf{R}_{bk} = \mathbf{R}_b$ and $\sigma_{ek}^2 = \sigma_e^2\}_{k=1}^K$

In this case, first, for a fixed transmitter, the $k$th user receiver is optimized by using the MAMSE method which yields

$$\widetilde{\mathbf{G}}_k \triangleq \mathbf{G}_k \boldsymbol{\alpha}_k \mathbf{Q}_k^{-1/2} = \boldsymbol{\Gamma} \widehat{\mathbf{H}}_k \mathbf{U}_k \mathbf{Q}_k^{1/2} \tag{2.18}$$

where $\boldsymbol{\Gamma} = [\sum_{i=1}^K (\widehat{\mathbf{H}}_i \mathbf{U}_i \mathbf{Q}_i \mathbf{U}_i^H \widehat{\mathbf{H}}_i^H + \widetilde{\sigma}_e^2 \mathrm{tr}\{\mathbf{Q}_i\}\mathbf{R}_b) + \sigma^2 \mathbf{I}_N]^{-1}$ and $\widetilde{\sigma}_e^2 = \sigma_e^2/(\sigma_e^2 + 1)$. Then, after substituting $\widetilde{\mathbf{G}}_k$ in $\overline{\boldsymbol{\xi}}_k^{UL}$, we get the $k$th user MAMSE matrix as $\widetilde{\boldsymbol{\xi}}_k^{UL} = \mathbf{I}_{S_k} - \mathbf{Q}_k^{1/2} \mathbf{U}_k^H \widehat{\mathbf{H}}_k^H \boldsymbol{\Gamma} \widehat{\mathbf{H}}_k \mathbf{U}_k \mathbf{Q}_k^{1/2}$. It follows that

$$\sum_{k=1}^{K} \mathrm{tr}\{\widetilde{\boldsymbol{\xi}}_k^{UL}\} = \sum_{k=1}^{K} \mathrm{tr}\{\mathbf{I}_{S_k} - \mathbf{Q}_k^{1/2} \mathbf{U}_k^H \widehat{\mathbf{H}}_k^H \boldsymbol{\Gamma} \widehat{\mathbf{H}}_k \mathbf{U}_k \mathbf{Q}_k^{1/2}\} \tag{2.19}$$
$$= S - N + \mathrm{tr}\{(\sum_{k=1}^{K} \widetilde{\sigma}_e^2 \mathrm{tr}\{\mathbf{Q}_k\} \mathbf{R}_b + \sigma^2 \mathbf{I}_N) \boldsymbol{\Gamma}\}$$

where the second equality is derived using the matrix inversion Lemma and the fact $(\mathbf{AB})^{-1} = \mathbf{B}^{-1}\mathbf{A}^{-1}$ [PP08]. Thus, (2.17) can be solved by applying a two step approach. First $\mathbf{U}_k$ and $\mathbf{Q}_k$ are optimized by

$$\begin{aligned}&\min_{\{\mathbf{U}_k,\mathbf{Q}_k\}_{k=1}^K} \operatorname{tr}\{(\sum_{k=1}^{K}\widetilde{\sigma}_e^2\operatorname{tr}\{\mathbf{Q}_k\}\mathbf{R}_b+\sigma^2\mathbf{I}_N)\mathbf{\Gamma}\}\\ &\text{s.t} \sum_{i=1}^{K}\operatorname{tr}\{\mathbf{Q}_\text{k}\}\leq P_{max},\ \{\mathbf{u}_{ki}^H\mathbf{u}_{ki}=1, q_{ki}>0\}_{i=1}^{S_k},\ \forall k,\end{aligned} \tag{2.20}$$

and then the optimum $\mathbf{G}_k$ and $\boldsymbol{\alpha}_k$ are computed by using (2.18). Note that in (2.18), $\mathbf{G}_k$ and $\boldsymbol{\alpha}_k$ are obtained such that $\{\mathbf{g}_{ki}^H\mathbf{g}_{ki}=1\}_{i=1}^{S_k}$, $\forall k$ and $\{\boldsymbol{\alpha}_k\}_{k=1}^K$ are diagonal matrices.

Using matrix inversion Lemma, (19) can also be written in terms of $\overline{\mathbf{Q}} \triangleq \mathbf{Q}/\operatorname{tr}\{\mathbf{Q}\}$ as $\sum_{k=1}^K \operatorname{tr}\{\widetilde{\boldsymbol{\xi}}_k^{UL}\} = \operatorname{tr}\{(\mathbf{I}_S + \overline{\mathbf{Q}}^{1/2}\mathbf{U}^H\widehat{\mathbf{H}}^H(\widetilde{\sigma}_e^2\mathbf{R}_b + \frac{\sigma^2}{\operatorname{tr}\{\mathbf{Q}\}}\mathbf{I}_N)^{-1}\widehat{\mathbf{H}}\mathbf{U}\overline{\mathbf{Q}}^{1/2})^{-1}\}$. According to [ZPO08], for the given $\overline{\mathbf{Q}}$, $\sum_{k=1}^K \operatorname{tr}\{\widetilde{\boldsymbol{\xi}}_k^{UL}\}$ is a non-increasing function of $\operatorname{tr}\{\mathbf{Q}\} = P_{\text{sum}}$. Since the difference between (2.19) and the objective function of (2.20) is only the constant term $S-N$, it is clear that the latter objective function is also non-increasing in $P_{\text{sum}}$. By defining $\{\overline{\mathbf{U}}_k \triangleq \mathbf{U}_k\mathbf{Q}_k\mathbf{U}_k^H\}_{k=1}^K$, problem (2.20) can thus be equivalently formulated as

$$\begin{aligned}&\min_{\{\overline{\mathbf{U}}_k\}_{k=1}^K} \operatorname{tr}\{(\widetilde{\sigma}_e^2 P_{\max}\mathbf{R}_b+\sigma^2\mathbf{I}_N)\widetilde{\mathbf{\Gamma}}\}\\ &\text{s.t} \sum_{i=1}^{K}\operatorname{tr}\{\overline{\mathbf{U}}_\text{k}\}= P_{max},\ \overline{\mathbf{U}}_k\succeq 0,\ \operatorname{rank}\{\overline{\mathbf{U}}_k\}=\min(M_k,S_k),\ \forall k\end{aligned} \tag{2.21}$$

where $\widetilde{\mathbf{\Gamma}} = [\sum_{i=1}^K \widehat{\mathbf{H}}_i\overline{\mathbf{U}}_i\widehat{\mathbf{H}}_i^H + \widetilde{\sigma}_e^2 P_{max}\mathbf{R}_b + \sigma^2\mathbf{I}_N]^{-1}$. If we ignore (relax) the rank-constraint of (2.21), the above problem can be formulated as a SDP problem for which global optimum is guaranteed [LDGW04], [SY04], [ZPO07]. Now, if the optimal solution of this SDP satisfies $\operatorname{rank}\{\overline{\mathbf{U}}_k\} = \min(M_k, S_k)$, the latter solution can be considered as a global minimizer of (2.21), otherwise, the solution is deemed as the lower bound solution of (2.21). After computing the solution of (2.21), the optimum $\{\mathbf{U}_k,\mathbf{Q}_k\}_{k=1}^K$ are determined from the eigenvalue decomposition of $\{\overline{\mathbf{U}}_k\}_{k=1}^K$ (see Table I of [SSB07]). It turns out that the optimum (either local or global) solution of (2.16) is computed by using our sum AMSE transfer (see Section 2.4.1.1).

For $\{M_k = S_k = L\}_{k=1}^K$, the approach of [SSB07] requires $O(K^3N^3)$ operations to transfer the powers from uplink to downlink channel (see appendix of [SSB07]) whereas our proposed method needs only $O(KL)$ operations. Thus, as claimed in Section 2.1, the proposed power transformation requires less computation than that of in [SSB07].

### Case 2: For any $\{\tau_k, \widetilde{\mathbf{R}}_{mk}, \mathbf{R}_{bk}$ and $\sigma_{ek}^2\}_{k=1}^K$

In such general case, (2.17) can not be formulated as an SDP problem. Thus, the solution method discussed for **Case 1** can not be applied. Due to this, here we first formulate the power allocation part of (2.17) as a GP for which global optimality is guaranteed. Then, based on the solution of GP, MAMSE receiver and AMSE duality results, we solve (2.16) using the alternating optimization method like in [ECV09]. To this end, we rewrite $\overline{\xi}_k^{UL}$ into a form which is suitable for the GP formulation. Using (2.6), we can express $\overline{\xi}_k^{UL}$ as

$$\overline{\xi}_k^{UL} = \lambda_k + \widetilde{q}_k^{-1} \sum_{i=1, i\neq k}^{K} \widetilde{q}_i v_{ki} + \sigma^2 \widetilde{q}_k^{-1} \vartheta_k \tag{2.22}$$

where $\mathbf{Q}_k = \widetilde{q}_k \widetilde{\mathbf{Q}}_k$, $\text{tr}\{\widetilde{\mathbf{Q}}_k\} = 1$, $\widetilde{\mathbf{U}}_k = \mathbf{U}_k \widetilde{\mathbf{Q}}_k \mathbf{U}_k^H$, $\lambda_k = \text{tr}\{\widetilde{\mathbf{Q}}_k^{-1} \boldsymbol{\alpha}_k^2 \mathbf{G}_k^H (\widehat{\mathbf{H}}_k \widetilde{\mathbf{U}}_k \widehat{\mathbf{H}}_k^H + \sigma_{ek}^2 \text{tr}\{\mathbf{R}_{mk} \widetilde{\mathbf{U}}_k\} \mathbf{R}_{bk}) \mathbf{G}_k - \boldsymbol{\alpha}_k \mathbf{G}_k^H \widehat{\mathbf{H}}_k \mathbf{U}_k - \mathbf{U}_k^H \widehat{\mathbf{H}}_k^H \mathbf{G}_k \boldsymbol{\alpha}_k\} + S_k$, $\vartheta_k = \text{tr}\{\widetilde{\mathbf{Q}}_k^{-1} \boldsymbol{\alpha}_k^2\}$ and $v_{ki} = \text{tr}\{\widetilde{\mathbf{Q}}_k^{-1} \boldsymbol{\alpha}_k^2 \mathbf{G}_k^H (\widehat{\mathbf{H}}_i \widetilde{\mathbf{U}}_i \widehat{\mathbf{H}}_i^H + \sigma_{ei}^2 \text{tr}\{\mathbf{R}_{mi} \widetilde{\mathbf{U}}_i\} \mathbf{R}_{bi}) \mathbf{G}_k\}$. Once again, we can rewrite $\overline{\xi}_k^{UL}$ as

$$\overline{\xi}_k^{UL} = \widetilde{q}_k^{-1} [\mathbf{Y}\widetilde{\mathbf{q}} + \sigma^2 \boldsymbol{\vartheta}]_k \tag{2.23}$$

where $\widetilde{\mathbf{q}} = [\widetilde{q}_1, \cdots, \widetilde{q}_K]$, $\boldsymbol{\vartheta} = [\vartheta_1, \cdots, \vartheta_K]^T$ and

$$[\mathbf{Y}]_{k,i} = \begin{cases} \lambda_k, & for\ \ k = i \\ v_{k,i}, & for\ \ k \neq i. \end{cases}$$

As can be seen from the above equation, for fixed $\{\widetilde{\mathbf{Q}}_k, \widetilde{\mathbf{U}}_k$ and $\boldsymbol{\alpha}_k\}_{k=1}^K$, (2.23) is a posynomial (where $\widetilde{\mathbf{q}} = [\widetilde{q}_1, \cdots, \widetilde{q}_K]$ are the variables). Thus, the power allocation part of (2.17) is formulated by the following GP problem.

$$\min_{\{\widetilde{q}_k\}_{k=1}^K} \sum_{k=1}^{K} \tau_k \overline{\xi}_k^{UL}, \text{ s.t } \|\widetilde{\mathbf{q}}\|_1 \leq P_{max}, \ \widetilde{q}_k > 0, \forall k. \tag{2.24}$$

Using the solution of (2.24) and the user wise AMSE duality results, we solve (2.16) by using the alternating optimization technique similar to that of [ECV09]. In general, we can optimize the powers and filters in many possible orders. In this paper we present a particular algorithm where optimization is started in the uplink channel.

### 1) Uplink channel

In the uplink channel first (2.24) is solved. With the solution $\{\widetilde{q}_k\}_{k=1}^K$, the powers $\{\mathbf{Q}_k = \widetilde{q}_k \widetilde{\mathbf{Q}}_k\}_{k=1}^K$ are computed and then $\{\mathbf{G}_k$ and $\boldsymbol{\alpha}_k\}_{k=1}^K$ are updated by the following uplink MAMSE receiver

$$\mathbf{G}_k \boldsymbol{\alpha}_k = \boldsymbol{\Gamma}_c^{-1} \widehat{\mathbf{H}}_k \mathbf{U}_k \mathbf{Q}_k. \tag{2.25}$$

### 2) Downlink channel

Now we switch the optimization to the downlink channel. Thus, we first ensure the same performance as in the uplink channel ($\{\overline{\xi}_k^{DL_1} = \overline{\xi}_k^{UL}\}_{k=1}^K$) by using our user wise AMSE transfer (2.9).

Then, for fixed $\{\mathbf{P}_k\}_{k=1}^K$, the matrices $\{\mathbf{U}_k$ and $\boldsymbol{\alpha}_k\}_{k=1}^K$ are updated by the downlink MAMSE receiver which is given as

$$\mathbf{U}_k \boldsymbol{\alpha}_k = (\boldsymbol{\Gamma}_k^{DL})^{-1} \widehat{\mathbf{H}}_k^H \mathbf{G}_k \mathbf{P}_k. \tag{2.26}$$

At this stage, the $k$th user achieves a new AMSE $\triangleq \overline{\xi}_k^{DL_2} \leq \overline{\xi}_k^{DL_1}$.

### 3) Uplink channel

Like in Step (2) above, we first ensure the same performance as in the downlink channel ($\{\overline{\xi}_k^{UL_1} = \overline{\xi}_k^{DL_2}\}_{k=1}^K$) and then we update $\{\mathbf{G}_k$ and $\boldsymbol{\alpha}_k\}_{k=1}^K$ by (2.25). We observe less overall computational time if the latter two steps are performed before proceeding to the next iteration. The detailed iterative steps to solve (2.16) are summarized in **Table 2.I** (**Algorithm 2.I**).

### 2.5.2 The robust weighted MSE min-max problem ($\mathcal{P}2.2$)

In the downlink channel, for given user wise AMSE weights $\{\eta_k\}_{k=1}^{K}$, $\mathcal{P}2.2$ can be formulated by

$$\min_{\mathbf{P}_k,\mathbf{G}_k,\mathbf{U}_k,\boldsymbol{\alpha}_k} \max \frac{\mathrm{E}_{\mathbf{E}_{wk}^{H}} \mathrm{tr}\{ ,_{k}^{DL} \}}{\eta_k} = \frac{\{\bar{\xi}_k^{DL}\}}{\eta_k}$$
$$\text{s.t} \sum_{k=1}^{K} \mathrm{tr}\{\mathbf{P}_k\} \leq P_{max}, \{\mathbf{g}_{ki}^{H}\mathbf{g}_{ki} = \mathbf{u}_{ki}^{H}\mathbf{u}_{ki} = 1, p_{ki} > 0\}_{i=1}^{S_k}, \forall k. \tag{2.27}$$

Here we first solve the power allocation part of (2.27), then we use the solution framework of $\mathcal{P}2.1$ (**Case 2**) to jointly optimize the transceivers. To this end, for fixed $\{\widetilde{\mathbf{Q}}_k, \boldsymbol{\alpha}_k, \mathbf{G}_k\}_{k=1}^{K}$, the uplink power allocation part of (2.27) reads

$$\mu^{UL} \triangleq \min_{\widetilde{q}_k} \max \frac{\bar{\xi}_k^{UL}}{\eta_k}, \text{ s.t } \|\widetilde{\mathbf{q}}\|_1 \leq P_{max}, \ \widetilde{q}_k > 0, \forall k. \tag{2.28}$$

The global optimal solution of the above optimization problem satisfies the following relations [SSB08a]

$$\mu^{UL} = \frac{\bar{\xi}_k^{UL}}{\eta_k}, \ \forall k \ \text{ and } \ \|\widetilde{\mathbf{q}}\|_1 = P_{max}. \tag{2.29}$$

Moreover, by defining $\boldsymbol{\eta} \triangleq \mathrm{diag}\{\eta_1, \eta_2, \cdots, \eta_K\}$ the following eigensystem can be formed from (2.23) and (2.29).

$$\boldsymbol{\Omega} \begin{bmatrix} \widetilde{\mathbf{q}} \\ 1 \end{bmatrix} = \mu^{UL} \begin{bmatrix} \widetilde{\mathbf{q}} \\ 1 \end{bmatrix}, \quad \text{where } \boldsymbol{\Omega} = \begin{bmatrix} \boldsymbol{\eta}^{-1}\mathbf{Y} & \sigma^2\boldsymbol{\eta}^{-1}\boldsymbol{\vartheta} \\ \frac{1}{P_{max}}\mathbf{1}_K^T\boldsymbol{\eta}^{-1}\mathbf{Y} & \frac{\sigma^2}{P_{max}}\mathbf{1}_K^T\boldsymbol{\eta}^{-1}\boldsymbol{\vartheta} \end{bmatrix}. \tag{2.30}$$

Therefore, the optimal solution of (2.28) is given by $\mu^{UL} = \lambda_{max}(\boldsymbol{\Omega})$ and $[\widetilde{\mathbf{q}}\ 1]^T$ is the eigenvector of $\boldsymbol{\Omega}$ corresponding to $\mu^{UL}$ [SSB08a]. By using the optimal $\widetilde{\mathbf{q}}$ of (2.28), MAMSE receiver and AMSE duality results, (2.27) can be solved as shown in **Table 2.I** (**Algorithm 2.II**).

## 2.6 Simulation results

In all simulations, we take $K = 2$, $\{M_k = S_k = 2\}_{k=1}^{K}$ and $N = 4$. We model $\{\mathbf{R}_{bk}, \widetilde{\mathbf{R}}_{mk}\}_{k=1}^{K}$ as $\{\mathbf{R}_{bk}\}_{k=1}^{K} = \mathbf{R}_b = \rho_b^{|i-j|}$, $\{\widetilde{\mathbf{R}}_{mk}\}_{k=1}^{K} = \widetilde{\mathbf{R}}_m = \rho_m^{|i-j|}$,

**Table 2.I**: Iterative solution for problems $\mathcal{P}2.1$ (2.16) and $\mathcal{P}2.2$ (2.27)
Initialization: Set equal powers for all symbols, i.e., $\{\mathbf{Q}_k = \frac{P_{max}}{S}\mathbf{I}_{S_k}\}_{k=1}^K$ and $\mathbf{U}_k \in \mathcal{C}^{M_k \times S_k}$ as the first $S_k$ right-hand singular value decomposition vectors of $\widehat{\mathbf{H}}_k$, $\forall k$ and then update $\{\mathbf{G}_k$ and $\boldsymbol{\alpha}_k\}_{k=1}^K$ by (2.25).

**Repeat Virtual uplink channel**

1. Set $\{\widetilde{\mathbf{Q}}_k = \mathbf{Q}_k/\mathrm{tr}\{\mathbf{Q}_k\}\}_{k=1}^K$.
2. For $\mathcal{P}2.1$, compute $\{\widetilde{q}_k\}_{k=1}^K$ using (2.24) (**Algorithm 2.I**).
3. For $\mathcal{P}2.2$, with the given $\boldsymbol{\eta} = \mathrm{diag}\{\eta_1, \eta_2, \cdots, \eta_K\}$, compute $\{\widetilde{q}_k\}_{k=1}^K$ and $\mu^{UL}$ using (2.30). In this step the power constraint $||\widetilde{\mathbf{q}}||_1 = P_{max}$ is ensured by scaling the eigenvector corresponding to $\lambda_{max}(\boldsymbol{\Omega})$ such that the last element equals 1. (**Algorithm 2.II**)
4. Update $\{\mathbf{Q}_k = \widetilde{q}_k\widetilde{\mathbf{Q}}_k\}_{k=1}^K$. Using the latter $\{\mathbf{Q}_k\}_{k=1}^K$, update $\{\mathbf{G}_k$ and $\boldsymbol{\alpha}_k\}_{k=1}^K$ by (2.25). Then, compute $\{\beta_k\}_{k=1}^K$ with (2.11).

**Downlink Channel**

5. Set $\{\mathbf{P}_k = \beta_k\boldsymbol{\alpha}_k^2\mathbf{Q}_k^{-1}\}_{k=1}^K$. Using this $\{\mathbf{P}_k\}_{k=1}^K$, update $\{\mathbf{U}_k$ and $\boldsymbol{\alpha}_k\}_{k=1}^K$ by (2.26). Then, compute $\{\widetilde{\beta}_k\}_{k=1}^K$ with (2.14).

**Virtual uplink Channel**.

6. Set $\{\mathbf{Q}_k = \widetilde{\beta}_k\boldsymbol{\alpha}_k^2\mathbf{P}_k^{-1}\}$. Using the latter $\{\mathbf{Q}_k\}_{k=1}^K$, update $\{\mathbf{G}_k$ and $\boldsymbol{\alpha}_k\}_{k=1}^K$ by (2.25).

**Until** convergence

**Convergence**: It can be shown that **Algorithm 2.I** and **Algorithm 2.II** are convergent [SSB08c], [SSB08a]. Different initializations affect the convergence speed of both algorithms. In most of our simulations ($> 95\%$), we observe fast convergence when the initialization is performed as in this table.
**Global optimality**: Since $\mathcal{P}2.1$ (**Case 2**) and $\mathcal{P}2.2$ are non-convex, we can not prove the global optimality of **Algorithm 2.I** and **Algorithm 2.II**. However, for $\mathcal{P}2.1$ (**Case 1** with $M_k = S_k$), simulation results show that **Algorithm 2.I** achieves global minimum (see the next Section).

where $0 \leq \rho_b$ $(\rho_m) < 1$. The SNR is defined as $P_{\max}/\sigma^2$ where $P_{\max}$ is the maximum total BS power and $\sigma^2$ is the noise variance. The SNR is controlled by varying $\sigma^2$ while $P_{\max}$ is set to $10mW$. All simulation results are obtained by averaging over 100 randomly chosen channel realizations.

### 2.6.1 Simulation results for problem $\mathcal{P}2.1$

#### 2.6.1.1 Simulation results for Case 1

For the aforementioned parameters, all our simulation results show that the optimal solution of the SDP problem (the problem (2.21) after rank relaxation) satisfy the rank constraint of (2.21)[3]. Consequently, for our setup, the SDP solution is considered as a GM for (2.17). Similar observation has also been made in the perfect CSI case [SSB07]. Now, we check whether **Algorithm 2.I** achieves the GM or not when $\mathbf{R}_b = \mathbf{I}$, $\sigma_{e1}^2 = \sigma_{e2}^2 = 0.0101$. Fig. 2.2.a shows that the GM can be achieved by **Algorithm 2.I** for the robust and perfect CSI designs. In the non-robust/naive design, which refers to the design in which the estimated channel is considered as perfect [HJU09], the gap between **Algorithm 2.I** and the GM is large in the high SNR zone.

#### 2.6.1.2 Comparison of robust and non-robust/naive designs

For **Case 1**, as can be seen from Fig. 2.2.a, the robust design has better performance than the non-robust design. Now, we compare the performance of our robust design with that of the non-robust design proposed in [SSB07] for **Case 2**. The comparison is based on the total sum AMSE and ASER[4] of all users versus the SNR when $\{\tau_k = 1\}_{k=1}^K$.

#### 2.1) When $\{\sigma_{ek}^2 \neq 0\}_{k=1}^K$, $\rho_b \neq 0$ and $\rho_m = 0$

Here we examine the joint effect of $\rho_b$ and $\{\sigma_{ek}^2\}_{k=1}^K$ on the system performance. To this end, we set $\sigma_{e1}^2 = 0.0101, \sigma_{e2}^2 = 0.0204$ and then we vary $\rho_b$ from

[3]We have noted that the SDP solution of this problem does not always satisfy its rank constraint when $S_k < M_k$. Simulation results for the case $S_k < M_k$ are not included for conciseness.

[4]For the sum AMSE design, ASER is also an appropriate metric for comparing the performance of the robust and non-robust designs [DB09]. QPSK modulation is utilized for each symbol stream.

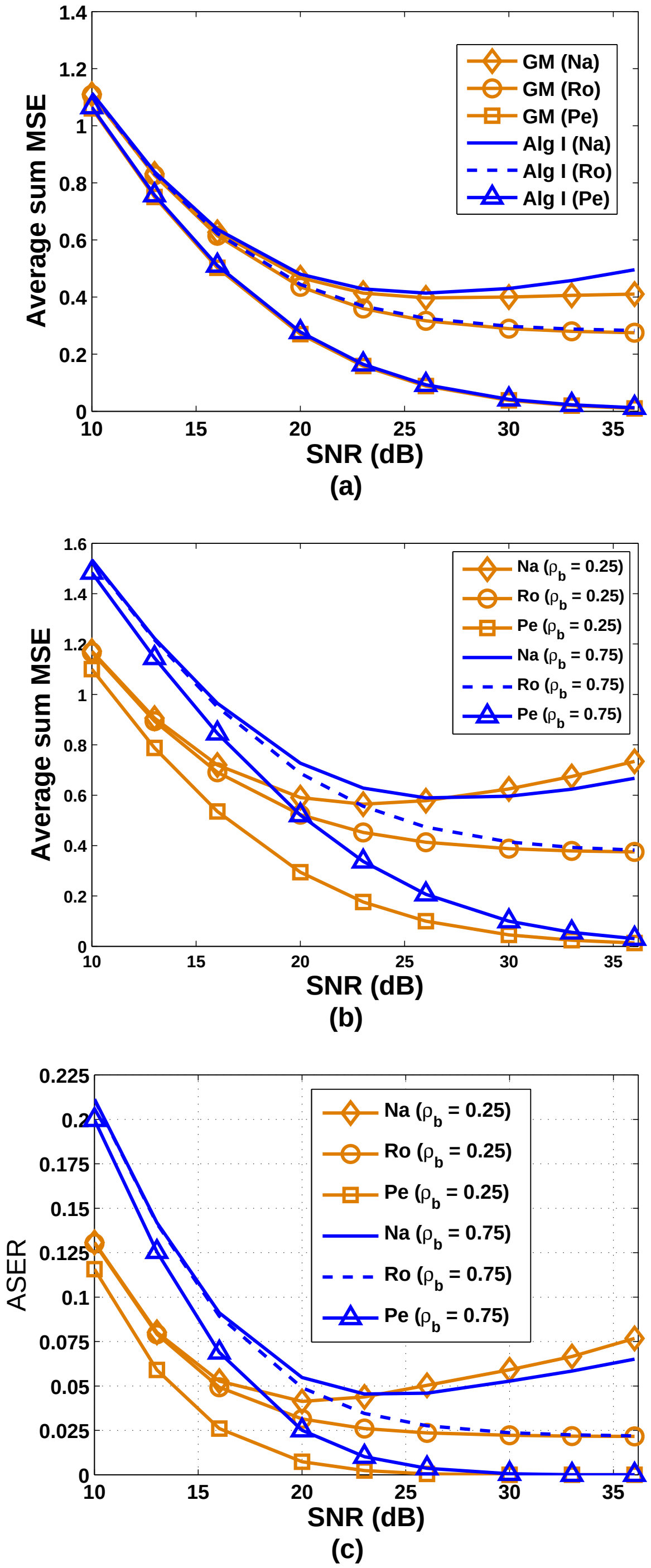

**Figure 2.2** (a) Comparison of the performance of **Algorithm 2.I** and the GM. (b)-(c) Comparison of the robust and non-robust/naive designs when $\rho_m = 0$, and $\rho_b = 0.25$ and 0.75. The non-robust/naive, robust and perfect CSI designs are denoted by 'Na', 'Ro' and 'Pe', respectively.

0.25 to 0.75. Figs. 2.2.(b-c) show that as the BS antenna correlation coefficient increases, the sum AMSE and ASER also increase.

**2.2) When** $\{\sigma_{ek}^2 \neq 0\}_{k=1}^K$, $\rho_b \neq 0, \rho_m \neq 0$

Now we discuss the effects of $\{\sigma_{ek}^2\}_{k=1}^K$, $\rho_b$ and $\rho_m$ on the system performance. We keep $\sigma_{e1}^2 = 0.0101, \sigma_{e2}^2 = 0.0204$, $\rho_b = 0.25$ and then we take $\rho_m$ as 0.25 and 0.75. Figs. 2.3.(a-b) show that the performance of the system degrades further as $\rho_m$ increases[5].

The results of Section 2.6.1.2 gracefully fit to that of [DB09] where (2.16) is examined for single user MIMO systems.

### 2.6.2 Simulation results for problem $\mathcal{P}2.2$

This simulation compares the performance of the robust design and the non-robust design proposed in [SSB08a]. Here we keep $\eta_1 = \eta_2 = 0.3$, $\sigma_{e1}^2 = 0.0101, \sigma_{e2}^2 = 0.0204$, $\rho_b = 0.25$ and then we take $\rho_m$ as 0.25 and 0.75. Fig. 2.3.c shows that the maximum AMSE of the robust design is less than that of the non-robust design proposed in [SSB08a]. Moreover, the performance gap between these designs increases as the SNR increases. This figure also illustrates the fact that large antenna correlation factor degrades the performance of the considered system.

**Effect of $\rho_b$ $(\rho_m)$ on the performances of the robust and non-robust designs for both $\mathcal{P}1$ and $\mathcal{P}1$:** As can be seen from the figures of this chapter, in all figures, the robust design outperforms the non-robust design. When $\rho_b$ $(\rho_m)$ increases, the system performance degrades. This is because as $\rho_b$ $(\rho_m)$ increases, the number of symbols with low channel gain increases (this can be easily seen from the eigenvalue decomposition of $\mathbf{R}_{bk}$ $(\mathbf{R}_{mk})$). Consequently, for fixed total BS power, the total sum AMSE ($\mathcal{P}2.1$) and maximum AMSE ($\mathcal{P}2.2$) also increase.

---

[5]**Remark:** When the SNR $\rightarrow \infty$ (i.e., $\sigma^2 \rightarrow 0$), a sum AMSE floor exists for our robust design. Such sum AMSE floor is observed in Figs. 2.2.(a-b) and Fig. 2.3.a. The analytical proof is given as follows: for any $\{\sigma_{ek}^2, \widetilde{\mathbf{R}}_{mk}, \mathbf{R}_{bk}\}_{k=1}^K$ and $\{\tau_k = 1\}_{k=1}^K$, after some mathematical manipulations the uplink sum MAMSE can be expressed as $\text{tr}\{(\mathbf{I}_S + \mathbf{Q}^{1/2}\mathbf{U}^H\widehat{\mathbf{H}}^H(\sigma^2\mathbf{I} + \sum_{i=1}^K \sigma_{ei}^2 \text{tr}\{\mathbf{R}_{mi}\mathbf{U}_i\mathbf{Q}_i\mathbf{U}_i^H\}\mathbf{R}_{bi})^{-1}\widehat{\mathbf{H}}\mathbf{U}\mathbf{Q}^{1/2})^{-1}\}$. Hence, when $\sigma^2 \rightarrow 0$ the sum MAMSE approaches to $\text{tr}\{(\mathbf{I}_S + \mathbf{Q}^{1/2}\mathbf{U}^H\widehat{\mathbf{H}}^H(\sum_{i=1}^K \sigma_{ei}^2 \text{tr}\{\mathbf{R}_{mi}\mathbf{U}_i\mathbf{Q}_i\mathbf{U}_i^H\}\mathbf{R}_{bi})^{-1}\widehat{\mathbf{H}}\mathbf{U}\mathbf{Q}^{1/2})^{-1}\} > 0$.

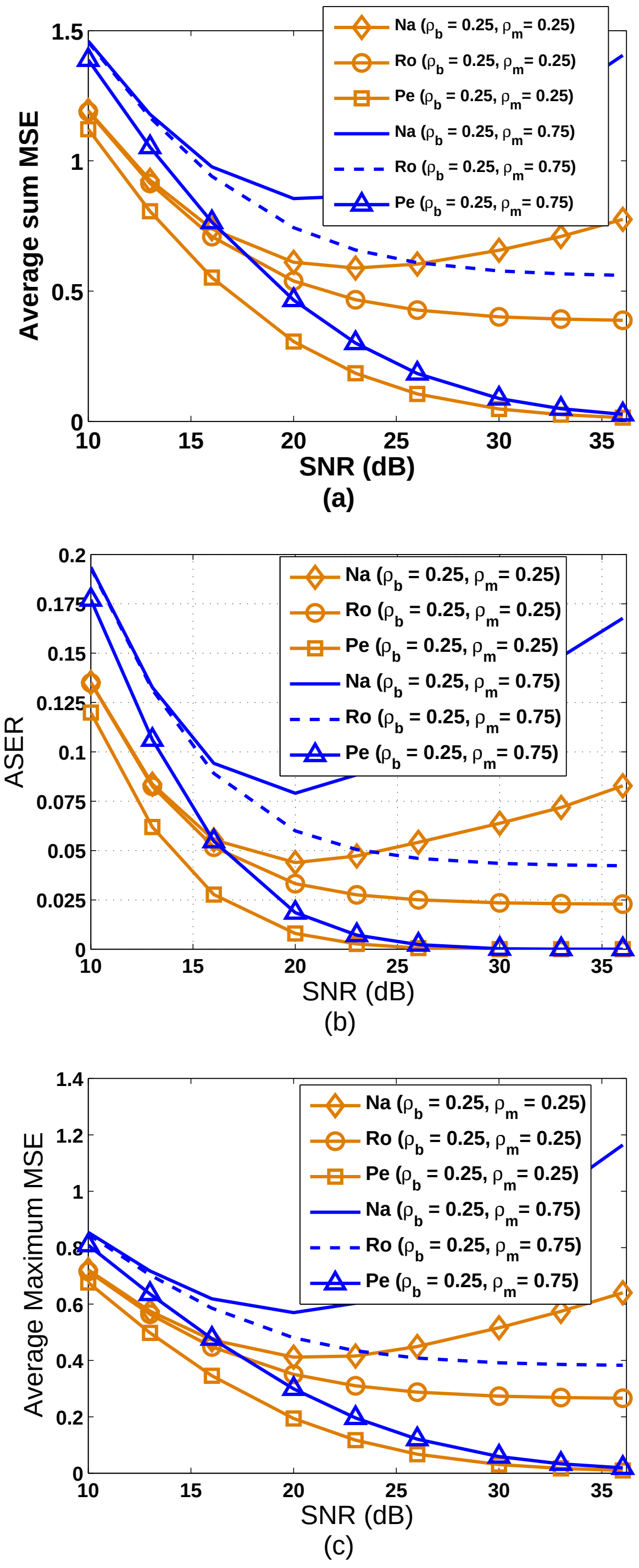

**Figure 2.3** Comparison of the robust and non-robust/naive designs when $\rho_b = 0.25$, and $\rho_m = 0.25$ and $0.75$. (a)-(b) For the robust sum MSE minimization problem ($\mathcal{P}2.1$). (c) For the robust Min-max problem ($\mathcal{P}2.2$).

**Effect of $\sigma_{ek}^2$ on the performances of the robust and non-robust designs for both $\mathcal{P}1$ and $\mathcal{P}1$**: For better explanation about the effect of $\sigma_{ek}^2$ on the performances of the robust and non-robust designs, let us consider the case where $\mathbf{R}_{mk} = \mathbf{R}_{bk} = \mathbf{I}$. For this setting, we will have $\tilde{\sigma}^2 \triangleq \sigma_{ek}^2 \text{tr}\{\mathbf{R}_{bk}\mathbf{G}\mathbf{P}\mathbf{G}^H\}\mathbf{R}_{mk} = \sigma_{ek}^2 P_{max}$ and the MSE of the $k$th user (2.3) can be equivalently expressed as

$$\bar{\xi}_k^{DL} = S_k + \text{tr}\{\mathbf{P}_k^{-1}\boldsymbol{\alpha}_k^2\mathbf{U}_k^H\tilde{\boldsymbol{\Gamma}}_k^{DL}\mathbf{U}_k - 2\Re\{\mathbf{G}_k^H\widehat{\mathbf{H}}_k\mathbf{U}_k\boldsymbol{\alpha}_k\}\}$$

where $\tilde{\boldsymbol{\Gamma}}_k^{DL} = (\widehat{\mathbf{H}}_k^H\mathbf{G}\mathbf{P}\mathbf{G}^H\widehat{\mathbf{H}}_k + \tilde{\tilde{\sigma}}^2)$ and $\tilde{\tilde{\sigma}}^2 = \tilde{\sigma}^2 + \sigma^2$.

Now let us consider the performance of the robust and non-robust designs for two extreme SNR values: very low SNR (i.e., $\sigma^2 \gg P_{max}$) and very high SNR (i.e., $\sigma^2 \rightarrow 0$).

At very low SNR, we will have $\tilde{\tilde{\sigma}}_{robust}^2 = \tilde{\sigma}^2 + \sigma^2 \approx \sigma^2 = \tilde{\tilde{\sigma}}_{non-robust}^2$. Thus, the performance gap between the robust and non-robust designs is insignificant. However, at very high SNR, we will have $\tilde{\tilde{\sigma}}_{robust}^2 = \tilde{\sigma}^2 + \sigma^2 \approx \tilde{\sigma}^2$(i.e., true $\tilde{\tilde{\sigma}}^2$) and $\tilde{\tilde{\sigma}}_{non-robust}^2 = \sigma^2 \approx 0$ (i.e., significantly deviated from the true $\tilde{\tilde{\sigma}}^2$). Thus, at very high SNR, since the non-robust design significantly deviated from the true $\tilde{\tilde{\sigma}}^2$, its performance is much worse compared to the robust design which employs the true $\tilde{\tilde{\sigma}}^2$. From this explanation, we can understand that the gap between the robust and non-robust designs will increase as the SNR increases. This scenario is clearly seen from all the figures of this section.

## 2.7 Conclusions

In this chapter, three kinds of MSE uplink-downlink duality are established by considering that the BS and MS antennas exhibit spatial correlations and the CSI at both the transmitter and receiver ends are imperfect. The duality are established by transforming only the power allocation diagonal matrices from uplink to downlink channel and vice versa. As application examples of the duality, two MSE-based robust transceiver design problems are considered. To solve the problems, we propose duality based iterative algorithms that perform optimization alternatively by switching between the uplink and downlink channels. Simulation results show the superior performance of the proposed robust design compared to the non-robust/naive design.

## 2.8 Looking ahead

The proposed duality based algorithm of this chapter can solve only total BS power constrained MSE-based problems[6]. Recently there is a growing interest to design precoders and decoders for coordinated BS systems. And in such systems, BSs are located at different places. For this reason, the maximum available power of each BS is limited. Clearly, one can see that the approach of this chapter can not handle such BS power constrained problems. Furthermore, for large networks, distributed precoder/design has an interest. In the next chapter, we will examine distributive transceiver designs for rate and MSE-based problems for coordinated BS systems.

[6]We would like to mention here that under perfect CSI, the duality of this chapter can also be applied to solve rate and SINR-base problems.

# Weighted Sum Rate Optimization for Downlink Multiuser MIMO Coordinated Base Station Systems: Centralized and Distributed Algorithms

# 3

This chapter considers the joint linear transceiver design problem for the downlink multiuser MIMO systems with fully coordinated BSs. We consider maximization of the weighted sum rate with per BS antenna power constraint problem. To solve the problem, new centralized and computationally efficient distributed iterative algorithms have been presented. These algorithms are described as follows. First, by introducing additional optimization variables, we reformulate the original problem into a new problem. Second, for the given precoder matrices of all users, the optimal receivers are computed using MMSE method and the optimal introduced variables are obtained in closed form expressions. Third, by keeping the introduced variables and receivers constant, the precoder matrices of all users are optimized by using SOCP and matrix fractional minimization approaches for the centralized and distributed algorithms, respectively. Finally, the second and third steps are repeated until these algorithms converge. It is shown that the proposed algorithms are guaranteed to converge. Also the proposed algorithms require less computational cost than that of the existing linear algorithm. All simulation results demonstrate that our distributed algorithm achieves the same performance as that of the centralized algorithm. Moreover, the proposed algorithms outperform the existing linear algorithm. In particular, when each of the users has single antenna, we have observed that the proposed algorithms achieve the global optimum. As will be clear later, the centralized and distributed algorithms of this chapter can be applied to solve per BS antenna power constrained MSE-based problems. The extension of this chapter to multicell systems has also been discussed.

## 3.1 Introduction

MIMO systems have been proven to enhance the spectral efficiency of wireless systems without additional transmission power [BL04]. In [VJG03, YC04, WSS04], the achievable sum rate of the BC channel obtained by dirty paper coding (DPC) technique of [Cos83] has been characterized for MIMO systems. In [VJG03,YC04,WSS04], it has been shown that DPC achieves the capacity region of a BC channel. However, since the DPC algorithm is nonlinear, implementing the DPC algorithm appears to be very complicated. Due to this reason, several modified but having less complexity transmission approaches have been proposed. Successive zero-forcing dirty paper coding (SZF-DPC) approach has been proposed in [TJBO13] (i.e., a combination of linear and non-linear transmitter design approach). This approach has also been applied in [TJBO13] to solve weighted sum rate maximization problem for the BC channel. But still all of these algorithms are not purely linear transceiver algorithms. As stated in the introduction chapter of this thesis, linear transceivers have elegant mathematical structure and are suitable for practical realization. In the following, we examine different linear processing algorithms to solve transmission rate based optimization problems[1].

In [SSH04], linear processing method that employs channel block-diagonalization is suggested. The latter method suffers from noise enhancement and has a restriction on the number of transmit and receive antennas.

[1]We are not aware of any linear processing scheme that achieves the capacity region of a BC channel. Thus, the development of a linear transceiver algorithm achieving the capacity region of a BC channel is still an ongoing research topic.

In [SSB08b], weighted sum rate maximization problem for the downlink multiuser MIMO system is formulated as the problem of minimizing the geometric product of MMSE. This paper solves its problem with a per BS antenna power constraint. The latter problem has also been examined in [SSB08c] with a total BS power constraint. To solve the optimization problem, an iterative approach which uses MSE uplink-downlink duality is suggested. Minimizing the product of all users MMSE matrix determinants is proposed as an equivalent formulation for the sum rate maximization problem of the downlink multiuser MIMO systems [TA08]. This problem is non-convex and it is solved by employing sequential quadratic programming. The work of [TA08] has been extended to the robust case in [ECV10]. The latter paper formulates the robust problem using the worst-case robust design approach, and utilizes MSE uplink-downlink duality approach to solve the sum rate maximization problem.

All of the aforementioned papers examine their problems for conventional downlink systems. In these systems, BSs from different cells communicate with their respective remote terminals independently. Hence, inter-cell interference is obliged to be considered as a background noise. Recently, it has been shown that BS coordination communication is a promising technique to significantly improve the capacity of wireless channels by mitigating (or possibly canceling) inter-cell interference [KFV06, DY10, BZGO10]. The BS coordination can be performed by two approaches: In the first approach, BSs are coordinated at the beamforming (precoder) level. In such kind of BS coordination, the system is termed as "multi-cell system" [DY10]. In the second approach, BS coordination takes place at both the signal and beamforming (precoder) levels. When BSs are coordinated in this approach, the system is termed as "network MIMO system" [KFV06], [BZGO10]. It is well know that the latter coordination approach has better performance gain compared to the former one [BZGO10], [BZGO09]. This performance improvement, however, requires additional signal coordination. In this chapter, we focus on the second BS coordination approach[2]. In [EV11], the joint optimization of the precoders to maximize the total sum rate with per BS antenna power constraint has been examined for the downlink multiuser systems with coordinated BSs.

[2]As will be clear later, the solution approach of this chapter can also be applied for multi-cell systems.

The latter paper assumes that the BSs are equipped with multiple antennas and MSs are equipped with single antenna. In [SSVB08], four MSE-based linear transceiver optimization problems have been considered for the downlink multiuser MIMO systems with coordinated BSs. These problems are examined by assuming that the total power of each BS or the individual power of each BS antenna (group of antennas) is constrained. The problems of [SSVB08] are solved as follows. First, by keeping the receivers constant, optimization of the precoder matrices are formulated as a SOCP problem [3]. Second, for the given BS precoders, the receiver of each user is optimized by MMSE technique. These two steps are repeated iteratively to jointly optimize the transmitters and receivers. Thus, in [SSVB08], the receiver of each user can be optimized independently and distributively. However, the joint optimization of the precoders of [SSVB08] has been carried out by a centralized algorithm. When the number of users and/or BSs increase, the computational cost of the joint precoder design also increases [TSC07]. Consequently, solving the precoder optimization problem in a centralized manner, especially for large-scale coordinated networks, is not a computationally efficient approach. This motivates us to develop distributed algorithms for transceiver design problems in the downlink coordinated BS systems with per BS antenna power constraint as in [BVC10]. This paper solves its optimization problems distributively by applying the Lagrangian dual decomposition, modified matrix fractional minimization and an iterative technique.

In this chapter, we design the transmitters and receivers of all users to maximize the weighted sum rate with per BS antenna power constraint problem[4]. We present novel centralized and computationally efficient distributed iterative algorithms that achieve local optimum to the latter problem. These algorithms are described as follows. First, by introducing additional optimization variables, we reformulate the original problem into a new problem. Second, for the given precoder matrices of all users, the optimal receivers are computed using MMSE method and the optimal introduced variables are obtained in closed form expressions. Third, by keeping the introduced variables and receivers constant, the precoder matrices of all users are optimized by using SOCP and matrix fractional minimization approaches for the centralized and distributed algorithms, respectively. Finally, the second and third steps are repeated until these algorithms converge. We have shown that the proposed algorithms are guaranteed to converge. All simulation results show that our proposed distributed algorithm achieve the same performance as that of the centralized algorithm. Moreover, the proposed centralized and distributed algorithms outperform the existing algorithm. In particular, when each of the users has single antenna, we have observed that the proposed algorithms achieve the global optimum. The contribution of this chapter is thus summarized as follows.

[3]Note that SOCP problems are convex and can be solved by using existing convex optimization tools

[4]According to [YL07], in a multi-antenna BS systems, each BS antenna has its own power amplifier and the maximum power of each BS antenna is limited by some value. This motivates us to consider the power constraint of each BS antenna. On the other hand, in some scenario, a per BS power constraint has practical interest. As will be clear later, our proposed algorithm can be extended straightforwardly to handle the latter power constraint and the sum power constraint of the whole network or groups of antennas.

1. New centralized and computationally efficient distributed iterative algorithms are presented to jointly optimize the transceivers of all users to maximize the weighted sum rate with a per BS antenna power constraint problem. Our proposed algorithms can be used for the case where the constraint of this problem is modified to sum power constraint of the whole network or groups of antennas. As will be clear later, we also show that the proposed algorithms can be applied to examine weighted sum rate optimization problem for multi-cell systems.

2. For the aforementioned problem, we have demonstrated that the proposed distributed algorithm has the same performance as that of the centralized algorithm.

3. As will be shown later, our problem has exactly the same mathematical structure as that of in [SSB08b] where weighted sum rate maximization with per antenna power constraint problem is considered for conventional downlink MIMO systems. The latter paper, however, solves the optimization problem by constraining that the power allocated for each symbol is always positive. In other words, the algorithm proposed in [SSB08b] can not handle inactive symbols. The algorithms of this chapter have four major advantages compared to the algorithm in [SSB08b].

First, the algorithms of this chapter constrained the powers of each symbol to be non-negative[5]. Second, simulation results show that our algorithm has better weighted sum rate compared to that of [SSB08b]. Third, as will be clear later in Section 3.4, the presented centralized and distributed algorithms require less computational cost compared to that of [SSB08b]. Fourth, the proposed algorithms have faster convergence speed than that of [SSB08b].

4. When each of the users has single antenna, the global optimal solution of weighted sum rate maximization problem can be obtained with the framework of MGO algorithm as in [BU09]. For this case, in all of our simulation results, we have observed that the proposed centralized and distributed algorithms achieve the global optimum.

The remaining part of this chapter is organized as follows. We present the downlink multiuser MIMO coordinated BS system model in Section 3.2. The problem formulation is discussed in Section 3.3. The existing centralized, and the proposed centralized and distributed algorithms are presented in Section 3.4. The extensions of the centralized and distributed algorithms for the robust versions of $\mathcal{P}3.1$ and $\mathcal{P}3.2$ are discussed in Section 3.5. In Section 3.6, computer simulations are used to compare the performance of the centralized and distributed algorithms, and the proposed algorithms with that of the other existing algorithms. Finally, conclusions are drawn in Section 3.7.

## 3.2 System model

We consider a downlink multiuser MIMO coordinated BS system as shown in Fig. 3.1 where $L$ BSs are serving $K$ decentralized multiantenna MSs. The $l$th BS and $k$th MS are equipped with $N_l$ and $M_k$ antennas, respectively. The total number of BS and MS antennas are thus $N = \sum_{l=1}^{L} N_l$ and $M = \sum_{k=1}^{K} M_k$, respectively. By denoting the symbol intended for the $k$th user as $\mathbf{d}_k \in \mathbb{C}^{S_k \times 1}$ and $S = \sum_{k=1}^{K} S_k$, the entire symbol can be written in a data vector $\mathbf{d} \in \mathbb{C}^{S \times 1}$ as $\mathbf{d} = [\mathbf{d}_1^T, \cdots, \mathbf{d}_K^T]^T$. The $l$th BS precodes $\mathbf{d}$ into an $N_l$ length vector by using its overall precoder matrix $\mathbf{B}_l = [\mathbf{b}_{l11}, \cdots, \mathbf{b}_{lKS_K}]$, where $\mathbf{b}_{lki} \in \mathbb{C}^{N_l \times 1}$ is the precoder vector of the $l$th BS for the $k$th MS $i$th symbol. The $i$th symbol of the $k$th MS employs a receiver $\mathbf{w}_{ki}$ to estimate its symbol $d_{ki}$. We follow the same channel matrix notations as in [BV11a]. The estimate of the $k$th MS $i$th symbol $(\hat{d}_{ki})$ is given by

$$\hat{d}_{ki} = \mathbf{w}_{ki}^H (\sum_{l=1}^{L} \mathbf{H}_{lk}^H \mathbf{B}_l \mathbf{d} + \mathbf{n}_k) = \mathbf{w}_{ki}^H (\mathbf{H}_k^H \mathbf{B} \mathbf{d} + \mathbf{n}_k) \tag{3.1}$$

where $\mathbf{H}_k^H = [\mathbf{H}_{1k}^H, \cdots, \mathbf{H}_{Lk}^H] \in \mathbb{C}^{M_k \times N}$, $\mathbf{B} = [\mathbf{B}_1; \cdots; \mathbf{B}_L] \in \mathbb{C}^{N \times S}$, $\mathbf{H}_{lk}^H \in \mathbb{C}^{M_k \times N_l}$ is the channel vector between the $l$th BS and the $k$th MS, and $\mathbf{n}_k$ is the additive noise at the $k$th MS. As can be seen from (3.1), the $k$th user decodes its symbol $d_{ki}$ independently with the receiver $\mathbf{w}_{ki}$. As will be clear later, this chapter applies MMSE approach to design $\mathbf{w}_{ki}$. On the other hand, the $k$th user can decode its symbol $d_{ki}$ by first canceling known interference (i.e., successive interference cancelation) and then applying MMSE receiver as in [ZCW05]. According to [ZCW05], the latter decoding approach achieves less BER compared to that of the former one. However, since the latter approach is non-linear [SSB08c], this chapter focuses on the former decoding approach which is linear.

[5]We would like to mention here that at optimality the powers of some of the symbols can be zero. This scenario happens, especially, for total sum rate maximization problems. This shows that the presented algorithms are more general than that of [SSB08b].

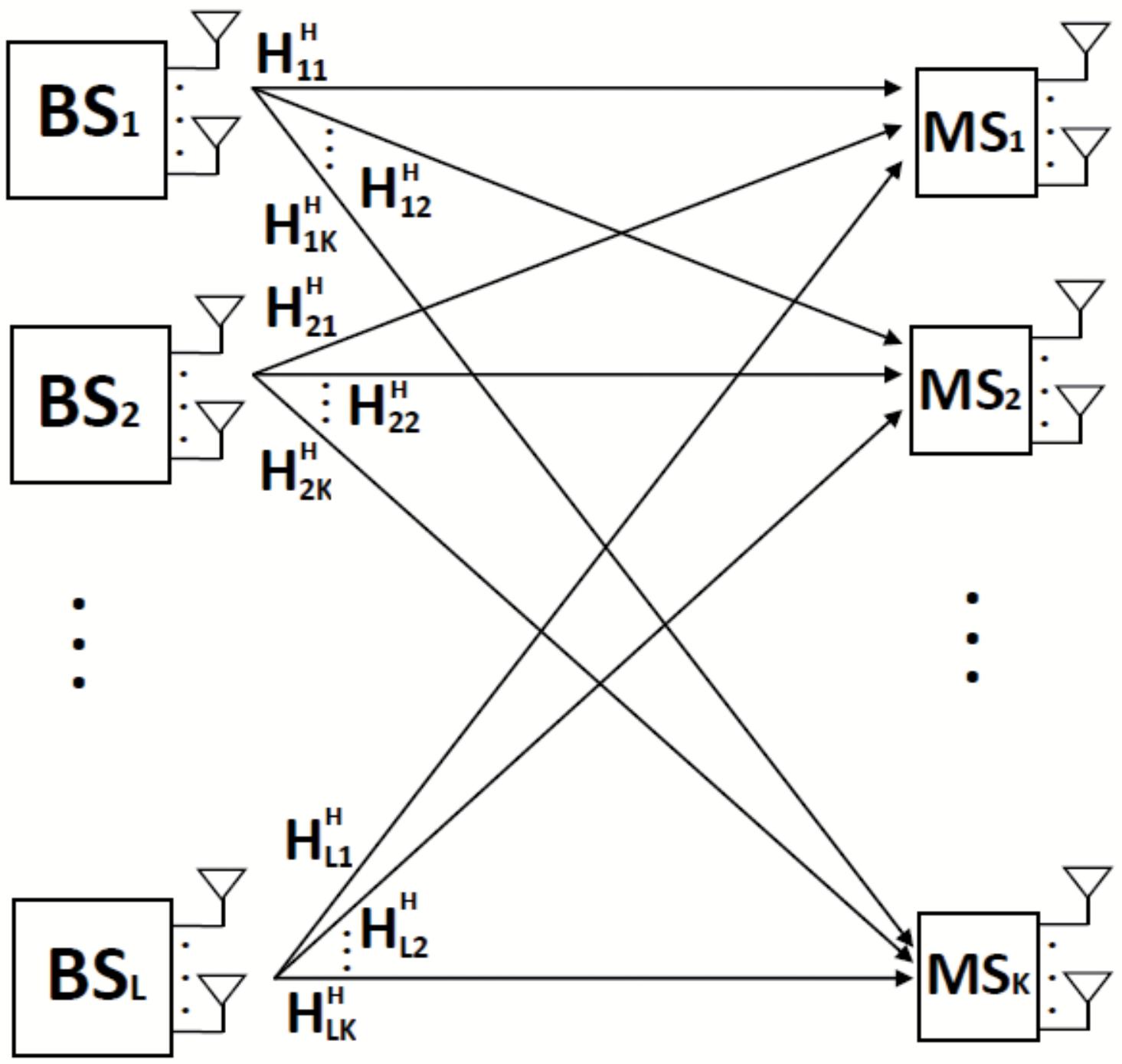

**Figure 3.1** MIMO Coordinated base station system model.

It is clearly seen that the last expression of (3.1) has exactly the same form as the estimate of $d_{ki}$ for the conventional downlink multiuser MIMO system where a BS equipped with $N$ transmit antennas is serving $K$ decentralized multiantenna MSs. Hence, we can interpret a coordinated BS system as one giant downlink system [SSVB08, TSC07]. It is assumed that the entries of $\mathbf{n}_k$ are i.i.d ZMCSCG random variables with the variance $\sigma_k^2$, i.e., $\mathbf{n}_k \sim \mathcal{NC}(\mathbf{0}, \sigma_k^2 \mathbf{I}_{M_k})$. We also assume that the symbol $\mathbf{d}_k$ consists of ZMCSCG random variables with unit variance and is independent of $\{\mathbf{d}_i\}_{i=1,i\neq k}^K$ and noise $\mathbf{n}_k$, i.e., $\mathrm{E}\{\mathbf{d}_k \mathbf{d}_k^H\} = \mathbf{I}_{S_k}$, $\mathrm{E}\{\mathbf{d}_k \mathbf{d}_i^H\} = \mathbf{0}$, $\forall i \neq k$ and $\mathrm{E}\{\mathbf{d}_k \mathbf{n}_k^H\} = \mathbf{0}$. For this system model, the MSE between $d_{ki}$ and $\hat{d}_{ki}$ is given by

$$\begin{aligned}\xi_{ki} =& \mathrm{E}_{\mathbf{d}}\{(\hat{d}_{ki} - d_{ki})(\hat{d}_{ki} - d_{ki})^H\} \\ =& \mathbf{w}_{ki}^H (\mathbf{H}_k^H \mathbf{B}\mathbf{B}^H \mathbf{H}_k + \sigma_k^2 \mathbf{I}_{M_k}) \mathbf{w}_{ki} - \mathbf{w}_{ki}^H \mathbf{H}_k^H \mathbf{b}_{ki} - \mathbf{b}_{ki}^H \mathbf{H}_k \mathbf{w}_{ki} + 1. \end{aligned} \quad (3.2)$$

For notational convenience, we represent $[\xi_{11}, \cdots, \xi_{1S_1}, \cdots, \xi_{K1}, \cdots, \xi_{KS_K}]$ by $[\xi_1, \xi_2, \cdots, \xi_S]$, $[\mathbf{w}_{11}, \cdots, \mathbf{w}_{1S_1}, \cdots, \mathbf{w}_{K1}, \cdots, \mathbf{w}_{KS_K}]$ by $[\mathbf{w}_1, \mathbf{w}_2, \cdots, \mathbf{w}_S]$, $\mathbf{B} = [\mathbf{b}_1, \cdots, \mathbf{b}_S]$, and the channel matrix and noise variance corresponding to the $s$th symbol is denoted by $\tilde{\mathbf{H}}_s$ and $\tilde{\sigma}_s^2$, respectively[6]. By doing so, the MSE of the $s$th symbol is given by

$$\xi_s = \mathbf{w}_s^H (\tilde{\mathbf{H}}_s^H \mathbf{B}\mathbf{B}^H \tilde{\mathbf{H}}_s + \tilde{\sigma}_s^2 \mathbf{I}) \mathbf{w}_s - \mathbf{w}_s^H \tilde{\mathbf{H}}_s^H \mathbf{b}_s - \mathbf{b}_s^H \tilde{\mathbf{H}}_s \mathbf{w}_s + 1. \quad (3.3)$$

When perfect CSI is available at the BS and MSs, the rate of the $s$th symbol is given by

$$R_s = \log_2 (1 + SINR_s) \quad (3.4)$$

where

$$SINR_s = \frac{\mathbf{w}_s^H \tilde{\mathbf{H}}_s^H \mathbf{b}_s \mathbf{b}_s \tilde{\mathbf{H}}_s \mathbf{w}_s}{\mathbf{w}_s^H (\tilde{\mathbf{H}}_s^H \sum_{i=1,i\neq s}^S \mathbf{b}_i \mathbf{b}_i \tilde{\mathbf{H}}_s + \tilde{\sigma}_s^2) \mathbf{w}_s}.$$

[6]Note that $\tilde{\mathbf{H}}_s$ and $\tilde{\sigma}_s^2$ are the same as the channel and noise variance of the MS associated with the $s$th symbol, respectively.

The MMSE receiver of the $s$th symbol is given as

$$\mathbf{w}_s = (\tilde{\mathbf{H}}_s^H \mathbf{B}\mathbf{B}^H \tilde{\mathbf{H}}_s + \tilde{\sigma}_s^2 \mathbf{I})^{-1} \tilde{\mathbf{H}}_s^H \mathbf{b}_s. \tag{3.5}$$

Plugging this equation into (3.3), we get the MMSE of the $s$th symbol as

$$\tilde{\tilde{\xi}}_s = 1 - \mathbf{b}_s^H \tilde{\mathbf{H}}_s (\tilde{\mathbf{H}}_s^H \mathbf{B}\mathbf{B}^H \tilde{\mathbf{H}}_s + \tilde{\sigma}_s^2 \mathbf{I})^{-1} \tilde{\mathbf{H}}_s^H \mathbf{b}_s. \tag{3.6}$$

When each of the symbols ($\{d_s\}_{s=1}^S$) is decoded individually independent of each other using a minimum Euclidean distance decoding rule, the achievable rate of the $s$th symbol can be expressed as [SSB08b, SSB08c, Lap96]

$$R_s = \log_2 (1 + SINR_s) = \log_2 (\tilde{\tilde{\xi}}_s)^{-1} \tag{3.7}$$

## 3.3 Problem formulation

Mathematically, the weighted sum rate maximization problem can be formulated as

$$\mathcal{P}3.1: \max_{\{\mathbf{b}_s\}_{s=1}^S} \sum_{s=1}^S \omega_s \log_2 (\tilde{\tilde{\xi}}_s)^{-1}, \text{ s.t } [\sum_{s=1}^S \mathbf{b}_s \mathbf{b}_s^H]_{n,n} \leq p_n, \forall n \tag{3.8}$$

where $p_n$ is the maximum allocated power to the $n$th BS antenna and $\omega_s \geq 0$ is the rate weighting factor of the $s$th symbol. The antenna numbers are assigned from the first antenna of $\text{BS}_1$ (which corresponds to antenna 1) to the last antenna of $\text{BS}_L$ (which corresponds to antenna $N$).

In a multimedia communication, different types of information (for example, audio and video information) can be sent to a user simultaneously [SSP01]. In such a case, for successful transmission, more priority could be given to symbols corresponding to the video information. Consequently, the symbols of a user (all users) can have different priorities. This motivates us to examine the joint transmitter and receiver design for symbol wise weighted sum rate maximization problem. However, as will be clear later, the proposed algorithms can be applied to get the suboptimal solution for user wise weighted (un-weighted) sum rate optimization problem. Without loss of generality, we assume that $\{0 < \omega_s < 1\}_{s=1}^S$. After straightforward mathematical manipulations, problem (3.8) can be equivalently expressed as

$$\min_{\{\mathbf{b}_s\}_{s=1}^S} \prod_{s=1}^S \tilde{\tilde{\xi}}_s^{\omega_s}, \text{ s.t } [\sum_{s=1}^S \mathbf{b}_s \mathbf{b}_s^H]_{n,n} \leq p_n, \ \forall n. \tag{3.9}$$

Note that although (3.8) and (3.9) are equivalent problems, the optimal (sub-optimal) values of these problems are not necessarily equal. Solving the latter problem in its current form has appeared to be intractable. Due to this, we introduce the receivers $\{\mathbf{w}_s\}_{s=1}^S$ and then reformulate the above problem as (see Appendix 3.A)

$$\min_{\{\mathbf{b}_s,\mathbf{w}_s\}_{s=1}^S} \prod_{s=1}^S \tilde{\xi}_s, \text{ s.t } [\sum_{s=1}^S \mathbf{b}_s\mathbf{b}_s^H]_{n,n} \leq p_n, \ \forall n \tag{3.10}$$

where $\tilde{\xi}_s = \xi_s^{\omega_s}$.

## 3.4 Existing and proposed solutions

The above optimization problem is non-convex. Thus, convex optimization techniques can not be applied. In this section, we present the exiting centralized, and the proposed centralized and distributed algorithms for (3.10). This problem has been examined in [SSB08b] for the case where the power of each symbol is strictly positive. The paper proposes an iterative algorithm that achieves a local optimum to (3.10). Assuming $\{M_k = \tilde{M}\}_{k=1}^K$, the complexity of each iteration is given by

$$\mathfrak{C}_e = O(\sqrt{(N+S)}(2NS+1)^2(2S^2+2NS+S)) + O(K\tilde{M}^{2.376}) + \mathfrak{C}_{GP} \tag{3.11}$$

where $\mathfrak{C}_{GP}$ is the complexity of the GP problem of [SSB08b]. In general, $\mathfrak{C}_{GP}$ depends on different parameter settings and solution methods (see Appendix 3.B for the details).

In the following, we present our novel centralized and distributed iterative algorithms that achieve local optimum to (3.10). The proposed algorithms require less computational cost per iteration than that of the algorithm in [SSB08b] (i.e., $\mathfrak{C}_e$). As will be clear later in Section VI, the proposed algorithms also have faster convergence speed than that of the algorithm in [SSB08b]. As a result, our algorithms require less overall computational complexities than that of the algorithm in [SSB08b]. To this end, we consider the following Lemma.

*Lemma 3.1*: The optimal/suboptimal $\{\mathbf{b}_s, \mathbf{w}_s\}_{s=1}^S$ of (3.10) can be obtained by solving the following problem.

$$\min_{\{\mathbf{b}_k, \mathbf{w}_s, \nu_s\}_{s=1}^S} \left(\frac{1}{S}\sum_{s=1}^S \nu_s \tilde{\tilde{\xi}}_s\right)^S, \ \text{s.t} \ [\sum_{s=1}^S \mathbf{b}_s \mathbf{b}_s^H]_{n,n} \leq p_n, \ \prod_{s=1}^S \nu_s = 1, \nu_s \geq 0, \ \forall s, n. \tag{3.12}$$

Moreover, for fixed $\{\mathbf{b}_s\}_{s=1}^S$, the optimal $\{\nu_s\}_{s=1}^S$ of this problem is given by

$$\nu_s^\star = \frac{\left[\prod_{i=1}^S \tilde{\tilde{\xi}}_i^{\omega_i}\right]^{\frac{1}{S}}}{\tilde{\tilde{\xi}}_s^{\omega_s}}, \ \forall s. \tag{3.13}$$

**Proof.** See Appendix 3.C. □

From the proof of this Lemma, one can realize that the optimal/suboptimal solution of (3.12) satisfies $\frac{1}{S}\sum_{s=1}^S \tilde{\xi}_s \nu_s > 0$ and $\{\nu_s > 0\}_{s=1}^S$. Thus, the objective function of (3.12) can be replaced by $\sum_{s=1}^S \nu_s \tilde{\tilde{\xi}}_s = \sum_{s=1}^S \nu_s \xi_s^{\omega_s}$. This is due to the fact that $\min_x (cf(x))^N$ is equivalent to $\min_x f(x)$ for any $c > 0$, $f(x) > 0, \forall x$ and positive integer $N$ [BV04]. Consequently, (3.12) can be equivalently expressed as

$$\min_{\{\nu_s, \mathbf{b}_s, \mathbf{w}_s\}_{s=1}^S} \sum_{s=1}^S \nu_s \xi_s^{\omega_s}, \ \text{s.t} \ [\sum_{s=1}^S \mathbf{b}_s \mathbf{b}_s^H]_{n,n} \leq p_n, \ \prod_{s=1}^S \nu_s = 1, \nu_s > 0, \ \forall n, s. \tag{3.14}$$

Due to $\{\xi_s^{\omega_s}\}_{s=1}^S$ terms in the objective function of (3.14), getting the suboptimal solution of this problem is not trivial. To simplify the latter problem, we present *Lemma 3.2*.

*Lemma 3.2*: For any strictly positive real numbers $a$ and $b$, and $0 < \omega < 1$, the following holds true

$$F = \min_{\tau > 0} \kappa\left(\frac{a^\gamma}{\tau} + b\tau^\mu\right) = ab^\omega \tag{3.15}$$

where $\gamma = \frac{1}{1-\omega}$, $\mu = \frac{1}{\omega} - 1$ and $\kappa = \omega\mu^{(1-\omega)}$.

**Proof.** The optimal $\tau$ of $F$ can be obtained by using the first order derivative of $F$ with respect to $\tau$ as

$$\frac{\mathrm{d}F}{\mathrm{d}\tau} = \kappa\left(-\frac{a^\gamma}{\tau^2} + \mu b \tau^{\frac{1}{\omega}-2}\right) = 0 \Rightarrow \tau = \frac{a^{\omega\gamma}}{(b\mu)^\omega}. \tag{3.16}$$

Substituting (3.16) in $F$, we get

$$F = \kappa(\frac{a^{\gamma}}{\tau} + b\tau^{\mu}) = \frac{\kappa}{\omega\mu}\left[a^{(1-\omega)\gamma}(b\mu)^{\omega}\right] = \frac{\kappa}{\omega\mu^{(1-\omega)}}\left[ab^{\omega}\right] = ab^{\omega}. \quad (3.17)$$

□

Following *Lemma 3.2*, it can be shown that $\{\mathbf{b}_s, \mathbf{w}_s, \nu_s\}_{s=1}^{S}$ of (3.14) can be optimized by solving (3.18).

$$\min_{\{\tau_s,\nu_s,\mathbf{b}_s,\mathbf{w}_s\}_{s=1}^{S}} \sum_{s=1}^{S} \kappa_s \left[\frac{1}{\tau_s}\nu_s^{\gamma_s} + \xi_s\tau_s^{\mu_s}\right]$$
$$\text{s.t } [\sum_{s=1}^{S} \mathbf{b}_s\mathbf{b}_s^H]_{n,n} \leq p_n, \ \prod_{s=1}^{S} \nu_s = 1, \nu_s > 0, \ \tau_s > 0 \ \forall s, n \quad (3.18)$$

where $\gamma_s = \frac{1}{1-\omega_s}$, $\mu_s = \frac{1}{\omega_s} - 1$ and $\kappa_s = \omega_s\mu_s^{(1-\omega_s)}$. The above problem is non-convex. Thus, convex optimization can not be applied. Next, we present our centralized and distributed iterative algorithms that achieve local optimum to this problem.

### 3.4.1 Proposed centralized algorithm

The key step of this centralized algorithm is the utilization of *Lemma 3.1* and *Lemma 3.2* which help us to transform the intractable problem (3.10) to a more convenient problem (3.18). In this subsection, we present our centralized iterative algorithm for the latter problem as follows. First, keeping the precoders of all symbols $\{\mathbf{b}_s\}_{s=1}^{S}$ constant, the optimal $\{\mathbf{w}_s\}_{s=1}^{S}$ can be obtained by using MMSE receiver approach (3.5) and $\{\nu_s, \tau_s\}_{s=1}^{S}$ are optimized by solving the following problem

$$\min_{\{\tau_s,\nu_s\}_{s=1}^{S}} \sum_{s=1}^{S} \kappa_s \left[\frac{1}{\tau_s}\nu_s^{\gamma_s} + \tilde{\tilde{\xi}}_s\tau_s^{\mu_s}\right], \ \text{s.t} \prod_{s=1}^{S} \nu_s = 1, \nu_s > 0, \ \tau_s > 0, \ \forall s. \quad (3.19)$$

The above optimization problem is a GP for which global optimal solution can be obtained by existing optimization tools [MKKB06]. However, here we provide closed form expressions for the optimal $\{\nu_s, \tau_s\}_{s=1}^{S}$ of this problem. For fixed $\{\nu_s\}_{s=1}^{S}$, the optimal $\{\tau_s\}_{s=1}^{S}$ of (3.19) can be obtained by applying the first order derivative of (3.19) with respect to $\{\tau_i\}_{i=1}^{S}$ and are given as

$$\tau_s^{\star} = \left[\frac{\nu_s^{\gamma_s}}{\mu_s\tilde{\tilde{\xi}}_s}\right]^{\frac{1}{\mu_s+1}}, \ \forall s. \quad (3.20)$$

Substituting these $\{\tau_s^\star\}_{s=1}^S$ back into the objective function of (3.19) we get

$$\min_{\{\nu_s\}_{s=1}^S} \sum_{s=1}^S \nu_s \xi_s^{\tilde{\tilde{\omega}}_s}, \quad \text{s.t} \prod_{s=1}^S \nu_s = 1, \nu_s > 0, \ \ \forall s. \tag{3.21}$$

It can be easily seen that the latter problem and (3.41) has the same optimal solution. Thus, the global optimal $\{\nu_s\}_{s=1}^S$ of (3.21) is given by (3.13). Second, for fixed $\{\mathbf{w}_s, \nu_s, \tau_s\}_{s=1}^S$, the optimal $\{\mathbf{b}_s\}_{s=1}^S$ of (3.18) can be obtained by solving the following problem

$$\min_{\{\mathbf{b}_s\}_{s=1}^S} \sum_{s=1}^S \eta_s \xi_s, \quad \text{s.t} \ [\sum_{s=1}^S \mathbf{b}_s \mathbf{b}_s^H]_{n,n} \leq p_n, \ \ \forall n \tag{3.22}$$

where $\eta_s = \kappa_s \tau_s^{\mu_s}$. The objective function of the above problem can be expressed as

$$\begin{aligned}
&\sum_{s=1}^S \eta_s \xi_s \\
=&\sum_{s=1}^S \left\{ \mathbf{b}_s^H (\sum_{i=1}^S \eta_i \tilde{\mathbf{H}}_i \mathbf{w}_i \mathbf{w}_i^H \tilde{\mathbf{H}}_i^H) \mathbf{b}_s - \eta_s \mathbf{w}_s^H \tilde{\mathbf{H}}_s^H \mathbf{b}_s - \eta_s \mathbf{b}_s^H \tilde{\mathbf{H}}_s \mathbf{w}_s + \tilde{\sigma}_s^2 \eta_s \mathbf{w}_s^H \mathbf{w}_s + \eta_s \right\} \\
=&\text{tr}\{(\sqrt{\boldsymbol{\eta}} \tilde{\mathbf{W}}^H \mathbf{H}^H \mathbf{B} - \sqrt{\boldsymbol{\eta}})^H (\sqrt{\boldsymbol{\eta}} \tilde{\mathbf{W}}^H \mathbf{H}^H \mathbf{B} - \sqrt{\boldsymbol{\eta}})\} + \text{tr}\{\boldsymbol{\eta} \tilde{\mathbf{W}}^H \boldsymbol{\sigma}^2 \tilde{\mathbf{W}}\}
\end{aligned} \tag{3.23}$$

where $\boldsymbol{\eta} = \text{diag}(\eta_1, \cdots, \eta_S)$, $\boldsymbol{\sigma}^2 = \text{blkdiag}(\sigma_1^2 \mathbf{I}_{M_1}, \cdots, \sigma_K^2 \mathbf{I}_{M_K})$, $\mathbf{H} = [\mathbf{H}_1, \cdots, \mathbf{H}_K]$, $\tilde{\mathbf{W}}_k$ as the decoder matrix of the $k$th MS and $\tilde{\mathbf{W}} = \text{blkdiag}(\tilde{\mathbf{W}}_1, \cdots, \tilde{\mathbf{W}}_K)$. By applying (3.23), problem (3.22) can be reexpressed as

$$\min_{\chi, \{\mathbf{b}_s\}_{s=1}^S} \chi \ \text{ s.t} \ \|[\text{vec}(\sqrt{\boldsymbol{\eta}} \tilde{\mathbf{W}}^H \mathbf{H}^H \mathbf{B} - \sqrt{\boldsymbol{\eta}})]\|_2 \leq \chi, \ \ \|\tilde{\mathbf{b}}_n\|_2 \leq \sqrt{p_n}, \ \forall n \tag{3.24}$$

where $\tilde{\mathbf{b}}_n^H$ as the $n$th row of $\mathbf{B}$. As we can see, (3.24) is a SOCP problem for which the global optimal solution is obtained by existing convex optimization tools [BV04]. Finally, the first and second steps are repeated iteratively until convergence is achieved. Our centralized iterative algorithm that achieves a local optimum to (3.8) is summarized as shown in **Algorithm 3.I**.

**Algorithm 3.I**: Centralized iterative algorithm for problem (3.8)

Initialization: Set $\{\mathbf{B}_k\}_{k=1}^K$ as the first $S_k$ vectors of $\{\mathbf{H}_k\}_{k=1}^K$ and the maximum number of iterations as $i_{max}$. Then, normalize $\{\mathbf{B}_k\}_{k=1}^K$ such that each BS antenna power constraint is satisfied with equality.

**Repeat**

1. With the current $\{\mathbf{b}_s\}_{s=1}^S$, $\{\mathbf{w}_s, \nu_s$ and $\tau_s\}_{s=1}^S$ are updated using (3.5), (3.13) and (3.20), respectively.
2. With the current $\{\nu_s, \mathbf{w}_s, \tau_s\}_{s=1}^S$, $\{\mathbf{b}_s\}_{s=1}^S$ are optimized by solving (3.24).
3. Compute the objective function of (3.8).

**Until** convergence.

**Convergence**:- For fixed $\{\mathbf{w}_s, \nu_s$ and $\tau_s\}_{s=1}^S$, the global minimum of (3.18) can be achieved by optimizing $\{\mathbf{b}_s\}_{s=1}^S$ with (3.24). Moreover, for fixed $\{\mathbf{b}_s\}_{s=1}^S$, the global minimum of (3.18) can be achieved by optimizing $\{\mathbf{w}_s, \nu_s$ and $\tau_s\}_{s=1}^S$ with (3.5), (3.13) and (3.20), respectively. As a result, $\mathcal{F}_n^2 \leq \mathcal{F}_n^1$ is satisfied, where $\mathcal{F}_n^i$ is the objective function of (3.18) at step $i$ of the $n$th iteration [SSB08c], [ECV10]. At the $(n+1)$th iteration, we achieve $\mathcal{F}_{(n+1)}^2 \leq \mathcal{F}_{(n+1)}^1 \leq \mathcal{F}_n^2$. These discussions reveal the fact that the objective function of (3.18) is non-increasing at each step. Which implies that the objective function of (3.8) also non-decreasing. On the other hand, the objective function of the latter problem is upper bounded by a positive finite value. These two facts show that the proposed iterative algorithm is always guaranteed to converge. However, since problem (3.8) is non-convex, we are not able to show the global optimality of **Algorithm 3.I** analytically.

**Initialization:-** In general, different initializations affect the convergence speed and optimal weighted sum rate of **Algorithm 3.I**. In most of our simulations, we observe faster convergence speed and better weighted sum rate when the initialization is performed as in **Algorithm 3.1**. Nonetheless, getting the best initialization that results the fastest convergence speed and best weighted sum rate of **Algorithm 3.I** is an open research topic.

**Computational complexity:-** The main computational load of **Algorithm 3.I** arises from solving (3.5) and (3.24). For the assumption discussed in Section 3.4, (3.5) can be performed with $O(K\tilde{M}^{2.376})$ [CW90]. It can be shown that problem (3.24) has $N$ SOC constraints where each of them consists of $2S$ real dimensions, one SOC constraint with $2S^2$ real dimensions and $2NS+1$ real optimization variables. According to [LVBL98] (see page 196 of [LVBL98]), the computational complexity of the latter problem in terms of number of iterations is upper bounded by $O(\sqrt{N+1})$ where the complexity of each iter-ation is on the order of $O((2SN+1)^2(2S^2+2SN))$. Thus, the total computational complexity of (3.24) is given by $O(\sqrt{N+1}(2SN+1)^2(2S^2+2SN))$. Therefore, in one iteration, **Algorithm 3.I** requires $\mathfrak{C}_c = O(\sqrt{N+1}(2SN+1)^2(2S^2+2SN)) + O(K\tilde{M}^{2.376})$ operations. This shows that our proposed centralized algorithm requires less computational cost per iteration than that of the algorithm in [SSB08b] (i.e., $\mathfrak{C}_c < \mathfrak{C}_e$).

However, although $\mathfrak{C}_c < \mathfrak{C}_e$, we still believe that for large-scale networks, $\mathfrak{C}_c$ is very large computational load and hence it is not suitable for practical realization. This motivates us to develop a distributed algorithm that achieves a local optimum to (3.8) with less computational cost than that of our centralized algorithm.

### 3.4.2 Proposed distributed algorithm

We have shown in the previous subsection that (3.8) can be solved equivalently by using (3.18). As can be seen from **Algorithm 3.I**, the optimal $\{\mathbf{w}_s, \nu_s, \tau_s\}_{s=1}^S$ of (3.18) can be solved independently and distributively. However, the optimal solution of (3.22) is computed using a centralized algorithm. In this subsection, we present our distributed algorithm for (3.22). The Lagrangian dual decomposition technique is applied to solve this problem distributively[7]. To this end, we first express the Lagrangian function associated with (3.22) as

$$\begin{aligned}
&L(\boldsymbol{\lambda}, \mathbf{B}) \\
&= \sum_{s=1}^{S} \eta_s \xi_s + \sum_{n=1}^{N} \lambda_n \left( [\sum_{i=1}^{S} \mathbf{b}_i \mathbf{b}_i^H]_{n,n} - p_n \right) \\
&= \sum_{s=1}^{S} \left\{ \mathbf{b}_s^H (\sum_{i=1}^{S} \eta_i \tilde{\mathbf{H}}_i \mathbf{w}_i \mathbf{w}_i^H \tilde{\mathbf{H}}_i^H) \mathbf{b}_s - \eta_s \mathbf{w}_s^H \tilde{\mathbf{H}}_s^H \mathbf{b}_s - \eta_s \mathbf{b}_s^H \tilde{\mathbf{H}}_s \mathbf{w}_s + \tilde{\sigma}_s^2 \eta_s \mathbf{w}_s^H \mathbf{w}_s + \eta_s \right\} \\
&\quad + \sum_{n=1}^{N} \lambda_n \left( [\sum_{i=1}^{S} \mathbf{b}_i \mathbf{b}_i^H]_{n,n} - p_n \right)
\end{aligned}$$

[7] Since (3.22) is convex and Slater's condition (i.e., the existence of strictly feasible points) is satisfied by choosing any $\{\mathbf{b}_s\}_{s=1}^S$ with $\{[\sum_{s=1}^S \mathbf{b}_s \mathbf{b}_s^H]_{n,n} < p_n\}_{n=1}^N$, the duality gap between the primal problem (3.22) and its dual problem is zero [BV04].

$$= \sum_{s=1}^{S} \left\{ \mathbf{b}_s^H \mathbf{A} \mathbf{b}_s - \eta_s \mathbf{w}_s^H \tilde{\mathbf{H}}_s^H \mathbf{b}_s - \eta_s \mathbf{b}_s^H \tilde{\mathbf{H}}_s \mathbf{w}_s + \tilde{\sigma}_s^2 \eta_s \mathbf{w}_s^H \mathbf{w}_s + \eta_s \right\} - \sum_{n=1}^{N} \lambda_n p_n \tag{3.25}$$

where $\boldsymbol{\lambda} = \text{diag}(\lambda_1, \cdots, \lambda_N)$ are the Lagrangian multipliers corresponding to the constraint sets of (3.22) and $\mathbf{A} = \sum_{i=1}^{S} \eta_i \tilde{\mathbf{H}}_i \mathbf{w}_i \mathbf{w}_i^H \tilde{\mathbf{H}}_i^H + \boldsymbol{\lambda}$. Thus, the dual function of (3.22) is

$$\begin{aligned} g(\boldsymbol{\lambda}) &= \min_{\{\mathbf{b}_s\}_{s=1}^{S}} L(\boldsymbol{\lambda}, \mathbf{B}) \\ &= \sum_{s=1}^{S} \left\{ \mathbf{b}_s^H \mathbf{A} \mathbf{b}_s - \eta_s \mathbf{w}_s^H \tilde{\mathbf{H}}_s^H \mathbf{b}_s - \eta_s \mathbf{b}_s^H \tilde{\mathbf{H}}_s \mathbf{w}_s + \tilde{\sigma}_s^2 \eta_s \mathbf{w}_s^H \mathbf{w}_s + \eta_s \right\} - \sum_{n=1}^{N} \lambda_n p_n \\ &= \sum_{s=1}^{S} \left\{ \eta_s (\tilde{\sigma}_s^2 \mathbf{w}_s^H \mathbf{w}_s + 1) - \eta_s^2 \mathbf{w}_s^H \tilde{\mathbf{H}}_s^H \mathbf{A}^{-1} \tilde{\mathbf{H}}_s \mathbf{w}_s \right\} - \sum_{n=1}^{N} \lambda_n p_n \end{aligned} \tag{3.26}$$

where the third equality is obtained after substituting the optimal $\mathbf{b}_s$ of (3.26) which is given by

$$\mathbf{b}_s^\star = \eta_s \mathbf{A}^{-1} \tilde{\mathbf{H}}_s \mathbf{w}_s, \Rightarrow \mathbf{b}_{ls}^\star = \eta_s [\mathbf{A}^{-1}]_l \tilde{\mathbf{H}}_s \mathbf{w}_s \; \forall l, s \tag{3.27}$$

where $[\mathbf{A}^{-1}]_l \in \mathbb{C}^{N_l \times N}$ is obtained by $[\mathbf{A}^{-1}]_{(F_l:F_l+N_l-1,:)}$ with $F_l = \sum_{i=0}^{l-1} N_i + 1$ and $N_0 = 0$. As can be seen from (3.27), for a given $\boldsymbol{\lambda}$, the precoder vector of each symbol can be optimized independently. The optimal $\boldsymbol{\lambda}$ of (3.25) can be obtained by solving the dual optimization problem of (3.22) which is given by

$$\begin{aligned} &\max_{\{\lambda_n \geq 0\}_{n=1}^{N}} g(\boldsymbol{\lambda}) = \\ &\max_{\{\lambda_n \geq 0\}_{n=1}^{N}} \sum_{s=1}^{S} \left\{ \eta_s (\tilde{\sigma}_s^2 \mathbf{w}_s^H \mathbf{w}_s + 1) - \eta_s^2 \mathbf{w}_s^H \mathbf{H}_s^H \mathbf{A}^{-1} \tilde{\mathbf{H}}_s \mathbf{w}_s \right\} - \sum_{n=1}^{N} \lambda_n p_n. \end{aligned} \tag{3.28}$$

By employing the eigenvalue decompositions of $\mathbf{H}\tilde{\mathbf{W}}\boldsymbol{\eta}^2\tilde{\mathbf{W}}^H\mathbf{H}^H \triangleq \tilde{\mathbf{V}}\tilde{\boldsymbol{\Lambda}}\tilde{\mathbf{V}}^H$ and $\mathbf{H}\tilde{\mathbf{W}}\boldsymbol{\eta}\tilde{\mathbf{W}}^H\mathbf{H}^H \triangleq \bar{\mathbf{V}}\bar{\boldsymbol{\Lambda}}\bar{\mathbf{V}}^H$, problem (3.28) can be written as

$$\min_{\{\lambda_n \geq 0\}_{n=1}^{N}} \text{tr}\left\{ \mathbf{F}^H (\mathbf{R}\mathbf{R}^H + \boldsymbol{\lambda})^{-1} \mathbf{F} \right\} + \sum_{n=1}^{N} \lambda_n p_n \tag{3.29}$$

where $\mathbf{F} = \tilde{\mathbf{V}}\sqrt{\tilde{\boldsymbol{\Lambda}}}$ and $\mathbf{R} = \bar{\mathbf{V}}\sqrt{\bar{\boldsymbol{\Lambda}}}$. The above optimization problem can be cast as an SDP problem where the global solution can be found by existing convex optimization tools [BV04]. The computational complexity of this problem is on the order of $O((2N^2 + N)^2 (4N)^{2.5})$ [LVBL98]. However, here our aim is to obtain the optimal values of $\{\lambda_n\}_{n=1}^N$ distributively with less computational load than that of the SDP method. In this regard, we present the following Lemma.

*Lemma 3.3*: The optimal $\{\lambda_n\}_{n=1}^N$ of the above optimization problem can be obtained by solving the following problem

$$\min_{\{\lambda_n,\mathbf{g}_n,\mathbf{t}_n\}_{n=1}^N} \sum_{n=1}^N \{\mathbf{g}_n^H \boldsymbol{\lambda}^{-1}\mathbf{g}_n + \mathbf{t}_n^H\mathbf{t}_n + \lambda_n p_n\} \triangleq \varphi$$
$$\text{s.t } \mathbf{R}\mathbf{t}_n + \mathbf{g}_n = \mathbf{f}_n, \ \forall n \tag{3.30}$$

where $\mathbf{f}_n$ is the $n$th column of $\mathbf{F}$.

**Proof.** By keeping $\boldsymbol{\lambda}$ constant, the Lagrangian function of (3.30) is given by

$$L = \sum_{n=1}^N \{\mathbf{g}_n^H \boldsymbol{\lambda}^{-1}\mathbf{g}_n + \mathbf{t}_n^H\mathbf{t}_n + \lambda_n p_n - \boldsymbol{\psi}_n^H(\mathbf{R}\mathbf{t}_n + \mathbf{g}_n - \mathbf{f}_n)\}$$

where $\boldsymbol{\psi}_n^H$ is the Lagrangian multiplier associated with the $n$th equality constraint of (3.30). Differentiation of $L$ with respect to $\{\mathbf{g}_i, \mathbf{t}_i\}_{i=1}^N$ yield $\{\mathbf{g}_i^\star = \boldsymbol{\lambda}\boldsymbol{\psi}_i\}_{i=1}^N$ and $\{\mathbf{t}_i^\star = \mathbf{R}^H\boldsymbol{\psi}_i\}_{i=1}^N$. By substituting these $\{\mathbf{g}_i^\star, \mathbf{t}_i^\star\}_{i=1}^N$ in the equality constraint of (3.30), we get $\{\boldsymbol{\psi}_i = (\mathbf{R}\mathbf{R}^H + \boldsymbol{\lambda})^{-1}\mathbf{f}_i\}_{i=1}^N$. It follows

$$\mathbf{g}_i^\star = \boldsymbol{\lambda}(\mathbf{R}\mathbf{R}^H + \boldsymbol{\lambda})^{-1}\mathbf{f}_i, \ \mathbf{t}_i^\star = \mathbf{R}^H(\mathbf{R}\mathbf{R}^H + \boldsymbol{\lambda})^{-1}\mathbf{f}_i, \ \forall i. \tag{3.31}$$

Plugging (3.31) into the objective function of (3.30) yields

$$\varphi = \sum_{i=1}^N \{\mathbf{g}_i^H \boldsymbol{\lambda}^{-1}\mathbf{g}_i + \mathbf{t}_i^H\mathbf{t}_i + \lambda_i p_i\} = \sum_{i=1}^N \{\mathbf{f}_i^H(\mathbf{R}\mathbf{R}^H + \boldsymbol{\lambda})^{-1}\mathbf{f}_i\} + \sum_{i=1}^N \lambda_i p_i$$
$$= \text{tr}\left\{\mathbf{F}^H(\mathbf{R}\mathbf{R}^H + \boldsymbol{\lambda})^{-1}\mathbf{F}\right\} + \sum_{i=1}^N \lambda_i p_i. \tag{3.32}$$

The above equation is the same as the objective function of the original optimization problem (3.29). It follows that (3.29) and (3.30) are equivalent problems. Note that *Lemma 3.3* is proved by modifying the idea of matrix fractional minimization (see [BVC10] and [BV04]). It can be shown that (3.30) is a convex optimization problem [BV04]. □

To develop distributed algorithm for (3.30), we reexpress $\mathbf{G} = [\mathbf{g}_1, \cdots, \mathbf{g}_N]$ as $\mathbf{G} = [\bar{\mathbf{g}}_1^H; \cdots; \bar{\mathbf{g}}_N^H]$, where $\bar{\mathbf{g}}_i^H$ is the $i$th row of $\mathbf{G}$. By doing so, $\mathbf{G}^\star =$ $[\mathbf{g}_1^\star, \cdots, \mathbf{g}_N^\star]$ of (3.31) can also be written as $\mathbf{G}^\star = [(\overline{\mathbf{g}}_1^\star)^H; \cdots; (\overline{\mathbf{g}}_N^\star)^H]$, where

$$\overline{\mathbf{g}}_i^\star = \lambda_i \mathbf{\Gamma}_i^H, \ \forall i \tag{3.33}$$

and $\mathbf{\Gamma}_i$ is the $i$th row of $\mathbf{\Gamma} = \mathbf{A}^{-1}\mathbf{F}$.

Now, problem (3.30) can be solved distributively as follows. First, keeping $\boldsymbol{\lambda}$ constant, the optimal $\overline{\mathbf{g}}_i^r$ can be computed independently using (3.33), i.e., $\overline{\mathbf{g}}_i^{r^\star} = \lambda_i^{r-1}(\mathbf{\Gamma}_i^{r-1})^H$, where the superscripts $(.)^r$ and $(.)^{r-1}$ represent the current and previous values, respectively. Then, $\lambda_i^{r^\star}$ is computed by

$$\frac{\partial \varphi}{\lambda_i} = -\frac{1}{\lambda_i^2}\beta_i + p_i = 0 \Rightarrow \lambda_i^{r^\star} = \sqrt{\beta_i^r / p_i}, \ \forall i \tag{3.34}$$

where $\beta_i^r = (\overline{\mathbf{g}}_i^{r^\star})^H \overline{\mathbf{g}}_i^{r^\star}$. As we can see from the above expression $\lambda_i^\star$ is always non-negative. Furthermore, from (3.33) and (3.34), one can observe that $\lambda_i^\star$ can be updated in parallel by using only $\overline{\mathbf{g}}_i^\star$. Thus, for our problem, the computation of $\{\mathbf{t}_i^\star, \mathbf{g}_i^\star\}_{i=1}^N$ is not required. To summarize, problem (3.29) can be solved iteratively in a distributed manner as shown in **Algorithm 3.II**.

**Algorithm 3.II**: Iterative algorithm to solve (3.29)

1. Initialization: Set $\{\lambda_n = 1\}_{n=1}^N$.

   **Repeat**

2. With the current $\{\lambda_n\}_{n=1}^N$, compute $\{\overline{\mathbf{g}}_n\}_{n=1}^N$ using (3.33) and update $\{\lambda_n\}_{n=1}^N$ with (3.34).

3. Share the latter $\{\lambda_n\}_{n=1}^N$ among all BSs/processors.

4. Calculate the objective function of (3.29).

   **Until** convergence.

**Convergence:** The convergence of this algorithm can be studied like that of **Algorithm 3.I**. Here, although we are not able to show the global optimality of **Algorithm 3.II** analytically, in all simulation results we observe that the optimal $\boldsymbol{\lambda}$ of (3.29) obtained by **Algorithm 3.II** and the SDP method are the same.

**Computational complexity:** The major computational task of **Algorithm 3.II** arises from matrix inversion which has a complexity on the order of $O(N^{2.376})$ [CW90]. Thus, **Algorithm 3.II** requires $O(N^{2.376})$ per iteration. As will be shown later in Section 3.6, in all our simulations, **Algorithm 3.II** converges to an optimal solution in less than 10 iterations. This shows that the proposed distributed algorithm significantly reduces the computational complexity of (3.22). Therefore, for (3.8), the distributed algorithm requires less overall computational cost than that of the centralized algorithm.

Using $\{\lambda_n\}_{n=1}^N$ of **Algorithm 3.III**, the suboptimal $\{\mathbf{b}_{ls}\}_{s=1}^S, \forall l$ of (3.8) can be computed by (3.27). With these $\{\mathbf{b}_{ls}\}_{s=1}^S, \forall l$, the introduced variables $\nu_s$ and $\tau_s$, and the receiver of the $s$th symbol $\mathbf{w}_s$ are updated by using (3.5), (3.13) and (3.20), respectively. In summary, the suboptimal solution of (3.8) can be obtained distributively as shown in **Algorithm 3.III**.

**Algorithm 3.III**: Distributed algorithm for problem (3.8).

Initialization: Set $\{\mathbf{b}_s\}_{s=1}^S$ like in **Algorithm 3.I** and the maximum number of iterations as $i_{max}$.

**Repeat**

1. With the current $\{\mathbf{b}_s\}_{s=1}^S$, optimize $\{\mathbf{w}_s, \nu_s, \text{ and } \tau_s\}_{s=1}^S$ using (3.5), (3.13) and (3.20), respectively.
2. Using the latter $\{\tau_s, \nu_s, \mathbf{w}_s\}_{s=1}^S$, compute the optimal $\{\lambda_n\}_{n=1}^N$ with **Algorithm 3.II**.
3. Solve for $\{\mathbf{b}_{ls}\}_{s=1}^S, \forall l$, using (3.27).
4. Compute the objective function of (3.8).

**Until** convergence.

**Convergence:** It can be shown that at each step the weighted sum rate of (3.8) is non-decreasing. Hence the algorithm is always convergent.

**Implementation of Algorithm 3.III:** This algorithm can be implemented distributively by two approaches. To be convenient for explanation, we assume $\{M_k = 1\}_{k=1}^K$ and $K = N = L$, i.e., $S = N$.

**First approach:** In this approach, it is assumed that the problem is examined in a central controller which has as many parallel processors as the number of optimization variables. **Algorithm 3.III** can be implemented distributively as follows.

Initialization: The $s$th processor sets $\mathbf{b}_s$ as in **Algorithm 3.III** and $\{\lambda_n = 1\}_{n=1}^N$.

1. The current $\{\mathbf{b}_s\}_{s=1}^S$ are shared among all processors. Once again, using these precoders, the $s$th processor computes its $\mathbf{w}_s, \nu_s$ and $\tau_s$ using (3.5), (3.13) and (3.20), respectively, and then $\{\mathbf{w}_s, \eta_s\}_{s=1}^S$ are shared to all processors.
2. With the current $\{\lambda_n\}_{n=1}^N$ and $\{\mathbf{w}_s, \eta_s\}_{s=1}^S$, the $n$th processor computes $\bar{\mathbf{g}}_n$ using (3.33) and updates its $\lambda_n$ by (3.34). Then, $\{\lambda_n\}_{n=1}^N$ are shared among all processors. The latter two steps are repeated until $\{\lambda_n\}_{n=1}^N$ are found to be optimal.
3. Using $\{\lambda_n\}_{n=1}^N$ of step 2, the $s$th processor computes the optimal $\mathbf{b}_s$ by (3.27).
4. Steps (1), (2) and (3) are repeated until **Algorithm 3.III** converges.
5. The controller finally sends the optimal precoders and decoders to the corresponding BSs and MSs, respectively.

**Second approach:** In this approach, we assume that each BS obtain the channel of all users trough the feedback channel prior to optimization. Here we do not consider any central controller. This is motivated by the fact that each BS is responsible to design its precoder matrix independently by exchanging limited information with the other BSs. In our case, each BS is allowed to exchange $\lambda_n$, $\mathbf{w}_s$ and $\eta_s$ (three scalars for the aforementioned assumption) with all other BSs to jointly design the transceivers of all users. In such approach, **Algorithm 3.III** is implemented distributively as given below.

Initialization: Each BS sets $\{\mathbf{b}_s\}_{s=1}^S$ as in **Algorithm 3.III** and $\{\lambda_n = 1\}_{n=1}^N$.

1. Using the current precoders, the $s$th BS computes its $\mathbf{w}_s, \nu_s$ and $\tau_s$ using (3.5), (3.13) and (3.20), respectively, and then $\{\mathbf{w}_s, \eta_s\}_{s=1}^S$ are shared to all BSs.

2. With the current $\{\lambda_n\}_{n=1}^N$ and $\{\mathbf{w}_s, \eta_s\}_{s=1}^S$, the $n$th BS computes $\bar{\mathbf{g}}_n$ using (3.33) and updates $\lambda_n$ with (3.34). Then, the latter $\{\lambda_n\}_{n=1}^N$ are shared among all BSs. These two steps are repeated until $\{\lambda_n\}_{n=1}^N$ are found to be optimal.

3. Using the current $\{\lambda_n\}_{n=1}^N$, the $s$th BS computes $\{\mathbf{b}_s\}_{s=1}^S$ using (3.27)[8].

4. Steps (1), (2) and (3) are repeated until **Algorithm 3.III** converges.

5. Once **Algorithm 3.III** converges, each BS uses its precoder matrix to precode the data symbols of all users, and also transmits $\{\tilde{\mathbf{W}}_k, \exists k\}$ to those users near to this BS[9].

**Note:** We would like to point out that the un-weighted sum rate optimization problem can be examined with our algorithms either by using (3.14) with $\{\omega_s = 1\}_{s=1}^S$ or employing (3.18) with $\{\omega_s = \tilde{\omega}\}_{s=1}^S$ and $0 < \tilde{\omega} < 1$. It can be clearly seen that our centralized and distributed algorithms are able to handle both of these cases. Furthermore, it is clearly seen that the proposed centralized and distributed algorithms can be extended straightforwardly for the case where the constraint of (3.8) is modified to sum power constraint of the whole network or groups of antennas.

The computational complexities per iteration of the proposed centralized and distributed algorithms, and the algorithm in [SSB08b] for problem (3.8) when $\{M_k = \tilde{M}\}_{k=1}^K$ are summarized in Table 3.1.

As we can see from (3.14), for fixed $\{\nu_s\}_{s=1}^S$ and $\{\omega_s = 1\}_{s=1}^S$, the problem (3.14) turns to the following weighted sum MSE minimization problem:

$$\mathcal{P}3.2: \min_{\{\mathbf{b}_s, \mathbf{w}_s\}_{s=1}^S} \sum_{s=1}^S \bar{\eta}_s \xi_s, \ \text{ s.t } [\sum_{s=1}^S \mathbf{b}_s \mathbf{b}_s^H]_{n,n} \leq p_n, \ \forall n, s. \tag{3.35}$$

[8] In this equation, since the precoders of all users depend on a common matrix inversion $\mathbf{A}^{-1}$, the precoders of all users can be obtained at each BS without significant additional cost.

[9] Note that in a practical scenario, the backhaul capacity is accurate and fast enough to exchange $\lambda_n$, $\mathbf{w}_s$ and $\eta_s$ between BSs (i.e., three scalars for our example setup since N=K=S). Moreover, since users are not expected to design their receivers, the knowledge of $\{\mathbf{b}_s\}_{s=1}^S$ is not required at the receiver side (this reduces the bandwidth requirement of the downlink channel).

**Table 3.1** Computational complexities of the proposed Algorithms and the algorithm in [SSB08b] for (3.8)

| Algorithm | Computational complexity per iteration |
|---|---|
| **Pro centralized** | $O(\sqrt{N+1}(2SN+1)^2(2S^2+2SN))+O(K\tilde{M}^{2.376})$ |
| **Pro distributed** | $O(N^{2.376})+O(K\tilde{M}^{2.376})$ |
| **Alg in [SSB08b]** | $O(\sqrt{(N+S)}(2NS+1)^2(2S^2+2NS+S))+O(K\tilde{M}^{2.376})+\mathfrak{C}_{GP}$ |

where $\{\bar{\eta}_s\}_{s=1}^S$ are the MSE weighting factors. It is clearly seen that the the proposed centralized and distributed algorithms can also be applied to solve the above weighted sum MSE minimization problem.

## 3.5 Extension of the proposed algorithms for the robust versions of $\mathcal{P}$3.1 and $\mathcal{P}$3.2

So far, we have examined the weighted sum rate maximization and MSE minimization problems by assuming that perfect CSI is available at the BS and MSs. In this section, the extension of the proposed centralized and distributed algorithms for robust weighted sum rate and MSE-based problems will be discussed. The robustness against imperfect CSI is incorporated into our designs using stochastic approach [BCV11]. To this end, the channel can be modeled as (see Chapter 2.3)

$$\tilde{\mathbf{H}}_s^H = \widehat{\mathbf{H}}_s^H + \mathbf{R}_{ms}^{1/2}\tilde{\mathbf{E}}_{ws}^H\mathbf{R}_{bs}^{1/2} = \widehat{\mathbf{H}}_s^H + \tilde{\mathbf{E}}_s^H, \ \forall s \tag{3.36}$$

where $\tilde{\mathbf{H}}_s^H$ ($\widehat{\mathbf{H}}_s^H$) is the true (estimated) channel, $\mathbf{R}_{bs} \in C^{N\times N}$ ($\widetilde{\mathbf{R}}_{ms} \in C^{M_k\times M_k}$) antenna correlation matrix at the BS ($s$th symbol), $\mathbf{R}_{ms} = (\mathbf{I}_{M_s} + \sigma_{es}^2\widetilde{\mathbf{R}}_{ms}^{-1})^{-1}$, $\mathbf{E}_s^H$ is the estimation error and the entries of $\mathbf{E}_{ws}^H$ are i.i.d with $\mathcal{CN}(0,\sigma_{es}^2)$.

Like in Chapter 2, we assume that each MS estimates its channel and feeds the estimated channel back to the BS without any error and delay. Thus, both the BS and MSs have the same channel imperfections. With these assumptions, the downlink AMSE of the $s$th symbol is given by

$$\bar{\xi}_s^{DL} = \mathrm{E}_{\mathbf{E}_{ws}^H}\{\xi_s^{DL}\} = \mathbf{w}_s^H\Gamma_s^{DL}\mathbf{w}_s - \mathbf{w}_s^H\widehat{\mathbf{H}}_s^H\mathbf{b}_s - \mathbf{b}_{ks}^H\widehat{\mathbf{H}}_s\mathbf{w}_s + 1 \tag{3.37}$$

where $\Gamma_s^{DL} = \widehat{\mathbf{H}}_s^H\mathbf{B}\mathbf{B}^H\widehat{\mathbf{H}}_s + \sigma_{es}^2\mathrm{tr}\{\mathbf{R}_{bs}\mathbf{B}\mathbf{B}^H\}\mathbf{R}_{ms} + \tilde{\sigma}_s^2\mathbf{I}$.

Using the above channel model, the robust design versions of $\mathcal{P}3.1$ and **P2** can be expressed as

$$\mathcal{P}3.1_{robust}: \max_{\{\mathbf{b}_s,\mathbf{w}_s\}_{s=1}^S} \sum_{s=1}^S \omega_s \mathrm{E}_{\mathbf{E}_{ws}^H}\{\log_2(1+SINR_s)\}$$
$$\text{s.t } [\sum_{s=1}^S \mathbf{b}_s\mathbf{b}_s^H]_{n,n} \leq p_n, \ \forall \tag{3.38}$$

$$\mathcal{P}3.2_{robust}: \min_{\{\mathbf{b}_s,\mathbf{w}_s\}_{s=1}^S} \sum_{s=1}^S \bar{\eta}_s \mathrm{E}_{\mathbf{E}_{ws}^H}\{\xi_s^{DL}\}, \ \text{s.t } [\sum_{s=1}^S \mathbf{b}_s\mathbf{b}_s^H]_{n,n} \leq p_n, \ \forall n,s.$$
$$\equiv \min_{\{\mathbf{b}_s,\mathbf{w}_s\}_{s=1}^S} \sum_{s=1}^S \bar{\eta}_s \bar{\xi}_s, \ \text{s.t } [\sum_{s=1}^S \mathbf{b}_s\mathbf{b}_s^H]_{n,n} \leq p_n, \ \forall n,s. \tag{3.39}$$

As we can see from (3.38), the robust version of $\mathcal{P}3.1$ (i.e., $\mathcal{P}3.1_{robust}$) contains the expectation term incorporating $\log_2(.)$ which is nonlinear. Due to this fact, solving $\mathcal{P}3.1_{robust}$ is not trivial and it is still an open problem. However, (3.35) and (3.39) have similar problem structures. This shows that one can apply the algorithms of this chapter to solve $\mathcal{P}3.2_{robust}$ (see also [BVC12] for more details).

Note that when BSs are coordinated only at the beamforming level (i.e., multi-cell systems), the overall precoder matrix **B** will have block-diagonal structure. Since the analysis of this chapter can be employed for block-diagonal matrix **B**, the algorithm of this chapter can be extended straightforwardly for multi-cell systems (see Section V of [BV11d] for further details).

## 3.6 Simulation results

In this section, we present the simulation results for problem (3.8) (i.e., $\mathcal{P}3.1$). All of our simulation results are averaged over 100 randomly chosen channel realizations. The channel between all BS and each MS consists of ZMCSCG entries with unit variance. It is assumed that the noise variances of all users are the same, i.e., $\{\sigma_k^2 = \sigma^2\}_{k=1}^K$. The SNR is defined as $P_{\text{sum}}/\sigma^2$ and it is controlled by varying $\sigma^2$, where $P_{\text{sum}}$ is the total sum power utilized by all antennas.

### 3.6.1 Comparison of the proposed algorithms and the algorithm in [SSB08b]

For the comparison of these three algorithms, we consider a system with $L = 2$ BSs where each BS has 4 antennas, and $K = 4$ MSs where each MS has 2 antennas. It is assumed that $\{p_n = 0.125mw\}_{n=1}^{8}$ and $\boldsymbol{\omega}_1 = [0.6\ 0.4\ 0.5\ 0.8\ 0.25\ 0.8\ 0.46\ 0.28]^T$. First, we compare these three algorithms based on their powers utilized at each antenna when $\sigma^2 = 0.1mw$. For this system setup, all of these three algorithms utilize the maximum available powers at each BS antenna[10]. Second, we compare the performance of the aforementioned algorithms based on their total achievable weighted sum rate. Fig. 3.2 shows that the proposed distributed algorithm achieves the same weighted sum rate as that of the centralized algorithm. Moreover, our proposed algorithms outperform the algorithm proposed in [SSB08b].

Next, we compare the performances of the proposed algorithms and the algorithm in [SSB08b] for different rate weighting factors. The comparison is based on the total weighted sum rate. For this purpose we use two sets of rate weighting factors $\boldsymbol{\omega}_2$ and $\boldsymbol{\omega}_3$ as $\boldsymbol{\omega}_2 = [0.9\ 0.2\ 0.5\ 0.95\ 0.1\ 0.9\ 0.2\ 0.05]$ and $\boldsymbol{\omega}_3 = [0.1\ 0.4\ 0.2\ 0.6\ 0.3\ 0.16\ 0.12\ 0.25]$. For these weighting factors, the weighted sum rates of the proposed algorithms and the algorithm in [SSB08b] are plotted in Fig. 3.3. As can be seen from this figure, the proposed algorithms outperform the algorithm in [SSB08b].

From Fig. 3.2 and Fig. 3.3, we can observe that the performance gap between the proposed algorithms and the algorithm in [SSB08b] depends on the weighting factors. Here, we would like to mention that for more than 90% of our channel realizations, we have noticed that the proposed algorithms achieve at least the same weighted sum rate as that of [SSB08b].

### 3.6.2 Convergence characteristics of the proposed algorithms and the algorithm in [SSB08b]

In Section 3.4, the computational complexities of the proposed centralized and distributed algorithms, and the algorithm in [SSB08b] is discussed for a single iteration only. Therefore, to compare the overall computational complexities of our algorithms and the algorithm in [SSB08b], the convergence speed of these algorithms should be examined. In this simulation, we examine the convergence speed of our algorithms and the algorithm proposed in [SSB08b] for the initialization as presented in **Algorithm 3.1**. We have used the same simulation parameters as in the first paragraph of Section 3.6.1. As can be seen from Fig. 3.4, the proposed algorithms have faster convergence speed and higher weighted sum rate than that of the algorithm proposed in [SSB08b].

[10]We would like to mention here that for problem (3.8) all antennas do not necessarily utilize their maximum powers to optimize the total weighted sum rate (see for example in [EV11]).

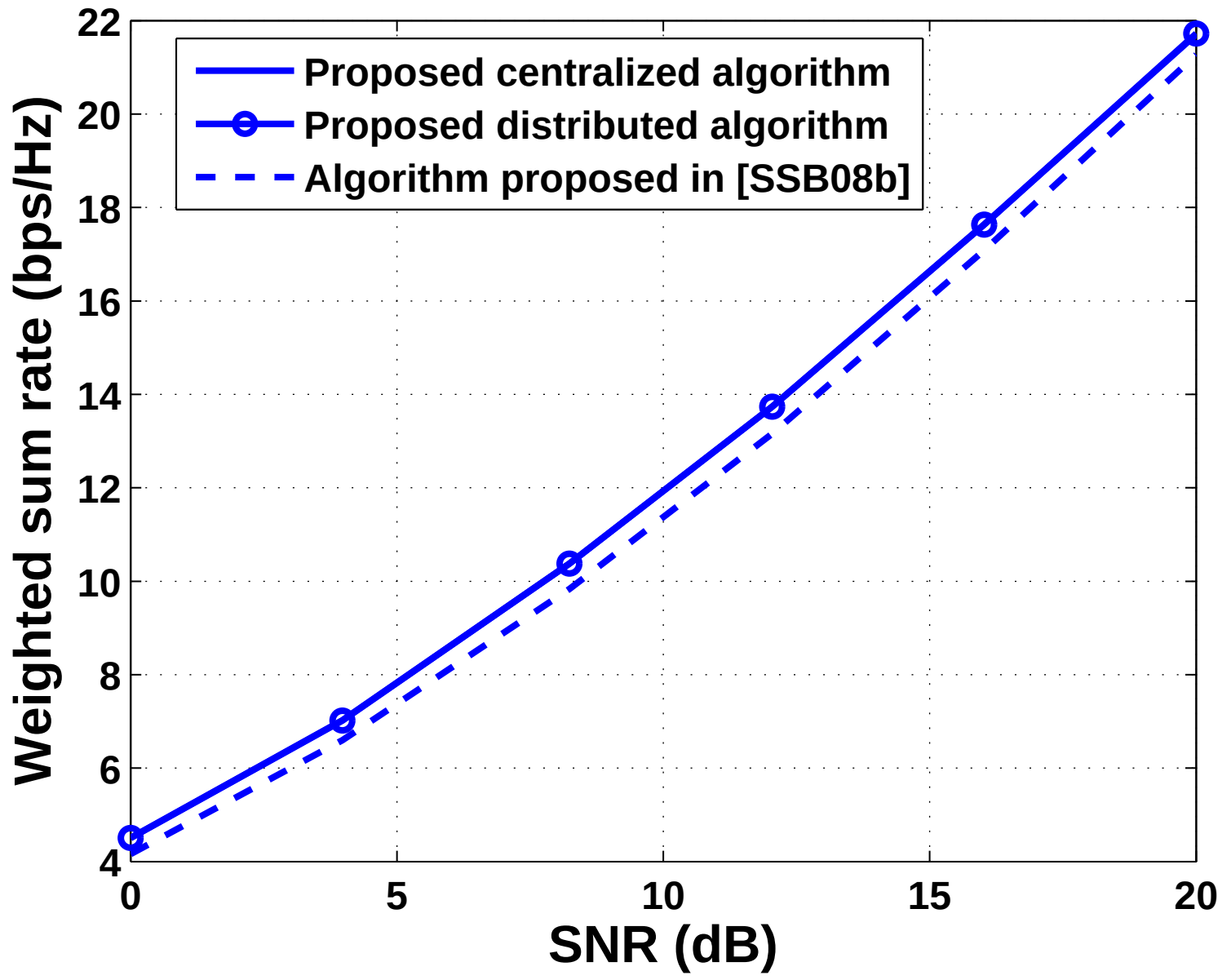

**Figure 3.2** Comparison of the proposed centralized and distributed algorithms, and the algorithm proposed in [SSB08b].

### 3.6.3 Convergence characteristics of Algorithm 3.II

To demonstrate the computational advantage of our distributed algorithm over the centralized algorithm, we examine the convergence characteristics of **Algorithm 3.II** for both small-scale and large-scale networks.

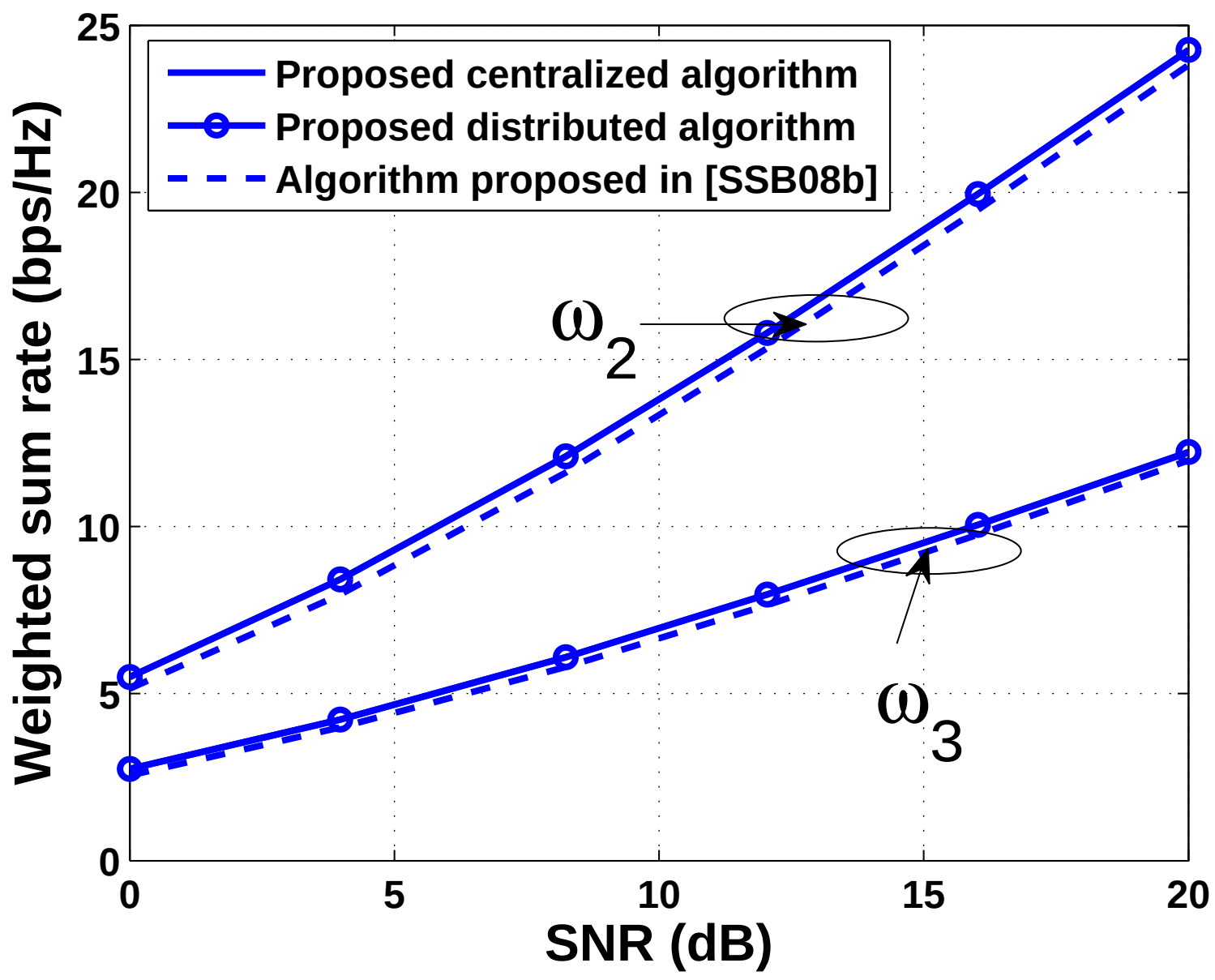

**Figure 3.3** Comparison of the proposed centralized and distributed algorithms, and the algorithm in [SSB08b].

#### 3.6.3.1 Small-scale network

In this simulation we demonstrate the convergence characteristics of **Algorithm 3.II** for a system with $L = 2$ coordinated BSs where each of them has two antennas and $K = 2$ MSs where each MS is equipped with 2 antennas. For this system Fig. 3.5 shows the convergence characteristics of **Algorithm 3.II** at different iterative stages of **Algorithm 3.III** (i.e., for different $\{\mathbf{w}_s, \eta_s\}_{s=1}^{S}$). As can be seen from this figure, **Algorithm 3.II** converges to an optimal solution in less than 10 iterations.

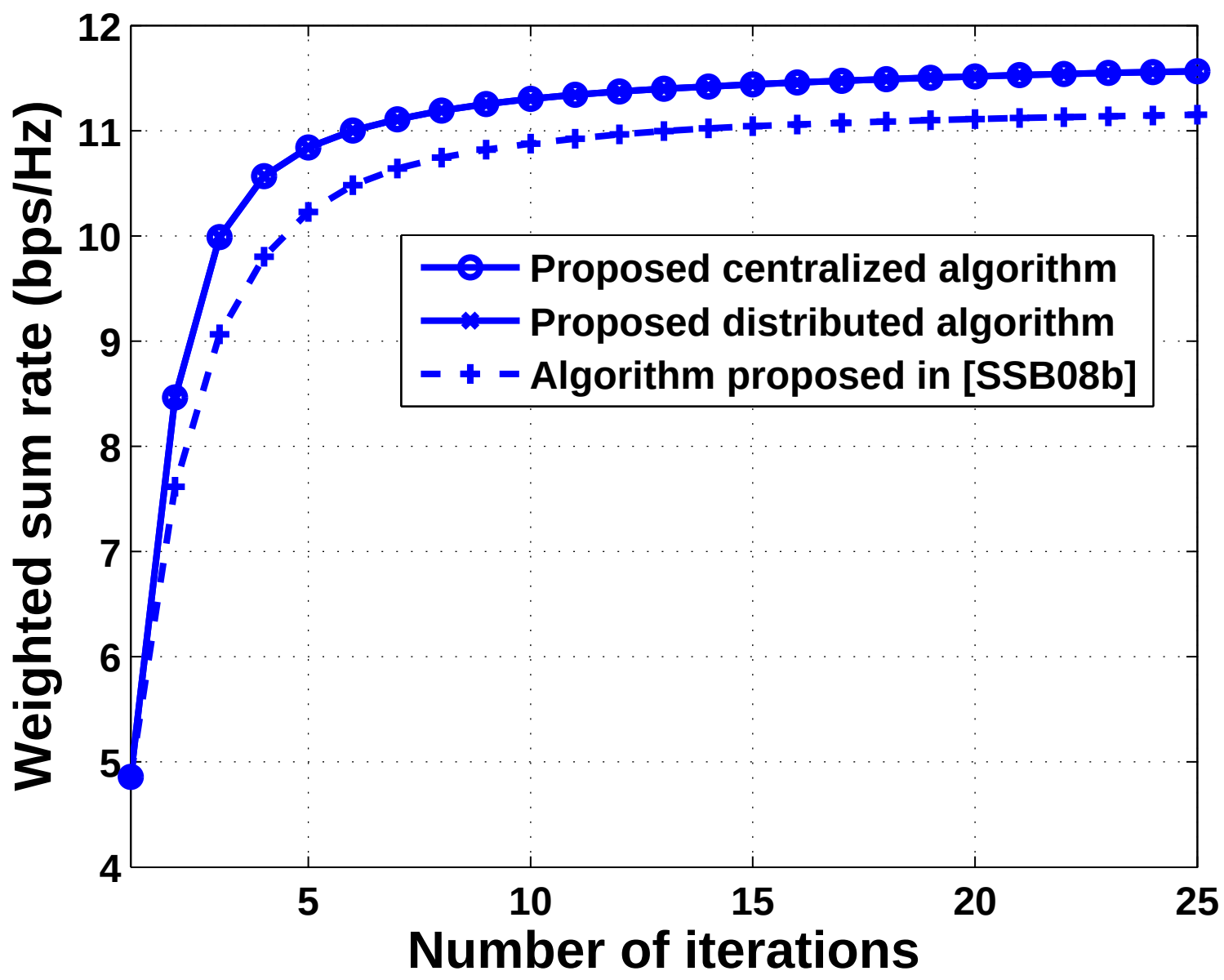

**Figure 3.4** Comparison of the convergence characteristics of our algorithms and the algorithm in [SSB08b] when $\sigma^2 = 0.1$.

$$\mathbf{Z} = \begin{bmatrix} 0.56-0.03i & 0.00+0.11i & 1.17 & 0.62-0.01i & 0.05+0.11i & 2.00 \\ -0.05-0.08i & 0.75+0.03i & 0.42 & -0.11+0.12i & 0.62+0.00i & 0.27 \\ 0.64-0.01i & -0.06+0.00i & 0.67 & 0.76-0.03i & -0.13+0.03i & 0.59 \\ -0.14+0.04i & 0.80+0.01i & 0.96 & -0.08-0.04i & 0.96+0.04i & 1.21 \\ 0.61+0.01i & 0.044+0.08i & 4.10 & 0.59+0.02i & 0.01+0.02i & 4.70 \\ -0.13+0.28i & 0.24-0.00i & 0.18 & -0.13+0.33i & 0.05-0.00i & 0.16 \\ 0.81-0.01i & -0.18-0.00i & 0.41 & 0.75+0.00i & -0.25-0.07i & 0.30 \\ 0.10-0.09i & 1.04+0.04i & 1.64 & 0.29-0.17i & 0.99+0.02i & 2.37 \end{bmatrix} \tag{3.40}$$

where $\mathbf{Z} = [\mathbf{Z}_1\ \mathbf{Z}_2; \mathbf{Z}_3\ \mathbf{Z}_4]$ with $\mathbf{Z}_k = [[\tilde{\mathbf{W}}_1; \cdots; \tilde{\mathbf{W}}_K]\ \tilde{\boldsymbol{\eta}}]$ and $\tilde{\boldsymbol{\eta}} = [\eta_1, \cdots, \eta_S]^T$.

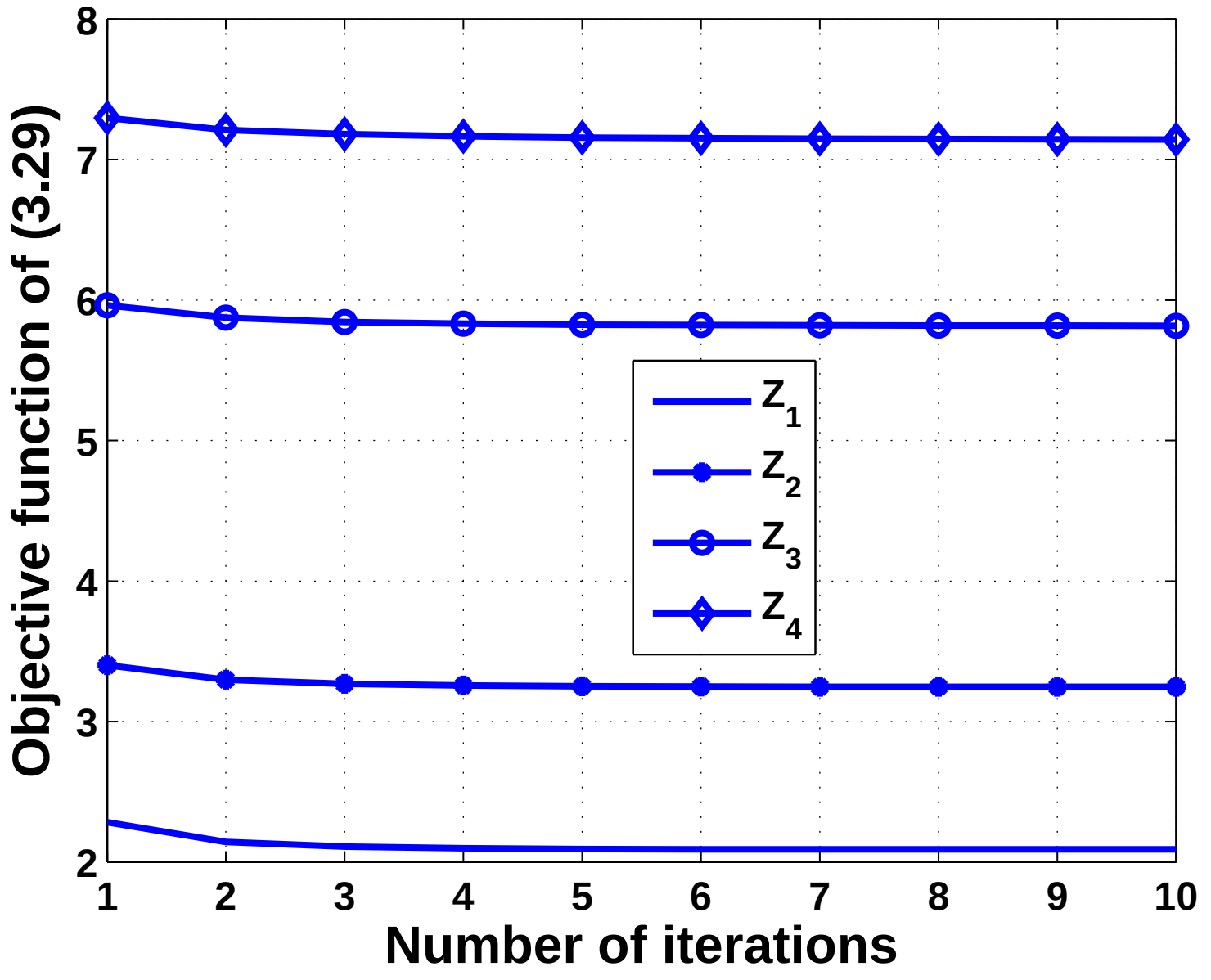

**Figure 3.5** Convergence characteristics of **Algorithm 3.II** at different iterative stages of **Algorithm 3.III** for small-scale network with the set **Z** as given in (3.40).

#### 3.6.3.2 Large-scale network

Next we examine the convergence characteristics of **Algorithm 3.II** for large-scale networks. We consider a system with $L = 25$ coordinated BSs where each of them has four antennas and $K = 50$ MSs where each MS is equipped with 2 antennas. For simplicity, we assume that $\{p_n = 0.25\}_{n=1}^{N}$ and $\sigma^2 = 0.1mw$[11]. The weighting factor of the $s$th symbol ($\omega_s$) is chosen from a uniform distribution with $\{0 < \omega_s < 1\}_{s=1}^{S}$. For these settings, we examine the convergence characteristics of **Algorithm 3.II** at different iterative stages of **Algorithm 3.III**. As can be seen from Fig. 3.6, **Algorithm 3.II** converges to an optimal solution within few iterations.

[11]Similar behavior is observed for the other $\sigma^2$.

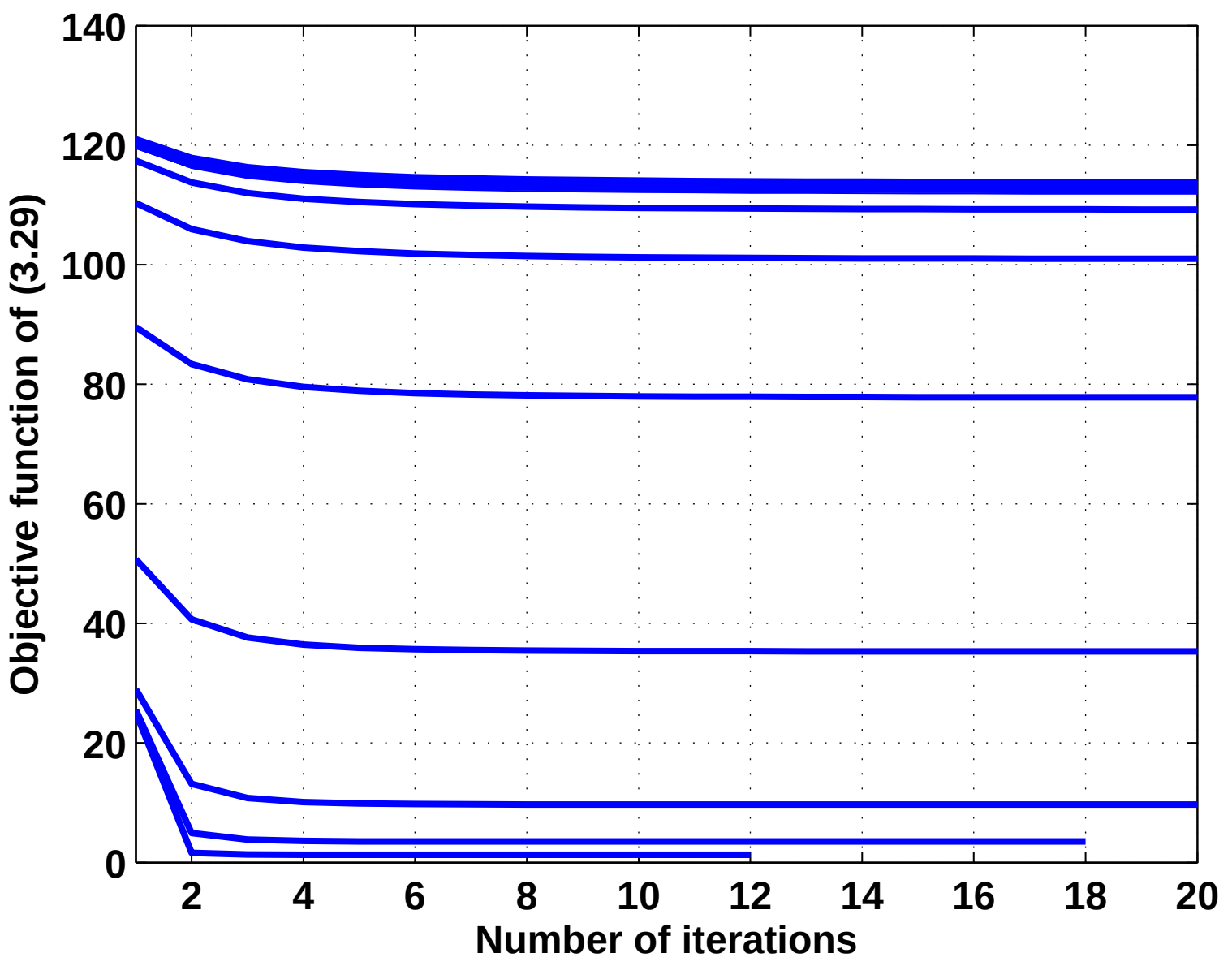

**Figure 3.6** Convergence characteristics of **Algorithm 3.II** at different iterative stages of **Algorithm 3.III** for large-scale network.

### 3.6.4 Simulation results for problem $\mathcal{P}$3.1 when $\{M_k = 1\}_{k=1}^{K}$

When each MS has single antenna, the global optimal solution of (3.8) can be obtained with the framework of MGO algorithm as discussed in [BU09]. The MGO algorithm requires solving a feasibility problem to get the upper boundary feasible points of a monotonic optimization problem (see also [RTM01] for more details about MGO and upper boundary feasible points of a monotonic optimization problem). For our case, this feasibility problem (i.e., rate feasibility problem) can be solved by the phase rotation technique of [YL07]. According to [BU09], the computational complexity of MGO algorithm grows quickly with the number of users. Thus, the MGO algorithm serves as a benchmark for suboptimal less complex algorithms. On the other hand, a simple IZF solution for (3.8) with $\{M_k = 1\}_{k=1}^{K}$ can be obtained by the approach of [WES08] (see Section V.B of [WES08]). These findings motivate us to compare our proposed algorithms with that of MGO, IZF and the algorithm of [SSB08b] when $\{M_k = 1\}_{k=1}^{K}$. The comparison of these algorithms is based on the total weighted sum rate of all users when $L = K = 3$, $\{N_l = 1\}_{l=1}^{L}$, $\{p_n = 1/N\}_{n=1}^{N}$, the rate weighting factors $\tilde{\omega}_1 = [0.46\ 0.83\ 0.79]^T$ and $\tilde{\omega}_2 = [0.9\ 0.54\ 0.1]^T$, and all the other settings are the same as the first paragraph of Section 3.6. For the MGO algorithm, we have used the following tolerance $(\epsilon = 0.001, \varrho = 0.01)$ which is analogous to $(\epsilon, \eta)$ of [RTM01]. Here, we have employed the weighting factors $\{\tilde{\omega}_i\}_{i=1}^{2}$ to differentiate from the weighting factors $\{\boldsymbol{\omega}_i\}_{i=1}^{3}$ which are used in Sections 3.6.1 - 3.6.3 for $\{M_k = 2\}_{k=1}^{K}$. As can be seen from Fig. 3.7.(a)-(b), the proposed algorithms achieve global optimum, whereas the algorithms in [SSB08b] and [WES08] do not achieve the global optimum. As expected, at high SNR regions, the weighted sum rate achieved by the IZF algorithm of [WES08] approaches the optimal weighted sum rate. However, the exact SNR value at which the weighted sum rate achieved by the latter algorithm approaches the optimal weighted sum rate is not necessarily the same for all rate weighting factors.

We would like to mention here that the performance characteristics of our proposed algorithms, the algorithm of [SSB08b], the IZF algorithm of [WES08] and the MGO algorithm of [BU09] for $\mathcal{P}3.2$ are like that of $\mathcal{P}3.1$. Due to this reason, we omit the simulation results of these algorithms for $\mathcal{P}3.2$.

## 3.7 Conclusions

This chapter considers the joint linear transceiver design problem for downlink multiuser MIMO systems with coordinated BSs. We examine maximization of the total weighted sum rate with per BS antenna power constraint problem. We propose novel centralized and computationally efficient distributed iterative algorithms that achieve local optimum to the latter problem. These algorithms are described as follows. First, by introducing additional optimization variables, we reformulate the original problem into a new problem. Second, for the given precoder matrices of all users, the optimal receivers are computed using MMSE method and the optimal introduced variables are obtained in closed form expressions. Third, by keeping the introduced variables and receivers constant, the precoder matrices of all users are optimized by using SOCP and matrix fractional minimization approaches for the centralized and distributed algorithms, respectively. Finally, the second and third steps are repeated until these algorithms converge. We have shown that the proposed algorithms require less computational cost than that of the existing algorithm. Moreover, the proposed algorithms achieve higher weighted sum rate than that of the existing linear algorithm. All simulation results show that the proposed distributed algorithm achieve the same performance as that of the centralized algorithm. In particular, when each of the users has single antenna, we have observed that the proposed algorithms achieve the global optimum.

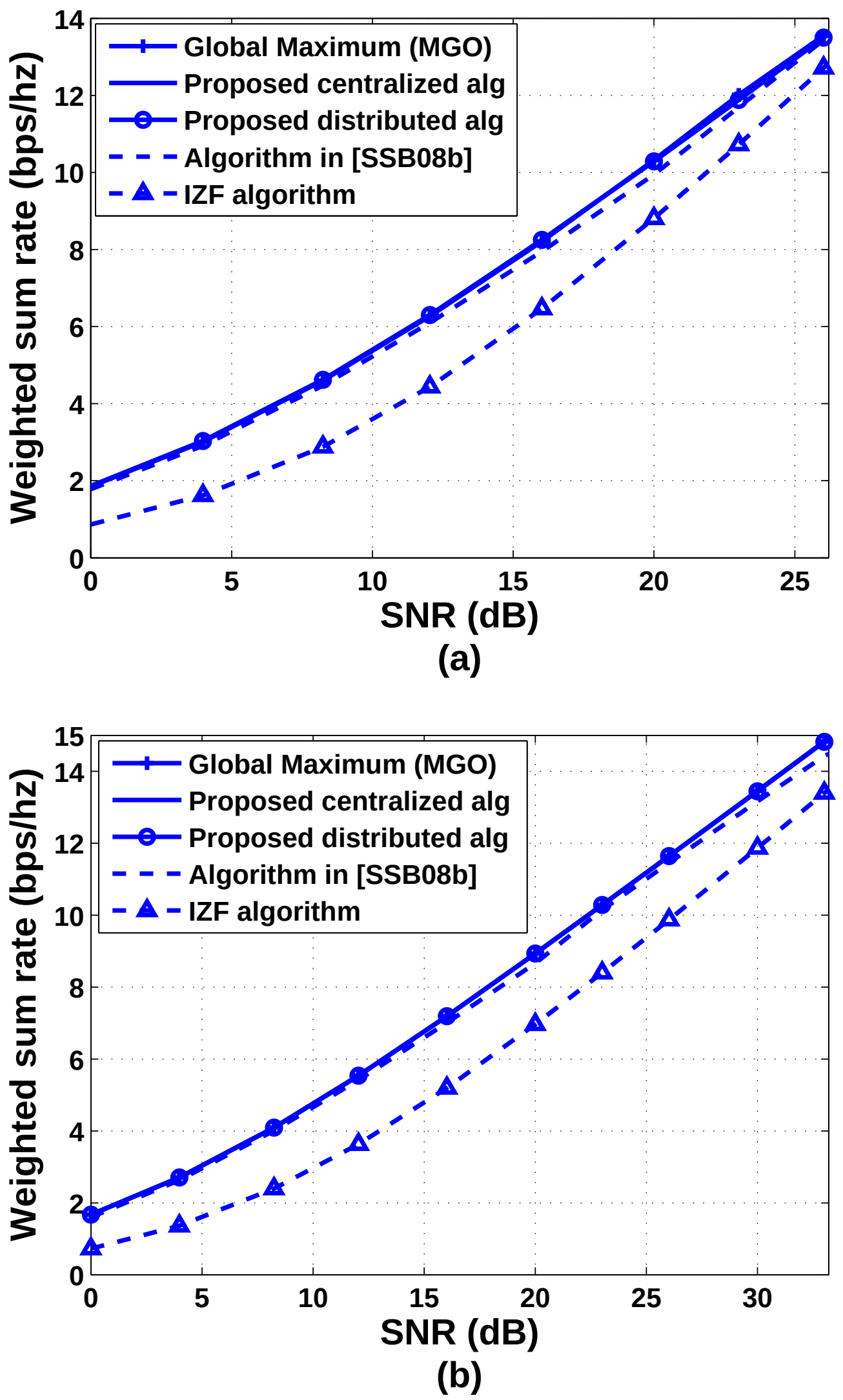

**Figure 3.7** Comparison of the proposed algorithms, MGO algorithm of [BU09], IZF algorithm of [WES08] and the algorithm in [SSB08b]. (a) For the rate weighting factor $\tilde{\omega}_1$. (b) For the rate weighting factor $\tilde{\omega}_2$.

## 3.8 Looking ahead

The proposed algorithms of this chapter employs Lagrangian dual decomposition for rate and MSE-based problems. And we know that the Lagrangian dual problem may have different structure for different primal problems. Due to this reason, the algorithms of this chapter may not be extended straightforwardly to solve other rate or MSE-based problems.

However, from this chapter we have learnt that the problem structure of coordinated BS systems is similar to that of the conventional multiuser MIMO systems. And, from Chapter 2, we have learnt that the duality approach of solving MSE-based problems can be extended straightforwardly to solve other total BS power constrained MSE-based problems. Also, the power allocation part of the MSE-based problems of Chapter 2 is formulated as a GP which can be solved by a central controller[12], whereas, the transmit and receive filters of each user can be computed independently which naturally leads to a distributive algorithm. These promising outcomes of duality algorithms motivate us to generalize the MSE duality of Chapter 2 for practically relevant power constraints (like per BS antenna, user etc) which is the topic of the next chapter.

[12] As mentioned in [SSB08c] (see Appendix A of [SSB08c]), a small desktop computer can solve a GP of 100 variables and 10000 constraints by standard interior point method under a minute. Thus, we believe that solving the GP problem centrally is feasible even for large scale network.

## 3.9 Appendix 3.A: Proof of the equivalence of (3.9) and (3.10)

Since the constraint functions of (3.10) do not depend on $\{\mathbf{w}_s\}_{s=1}^S$, the receivers $\{\mathbf{w}_s\}_{s=1}^S$ of (3.10) can be optimized by applying standard first order differentiation of $\prod_{s=1}^S \tilde{\xi}_s$ with respect to $\{\mathbf{w}_s^H\}_{s=1}^S$. By doing so, we get

$$\frac{\partial(\prod_{i=1}^S \tilde{\xi}_i)}{\partial \mathbf{w}_s^H} = \omega_s \prod_{i=1, i\neq s}^S \tilde{\xi}_i \xi_s^{(\omega_s - 1)} \cdot \frac{\partial \xi_s}{\partial \mathbf{w}_s^H} = 0$$

$$\Rightarrow \mathbf{w}_s^\star = (\tilde{\mathbf{H}}_s^H \mathbf{B}\mathbf{B}^H \tilde{\mathbf{H}}_s + \tilde{\sigma}_s^2 \mathbf{I})^{-1} \tilde{\mathbf{H}}_s^H \mathbf{b}_s$$

where the last equality follows from the fact that $\omega_s \prod_{i=1, i\neq s}^S \tilde{\xi}_i \xi_s^{(\omega_s - 1)}$ is always positive for any $\{\mathbf{b}_s, \omega_s\}_{s=1}^S$ with $0 < \omega_s < 1$. Now, by substituting the above $\mathbf{w}_s^\star$ into (3.10), we get (3.9).

## 3.10 Appendix 3.B: Computation of $\mathfrak{C}_e$

Problem (3.10) has been examined in [SSB08b] for the case where the power of each symbol is strictly positive. In [SSB08b], the transmitters are decomposed into a product of unity norm filter and square root of power allocation matrices, and the receiver matrix of each user is decomposed as a product of the inverse of the square root of power allocation, unity norm filter and diagonal scaling factor matrices (see Section 2 of [SSB08b]). Upon doing so, the weighted sum rate maximization problem is formulated as in (2) of [SSB08b]. Then, for (2) of [SSB08b], this paper utilizes Algorithm 1 of **Table 1**. Here, we summarize the computational cost required to perform one iteration of Algorithm 1 in [SSB08b] by using the system model parameter settings as discussed in Section 3.2 of our paper. For simplicity, we assume that $\{M_k = \tilde{M}\}_{k=1}^K$. The major computational cost of Algorithm 1 of [SSB08b] comes from the steps 5, 6, 7 and 8. The steps 5 and 7 of the latter algorithm contain matrix inversion. According to [CW90], matrix inversion can be performed with a complexity of $O(K\tilde{M}^{2.376})$. The computational load of step 8 is on the order of $O(\sqrt{(N+S)}(2NS+1)^2(2S^2+2NS+S))$ [LVBL98] (see page 196 of [LVBL98] for the details). In general, the computational complexity of step 6 depends on different parameter settings and solution methods. The detail analysis on the computational complexity of GP problems (i.e., step 6) can be found in [BV04], [Chi05] (see page 36 of [Chi05] for the Barrier method of solving GP problems). Therefore, the computational complexity of Algorithm 1 in [SSB08b] per iteration is given by $O(\sqrt{(N+S)}(2NS+1)^2(2S^2+2NS+S)) + O(K\tilde{M}^{2.376}) + \mathfrak{C}_{GP}$, where $\mathfrak{C}_{GP}$ is the computational cost of the GP in [SSB08b].

## 3.11 Appendix 3.C: Proof of Lemma 3.1

**Proof.** For fixed $\{\mathbf{b}_s, \mathbf{w}_s\}_{s=1}^S$, optimizing $\{\nu_s\}_{s=1}^S$ of (3.12) can be expressed as

$$\min_{\{\nu_s\}_{s=1}^S} \left(\frac{1}{S}\sum_{s=1}^S \tilde{\xi}_s \nu_s\right)^S, \ \text{ s.t} \prod_{s=1}^S \nu_s = 1, \nu_s \geq 0, \ \forall s. \tag{3.41}$$

The above problem is GP for which global optimality is guaranteed. Clearly, the optimal solution of (3.41) satisfy $\{\nu_s > 0\}_{s=1}^S$, and the objective and constraint functions of this problem are continuously differentiable. Moreover, by replacing $\nu_1 = (\prod_{s=2}^S \nu_s)^{-1}$, the equality constraint of the latter problem can be removed. These two facts show that the optimal $\{\nu_s\}_{s=1}^S$ of the above problem are regular [DB09], [JMT97][13]. Thus, the global optimal solution of (3.41) can be obtained by choosing $\{\nu_s\}_{s=1}^S$ that satisfy the Karush-Kuhn-Tucker (KKT) optimality conditions which are given by [BV04]

$$\tilde{\xi}_s\left(\frac{1}{S}\sum_{i=1}^S \tilde{\xi}_i \nu_i\right)^{S-1} - \gamma \prod_{i=1, i\neq s}^S \nu_i - \lambda_s = 0 \tag{3.42}$$

$$\lambda_s \nu_s = 0 \tag{3.43}$$

$$\lambda_s \geq 0, \ \forall s \tag{3.44}$$

where $\gamma$ and $\{\lambda_s\}_{s=1}^S$ are the Lagrangian multipliers corresponding to the constraints $\prod_{s=1}^S \nu_s = 1$ and $\{\nu_s \geq 0\}_{s=1}^S$, respectively. Multiplying (3.42) by $\nu_s$, and employing $\prod_{s=1}^S \nu_s = 1$ and (3.43) results

$$\tilde{\xi}_s \nu_s \left(\frac{1}{S}\sum_{i=1}^S \tilde{\xi}_i \nu_i\right)^{S-1} - \gamma \nu_s \prod_{i=1, i\neq s}^S \nu_i - \lambda_s \nu_s = 0$$

$$\Rightarrow \tilde{\xi}_s \nu_s \left(\frac{1}{S}\sum_{i=1}^S \tilde{\xi}_i \nu_i\right)^{S-1} = \gamma \prod_{i=1}^S \nu_i = \gamma. \tag{3.45}$$

[13] For the inequality constrained optimization problems, a feasible point is said to be regular if all the inequality constraints are inactive at this point [DB09], [JMT97].

By summing the $S$ equalities of (3.45), $\gamma$ can be determined by

$$\gamma = \left( \frac{1}{S} \sum_{i=1}^{S} \tilde{\xi}_i \nu_i \right)^S. \tag{3.46}$$

Substituting $\gamma$ of (3.46) into (3.45), and noting that $\frac{1}{S}\sum_{i=1}^{S} \tilde{\xi}_i \nu_i > 0$ we obtain

$$\tilde{\xi}_s \nu_s = \frac{1}{S} \sum_{i=1}^{S} \tilde{\xi}_i \nu_i. \tag{3.47}$$

Multiplying the $S$ equalities of (3.47) and utilizing $\prod_{s=1}^{S} \nu_s = 1$ yields

$$\prod_{s=1}^{S} \tilde{\xi}_s = (\frac{1}{S} \sum_{i=1}^{S} \tilde{\xi}_i \nu_i)^S. \tag{3.48}$$

The above expression shows that the optimal/suboptimal solution of (3.10) can be equivalently obtained by solving (3.12). By employing (3.47) and (3.48), it can be shown that the optimal $\{\nu_s\}_{s=1}^{S}$ of (3.41) can be expressed as in (3.13) [EV11], [JMT97]. □

# Linear Transceiver Design for Downlink Multiuser MIMO Systems: Downlink-Interference Duality Approach

4

This chapter considers linear transceiver design for downlink multiuser MIMO systems. We examine different transceiver design problems. We focus on two groups of design problems. The first group is the WSMSE (i.e., symbol-wise or user-wise WSMSE) minimization problems and the second group is the minimization of the maximum WMSE (symbol-wise or user-wise WMSE) problems. The problems are examined for the practically relevant scenario where the power constraint is a combination of per BS-antenna and per symbol (user), and the noise vector of each MS is a ZMCSCG random variable with arbitrary covariance matrix. For each of these problems, we propose a novel downlink-interference duality based iterative solution. Each of these problems is solved as follows. First, we establish a new MSE downlink-interference duality. Second, we formulate the power allocation part of the problem in the downlink channel as a GP. Third, using the duality result and the solution of GP, we utilize alternating optimization technique to solve the original downlink problem. For the first group of problems, we have established symbol-wise and user-wise WSMSE downlink-interference duality. These duality are established by formulating the noise covariance matrices of the interference channels as fixed point functions. On the other hand, for the second group of problems, we have established symbol-wise and user-wise MSE downlink-interference duality. These duality are established by formulating the noise covariance matrices of the interference channels as marginally stable (convergent) discrete-time-switched systems. The proposed duality based iterative solutions can be extended straightforwardly to solve many other linear transceiver design problems. We also show that our MSE downlink-interference duality unify all existing MSE duality. In our simulation results, we have observed that the proposed duality based iterative algorithms utilize less total BS power than that of the existing algorithms.

## 4.1 Introduction

In the downlink multiuser MIMO systems, most practically relevant transceiver design problems such as weighted sum rate maximization, rate or SINR balancing and rate or SINR constrained power minimization can be equivalently expressed as MSE-based problems (see for example [SSB08c] and Chapter 3 of this thesis). Because of this, the current chapter examines MSE-based problems in the downlink channel. The downlink MSE-based problems can be solved by direct approach as in [SSVB08, UC08] or by uplink-downlink duality approach as in [SSB07, HJU09, SSJB05].

Several MSE-based problems have been examined by duality approach [SSB07, SSJB05, HJU09, BV11a, BCV11] (see also Chapter 2 of this thesis). However, the duality of these works are able to solve total BS power constrained MSE-based problems only. In a practical multi-antenna BS system, the maximum power of each BS antenna is limited [YL07]. Also, in some scenario allocating different powers to different users (symbols) according to their priority or protection level has some interest. Also as mentioned in the previous chapter, in a multimedia communication, different types of information (for example, audio and video information) can be sent to a user (all users) simultaneously [BV11d, SSP01]. In such a scenario, for successful transmission, more priority (power) could be given to symbols (users) corresponding to the video information. Thus, for this scenario, the design criteria may incorporate fairness/priority and power constraints for each symbol (user). On the other hand, examining combined (per antenna and symbol (user)) power con-straints may have practical interest (for example in network MIMO). For these reasons, the current chapter generalizes the existing MSE duality for more general classes of power constraints while incorporating both symbol-wise and user-wise MSE fairness/priority, and combined (i.e., per antenna and symbol (user)) power constraints. We examine the following problems: Minimization of symbol-wise WSMSE constrained with per BS antenna and symbol powers ($\mathcal{P}4.1$), minimization of user-wise WSMSE constrained with per BS antenna and user powers ($\mathcal{P}4.2$), minimization of the maximum symbol-wise WMSE constrained with per BS antenna and symbol powers ($\mathcal{P}4.3$) and minimization of the maximum user-wise WMSE constrained with per BS antenna and user powers ($\mathcal{P}4.4$). Each of these problems is examined for the scenario where the noise vector of each MS is a ZMCSCG random variable with arbitrary covariance matrix.

To the best of our knowledge, the problems $\mathcal{P}4.1$ - $\mathcal{P}4.4$ are non-convex. Furthermore, duality based solutions for these problems with our noise covariance matrix assumptions are not known. In the current chapter, we propose duality based iterative solutions to solve the problems. Each of these problems is solved as follows. First, we establish a new MSE downlink-interference duality. Second, we formulate the power allocation part of the problem in the downlink channel as a GP. Third, using the duality result and the solution of GP, we utilize alternating optimization technique to solve the original downlink problem. For the problems $\mathcal{P}4.1$ and $\mathcal{P}4.2$, the duality are established by formulating the noise covariance matrices of the interference channels as fixed point functions. For these two problems, the noise covariance matrices of the dual interference channels are computed by modifying the approach of [BV11b] to $\mathcal{P}4.1$ and $\mathcal{P}4.2$ of the current paper. On the other hand, for the problems $\mathcal{P}4.3$ and $\mathcal{P}4.4$, the duality are established by formulating the noise covariance matrices of the interference channels as new marginally stable (convergent) discrete-time-switched systems. The proposed duality based iterative solutions can be extended straightforwardly to solve many other linear transceiver design problems. We also show that our MSE downlink-interference duality unify all existing MSE duality. In our simulation results, we have observed that the proposed duality based iterative algorithms utilize less total BS power than that of the existing algorithms. The main contributions of the current chapter is summarized as follows:

1. To solve the problems $\mathcal{P}4.1$ and $\mathcal{P}4.2$, we have established new WSMSE downlink-interference duality by formulating the noise covariance matrices of the interference channels as fixed point functions. As will be clear later, for WSMSE-based problems with a total BS power constraint function, the proposed duality based algorithm requires less computation than that of the existing duality based algorithms.

2. To solve the problems $\mathcal{P}4.3$ and $\mathcal{P}4.4$, we have established novel MSE (symbol-wise and user-wise) downlink-interference duality by formulating the noise covariance matrices of the interference channels as marginally stable (convergent) discrete-time-switched systems. Furthermore, as will be shown later in Section 4.9, the proposed duality based iterative solutions can be extended straightforwardly to solve many other linear transceiver design problems. We also show that the MSE downlink-interference duality of the current paper is also applicable to solve total BS power based linear transceiver design problems. Thus, the current duality unify all existing MSE duality[1].

3. By employing the system model of [SSB08c] and [BV11a], we formulate the power allocation parts of $\mathcal{P}4.1$ - $\mathcal{P}4.4$ as GPs. The GPs are formulated by applying the GP formulation approach of [SSB08c]. Consequently, we are able to solve our problems by alternating optimization technique [SSB07, BV11a, BCV11, BV11b] (i.e., duality based iterative algorithm).

4. In our simulation results, we have observed that the proposed duality based iterative algorithms utilize less total BS power than that of the existing algorithms.

This chapter is organized as follows. In Section 3.2, multiuser MIMO downlink and virtual interference channel system models are presented. In Section 4.3, we formulate our problems $\mathcal{P}1$ - $\mathcal{P}4.4$ and discuss the general framework of our duality based iterative solutions. Sections 4.4 - 4.8 present the proposed duality based iterative solutions for solving these problems. The extension of our duality based iterative algorithms to other non-robust and robust design problems is discussed in Sections 4.9 and 4.10. In Section 4.11, computer simulations are used to compare the performance of the proposed duality algorithms with that of existing algorithms. Finally, conclusions are drawn in Section 4.12.

[1] Note that the existing MSE duality are established for a total BS power based linear transceiver design problems (see [HJU09, SSB07, SSB08c]).

## 4.2 System model

In this section, multiuser MIMO downlink and virtual interference channel system models are discussed which are shown in Fig. 4.1. In the downlink channel, the BS and $k$th MS are equipped with $N$ and $M_k$ antennas, respectively. The total number of MS antennas are thus $M = \sum_{k=1}^{K} M_k$. By denoting the symbol intended for the $k$th user as $\mathbf{d}_k \in \mathbb{C}^{S_k \times 1}$ and $S = \sum_{k=1}^{K} S_k$, the entire symbol can be written in a data vector $\mathbf{d} \in \mathbb{C}^{S \times 1}$ as $\mathbf{d} = [\mathbf{d}_1^T, \cdots, \mathbf{d}_K^T]^T$. The BS precodes $\mathbf{d}$ into an $N$ length vector by using its overall precoder matrix $\mathbf{B} = [\mathbf{b}_{11}, \cdots, \mathbf{b}_{KS_K}]$, where $\mathbf{b}_{ks} \in \mathbb{C}^{N \times 1}$ is the precoder vector of the BS for the $k$th MS $s$th symbol. The $k$th MS employs a receiver $\mathbf{w}_{ks}$ to estimate the symbol $d_{ks}$. We follow the same channel matrix notations as in [BV11a]. The estimates of the $k$th MS $s$th symbol ($\hat{d}_{ks}$) and $k$th user ($\hat{\mathbf{d}}_k$) are given by

$$\hat{d}_{ks} = \mathbf{w}_{ks}^H(\mathbf{H}_k^H \sum_{i=1}^{K} \mathbf{B}_i \mathbf{d}_i + \mathbf{n}_k) = \mathbf{w}_{ks}^H(\mathbf{H}_k^H \mathbf{B}\mathbf{d} + \mathbf{n}_k) \quad (4.1)$$

$$\hat{\mathbf{d}}_k = \mathbf{W}_k^H(\mathbf{H}_k^H \mathbf{B}\mathbf{d} + \mathbf{n}_k) \quad (4.2)$$

where $\mathbf{H}_k^H \in \mathbb{C}^{M_k \times N}$ is the channel matrix between the BS and $k$th MS, $\mathbf{W}_k = [\mathbf{w}_{k1}, \cdots \mathbf{w}_{kS_k}]$, $\mathbf{B}_k = [\mathbf{b}_{k1}, \cdots \mathbf{b}_{kS_k}]$ and $\mathbf{n}_k$ is the $k$th MS additive noise. Without loss of generality, we can assume that the entries of $\mathbf{d}_k$ are independent and identically distributed (i.i.d) ZMCSCG random variables all with unit variance, i.e., $\mathrm{E}\{\mathbf{d}_k \mathbf{d}_k^H\} = \mathbf{I}_{S_k}$, $\mathrm{E}\{\mathbf{d}_k \mathbf{d}_i^H\} = \mathbf{0}$, $\forall i \neq k$, and $\mathrm{E}\{\mathbf{d}_k \mathbf{n}_i^H\} = \mathbf{0}$, $\forall i, k$. The $k$th MS noise vector is a ZMCSCG random variable with covariance matrix $\mathbf{R}_{nk} \in \mathbb{C}^{M_k \times M_k}$.

To establish our MSE downlink-interference duality, we model the virtual interference channel (Fig. 4.1.[lower]) is modeled by introducing precoders $\{\mathbf{V}_k = [\mathbf{v}_{k1}, \cdots, \mathbf{v}_{kS_k}]\}_{k=1}^{K}$ and decoders $\{\mathbf{T}_k = [\mathbf{t}_{k1}, \cdots, \mathbf{t}_{kS_k}]\}_{k=1}^{K}$, where $\mathbf{v}_{ks} \in \mathbb{C}^{M_k \times 1}$ and $\mathbf{t}_{ks} \in \mathbb{C}^{N \times 1}, \forall k, s$. In this channel, it is assumed that the $k$th user's $s$th symbol ($d_{ks}$) is an i.i.d ZMCSCG random variable with variance $\zeta_{ks}$ and estimated independently by $\mathbf{t}_{ks} \in \mathcal{C}^{N \times 1}$, i.e., $\mathrm{E}\{d_{ks} d_{ks}^H\} = \zeta_{ks}$, $\mathrm{E}\{d_{ks}d_{ij}^H\} = 0, \forall(i,j) \neq (k,s)$, and $\mathrm{E}\{\mathbf{d}_k\mathbf{n}_i^H\} = \mathbf{0}, \forall i,k$. Moreover, $\{\mathbf{n}_{ks}^I, \forall s\}_{k=1}^K$ (Fig. 4.1.[lower]) are also ZMCSCG random variables with covariance matrices $\{\mathbf{\Delta}_{ks} \in \Re^{N\times N} = \mathrm{diag}(\delta_{ks1}, \cdots, \delta_{ksN}), \forall s\}_{k=1}^K$ and the channels between the $k$th transmitter and all receivers are the same (i.e., $\{\mathbf{H}_{kjs} = \mathbf{H}_k, \forall j,s\}_{k=1}^K$).

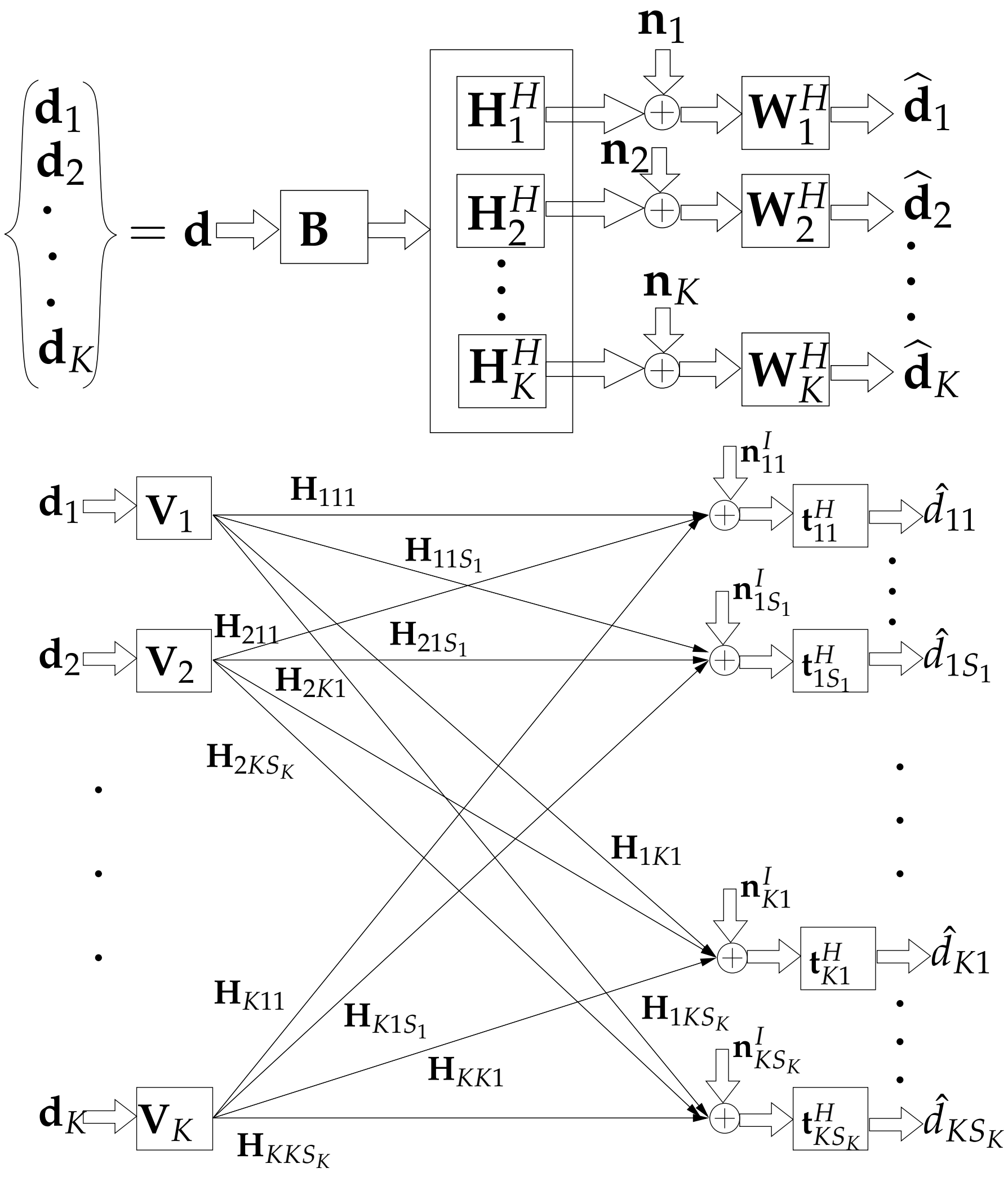

**Figure 4.1** Multiuser MIMO system model. [upper] Downlink channel. [lower] Virtual interference channel.

As can be seen from Fig. 4.1, the outputs of Fig. 4.1.[upper] and Fig. 4.1.[lower] are not the same. However, since Fig. 4.1.[lower] is a "virtual" interference channel which is introduced just to solve the downlink MSE-based problems by duality approach, the output of Fig. 4.1.[lower] is not required in practice. For this reason, the difference in the outputs of the downlink and interference channels of Fig. 4.1 will not affect the downlink MSE-based problem formulations and the duality based solutions.

For the downlink system model of Fig. 4.1, the symbol-wise and user-wise MSEs can be expressed as

$$\begin{aligned}\xi_{ks}^{DL} =& \mathrm{E}_{\mathbf{d}}\{(\hat{d}_{ks} - d_{ks})(\hat{d}_{ks} - d_{ks})^H\}\\ =& \mathbf{w}_{ks}^H(\mathbf{H}_k^H\mathbf{B}\mathbf{B}^H\mathbf{H}_k + \mathbf{R}_{nk})\mathbf{w}_{ks} - \mathbf{w}_{ks}^H\mathbf{H}_k^H\mathbf{b}_{ks} - \mathbf{b}_{ks}^H\mathbf{H}_k\mathbf{w}_{ks} + 1 \qquad (4.3)\\ \xi_k^{DL} =& \mathrm{E}_{\mathbf{d}}\{(\hat{\mathbf{d}}_k - \mathbf{d}_k)(\hat{\mathbf{d}}_k - \mathbf{d}_k)^H\}\\ =& \mathrm{tr}\{\mathbf{I}_{S_k} + \mathbf{W}_k^H(\mathbf{H}_k^H\mathbf{B}\mathbf{B}^H\mathbf{H}_k + \mathbf{R}_{nk})\mathbf{W}_k - \mathbf{W}_k^H\mathbf{H}_k^H\mathbf{B}_k - \mathbf{B}_k^H\mathbf{H}_k\mathbf{W}_k\}. \qquad (4.4)\end{aligned}$$

Using these two equations, the symbol-wise and user-wise WSMSEs can be expressed as

$$\begin{aligned}\xi_{ws}^{DL} =& \sum_{k=1}^K\sum_{s=1}^{S_k}\eta_{ks}\xi_{ks}^{DL}\\ =& \mathrm{tr}\{\boldsymbol{\eta} + \boldsymbol{\eta}\mathbf{W}^H\mathbf{H}^H\mathbf{B}\mathbf{B}^H\mathbf{H}\mathbf{W} + \boldsymbol{\eta}\mathbf{W}^H\mathbf{R}_n\mathbf{W} - \boldsymbol{\eta}\mathbf{W}^H\mathbf{H}^H\mathbf{B} - \boldsymbol{\eta}\mathbf{B}^H\mathbf{H}\mathbf{W}\} \qquad (4.5)\\ \xi_{wu}^{DL} =& \sum_{k=1}^K\tilde{\eta}_k\xi_k^{DL}\\ =& \mathrm{tr}\{\tilde{\boldsymbol{\eta}} + \tilde{\boldsymbol{\eta}}\mathbf{W}^H\mathbf{H}^H\mathbf{B}\mathbf{B}^H\mathbf{H}\mathbf{W} + \tilde{\boldsymbol{\eta}}\mathbf{W}^H\mathbf{R}_n\mathbf{W} - \tilde{\boldsymbol{\eta}}\mathbf{W}^H\mathbf{H}^H\mathbf{B} - \tilde{\boldsymbol{\eta}}\mathbf{B}^H\mathbf{H}\mathbf{W}\} \qquad (4.6)\end{aligned}$$

where $\mathbf{R}_n = \mathrm{blkdiag}(\mathbf{R}_{n1}, \cdots, \mathbf{R}_{nK})$, $\boldsymbol{\eta} = \mathrm{diag}(\eta_{11}, \cdots, \eta_{1S_1}, \cdots, \eta_{K1}, \cdots, \eta_{KS_K})$ and $\tilde{\boldsymbol{\eta}} = \mathrm{blkdiag}(\tilde{\eta}_1\mathbf{I}_{S_1}, \cdots, \tilde{\eta}_K\mathbf{I}_{S_K})$ with $\eta_{ks}$ and $\tilde{\eta}_k$ are the MSE weights of the $k$th user $s$th symbol and $k$th user, respectively. Like in the downlink channel, the interference channel symbol, user MSE and WSMSEs are expressed as

$$\xi_{ks}^{I} = \mathbf{t}_{ks}^{H}\mathbf{\Gamma}_c\mathbf{t}_{ks} + \mathbf{t}_{ks}^{H}\mathbf{\Delta}_{ks}\mathbf{t}_{ks} - \mathbf{t}_{ks}^{H}\mathbf{H}_k\mathbf{v}_{ks}\zeta_{ks} - \zeta_{ks}\mathbf{v}_{ks}^{H}\mathbf{H}_k^{H}\mathbf{t}_{ks} + \zeta_{ks} \tag{4.7}$$

$$\xi_{k}^{I} = \mathrm{tr}\{\mathbf{T}_k^{H}\mathbf{\Gamma}_c\mathbf{T}_k - \mathbf{T}_k^{H}\mathbf{H}_k\mathbf{V}_k\boldsymbol{\zeta}_k - \boldsymbol{\zeta}_k\mathbf{V}_k^{H}\mathbf{H}_k^{H}\mathbf{T}_k + \boldsymbol{\zeta}_k\} + \sum_{s=1}^{S_k}\mathbf{t}_{ks}^{H}\mathbf{\Delta}_{ks}\mathbf{t}_{ks} \tag{4.8}$$

$$\begin{aligned}\xi_{ws}^{I} &= \sum_{k=1}^{K}\sum_{s=1}^{S_k}\lambda_{ks}\xi_{ks}^{I} = \mathrm{tr}\{\boldsymbol{\lambda}\mathbf{T}^{H}\mathbf{\Gamma}_c\mathbf{T} - \boldsymbol{\lambda}\mathbf{T}^{H}\mathbf{H}\mathbf{V}\boldsymbol{\zeta} - \boldsymbol{\lambda}\boldsymbol{\zeta}\mathbf{V}^{H}\mathbf{H}^{H}\mathbf{T} + \boldsymbol{\lambda}\boldsymbol{\zeta}\} \\ &\quad + \sum_{k=1}^{K}\sum_{s=1}^{S_k}\lambda_{ks}\mathbf{t}_{ks}^{H}\mathbf{\Delta}_{ks}\mathbf{t}_{ks}\end{aligned} \tag{4.9}$$

$$\begin{aligned}\xi_{wu}^{I} &= \sum_{k=1}^{K}\tilde{\lambda}_k\xi_k^{I} = \mathrm{tr}\{\tilde{\boldsymbol{\lambda}}\mathbf{T}^{H}\mathbf{\Gamma}_c\mathbf{T} - \tilde{\boldsymbol{\lambda}}\mathbf{T}^{H}\mathbf{H}\mathbf{V}\boldsymbol{\zeta} - \tilde{\boldsymbol{\lambda}}\boldsymbol{\zeta}\mathbf{V}^{H}\mathbf{H}^{H}\mathbf{T} + \tilde{\boldsymbol{\lambda}}\boldsymbol{\zeta}\} \\ &\quad + \sum_{k=1}^{K}\sum_{s=1}^{S_k}\tilde{\lambda}_k\mathbf{t}_{ks}^{H}\mathbf{\Delta}_{ks}\mathbf{t}_{ks}\end{aligned} \tag{4.10}$$

where $\boldsymbol{\zeta}_k = \mathrm{diag}(\zeta_{k1}, \cdots, \zeta_{kS_k})$, $\boldsymbol{\zeta} = \mathrm{blkdiag}(\boldsymbol{\zeta}_1, \cdots, \boldsymbol{\zeta}_K)$, $\boldsymbol{\lambda} = \mathrm{diag}(\lambda_{11}, \cdots, \lambda_{1S_1}, \cdots, \lambda_{K1}, \cdots, \lambda_{KS_K})$, $\tilde{\boldsymbol{\lambda}} = \mathrm{blkdiag}(\tilde{\lambda}_1\mathbf{I}_{S_1}, \cdots, \tilde{\lambda}_K\mathbf{I}_{S_K})$ and $\mathbf{\Gamma}_c = \sum_{i=1}^{K}\sum_{j=1}^{S_i}\zeta_{ij}\mathbf{H}_i\mathbf{v}_{ij}\mathbf{v}_{ij}^{H}\mathbf{H}_i^{H}$ with $\lambda_{ks}$ and $\tilde{\lambda}_k$ are the MSE weights of the $k$th user $s$th symbol and $k$th user, respectively.

## 4.3 Problem formulation

The aforementioned MSE-based optimization problems can be formulated as

$$\begin{aligned}\mathcal{P}4.1: &\min_{\{\mathbf{B}_k,\mathbf{W}_k\}_{k=1}^{K}} \sum_{k=1}^{K}\sum_{s=1}^{S_k}\eta_{ks}\xi_{ks}^{DL}, \\ &\text{s.t } [\mathbf{B}\mathbf{B}^{H}]_{(n,n)} \leq \check{p}_n,\ \mathbf{b}_{ks}^{H}\mathbf{b}_{ks} \leq \check{p}_{ks},\ \forall n,k,s\end{aligned} \tag{4.11}$$

$$\begin{aligned}\mathcal{P}4.2: &\min_{\{\mathbf{B}_k,\mathbf{W}_k\}_{k=1}^{K}} \sum_{k=1}^{K}\tilde{\eta}_k\xi_k^{DL}, \\ &\text{s.t } [\mathbf{B}\mathbf{B}^{H}]_{(n,n)} \leq \check{p}_n,\ \mathrm{tr}\{\mathbf{B}_k^{H}\mathbf{B}_k\} \leq \hat{p}_k,\ \forall n,k\end{aligned} \tag{4.12}$$

$$\begin{aligned}\mathcal{P}4.3: &\min_{\{\mathbf{B}_k,\mathbf{W}_k\}_{k=1}^{K}} \max\ \rho_{ks}\xi_{ks}^{DL}, \\ &\text{s.t } [\mathbf{B}\mathbf{B}^{H}]_{(n,n)} \leq \check{p}_n,\ \mathbf{b}_{ks}^{H}\mathbf{b}_{ks} \leq \check{p}_{ks},\ \forall n,k,s\end{aligned} \tag{4.13}$$

$$\mathcal{P}4.4: \min_{\{\mathbf{B}_k,\mathbf{W}_k\}_{k=1}^{K}} \max\ \tilde{\rho}_k\xi_k^{DL},$$

$$\text{s.t } [\mathbf{B}\mathbf{B}^H]_{(n,n)} \leq \check{p}_n, \ \text{tr}\{\mathbf{B}_k^H\mathbf{B}_k\} \leq \hat{p}_k, \ \forall n, k \tag{4.14}$$

where $\tilde{\rho}_k(\hat{p}_k)$ and $\rho_{ks}(\check{p}_{ks})$ are the MSE balancing weights (maximum available power) of the $k$th user and $k$th user $s$th symbol, respectively, and $\check{p}_n$ denotes the maximum transmitted power by the $n$th antenna.

For both the WSMSE minimization and min max WMSE problems, different weights are given to different symbols (users). However, at optimality the solutions of these two problems are not necessarily the same. This is due to the fact that the aim of the WSMSE minimization problem is just to minimize the WSMSE of all symbols (users) (i.e., in such a problem the minimized WMSE of each symbol (user) depends on its corresponding channel gain), whereas the aim of min max WMSE problem is to minimize and balance the WMSE of each symbol (user) simultaneously (i.e., in such a problem all symbols (users) achieve the same minimized WMSE [SSB08a]). Moreover, as will be clear later, the solution approach of WSMSE minimization problem can not be extended straightforwardly to solve the min max WMSE problem. Due to these facts, we examine the WSMSE minimization and min max WMSE problems separately.

Since, the problems $\mathcal{P}4.1$ - $\mathcal{P}4.4$ are not convex, convex optimization framework can not be applied to solve them. To the best of our knowledge, duality based solutions for these problems are not known. In the following, we present an MSE downlink-interference duality based approach for solving each of these problems which is shown in **Algorithm 4.I**[2].

**Algorithm 4.I**

Initialization: For each problem, initialize $\{\mathbf{B}_k \neq \mathbf{0}\}_{k=1}^K$ such that the power constraint functions are satisfied[3]. Then, update $\{\mathbf{W}_k\}_{k=1}^K$ by using MMSE receiver approach, i.e.,

$$\mathbf{W}_k = (\mathbf{H}_k^H\mathbf{B}\mathbf{B}^H\mathbf{H}_k + \mathbf{R}_{nk})^{-1}\mathbf{H}_k^H\mathbf{B}_k, \ \forall k. \tag{4.15}$$

[2]As will be clear later in Section 4.8, to solve $\mathcal{P}4.3$ and $\mathcal{P}4.4$ (and more general MSE-based problems), an additional power allocation step is required. In **Algorithm I**, this step is omitted for clarity of presentation.

[3]For the simulation, we use $\{\mathbf{B}_k = [\mathbf{H}_k]_{(:,1:S_k)}\}_{k=1}^K$ followed by the appropriate normalization of $\{\mathbf{B}_k\}_{k=1}^K$ to ensure the power constraints.

**Repeat**: **Interference channel**

1. Transfer the symbol-wise (user-wise) WSMSE or WMSE from downlink to interference channel.

2. Update the receivers of the interference channel $\{\mathbf{t}_{ks}, \forall s\}_{k=1}^{K}$ using MMSE receiver technique.

**Downlink channel**

3. Transfer the symbol-wise (user-wise) WSMSE or WMSE from interference to downlink channel.

4. Update the receivers of the downlink channel $\{\mathbf{W}_k\}_{k=1}^{K}$ by MMSE receiver approach (4.15).

**Until** convergence.

The above iterative algorithm is already known in [HJU09], [BV11a] and [BV11b]. However, the approaches of these papers can not ensure the power constraints of $\mathcal{P}4.1$ - $\mathcal{P}4.4$ at step 3 of **Algorithm 4.I**. Hence, one can not apply the approaches of these papers to solve $\mathcal{P}4.1$ - $\mathcal{P}4.4$. In the following sections, we establish our MSE downlink-interference duality.

## 4.4 Symbol-wise WSMSE downlink-interference duality

This duality is established to solve symbol-wise WSMSE-based problems (for example $\mathcal{P}4.1$).

### 4.4.1 Symbol-wise WSMSE transfer (From downlink to interference channel)

In order to use this WSMSE transfer for solving $\mathcal{P}4.1$, we set the interference channel precoder, decoder, noise covariance, input covariance and MSE weight matrices as

$$\mathbf{V} = \bar{\beta}\mathbf{W},\ \mathbf{T} = \mathbf{B}/\bar{\beta},\ \boldsymbol{\zeta} = \boldsymbol{\eta},\ \boldsymbol{\lambda} = \mathbf{I},\ \boldsymbol{\Delta}_{ks} = \boldsymbol{\Psi} + \mu_{ks}\mathbf{I} \tag{4.16}$$

where $\bar{\beta}$, $\{\psi_n\}_{n=1}^N$ and $\{\mu_{ks}, \forall s\}_{k=1}^K$ are positive real scalars that will be determined in the sequel and $\mathbf{\Psi} = \text{diag}(\psi_1, \cdots, \psi_N)$. Substituting (4.16) into (4.9) and equating $\xi_{ws}^I = \xi_{ws}^{DL}$ yields

$$\text{tr}\{\mathbf{B}^H\mathbf{H}\mathbf{W}\boldsymbol{\eta}\mathbf{W}^H\mathbf{H}^H\mathbf{B} - \mathbf{B}^H\mathbf{H}\mathbf{W}\boldsymbol{\eta} - \boldsymbol{\eta}\mathbf{W}^H\mathbf{H}^H\mathbf{B} + \boldsymbol{\eta}\}+$$
$$\frac{1}{\bar{\beta}^2}\sum_{k=1}^K\sum_{s=1}^{S_k}\mathbf{b}_{ks}^H(\mathbf{\Psi} + \mu_{ks}\mathbf{I}_N)\mathbf{b}_{ks} = \text{tr}\{\boldsymbol{\eta}\mathbf{W}^H\mathbf{H}^H\mathbf{B}\mathbf{B}^H\mathbf{H}\mathbf{W}$$
$$+ \boldsymbol{\eta}\mathbf{W}^H\mathbf{R}_n\mathbf{W} - \boldsymbol{\eta}\mathbf{W}^H\mathbf{H}^H\mathbf{B} - \boldsymbol{\eta}\mathbf{B}^H\mathbf{H}\mathbf{W} + \boldsymbol{\eta}\}.$$

It follows

$$\bar{\beta}^2\tau = \sum_{n=1}^N \psi_n\check{p}_n + \sum_{k=1}^K\sum_{s=1}^{S_k}\mu_{ks}\bar{p}_{ks} = \check{\mathbf{p}}^T\boldsymbol{\psi} + \bar{\mathbf{p}}^T\boldsymbol{\mu} \tag{4.17}$$

where $\tau = \text{tr}\{\boldsymbol{\eta}\mathbf{W}^H\mathbf{R}_n\mathbf{W}\}$, $\boldsymbol{\psi} = [\psi_1, \cdots, \psi_N]^T$, $\boldsymbol{\mu} = [\mu_{11}, \cdots, \mu_{1S_1}, \cdots, \mu_{K1}, \cdots, \mu_{KS_K}]^T$, $\check{\mathbf{p}} = [\check{p}_1, \cdots, \check{p}_N]^T$ and $\bar{\mathbf{p}} = [\bar{p}_{11}, \cdots, \bar{p}_{1S_1}, \cdots, \bar{p}_{K1}, \cdots, \bar{p}_{KS_K}]^T$ with $\bar{p}_{ks} = \mathbf{b}_{ks}^H\mathbf{b}_{ks}$, $\check{p}_n = \check{\mathbf{b}}_n^H\check{\mathbf{b}}_n$ and $\check{\mathbf{b}}_n^H$ is the $n$th row of $\mathbf{B}$.

The above equation shows that by choosing any $\{\psi_n\}_{n=1}^N$ and $\{\mu_{ks}, \forall s\}_{k=1}^K$ that satisfy (4.17), one can transfer the downlink channel precoder/decoder to the interference channel decoder/precoder ensuring $\xi_{ws}^{DL} = \xi_{ws}^{I_1}$, where $\xi_w^{I_1}$ is the interference WSMSE at step 1 of **Algorithm 4.I**. However, here $\{\psi_n\}_{n=1}^N$ and $\{\mu_{ks}, \forall s\}_{k=1}^K$ should be selected in a way that $\mathcal{P}4.1$ can be solved by **Algorithm 4.I**. To this end, we choose $\boldsymbol{\psi}$ and $\boldsymbol{\mu}$ as

$$\bar{\beta}^2\tau \geq \check{\mathbf{p}}^T\boldsymbol{\psi} + \bar{\mathbf{p}}^T\boldsymbol{\mu}. \tag{4.18}$$

By doing so, the interference channel symbol-wise WSMSE is upper bounded by that of the downlink channel (i.e., $\xi_{ws}^{I_1} \leq \xi_{ws}^{DL}$). As will be clear later, to solve (4.11) with **Algorithm 4.I**, $\bar{\beta}$, $\boldsymbol{\psi}$ and $\boldsymbol{\mu}$ should be selected as in (4.18). This shows that step 1 of **Algorithm 4.I** can be carried out with (4.16). To perform step 2 of **Algorithm 4.I**, we update $\mathbf{t}_{ks}$ of (4.16) by using the interference channel MMSE receiver approach which is expressed as

$$\begin{aligned}\mathbf{t}_{ks} =& (\mathbf{\Gamma}_c + \mathbf{\Delta}_{ks})^{-1}\mathbf{H}_k\mathbf{v}_{ks}\zeta_{ks}\\ =& \bar{\beta}(\mathbf{H}\mathbf{W}\boldsymbol{\eta}\mathbf{W}^H\mathbf{H}^H + \mathbf{\Psi} + \mu_{ks}\mathbf{I})^{-1}\mathbf{H}_k\mathbf{w}_{ks}\eta_{ks}\end{aligned} \tag{4.19}$$

where the second equality is obtained from (4.16). The above expression shows that by choosing $\{\mu_{ks} > 0, \forall s\}_{k=1}^{K}$, $\{\psi_n > 0\}_{n=1}^{N}$, we ensure $(\mathbf{H}\mathbf{W}\boldsymbol{\eta}\mathbf{W}^H\mathbf{H}^H + \boldsymbol{\Psi} + \mu_{ks}\mathbf{I})^{-1}$ exists. Next, we transfer the symbol-wise WSMSE from interference to downlink channel by ensuring the power constraint of $\mathcal{P}4.1$ (i.e., we perform step 3 of **Algorithm 4.I**).

### 4.4.2 Symbol-wise WSMSE transfer (From interference to downlink channel)

For a given symbol-wise WSMSE in the interference channel with $\boldsymbol{\zeta} = \boldsymbol{\eta}$ and $\boldsymbol{\lambda} = \mathbf{I}$, we can achieve the same WSMSE in the downlink channel (with the MSE weighting matrix $\boldsymbol{\eta}$) using a nonzero scaling factor ($\beta$) satisfying

$$\widetilde{\mathbf{B}} = \beta\mathbf{T}, \ \ \widetilde{\mathbf{W}} = \mathbf{V}/\beta. \tag{4.20}$$

In this precoder/decoder transformation, we use the notations $\widetilde{\mathbf{B}}$ and $\widetilde{\mathbf{W}}$ to differentiate from the precoder and decoder matrices used in Section 4.4.1. By substituting (4.20) into $\xi_{ws}^{DL}$ (with $\widetilde{\mathbf{B}}$=$\mathbf{B}$, $\widetilde{\mathbf{W}}$=$\mathbf{W}$), equating the resulting symbol-wise WSMSE with that of the interference channel (4.9) and after some simple manipulations, we get

$$\begin{aligned}
&\sum_{k=1}^{K}\sum_{s=1}^{S_k}\mathbf{t}_{ks}^H(\boldsymbol{\Psi}+\mu_{ks}\mathbf{I}_N)\mathbf{t}_{ks} = \frac{1}{\beta^2}\mathrm{tr}\{\boldsymbol{\eta}\mathbf{V}^H\mathbf{R}_n\mathbf{V}\}\\
\Rightarrow \beta^2 =& \frac{\mathrm{tr}\{\boldsymbol{\eta}\mathbf{V}^H\mathbf{R}_n\mathbf{V}\}}{\sum_{k=1}^{K}\sum_{s=1}^{S_k}\mathbf{t}_{ks}^H(\boldsymbol{\Psi}+\mu_{ks}\mathbf{I}_N)\mathbf{t}_{ks}}\\
=& \frac{\bar{\beta}^2\tau}{\sum_{n=1}^{N}\psi_n\check{\mathbf{t}}_n^H\check{\mathbf{t}}_n + \sum_{k=1}^{K}\sum_{i=1}^{S_i}\mu_{ki}\mathbf{t}_{ki}^H\mathbf{t}_{ki}}
\end{aligned} \tag{4.21}$$

where $\check{\mathbf{t}}_n^H$ is the $n$th row of the MMSE matrix $\mathbf{T}$ (4.19) and the third equality follows from (4.16). The power constraints of each BS antenna and symbol in the downlink channel are thus given by

$$\begin{aligned}
\check{\widetilde{\mathbf{b}}}_n^H\check{\widetilde{\mathbf{b}}}_n =& \beta^2\check{\mathbf{t}}_n^H\check{\mathbf{t}}_n && (4.22)\\
=& \frac{\bar{\beta}^2\tau\check{\mathbf{t}}_n^H\check{\mathbf{t}}_n}{\sum_{i=1}^{N}\psi_i\check{\mathbf{t}}_i^H\check{\mathbf{t}}_i + \sum_{i=1}^{K}\sum_{j=1}^{S_i}\mu_{ij}\mathbf{t}_{ij}^H\mathbf{t}_{ij}} \leq \check{p}_n, \ \forall n\\
\widetilde{\mathbf{b}}_{ks}^H\widetilde{\mathbf{b}}_{ks} =& \beta^2\mathbf{t}_{ks}^H\mathbf{t}_{ks} && (4.23)
\end{aligned}$$

$$=\frac{\bar{\beta}^2\tau\mathbf{t}_{ks}^H\mathbf{t}_{ks}}{\sum_{i=1}^N \psi_i\check{\mathbf{t}}_i^H\check{\mathbf{t}}_i+\sum_{i=1}^K\sum_{j=1}^{S_i}\mu_{ij}\mathbf{t}_{ij}^H\mathbf{t}_{ij}} \leq \check{p}_{ks}, \forall k,s$$

where $\check{\tilde{\mathbf{b}}}_n^H$ is the $n$th row of $\widetilde{\mathbf{B}}$. By multiplying both sides of (4.22) and (4.23) with $\psi_n, \forall n$ and $\mu_{ks}, \forall k,s$, we get

$$\psi_n \geq \check{f}_n \text{ and } \mu_{ks} \geq f_{ks}, \forall n,k,s \tag{4.24}$$

where $\check{f}_n = \frac{\bar{\beta}^2\tau}{\check{p}_n}\frac{\psi_n\check{\mathbf{t}}_n^H\check{\mathbf{t}}_n}{\sum_{i=1}^N \psi_i\check{\mathbf{t}}_i^H\check{\mathbf{t}}_i+\sum_{i=1}^K\sum_{j=1}^{S_i}\mu_{ij}\mathbf{t}_{ij}^H\mathbf{t}_{ij}}$ and $f_{ks} = \frac{\bar{\beta}^2\tau}{\check{p}_{ks}}\frac{\mu_{ks}\mathbf{t}_{ks}^H\mathbf{t}_{ks}}{\sum_{i=1}^N \psi_i\check{\mathbf{t}}_i^H\check{\mathbf{t}}_i+\sum_{i=1}^K\sum_{j=1}^{S_i}\mu_{ij}\mathbf{t}_{ij}^H\mathbf{t}_{ij}}$. Now, for any given $\bar{\beta}$, $\{\check{\mathbf{t}}_n^H\check{\mathbf{t}}_n\}_{n=1}^N$ and $\{\mathbf{t}_{ks}^H\mathbf{t}_{ks}, \forall s\}_{k=1}^K$, suppose that there exist $\{\psi_n > 0\}_{n=1}^N$ and $\{\mu_{ks} > 0, \forall s\}_{k=1}^K$ that satisfy

$$\psi_n = \check{f}_n \text{ and } \mu_{ks} = f_{ks}, \forall n,k,s. \tag{4.25}$$

From the above equation one can also achieve $\psi_n\check{p}_n = \check{f}_n\check{p}_n$, $\mu_{ks}\check{p}_{ks} = f_{ks}\check{p}_{ks}$ $\forall n,k,s$. Summing up these expressions for all $n,k$ and $s$ results

$$\begin{aligned}\sum_{n=1}^N \psi_n\check{p}_n + \sum_{k=1}^K\sum_{s=1}^{S_k}\mu_{ks}\check{p}_{ks} &= \sum_{n=1}^N \check{f}_n\check{p}_n + \sum_{k=1}^K\sum_{s=1}^{S_k} f_{ks}\check{p}_{ks} \\ &= \bar{\beta}^2\tau. \end{aligned}\tag{4.26}$$

This equation shows that the solution of (4.25) satisfies (4.26). Moreover, as $\{\check{p}_n \geq \check{p}_n\}_{n=1}^N$ and $\{\check{p}_{ks} \geq \bar{p}_{ks}, \forall s\}_{k=1}^K$, the latter solution also ensures (4.18). Therefore, by choosing $\{\psi_n\}_{n=1}^N$ and $\{\mu_{ks}, \forall s\}_{k=1}^K$ such that (4.25) is satisfied, step 3 of **Algorithm 4.I** can be performed. Furthermore, one can notice from (4.26) that $\bar{\beta}^2$ can be any positive value.

Next, we show that there exists at least a set of feasible $\{\psi_n > 0\}_{n=1}^N$ and $\{\mu_{ks} > 0, \forall s\}_{k=1}^K$ that satisfy (4.25). To this end, we consider the following Theorem [Ber08].

*Theorem 4.1:* Let $(\mathbf{X}, \|.\|_2)$ be a complete metric space. We say that $\mathcal{F} : \mathbf{X} \rightarrow \mathbf{X}$ is an almost contraction, if there exist $\kappa(\tilde{\kappa}) \in [0,1)$ and $\chi(\tilde{\chi}) \geq 0$ such that

$$\begin{aligned}&\|\mathcal{F}(\mathbf{x}) - \mathcal{F}(\mathbf{y})\|_2 \leq \kappa\|\mathbf{x}-\mathbf{y}\|_2 + \chi\|\mathbf{y}-\mathcal{F}(\mathbf{x})\|_2, \text{ or} \\ &\|\mathcal{F}(\mathbf{x}) - \mathcal{F}(\mathbf{y})\|_2 \leq \tilde{\kappa}\|\mathbf{x}-\mathbf{y}\|_2 + \tilde{\chi}\|\mathbf{x}-\mathcal{F}(\mathbf{y})\|_2, \forall \mathbf{x},\mathbf{y} \in \mathbf{X}.\end{aligned}\tag{4.27}$$

If $\mathcal{F}$ satisfies (4.27), then the following holds true:

1. $\exists \mathbf{x} \in \mathbf{X} : \mathbf{x} = \mathcal{F}(\mathbf{x})$.
2. For any initial $\mathbf{x}_0 \in \mathbf{X}$, the iteration $\mathbf{x}_{n+1} = \mathcal{F}(\mathbf{x}_n)$ for n = 0, 1, 2, $\cdots$ converges to some $\mathbf{x}^\star \in \mathbf{X}$.
3. The solution $\mathbf{x}^\star$ is not necessarily unique.

**Proof.** See *Theorem 1.1* of [Ber08]. □

Note that according to [Ber04] (see (1.1) and (1.2) of [Ber04]), the two inequalities of (4.27) are dual to each other.

Define $\mathbf{x}$ and $\mathcal{F}$ as $\mathbf{x} \triangleq [x_1, \cdots, x_{S+N}]^T = [\psi_1, \cdots, \psi_N, \mu_{11}, \cdots, \mu_{1S_1} \cdots, \mu_{K1}, \cdots, \mu_{KS_K}]^T$, $\mathcal{F}(\mathbf{x}) \triangleq [\check{f}_1, \cdots, \check{f}_N, f_{11}, \cdots, f_{1S_1}, \cdots, f_{K1}, \cdots, f_{KS_K}]$ with $\{x_n = \psi_n \in [\epsilon, (\bar{\beta}^2\tau - \epsilon\sum_{i=1,\, i\neq n}^{N} p_{im})/p_{nm}]\}_{n=1}^{N}$ and $\{x_r\}_{r=N+1}^{S+N} = \{\mu_{ks} =\in [\epsilon, (\bar{\beta}^2\tau - \epsilon\sum_{i=1}^{K}\sum_{j=1,(i,j)\neq(k,s)}^{S_i} p_{ijm})/p_{ksm}], \forall s\}_{k=1}^{K}$[4]. As we can see from (4.27), when $\|\mathcal{F}(\mathbf{x}_1) - \mathcal{F}(\mathbf{x}_2)\|_2 = 0$ with $\mathbf{x}_1 = \mathbf{x}_2$ or $\mathbf{x}_1 \neq \mathbf{x}_2$, one can set $\kappa(\tilde{\kappa}) = 0$ and $\chi(\tilde{\chi}) = 0$ to satisfy this inequality. And when $\|\mathcal{F}(\mathbf{x}_1) - \mathcal{F}(\mathbf{x}_2)\|_2 > 0$ (i.e., $\mathbf{x}_1 \neq \mathbf{x}_2$), one can select appropriate $\kappa(\tilde{\kappa}) \in [0, 1)$ and $\chi(\tilde{\chi}) \geq 0$ such that (4.27) is satisfied. This is due to the fact that in the latter case, $\|\mathbf{x}_2 - \mathcal{F}(\mathbf{x}_1)\|_2 > 0$ and/or $\|\mathbf{x}_1 - \mathcal{F}(\mathbf{x}_2)\|_2 > 0$ and $\|\mathbf{x}_1 - \mathbf{x}_2\|_2 > 0$ are positive and bounded. This explanation shows the existence of $\kappa(\tilde{\kappa}) \in [0, 1)$ and $\chi(\tilde{\chi}) \geq 0$ ensuring (4.27) for any $\|\mathcal{F}(\mathbf{x}_1) - \mathcal{F}(\mathbf{x}_2)\|_2$, $\mathbf{x}_1, \mathbf{x}_2 \in \mathbf{X}$ . Consequently, $\mathcal{F}(\mathbf{x})$ is an almost contraction which implies

$$\mathbf{x}_{n+1} = \mathcal{F}(\mathbf{x}_n),\ \mathbf{x}_0 = [x_{01}, x_{02}, \cdots, x_{0(S+N)}]^T \geq \epsilon\mathbf{1}_{N+S}, \quad \text{for } n = 0, 1, 2, \cdots \text{ converges} \tag{4.28}$$

where $\mathbf{1}_{N+S}$ is an $N+S$ length vector with each element equal to unity. Thus, there exist $\{\psi_n \geq \epsilon\}_{n=1}^{N}$ and $\{\mu_{ks} \geq \epsilon, \forall s\}_{k=1}^{K}$ that satisfy (4.25) and can be computed using (4.28). For numerical simulation we initialize $\mathbf{x}_0$ as $x_{01} = x_{02} = \cdots = x_{0(S+N)}$. However, finding the optimal initialization strategy is still an open research topic.

Once the appropriate $\{\psi_n\}_{n=1}^{N}$ and $\{\mu_{ks}\forall s\}_{k=1}^{K}$ are obtained, step 4 of **Algorithm 4.I** is immediate and hence $\mathcal{P}4.1$ can be solved iteratively using this algorithm.

[4]For our simulation, we use $\epsilon = \min(10^{-6}, \{\bar{\beta}\tau/p_{nm}\}_{n=1}^{N}, \{\bar{\beta}\tau/p_{ksm}, \forall s\}_{k=1}^{K})$.

### 4.4.3 Extension of the current duality for $\mathcal{P}4.1$ with a total BS power constraint

If the constraints of $\mathcal{P}4.1$ are modified to a total BS power, the power constraint at step 3 of **Algorithm 4.I** can be ensured by applying the precoder/decoder transformation expression of [HJU09]. The precoder/decoder transformation of [HJU09] is performed by computing $S$ scaling factors. These scaling factors are obtained by solving $S$ systems of equations which require matrix inversion with complexity $O(S^3)$ (see (23) of [HJU09]).

In the current paper, if the constraints of $\mathcal{P}4.1$ are modified to a total BS power, one can ensure the power constraint at step 3 of **Algorithm 4.I** just by assigning $\Delta_{ks}$ of (4.16) as $\Delta_{ks} = \mathbf{I}$. By doing so, $\bar{\beta}^2$ of (4.17) and $\beta^2$ of (4.21) can be expressed as $\bar{\beta}^2 = \frac{\sum_{k=1}^{K}\sum_{s=1}^{S_k}\mathbf{b}_{ks}^H\mathbf{b}_{ks}}{\tau} = \frac{P_{max}}{\tau}$ and $\beta^2 = \frac{\text{tr}\{\boldsymbol{\eta}\mathbf{V}^H\mathbf{R}_n\mathbf{V}\}}{\sum_{k=1}^{K}\sum_{s=1}^{S_k}\mathbf{t}_{ks}^H\mathbf{t}_{ks}}$, where $P_{max}$ is the total BS power. Now by employing (4.20), the total BS power at step 3 of **Algorithm 4.I** can thus be given as $\text{tr}\{\widetilde{\mathbf{B}}\widetilde{\mathbf{B}}^H\} = \beta^2\text{tr}\{\mathbf{T}\mathbf{T}^H\} = \bar{\beta}^2\tau = P_{max}$ (i.e., the total BS power constraint is satisfied). Thus, for $\mathcal{P}4.1$ (with a total BS power constraint), we do not need to use *Theorem I*. Moreover, our duality requires only one scaling factor to perform the precoder/decoder transformation (i.e., $\beta^2(\bar{\beta}^2)$). This shows that for this problem, the proposed duality based algorithm requires less computation compared to that of [HJU09]. Note that the duality algorithm of [HJU09] requires the same computation as that of [BV11a] and less computation than that of [SSB08c] and [SSB07]. Thus, it is sufficient to compare the current duality algorithm with the duality algorithm of [HJU09].

For other WSMSE-based problems with a total BS power constraint function, the computational advantage of the current duality based algorithm over that of [HJU09] can be analysed like in this subsection.

## 4.5 User-wise WSMSE downlink-interference duality

This duality is established to solve user-wise WSMSE-based problems (for example $\mathcal{P}4.2$).

### 4.5.1 User-wise WSMSE transfer (From downlink to interference channel)

To apply this WSMSE transfer for solving $\mathcal{P}4.2$, we set the precoder, decoder and noise covariance matrices as

$$\mathbf{V} = \tilde{\beta}\mathbf{W},\ \mathbf{T} = \mathbf{B}/\tilde{\beta},\ \boldsymbol{\zeta} = \tilde{\boldsymbol{\eta}},\ \tilde{\boldsymbol{\lambda}} = \mathbf{I}, \boldsymbol{\Delta}_{ks} = \boldsymbol{\Psi} + \mu_k \mathbf{I} \tag{4.29}$$

where $\tilde{\beta}$, $\{\psi_n\}_{n=1}^{N}$ and $\{\mu_k\}_{k=1}^{K}$ are real positive scalars. Substituting (4.29) into (4.10) and equating $\xi_{wu}^{I} = \xi_{wu}^{DL}$ yields

$$\tilde{\beta}^2\tilde{\tau} = \sum_{n=1}^{N} \psi_n \check{p}_n + \sum_{k=1}^{K} \mu_k p_k = \check{\mathbf{p}}^T \boldsymbol{\psi} + \tilde{\mathbf{p}}^T \tilde{\boldsymbol{\mu}} \tag{4.30}$$

where $\tilde{\tau} = \text{tr}\{\tilde{\boldsymbol{\eta}}\mathbf{W}^H\mathbf{R}_n\mathbf{W}\}$, $\tilde{\boldsymbol{\mu}} = [\mu_1, \cdots, \mu_K]^T$, $\tilde{\mathbf{p}} = [\tilde{p}_1, \cdots, \tilde{p}_K]^T$ with $\tilde{p}_k = \text{tr}\{\mathbf{B}_k\mathbf{B}_k^H\}$. Like in Section 4.4.1, we perform step 1 of **Algorithm 4.I** by choosing $\tilde{\beta}^2$, $\boldsymbol{\psi}$ and $\tilde{\boldsymbol{\mu}}$ as

$$\tilde{\beta}\tilde{\tau} \geq \check{\mathbf{p}}^T \boldsymbol{\psi} + \tilde{\mathbf{p}}^T \tilde{\boldsymbol{\mu}}. \tag{4.31}$$

To perform step 2 of **Algorithm 4.I**, we update $\mathbf{t}_{ks}$ of (4.29) using the interference channel MMSE receiver as

$$\mathbf{t}_{ks} = \tilde{\beta}(\mathbf{H}\mathbf{W}\tilde{\boldsymbol{\eta}}\mathbf{W}^H\mathbf{H}^H + \boldsymbol{\Psi} + \mu_k\mathbf{I})^{-1}\mathbf{H}_k\mathbf{w}_{ks}\tilde{\eta}_k. \tag{4.32}$$

This expression shows that by choosing $\{\mu_k > 0\}_{k=1}^{K}$, $\{\psi_n > 0\}_{n=1}^{N}$, we ensure that $(\mathbf{H}\mathbf{W}\boldsymbol{\eta}\mathbf{W}^H\mathbf{H}^H + \boldsymbol{\Psi} + \mu_k\mathbf{I})^{-1}$ exists.

### 4.5.2 User-wise WSMSE transfer (From interference to downlink channel)

For a given user-wise WSMSE in the interference channel with $\boldsymbol{\zeta} = \tilde{\boldsymbol{\eta}}$ and $\tilde{\boldsymbol{\lambda}} = \mathbf{I}$, we can achieve the same WSMSE in the downlink channel (with the weighting matrix $\tilde{\boldsymbol{\eta}}$) by using a nonzero scaling factor ($\tilde{\tilde{\beta}}$) which satisfies

$$\widetilde{\mathbf{B}} = \tilde{\tilde{\beta}}\mathbf{T},\ \ \widetilde{\mathbf{W}} = \mathbf{V}/\tilde{\tilde{\beta}}. \tag{4.33}$$

In this precoder/decoder transformation, we use the notations $\widetilde{\mathbf{B}}$ and $\widetilde{\mathbf{W}}$ to differentiate from the precoder and decoder matrices used in Section 4.5.1.

By substituting (4.33) into $\xi_{wu}^{DL}$ (with $\widetilde{\mathbf{B}}$=$\mathbf{B}$, $\widetilde{\mathbf{W}}$=$\mathbf{W}$), then equating the resulting user-wise WSMSE with that of the interference channel ($\xi_{wu}^{I}$) and after simple manipulations, we get

$$\tilde{\tilde{\beta}}^2 = \frac{\tilde{\beta}^2\tilde{\tau}}{\sum_{n=1}^{N}\psi_n \check{\mathbf{t}}_n^H\check{\mathbf{t}}_n + \sum_{k=1}^{K}\mu_k \mathrm{tr}\{\mathbf{T}_k^H\mathbf{T}_k\}} \tag{4.34}$$

where $\check{\mathbf{t}}_n^H$ is the $n$th row of the MMSE matrix $\mathbf{T}$ (4.32). The power constraints of each BS antenna and user (i.e., step 3 of **Algorithm 4.I**) in the downlink channel can be expressed as

$$\psi_n \geq \check{\check{f}}_n \ \text{ and } \ \mu_k \geq \tilde{f}_k, \ \forall k \tag{4.35}$$

where

$$\check{\check{f}}_n = \frac{\tilde{\beta}^2\tilde{\tau}}{\check{p}_n}\frac{\psi_n\check{\mathbf{t}}_n^H\check{\mathbf{t}}_n}{\sum_{i=1}^{N}\psi_i\check{\mathbf{t}}_i^H\check{\mathbf{t}}_i + \sum_{i=1}^{K}\mu_i \mathrm{tr}\{\mathbf{T}_i^H\mathbf{T}_i\}} \tag{4.36}$$

$$\tilde{f}_k = \frac{\tilde{\beta}^2\tilde{\tau}}{\hat{p}_k}\frac{\mu_k \mathrm{tr}\{\mathbf{T}_k^H\mathbf{T}_k\}}{\sum_{i=1}^{N}\psi_i\check{\mathbf{t}}_i^H\check{\mathbf{t}}_i + \sum_{i=1}^{K}\mu_i \mathrm{tr}\{\mathbf{T}_i^H\mathbf{T}_i\}}. \tag{4.37}$$

For given $\tilde{\beta}$, $\{\check{\mathbf{t}}_n^H\check{\mathbf{t}}_n\}_{n=1}^N$ and $\{\mathrm{tr}\{\mathbf{T}_k^H\mathbf{T}_k\}\}_{k=1}^K$, one can show that there exist $\{\psi_n\}_{n=1}^N$ and $\{\mu_k\}_{k=1}^K$ which satisfy

$$\psi_n = \check{\check{f}}_n \ \text{ and } \ \mu_k = \tilde{f}_k, \ \forall n, k. \tag{4.38}$$

The solution of (4.38) can be obtained exactly like that of (4.25). As $\{\check{p}_n \geq \check{p}_n\}_{n=1}^N$ and $\{\hat{p}_k \geq \tilde{p}_k\}_{k=1}^K$, the latter solution also satisfies (4.31). Thus, $\mathcal{P}4.2$ can be solved using **Algorithm 4.I**.

## 4.6 Symbol-wise MSE downlink-interference duality

In this section, we establish the symbol-wise MSE duality between downlink and interference channels. If all symbols are active, this duality can be applied to solve MSE based problems. However, as will be clear later, this duality requires more computation compared to the duality of Sections 4.4 and 4.5. Thus, we propose this duality to be employed for problems like in $\mathcal{P}4.3$ since this problem maintains all symbols active and can not be solved by the duality in Sections 4.4 and 4.5.

### 4.6.1 Symbol-wise MSE transfer (From downlink to interference channel)

To apply this duality for $\mathcal{P}4.3$, we set the interference channel precoder, decoder, noise covariance, input covariance and MSE weight matrices as

$$\mathbf{v}_{ks} = \bar{\beta}_{ks}\mathbf{w}_{ks},\ \mathbf{t}_{ks} = \mathbf{b}_{ks}/\bar{\beta}_{ks}, \boldsymbol{\zeta} = \mathbf{I},\ \boldsymbol{\Delta}_{ks} = \boldsymbol{\Psi} + \mu_{ks}\mathbf{I}_N, \forall k, s. \tag{4.39}$$

Substituting (4.39) into (4.7) and $\{\xi^{DL}_{ks} = \xi^{I}_{ks}, \forall s\}_{k=1}^{K}$ yields

$$\mathbf{w}^H_{ks}(\mathbf{H}^H_k \sum_{i=1}^{K}\sum_{j=1}^{S_i} \mathbf{b}_{i,j}\mathbf{b}^H_{i,j}\mathbf{H}_k + \mathbf{R}_{nk})\mathbf{w}_{ks} - \mathbf{w}^H_{ks}\mathbf{H}^H_k\mathbf{b}_{ks}$$
$$- \mathbf{b}^H_{ks}\mathbf{H}_k\mathbf{w}_{ks} + 1 = \frac{1}{\bar{\beta}^2_{ks}}\mathbf{b}^H_{ks}(\sum_{i=1}^{K}\sum_{j=1}^{S_i} \bar{\beta}^2_{ij}\mathbf{H}_i\mathbf{w}_{ij}\mathbf{w}^H_{ij}\mathbf{H}^H_i +$$
$$\boldsymbol{\Psi} + \mu_{ks}\mathbf{I}_N)\mathbf{b}_{ks} - \mathbf{b}^H_{ks}\mathbf{H}_k\mathbf{w}_{ks} - \mathbf{w}^H_{ks}\mathbf{H}^H_k\mathbf{b}_{ks} + 1,\ \forall k, s.$$

It implies

$$\mathbf{w}^H_{ks}(\mathbf{H}^H_k \sum_{i=1}^{K}\sum_{j=1,(i,j)\neq(k,s)}^{S_i} \mathbf{b}_{i,j}\mathbf{b}^H_{i,j}\mathbf{H}_k + \mathbf{R}_{nk})\mathbf{w}_{ks} =$$
$$\frac{1}{\bar{\beta}^2_{ks}}\mathbf{b}^H_{ks}(\sum_{i=1}^{K}\sum_{j=1,(i,j)\neq(k,s)}^{S_i} \bar{\beta}^2_{ij}\mathbf{H}_i\mathbf{w}_{ij}\mathbf{w}^H_{ij}\mathbf{H}^H_i + \boldsymbol{\Psi} + \mu_{ks}\mathbf{I}_N)\mathbf{b}_{ks},\ \forall k, s. \tag{4.40}$$

Collecting the above expression for all $k$ and $s$ gives

$$(\bar{\mathbf{Y}} + \boldsymbol{\Theta})\bar{\boldsymbol{\beta}}^2 = [a_{11}, \cdots, a_{1S_1}, \cdots, a_{K1}, \cdots, a_{KS_K}]^T = \tilde{\bar{\mathbf{P}}}\mathbf{x}$$
$$\Rightarrow \bar{\boldsymbol{\beta}}^2 = \boldsymbol{\Theta}^{-1}(\mathbf{I} + \bar{\mathbf{Y}}\boldsymbol{\Theta}^{-1})^{-1}\tilde{\bar{\mathbf{P}}}\mathbf{x} \tag{4.41}$$

where $\bar{\boldsymbol{\beta}}^2 = [\bar{\beta}^2_{11}, \cdots, \bar{\beta}^2_{1S_1}, \cdots, \bar{\beta}^2_{K1}, \cdots, \bar{\beta}^2_{KS_K}]^T$, $\boldsymbol{\Theta} = \mathrm{diag}(\theta_{11}, \cdots, \theta_{1K_1}, \cdots, \theta_{K1}, \cdots, \theta_{KS_K})$, $a_{ks} = \mathbf{b}^H_{ks}\boldsymbol{\Psi}\mathbf{b}_{ks} + \mu_{ks}\mathbf{b}^H_{ks}\mathbf{b}_{ks}$, $\tilde{\bar{\mathbf{P}}} = [\bar{\bar{\mathbf{P}}}, \bar{\mathbf{P}}]$ and $\bar{\mathbf{Y}} = [\bar{\mathbf{y}}_{11}, \cdots, \bar{\mathbf{y}}_{1S_1}, \cdots, \bar{\mathbf{y}}_{K1} \cdots, \bar{\mathbf{y}}_{KS_K}]^T$ with $\theta_{ks} = \mathbf{w}^H_{ks}\mathbf{R}_{nk}\mathbf{w}_{ks}$, $\bar{\bar{\mathbf{P}}} \in \Re^{S\times N} = |\mathbf{B}^H|^2$, $\bar{\mathbf{P}} = \mathrm{diag}(\bar{p}_{11}, \cdots, \bar{p}_{1S_1}, \cdots, \bar{p}_{K1}, \cdots, \bar{p}_{KS_K})$, $\bar{\mathbf{y}}_{ks} = [-|\mathbf{b}^H_{ks}\mathbf{H}_1\mathbf{w}_{11}|^2, \cdots, \bar{z}_{ks}, \cdots, -|\mathbf{b}^H_{ks}\mathbf{H}_K\mathbf{w}_{K1}|, \cdots, -|\mathbf{b}^H_{ks}\mathbf{H}_K\mathbf{w}_{KS_K}|]^T$ and $\bar{z}_{ks} = \mathbf{w}^H_{ks}\mathbf{H}^H_k\sum_{i=1}^{K}\sum_{j=1,(i,j)\neq(k,s)}^{S_i}\mathbf{b}_{i,j}\mathbf{b}^H_{i,j}\mathbf{H}_k\mathbf{w}_{ks}$. Next we examine two important properties of $(\mathbf{I} + \bar{\mathbf{Y}}\boldsymbol{\Theta}^{-1})^{-1}$. To this end, we examine the following Theorem.

*Theorem 4.2:* Let $\mathbf{A} \in \Re^{n\times n}$ and $\mathbf{A}_{(i,j),(i\neq j)} \leq 0, 1 \leq i(j) \leq n$. If the diagonal elements of $\mathbf{A}$ are $\mathbf{A}_{(i,i)} = 1 - \sum_{j=1,j\neq i}^{n} \mathbf{A}_{(j,i)}$, then

$$\text{Property 1} : \mathbf{A}^{-1} \geq 0 \tag{4.42}$$

$$\text{Property 2} : |||\mathbf{A}^{-1}|||_1 = 1 \tag{4.43}$$

where $(.) \geq 0$ and $|||.|||_1$ denote matrix non-negativity and one norm, respectively.

**Proof.** See Appendix 4.A. □

According to the first property of *Theorem 4.2:*, if $\{\theta_{ks} > 0, \forall s\}_{k=1}^{K}$[5], the inverse of $(\mathbf{I} + \bar{\mathbf{Y}}\mathbf{\Theta}^{-1})$ exists and it has nonnegative entries. Consequently, for any positive $\{\psi_n\}_{n=1}^{N}$ and $\{\mu_{ks}, \forall s\}_{k=1}^{K}$, $\{\bar{\beta}_{ks}, \forall s\}_{k=1}^{K}$ of (4.41) are strictly positive[6]. Now, by selecting $\{\psi_n\}_{n=1}^{N}$ and $\{\mu_{ks}, \forall s\}_{k=1}^{K}$ such that (4.41) is fulfilled, we can transfer the MSE of each symbol from downlink to interference channel ensuring $\{\xi_{ks}^{DL} = \xi_{ks}^{I_1}, \forall s\}_{k=1}^{K}$, where $\xi_{ks}^{I_1}$ is the MSE of the $k$th user $s$th symbol at step 1 of **Algorithm 4.I**. Here we should also select $\{\psi_n\}_{n=1}^{N}$ and $\{\mu_{ks}, \forall s\}_{k=1}^{K}$ such that the power constraint of $\mathcal{P}4.3$ at step 3 of **Algorithm 4.I** is satisfied. To this end, we examine the steps (2) and (3) of this algorithm.

Like in Section 4.4, we perform step 2 of **Algorithm 4.I** by updating $\mathbf{t}_{ks}$ using MMSE receiver as

$$\begin{aligned}\mathbf{t}_{ks} =& (\mathbf{\Gamma}_c + \mathbf{\Delta}_{ks})^{-1}\mathbf{H}_k\mathbf{v}_{ks}\zeta_{ks} \\ =& (\sum_{i=1}^{K}\sum_{j=1}^{S_i}\bar{\beta}_{ij}\mathbf{H}_i\mathbf{w}_{ij}\mathbf{w}_{ij}^H\mathbf{H}_i^H + \mathbf{\Psi} + \mu_{ks}\mathbf{I})^{-1}\mathbf{H}_k\mathbf{w}_{ks}\bar{\beta}_{ks}\end{aligned} \tag{4.44}$$

where the second equality is obtained from (4.39). The above expression shows that by choosing $\{\mu_{ks} > 0, \forall s\}_{k=1}^{K}$, $\{\psi_n > 0\}_{n=1}^{N}$, we ensure $(\sum_{i=1}^{K}\sum_{j=1}^{S_i}\bar{\beta}_{ij}\mathbf{H}_i\mathbf{w}_{ij}\mathbf{w}_{ij}^H\mathbf{H}_i^H + \mathbf{\Psi} + \mu_{ks}\mathbf{I})^{-1}$ exists. Next, we transfer the symbol-wise MSE from interference to downlink channel by satisfying the power constraint of $\mathcal{P}4.3$ (i.e., we perform step 3).

[5]For $\mathcal{P}3$, $\{\mathbf{w}_{ks}^H\mathbf{R}_{nk}\mathbf{w}_{ks} > 0, \forall s\}_{k=1}^{K}$ is always true.

[6]Note that the application of (4.43) will be clear in the sequel (see (4.55)).

### 4.6.2 Symbol-wise MSE transfer (From interference to downlink channel)

For a given symbol MSE in the interference channel with $\boldsymbol{\zeta} = \mathbf{I}$, we can achieve the same symbol MSE in the downlink channel by using a nonzero scaling factor ($\beta_{ks}$) which satisfies

$$\widetilde{\mathbf{b}}_{ks} = \beta_{ks}\mathbf{t}_{ks}, \;\; \widetilde{\mathbf{w}}_{ks} = \mathbf{v}_{ks}/\beta_{ks}. \tag{4.45}$$

Here we use the notations $\widetilde{\mathbf{B}}$ and $\widetilde{\mathbf{W}}$ to differentiate with the precoder and decoder matrices used in Section 4.6.1. By substituting (4.45) into $\xi_{ks}^{DL}$ (with $\widetilde{\mathbf{B}}$=$\mathbf{B}$, $\widetilde{\mathbf{W}}$=$\mathbf{W}$), then equating the resulting symbol MSE with that of the interference channel (4.7) and after some straightforward steps, we get

$$\begin{aligned}&\frac{1}{\beta_{ks}^2}\mathbf{v}_{ks}^H(\mathbf{H}_k^H\sum_{i=1}^{K}\sum_{j=1,(i,j)\neq(k,s)}^{S_i}\beta_{ij}^2\mathbf{t}_{ij}\mathbf{t}_{ij}^H\mathbf{H}_k+\mathbf{R}_{nk})\mathbf{v}_{ks}=\\&\mathbf{t}_{ks}^H(\sum_{i=1}^{K}\sum_{j=1,(i,j)\neq(k,s)}^{S_i}\mathbf{H}_i\mathbf{v}_{ij}\mathbf{v}_{ij}^H\mathbf{H}_i^H+\boldsymbol{\Psi}+\mu_{ks}\mathbf{I})\mathbf{t}_{ks}, \forall k,s.\end{aligned}$$

By collecting the above equalities for all $k$ and $s$, $\{\beta_{ks},\forall s\}_{k=1}^K$ can be determined by

$$\begin{aligned}(\check{\mathbf{Y}}+\check{\boldsymbol{\Omega}})\boldsymbol{\beta}^2 =&[\mathbf{v}_{11}^H\mathbf{R}_{n1}\mathbf{v}_{11},\cdots,\mathbf{v}_{1S_1}^H\mathbf{R}_{n1}\mathbf{v}_{1S_1},\cdots,\mathbf{v}_{K1}^H\mathbf{R}_{nK}\mathbf{v}_{K1},\cdots,\mathbf{v}_{KS_K}^H\mathbf{R}_{nK}\mathbf{v}_{KS_K}]^T\\=&\boldsymbol{\Theta}\bar{\boldsymbol{\beta}}^2=\boldsymbol{\Theta}\boldsymbol{\Theta}^{-1}(\mathbf{I}+\bar{\mathbf{Y}}\boldsymbol{\Theta}^{-1})^{-1}\tilde{\bar{\mathbf{P}}}\mathbf{x}\\\Rightarrow\boldsymbol{\beta}^2=&(\check{\mathbf{Y}}+\check{\boldsymbol{\Omega}})^{-1}(\mathbf{I}+\bar{\mathbf{Y}}\boldsymbol{\Theta}^{-1})^{-1}\tilde{\bar{\mathbf{P}}}\mathbf{x}\\=&\check{\boldsymbol{\Omega}}^{-1}(\mathbf{I}+\check{\mathbf{Y}}\check{\boldsymbol{\Omega}}^{-1})^{-1}(\mathbf{I}+\bar{\mathbf{Y}}\boldsymbol{\Theta}^{-1})^{-1}\tilde{\bar{\mathbf{P}}}\mathbf{x}\end{aligned} \tag{4.46}$$

where the third equality follows from (4.41), $\boldsymbol{\beta}^2 = [\beta_{11}^2,\cdots,\beta_{1S_1}^2,\cdots,\beta_{K1}^2,\cdots,\beta_{KS_K}^2]^T$, $\boldsymbol{\Omega} = \mathrm{diag}(\mathbf{t}_{11}^H\boldsymbol{\Psi}\mathbf{t}_{11},\cdots,\mathbf{t}_{1S_1}^H\boldsymbol{\Psi}\mathbf{t}_{1S_1},\cdots,\mathbf{t}_{K1}^H\boldsymbol{\Psi}\mathbf{t}_{K1},\cdots,\mathbf{t}_{KS_K}^H\boldsymbol{\Psi}\mathbf{t}_{KS_K})$, $\bar{\boldsymbol{\Omega}} = \mathrm{diag}(\mu_{11}\mathbf{t}_{11}^H\mathbf{t}_{11},\cdots,\mu_{1S_1}\mathbf{t}_{1S_1}^H\mathbf{t}_{1S_1},\cdots,\mu_{K1}\mathbf{t}_{K1}^H\mathbf{t}_{K1},\cdots,\mu_{KS_K}\mathbf{t}_{KS_K}^H\mathbf{t}_{KS_K})$, $\check{\boldsymbol{\Omega}} = \boldsymbol{\Omega}+\bar{\boldsymbol{\Omega}}$ and $\check{\mathbf{Y}} = [\check{\mathbf{y}}_{11},\cdots,\check{\mathbf{y}}_{1S_1},\cdots,\check{\mathbf{y}}_{K1}\cdots,\check{\mathbf{y}}_{KS_K}]^T$ with $\check{\mathbf{y}}_{ks} = [-|\mathbf{t}_{11}^H\mathbf{H}_1\mathbf{v}_{ks}|^2,\cdots,\check{z}_{ks},\cdots,-|\mathbf{t}_{K1}^H\mathbf{H}_K\mathbf{v}_{ks}|^2,\cdots,-|\mathbf{t}_{KS_K}^H\mathbf{H}_K\mathbf{v}_{ks}|^2]^T$ and $\check{z}_{ks} = \mathbf{t}_{ks}^H\sum_{i=1}^K\sum_{j=1,(i,j)\neq(k,s)}^{S_i}\mathbf{H}_i\mathbf{v}_{ij}\mathbf{v}_{ij}^H\mathbf{H}_i^H\mathbf{t}_{ks}$. By applying *Theorem 4.2:*, it can be shown that $\{\beta_{ks},\forall s\}_{k=1}^K$ are strictly positive for $\{\psi_n > 0\}_{n=1}^N$ and $\{\mu_{ks} > 0,\forall s\}_{k=1}^K$. The power constraints of the $n$th BS antenna and $k$th user

$s$th symbol are given by

$$\check{\tilde{\mathbf{b}}}_n^H \check{\tilde{\mathbf{b}}}_n = \check{\mathbf{t}}_n^H \mathbf{Y} \check{\mathbf{t}}_n \leq \check{p}_n, \forall n \quad (4.47)$$

$$\tilde{\mathbf{b}}_{ks}^H \tilde{\mathbf{b}}_{ks} = \beta_{ks}^2 \mathbf{t}_{ks}^H \mathbf{t}_{ks} \leq \check{p}_{ks}, \ \forall k, s \quad (4.48)$$

where $\mathbf{Y} = \text{diag}(\beta_{11}^2, \cdots, \beta_{1S_1}, \cdots, \beta_{K1}^2, \cdots, \beta_{KS_K})$. Multiplying both sides of (4.47) by $\psi_n$ and stacking the resulting inequality for all $n$ yields

$$\check{\mathbf{P}}\boldsymbol{\psi} \geq \tilde{\boldsymbol{\Omega}}\boldsymbol{\beta}^2 \quad (4.49)$$

where $\check{\mathbf{P}} = \text{diag}(\check{p}_1, \cdots, \check{p}_N)$ and $\tilde{\boldsymbol{\Omega}} = \boldsymbol{\Psi}|\mathbf{T}|^2$. Like in the above expression, by multiplying both sides of (4.48) with $\mu_{ks}$ and collecting the resulting inequality for all $k$ and $s$, the power constraints (4.48) can be expressed as

$$\bar{\check{\mathbf{P}}}\boldsymbol{\mu} \geq \bar{\boldsymbol{\Omega}}\boldsymbol{\beta}^2 \quad (4.50)$$

where $\bar{\check{\mathbf{P}}} = \text{diag}(\check{p}_{11}, \cdots, \check{p}_{1S_1}, \cdots, \check{p}_{K1}, \cdots, \check{p}_{KS_K})$. By employing $\boldsymbol{\beta}^2$ of (4.46), (4.49) and (4.50) can be combined as

$$\begin{aligned} \mathbf{x}' \geq \tilde{\bar{\boldsymbol{\Omega}}}\boldsymbol{\beta}^2 =& \tilde{\bar{\boldsymbol{\Omega}}}\check{\boldsymbol{\Omega}}^{-1}(\mathbf{I} + \check{\mathbf{Y}}\check{\boldsymbol{\Omega}}^{-1})^{-1}(\mathbf{I} + \bar{\mathbf{Y}}\boldsymbol{\Theta}^{-1})^{-1}\tilde{\bar{\mathbf{P}}}(\tilde{\bar{\mathbf{P}}})^{-1}\mathbf{x}' \\ =& \daleth(\mathbf{x}')\mathbf{x}' \end{aligned} \quad (4.51)$$

where $\tilde{\bar{\mathbf{P}}} = \text{blkdiag}(\check{\mathbf{P}}, \bar{\check{\mathbf{P}}})$, $\tilde{\bar{\boldsymbol{\Omega}}} = [\tilde{\boldsymbol{\Omega}}^T, \bar{\boldsymbol{\Omega}}^T]^T$, $\mathbf{x}' = \tilde{\bar{\mathbf{P}}}[\boldsymbol{\psi} \ \boldsymbol{\mu}]^T$ and $\daleth(\mathbf{x}') = \tilde{\bar{\boldsymbol{\Omega}}}\check{\boldsymbol{\Omega}}^{-1}(\mathbf{I} + \check{\mathbf{Y}}\check{\boldsymbol{\Omega}}^{-1})^{-1}(\mathbf{I} + \bar{\mathbf{Y}}\boldsymbol{\Theta}^{-1})^{-1}\tilde{\bar{\mathbf{P}}}(\tilde{\bar{\mathbf{P}}})^{-1}$. Next we show that there exists $\mathbf{x}' > 0$ such that (4.51) is satisfied. Towards this end, we consider the following discrete-time switched system [Sun08].

$$\bar{\mathbf{x}}_{n+1} = \mathbf{F}_{\sigma_n}\bar{\mathbf{x}}_n \text{ for } n = 0, 1, 2, \cdots \quad (4.52)$$

where $\bar{\mathbf{x}} \in \Re^{m\times 1}$ is a state, $\mathbf{F}_{\sigma_n} \in \Re^{m\times m}$ is a switching matrix and $\sigma_n \in \{0, 1, 2, \cdots\}$. According to [Sun08] (Remark 2 of [Sun08]), the above system is marginally stable (convergent) if

$$\max_{\sigma_n} \|\mathbf{F}_{\sigma_n}\|_\star = 1 \text{ for } n = 0, 1, \cdots \quad (4.53)$$

where $\|.\|_\star$ denotes an induced matrix norm.

Let us consider the following iteration

$$\mathbf{x}'_{n+1} = \daleth(\mathbf{x}'_n)\mathbf{x}'_n, \ \text{ for } n = 0, 1, 2, \cdots. \quad (4.54)$$

Now if we assume $\gimel(\mathbf{x}'_n) = \mathbf{F}_{\sigma_n}, \forall n$[7], we can interpret (4.54) as a discrete time switched system. Consequently, the above iteration is guaranteed to converge if $\max_n \|\gimel(\mathbf{x}'_n)\|_\star = 1$. It is known that $|||.|||_1$ is an induced matrix norm [HJ85]. For any $\mathbf{x}'$, the matrix one norm of $\gimel(\mathbf{x}')$ is given by

$$\begin{aligned} |||\gimel(\mathbf{x}')|||_1 =& |||\tilde{\bar{\mathbf{\Omega}}}\check{\mathbf{\Omega}}^{-1}(\mathbf{I}+\check{\mathbf{Y}}\check{\mathbf{\Omega}}^{-1})^{-1}(\mathbf{I}+\bar{\mathbf{Y}}\mathbf{\Theta}^{-1})^{-1}\tilde{\bar{\mathbf{P}}}(\tilde{\check{\mathbf{P}}})^{-1}|||_1 \\ \leq& |||\check{\check{\mathbf{\Omega}}}|||_1|||(\mathbf{I}+\check{\mathbf{Y}}\check{\mathbf{\Omega}}^{-1})^{-1}|||_1|||(\mathbf{I}+\bar{\mathbf{Y}}\mathbf{\Theta}^{-1})^{-1}|||_1|||\check{\check{\mathbf{P}}}|||_1 \\ =& |||\check{\check{\mathbf{\Omega}}}|||_1|||\check{\check{\mathbf{P}}}|||_1 \leq 1 \end{aligned} \tag{4.55}$$

where $\check{\check{\mathbf{\Omega}}} = [\tilde{\bar{\mathbf{\Omega}}}\check{\mathbf{\Omega}}^{-1} \; \mathbf{0}_{(N+S)\times N}]$, $\check{\check{\mathbf{P}}} = [\tilde{\bar{\mathbf{P}}}(\tilde{\check{\mathbf{P}}})^{-1}; \; \mathbf{0}_{N\times(N+S)}]$, the second inequality is due to the fact that $|||\mathbf{XY}|||_1 \leq |||\mathbf{X}|||_1|||\mathbf{Y}|||_1$ [HJ85] (page 290), the third equality is obtained by applying *Theorem 4.2:* and the last inequality employs the following facts. Using the definition (4.79) (see Appendix 4.A), one can get $|||\check{\check{\mathbf{\Omega}}}|||_1 \leq 1$ by applying (4.46) and (4.51), and $|||\check{\check{\mathbf{P}}}|||_1 \leq 1$ by applying (4.13), (4.41) and (4.51).

Thus, $\max_n \|\gimel(\mathbf{x}'_n)\|_\star = 1$ holds true and (4.54) is guaranteed to converge. As we can see (4.54) is derived by using (4.41) and (4.46). Thus, the solution of (4.54) also satisfies (4.41) and (4.46). Moreover, for any initial $\mathbf{x}'_0 > 0$, since $\gimel(\mathbf{x}'_n), \forall n$ is positive, the solution of (4.54) is strictly positive and $[\boldsymbol{\psi} \; \boldsymbol{\mu}]^T = (\tilde{\check{\mathbf{P}}})^{-1}\mathbf{x}' > 0$ which is the desired result.

Once the feasible $\{\mu_{ks}, \forall s\}_{k=1}^K$ and $\{\lambda_n\}_{n=1}^N$ are obtained, step 4 of **Algorithm 4.I** is immediate. As a result, $\mathcal{P}4.3$ can be solved using **Algorithm 4.I** with an additional power allocation step which will be detailed in Section 4.8.

### 4.6.3 Extension of the current duality for $\mathcal{P}4.3$ with a total BS power constraint

In this subsection, we show the extension of the current duality for $\mathcal{P}4.3$ with a total BS power constraint. For this problem, we set $\Delta_{ks}$ of (4.39) as $\Delta_{ks} = \mathbf{I}, \forall k, s$ (i.e., like in Section 4.4.2). Upon doing so, $\boldsymbol{\beta}^2$ is computed directly from the first equality of (4.46) (i.e., the bound (4.55) is not needed). By summing the left and right hand sides of this equality, one can get $\text{tr}\{\widetilde{\mathbf{B}}\widetilde{\mathbf{B}}^H\} = P_{max}$.

[7]Since $\gimel(\mathbf{x}'_n)$ is the products of stochastic matrices (see the proof of *Theorem 4.2:*), $\gimel(\mathbf{x}'_n)$ is a bounded matrix for any $\mathbf{x}' > 0$. Thus, the assumption $\gimel(\mathbf{x}'_n) = \mathbf{F}_{\sigma_n}, \forall n$ holds true.

This shows that for $\mathcal{P}4.3$ with a total BS power constraint problem, the total BS power at step 3 of **Algorithm 4.I** is satisfied. Thus, one can apply **Algorithm 4.I** (with the additional power allocation step) to solve the latter problem by setting $\Delta_{ks}$ of (4.39) as $\Delta_{ks} = \mathbf{I}, \forall k, s$.

For other total BS power constrained WMSE-based problems, the current duality based algorithm can be applied like in this subsection. The details are omitted for conciseness. Note that for such problem types, the duality algorithm of the current paper has the same complexity as that of [HJU09].

## 4.7 User-wise MSE downlink-interference duality

This section establishes user-wise MSE duality between downlink and interference channels. This duality is established to solve the problems of type $\mathcal{P}4.4$.

### 4.7.1 User-wise MSE transfer (From downlink to interference channel)

To apply this MSE transfer for $\mathcal{P}4.4$, we set the interference channel precoder, decoder, noise covariance, input covariance and MSE weight matrices as

$$\mathbf{V}_k = \tilde{\beta}_k \mathbf{W}_k,\ \mathbf{T}_k = \mathbf{B}_k/\tilde{\beta}_k, \boldsymbol{\zeta} = \mathbf{I},\ \boldsymbol{\Delta}_{ks} = \boldsymbol{\Psi} + \mu_k \mathbf{I}. \tag{4.56}$$

Like in Section 4.6, substituting (4.56) into (4.8), equating $\{\xi_k^{DL} = \xi_k^{I}\}_{k=1}^K$ and after some straightforward steps, we get the following system of equations

$$(\tilde{\mathbf{Y}} + \tilde{\boldsymbol{\Theta}})\tilde{\boldsymbol{\beta}}^2 = \tilde{\bar{\mathbf{P}}}\mathbf{x},\ \Rightarrow \tilde{\boldsymbol{\beta}}^2 = \tilde{\boldsymbol{\Theta}}^{-1}(\mathbf{I} + \tilde{\mathbf{Y}}\tilde{\boldsymbol{\Theta}}^{-1})^{-1}\tilde{\bar{\mathbf{P}}}\mathbf{x} \tag{4.57}$$

where

$$\tilde{\mathbf{Y}}_{(k,l)} = \begin{cases} \sum_{i=1, i\neq k}^K \|\mathbf{W}_k^H \mathbf{H}_k^H \mathbf{B}_i\|_F^2 & \text{for } k = l \\ -\|\mathbf{W}_l^H \mathbf{H}_l^H \mathbf{B}_k\|_F^2, & \text{for } k \neq l \end{cases} \tag{4.58}$$

$\tilde{\boldsymbol{\Theta}} = \mathrm{diag}(\theta_1, \cdots, \theta_K)$, $\tilde{\boldsymbol{\beta}}^2 = [\tilde{\beta}_1^2, \cdots, \tilde{\beta}_K^2]^T$, $\tilde{\bar{\mathbf{P}}} = [\bar{\bar{\mathbf{P}}}, \tilde{\mathbf{P}}]$ with $\tilde{\theta}_k = \mathrm{tr}\{\mathbf{W}_k^H \mathbf{R}_{nk} \mathbf{W}_k\}$, $\bar{\bar{\mathbf{P}}} \in \Re^{S\times N} = |\mathbf{B}^H|^2$, $\tilde{\mathbf{P}} = \mathrm{diag}(\bar{p}_1, \cdots, \bar{p}_K)$. By applying

*Theorem 2*, it can be shown that $\{\tilde{\beta}_k\}_{k=1}^K$ of (4.57) are strictly positive. Thus, step 1 of **Algorithm 4.I** can be performed using (4.57). We perform step 2 of **Algorithm 4.I** by updating $\mathbf{t}_{ks}$ using MMSE receiver as

$$\mathbf{t}_{ks} =(\sum_{i=1}^{K} \tilde{\beta}_i \mathbf{H}_i \mathbf{W}_i \mathbf{W}_i^H \mathbf{H}_i^H + \boldsymbol{\Psi} + \mu_k \mathbf{I})^{-1} \mathbf{H}_k \mathbf{w}_{ks} \tilde{\beta}_k. \tag{4.59}$$

### 4.7.2 User-wise MSE transfer (From interference to downlink channel)

For a given user MSE in the interference channel with $\boldsymbol{\zeta} = \mathbf{I}$, we can achieve the same MSE in the downlink channel by using nonzero scaling factors ($\{\tilde{\tilde{\beta}}_k\}_{k=1}^K$) that satisfy

$$\widetilde{\mathbf{B}}_k = \tilde{\tilde{\beta}}_k \mathbf{T}_k, \;\; \widetilde{\mathbf{W}}_k = \mathbf{V}_k / \tilde{\tilde{\beta}}_k. \tag{4.60}$$

Here we also use the notations $\widetilde{\mathbf{B}}$ and $\widetilde{\mathbf{W}}$ to differentiate with the precoder and decoder matrices used in Section 4.7.1. By substituting (4.60) into $\xi_k^{DL}$ (with $\widetilde{\mathbf{B}}$=$\mathbf{B}$, $\widetilde{\mathbf{W}}$=$\mathbf{W}$) and then equating the resulting user-wise MSE with that of the interference channel (4.7) and after some steps, $\{\tilde{\tilde{\beta}}_k\}_{k=1}^K$ are determined as

$$\begin{aligned}(\mathbf{I} + \check{\tilde{\mathbf{Y}}} \check{\tilde{\boldsymbol{\Omega}}}^{-1}) \check{\tilde{\boldsymbol{\Omega}}} \tilde{\tilde{\boldsymbol{\beta}}}^2 &= [\mathrm{tr}\{\mathbf{V}_1^H \mathbf{R}_{n1} \mathbf{V}_1\}, \cdots, \mathrm{tr}\{\mathbf{V}_K^H \mathbf{R}_{nK} \mathbf{V}_K\}]^T \\ \Rightarrow \tilde{\tilde{\boldsymbol{\beta}}}^2 &= \check{\tilde{\boldsymbol{\Omega}}}^{-1} (\mathbf{I} + \check{\tilde{\mathbf{Y}}} \check{\tilde{\boldsymbol{\Omega}}}^{-1})^{-1} (\mathbf{I} + \tilde{\mathbf{Y}} \tilde{\boldsymbol{\Theta}}^{-1})^{-1} \tilde{\tilde{\mathbf{P}}} \tilde{\mathbf{x}} \end{aligned} \tag{4.61}$$

where the second equality follows from (4.57), $\tilde{\tilde{\boldsymbol{\beta}}}^2 = [\tilde{\tilde{\beta}}_1^2, \cdots, \tilde{\tilde{\beta}}_K^2]^T$, $\boldsymbol{\Omega}' = \mathrm{diag}(\mathrm{tr}\{\mathbf{T}_1^H \boldsymbol{\Psi} \mathbf{T}_1\}, \cdots, \mathrm{tr}\{\mathbf{T}_K^H \boldsymbol{\Psi} \mathbf{T}_K\})$, $\hat{\boldsymbol{\Omega}} = \mathrm{diag}(\mu_1 \mathrm{tr}\{\mathbf{T}_1^H \mathbf{T}_1\}, \cdots, \mu_K \mathrm{tr}\{\mathbf{T}_K^H \mathbf{T}_K\})$, $\check{\tilde{\boldsymbol{\Omega}}} = \boldsymbol{\Omega}' + \hat{\boldsymbol{\Omega}}$. By applying *Theorem 4.2:*, it can be shown that $\{\tilde{\tilde{\beta}}_k\}_{k=1}^K$ are strictly positive for $\{\psi_n > 0\}_{n=1}^N$ and $\{\mu_k > 0\}_{k=1}^K$. Like in Section 4.6.2, the power constraint of the $n$th BS antenna and $k$th user can thus be expressed as

$$\begin{aligned}\tilde{\mathbf{x}}' \geq \hat{\hat{\boldsymbol{\Omega}}} \tilde{\tilde{\boldsymbol{\beta}}}^2 =& \hat{\hat{\boldsymbol{\Omega}}} \check{\tilde{\boldsymbol{\Omega}}}^{-1} (\mathbf{I} + \check{\tilde{\mathbf{Y}}} \check{\tilde{\boldsymbol{\Omega}}}^{-1})^{-1} (\mathbf{I} + \tilde{\mathbf{Y}} \tilde{\boldsymbol{\Theta}}^{-1})^{-1} \tilde{\tilde{\mathbf{P}}} \tilde{\mathbf{x}} \\ =& \tilde{\daleth}(\tilde{\mathbf{x}}') \tilde{\mathbf{x}}' \end{aligned} \tag{4.62}$$

where $\hat{\mathbf{P}} = \mathrm{diag}(\hat{p}_1, \cdots, \hat{p}_K)$, $\mathring{\hat{\mathbf{P}}} = \mathrm{blkdiag}(\check{\mathbf{P}}, \hat{\mathbf{P}})$, $\hat{\hat{\boldsymbol{\Omega}}} = [\tilde{\boldsymbol{\Omega}}^T, \hat{\boldsymbol{\Omega}}^T]^T$, $\tilde{\mathbf{x}}' = \tilde{\tilde{\mathbf{P}}}[\boldsymbol{\psi} \, \tilde{\boldsymbol{\mu}}]^T$ and $\tilde{\daleth}(\mathbf{x}') = \hat{\hat{\boldsymbol{\Omega}}} \check{\tilde{\boldsymbol{\Omega}}}^{-1} (\mathbf{I} + \check{\tilde{\mathbf{Y}}} \check{\tilde{\boldsymbol{\Omega}}}^{-1})^{-1} (\mathbf{I} + \tilde{\mathbf{Y}} \tilde{\boldsymbol{\Theta}}^{-1})^{-1} \tilde{\tilde{\mathbf{P}}} (\mathring{\hat{\mathbf{P}}})^{-1}$. Like in Section 4.6.2, it can be shown that there exists a feasible $\tilde{\mathbf{x}}' > 0$ that satisfy (4.62) and can be obtained iteratively by

$$\tilde{\mathbf{x}}'_{n+1} = \tilde{\mathfrak{I}}(\tilde{\mathbf{x}}'_n)\tilde{\mathbf{x}}'_n, \quad \text{for } n = 0, 1, 2, \cdots. \tag{4.63}$$

By initializing $\mathbf{x}'_0 > 0$, the solution of the above iteration is always positive. Consequently, $\{\lambda_n > 0\}_{n=1}^N$ and $\{\mu_k > 0\}_{k=1}^K$ holds true. Once the feasible $\{\mu_k\}_{k=1}^K$ and $\{\lambda_n\}_{n=1}^N$ are obtained, step 4 of **Algorithm 4.I** is straightforward. As a result, $\mathcal{P}4.4$ can be solved using **Algorithm 4.I** with the additional power allocation step of Section 4.8.

## 4.8 Generalized and improved version of Algorithm 4.I

From the discussions of Sections 4.4 - 4.7, one can understand that each iteration of **Algorithm 4.I** gives a non increasing sequence of symbol (user) WMSE/WSMSE. As can be seen from Section 4.3, the objective of $\mathcal{P}1$ ($\mathcal{P}4.2$) is just to minimize the total WSMSE of all symbols (users), whereas the objective of $\mathcal{P}4.3$ ($\mathcal{P}4$) is to simultaneously minimize and balance the WMSE of all symbols (users). Thus, **Algorithm 4.I** is appropriate to solve $\mathcal{P}4.1$($\mathcal{P}4.2$) of the current paper. For $\mathcal{P}4.3$($\mathcal{P}4.4$), although each iteration of **Algorithm 4.I** is able to provide a non increasing sequence of symbol (user) WMSE (i.e., minimizes the maximum WMSE of all symbols (users)), each iteration of this algorithm is not able to guarantee balanced WMSEs of all symbols (users). On the other hand, for an MSE constrained total BS power minimization problem (for example $\mathcal{P}4.7$ in Section 4.9), an iterative algorithm that can provide a non increasing sequence of total BS power is required. This shows that **Algorithm 4.I** also can not solve the latter problem. In the following we address the drawbacks of **Algorithm 4.I** just by including a power allocation step into **Algorithm 4.I** as explained below.

In [BV11a], for fixed transmit and receive filters, the power allocation parts of total BS power constrained MSE-based problems have been formulated as GPs by employing the approach and system model of [SSB08c] under the as-sumption that all symbols are strictly active[8]. For this assumption, in [BV11a], we show that the system model of [SSB08c] is appropriate to solve any kind of total BS power constrained MSE-based problems using duality approach (alternating optimization). This motivates us to utilize the system model of [SSB08c] in the downlink channel only and then include the power allocation step (i.e., GP) into **Algorithm 4.I**. Towards this end, we decompose the precoders and decoders of the downlink channel as

$$\mathbf{B}_k = \mathbf{G}_k \mathbf{P}_k^{1/2}, \;\; \mathbf{W}_k = \mathbf{U}_k \boldsymbol{\alpha}_k \mathbf{P}_k^{-1/2}, \; \forall k \tag{4.64}$$

where $\mathbf{P}_k = \mathrm{diag}(p_{k1}, \cdots, p_{kS_k}) \in \Re^{S_k \times S_k}$, $\mathbf{G}_k = [\mathbf{g}_{k1} \; \cdots \; \mathbf{g}_{kS_k}] \in \mathbb{C}^{N \times S_k}$, $\mathbf{U}_k = [\mathbf{u}_{k1} \; \cdots \; \mathbf{u}_{kS_k}] \in \mathbb{C}^{M_k \times S_k}$ and $\boldsymbol{\alpha}_k = \mathrm{diag}(\alpha_{k1}, \cdots, \alpha_{kS_k}) \in \Re^{S_k \times S_k}$ are the transmit power, unity norm transmit filter, unity norm receive filter and receiver scaling factor matrices of the $k$th user, respectively, i.e., $\{\mathbf{g}_{ks}^H \mathbf{g}_{ks} = \mathbf{u}_{ks}^H \mathbf{u}_{ks} = 1, \forall s\}_{k=1}^K$.

By employing (4.64) and stacking $\xi = [\xi_{1,1}^{DL}, \cdots, \xi_{K,S_K}^{DL}]^T = [\xi_1^{DL}, \cdots, \xi_S^{DL}]^T = [\{\xi_l^{DL}\}_{l=1}^S]^T$, the $l$th downlink symbol MSE can be expressed as (see [SSB08c] and [BV11b] for more details about (4.64) and the above descriptions)

$$\xi_l^{DL} = p_l^{-1}[(\mathbf{D} + \boldsymbol{\alpha}^2 \boldsymbol{\Phi}^T)\mathbf{p}]_l + p_l^{-1}\alpha_l^2 \mathbf{u}_l^H \mathbf{R}_n \mathbf{u}_l \tag{4.65}$$

where

$$\boldsymbol{\Phi}_{(l,j)} = \begin{cases} |\mathbf{g}_l^H \mathbf{H} \mathbf{u}_j|^2, & \text{for } l \neq j \\ 0, & \text{for } l = j \end{cases} \tag{4.66}$$

$$\mathbf{D}_{(l,l)} = \alpha_l^2 |\mathbf{g}_l^H \mathbf{H} \mathbf{u}_l|^2 - 2\alpha_l \Re(\mathbf{u}_l^H \mathbf{H}^H \mathbf{g}_l) + 1, \tag{4.67}$$

$1 \leq l(j) \leq S$, $\mathbf{P} = \mathrm{blkdiag}(\mathbf{P}_1, \cdots, \mathbf{P}_K) = \mathrm{diag}(p_1, \cdots, p_S)$, $\mathbf{p} = [p_1, \cdots, p_S]^T$, $\mathbf{G} = [\mathbf{G}_1, \cdots, \mathbf{G}_K] = [\mathbf{g}_1, \cdots, \mathbf{g}_S]$, $\mathbf{U} = \mathrm{blkdiag}(\mathbf{U}_1, \cdots, \mathbf{U}_K) = [\mathbf{u}_1, \cdots, \mathbf{u}_S]$ and $\boldsymbol{\alpha} = \mathrm{blkdiag}(\mathbf{ff}_1, \cdots, \mathbf{ff}_K) = \mathrm{diag}(\alpha_1, \cdots, \alpha_S)$ with $\|\mathbf{g}_l\|_2 = \|\mathbf{u}_l\|_2 = 1$. Using (4.65), for fixed $\mathbf{G}, \mathbf{U}$ and $\boldsymbol{\alpha}$, the power allocation part of $\mathcal{P}4.1$ can be formulated as

$$\min_{\{p_l\}_{l=1}^S} \sum_{l=1}^S \eta_l \xi_l^{DL}, \text{ s.t } \boldsymbol{\varsigma}_n^T \mathbf{p} \leq \breve{p}_n, \; p_l \leq \breve{p}_l \; \forall n, l \tag{4.68}$$

[8]Note that this assumption is not always true for all MSE-based problems. However, as mentioned in [SSB08c], in practice replacing zero powers by a small value will not affect the overall optimization. Due to this reason, we replace zero powers by $10^{-6}$ in the simulation section.

where $\boldsymbol{\varsigma}_n^T \in \Re^{1\times S} = \{|[\mathbf{G}_{(n,i)}|^2\}_{i=1}^S$, $[\eta_1, \cdots, \eta_S]^T = [\eta_{11}, \cdots, \eta_{KS_K}]^T$ and $[\breve{p}_l, \cdots, \breve{p}_S]^T = [\breve{p}_{11}, \cdots, \breve{p}_{KS_K}]^T$. As $\xi_l^{DL}$ is a posynomial (where $\{p_l\}_{l=1}^S$ are the variables), (4.68) is a GP for which global optimality is guaranteed. Thus, it can be efficiently solved using interior point methods with a worst-case polynomial-time complexity [BV04].

For fixed $\mathbf{G}, \mathbf{U}$ and $\boldsymbol{\alpha}$, the power allocation parts of $\mathcal{P}4.2$ - $\mathcal{P}4.4$ can be formulated as GPs like in $\mathcal{P}4.1$. Our duality based algorithm for each of these problems including the power allocation step is summarized in **Algorithm 4.II**.

**Algorithm 4.II**

Initialization: Like in **Algorithm 4.I**.

**Repeat**

**Interference channel**

1. For $\mathcal{P}4.1$ and $\mathcal{P}4.2$, set $\mathbf{V} = \mathbf{W}, \mathbf{T} = \mathbf{B}$ (i.e., $\bar{\beta} = \tilde{\beta} = 1$), then compute $\{\psi_n, \mu_{ks}, \forall k, s, n\}$ and $\{\psi_n, \mu_k, \forall k, n\}$ using (4.25) and (4.38), respectively. For $\mathcal{P}4.3$ and $\mathcal{P}4.4$, first compute $\{\psi_n, \mu_{ks}, \forall k, s, n\}$ and $\{\psi_n, \mu_k, \forall k, n\}$ using (4.54) and (4.63), respectively, then transfer each symbol and user MSE from downlink to interference channels by (4.39) and (4.56), respectively.

2. Update the MMSE receivers of the interference channel for $\mathcal{P}1$, $\mathcal{P}4.2$, $\mathcal{P}4.3$ and $\mathcal{P}4.4$ using (4.19), (4.32), (4.44) and (4.59), respectively.

**Downlink channel**

3. Transfer the MSE (weighted sum, user or symbol MSE) from interference to downlink channel using (4.20), (4.33), (4.45) and (4.60) for $\mathcal{P}4.1$, $\mathcal{P}4.2$, $\mathcal{P}4.3$ and $\mathcal{P}4.4$, respectively.

4. For each of the problems $\mathcal{P}4.1$ - $\mathcal{P}4.4$, decompose the precoder and decoder matrices of each user as in (4.64). Then, formulate and solve the GP power allocation part. For example, the power allocation part of $\mathcal{P}4.1$ can be expressed in GP form as (4.68).

5. For each of the problems $\mathcal{P}4.1$ - $\mathcal{P}4.4$, by keeping $\{\mathbf{P}_k\}_{k=1}^K$ constant, update the receive filters $\{\mathbf{U}_k\}_{k=1}^K$ and scaling factors $\{\boldsymbol{\alpha}_k\}_{k=1}^K$ by applying downlink MMSE receiver approach i.e., $\{\mathbf{U}_k\boldsymbol{\alpha}_k = (\mathbf{H}_k^H\mathbf{GPG}^H\mathbf{H}_k + \mathbf{R}_{nk})^{-1}\mathbf{H}_k^H\mathbf{G}_k\mathbf{P}_k\}_{k=1}^K$. Note that in these expressions, $\{\boldsymbol{\alpha}_k\}_{k=1}^K$ are chosen such that each column of $\{\mathbf{U}_k\}_{k=1}^K$ has unity norm. Then, compute $\{\mathbf{B}_k, \mathbf{W}_k\}_{k=1}^K$ by (4.64).

   **Until** convergence.

   **Convergence**: It can be shown that at each iteration of this algorithm, the objective function of each of the problems $\mathcal{P}4.1$ - $\mathcal{P}4.4$ is non-increasing [SSB07], [BCV11], [BV11c]. Thus, the above iterative algorithm is convergent. However, since $\mathcal{P}4.1$ - $\mathcal{P}4.4$ are non-convex, this iterative algorithm is not guaranteed to converge to the global optimum.
   In this algorithm, we stop iteration (i.e., our convergence condition) when the difference between the objective functions in two consecutive iterations is smaller than some small value $\tilde{\epsilon}$ (we use $\tilde{\epsilon} = 10^{-6}$ for the simulation).
   **Computational complexity**: As can be seen from this algorithm, when we increase the number of users and/or (BS and/or MS antennas), the number and size of optimization variables increase. Because of this, the computational complexity of **Algorithm 4.II** increases as $K$ and/or $N$ and/or $M$ increases. However, studying the complexity of this algorithm as a function of $K, N$ and $M$ needs effort and time. And such a task is beyond the scope of this work and is an open research topic.

The power allocation step of **Algorithm 4.II** has thus the following benefits: (1) For BS power constrained WSMSE minimization problems, this step improves the convergence speed of **Algorithm 4.II** compared to that of **Algorithm 4.I**[9] (for example in $\mathcal{P}1 - P2$). The degree of improvement depends on different parameters (for example $\mathbf{H}_k$, $\boldsymbol{\Delta}_{ks}, \forall k, s$ etc). Thus, the theoretical comparison of these two algorithms in terms of convergence speed requires time and effort. And this task is beyond the scope of this work and it is an open research topic. (2) For symbol-wise (user-wise) WMSE balancing problems, this step helps to balance the WMSE of all symbols (users) (for example in $\mathcal{P}3 - P4$). (3) For MSE constrained total BS power minimization problems, this step ensures a non increasing total BS power at each iteration of **Algorithm 4.II**.

[9]This is at the expense of additional computation. However, as mentioned in [SSB08c] (see Appendix A of [SSB08c]), a small desktop computer can solve a GP of 100 variables and 10000 constraints by standard interior point method under a minute. Thus, we believe that the complexity of **Algorithm 4.I** and **Algorithm 4.II** are almost the same.

## 4.9 Application of the proposed duality based algorithm for other problems

### 4.9.1 MSE based problem with entry-wise power constraint

The symbol-wise WSMSE minimization constrained with entry wise power i.e,. $b_{ksn}^H b_{ksn} \leq \bar{\bar{p}}_{ksn}, \forall k, s, n$ problem is formulated as

$$\mathcal{P}4.5: \min_{\{\mathbf{B}_k, \mathbf{W}_k\}_{k=1}^K} \sum_{k=1}^K \sum_{s=1}^{S_k} \eta_{ks} \xi_{ks}^{DL}, \text{ s.t } b_{ksn}^H b_{ksn} \leq \bar{\bar{p}}_{ksn}. \tag{4.69}$$

It can be shown that this problem can be solved by **Algorithm 4.II** with $\{\mathbf{\Delta}_{ks} = \text{diag}(\delta_{ks1}, \cdots, \delta_{ksN}), \forall s\}_{k=1}^K$.

### 4.9.2 Weighted sum rate optimization constrained with per antenna and symbol power problem

By employing the approach of [BV11d] (see (16) of [BV11d]), one can equivalently express the weighted sum rate maximization constrained with per antenna and symbol power problem as

$$\begin{aligned} \mathcal{P}4.6: & \min_{\{\bar{\tau}_{ks}, \bar{v}_{ks}, \mathbf{b}_{ks}, \mathbf{w}_{ks}, \forall s\}_{k=1}^K} \sum_{k=1}^K \sum_{s=1}^{S_k} \bar{\theta}_{ks} \frac{1}{\bar{\tau}_{ks}} \bar{v}_{ks}^{\bar{\gamma}_{ks}} + \sum_{k=1}^K \sum_{s=1}^{S_k} \bar{\eta}_{ks} \xi_{ks}^{DL}, \\ & \text{s.t } [\mathbf{B}\mathbf{B}^H]_{(n,n)} \leq \breve{p}_n, \mathbf{b}_{ks}^H \mathbf{b}_{ks} \leq \breve{p}_{ks}, \prod_{k=1}^K \prod_{s=1}^{S_k} \bar{v}_{ks} = 1, \bar{\tau}_{ks} > 0 \end{aligned} \tag{4.70}$$

where $\{0 < \bar{\omega}_{ks} < 1, \forall s\}_{k=1}^K$ are the rate weighting factors for all symbols, $\bar{\eta}_{ks} = \bar{\tau}_{ks}^{\bar{\mu}_{ks}}$, $\bar{\gamma}_{ks} = \frac{1}{1-\bar{\omega}_{ks}}$ $\bar{\mu}_{ks} = \frac{1}{\bar{\omega}_{ks}} - 1$ and $\bar{\theta}_{ks} = \bar{\omega}_{ks} \bar{\mu}_{ks}^{(1-\bar{\omega}_{ks})}$. For fixed $\{\bar{\tau}_{ks}, \bar{v}_{ks}, \forall s\}_{k=1}^K$, the above optimization problem has the same mathematical structure as that of $\mathcal{P}4.1$. Thus, by keeping $\{\bar{\tau}_{ks}, \bar{v}_{ks}, \forall s\}_{k=1}^K$ constant, $\{\mathbf{b}_{ks}, \mathbf{w}_{ks}, \forall s\}_{k=1}^{K}$ can be optimized by applying the MSE duality discussed in Section 4.4. Moreover, $\{\bar{\tau}_{ks}, \bar{\upsilon}_{ks}, \forall s\}_{k=1}^{K}$ and the power allocation part of the above problem can be optimized by a GP method like in (25) of [BV12]. Consequently, we can apply **Algorithm 4.II** to solve (4.70). The detailed explanations are omitted for conciseness. The following problems can also be solved by simple modification of **Algorithm 4.II**

$$\begin{aligned}
\mathcal{P}4.7: & \min_{\{\mathbf{B}_k, \mathbf{W}_k\}_{k=1}^{K}} \sum_{k=1}^{K} \mathrm{tr}\{\mathbf{B}_k\mathbf{B}_k^H\}, \\
& \text{s.t } [\mathbf{B}\mathbf{B}^H]_{(n,n)} \leq \check{p}_n,\ \mathrm{tr}\{\mathbf{B}_k^H\mathbf{B}_k\} \leq \hat{p}_k,\ SINR_{ks} \geq \varrho_{ks}, \qquad \forall n,k,s \\
\equiv: & \min_{\{\mathbf{B}_k, \mathbf{W}_k\}_{k=1}^{K}} \sum_{k=1}^{K} \mathrm{tr}\{\mathbf{B}_k\mathbf{B}_k^H\}, \\
& \text{s.t } [\mathbf{B}\mathbf{B}^H]_{(n,n)} \leq \check{p}_n,\ \mathrm{tr}\{\mathbf{B}_k^H\mathbf{B}_k\} \leq \hat{p}_k,\ \xi_{ks}^{DL} \leq (1+\varrho_{ks})^{-1}, \qquad \forall n,k,s \\
\mathcal{P}4.8: & \max_{\{\mathbf{B}_k, \mathbf{W}_k\}_{k=1}^{K}} \min R_{ks} \\
& \text{s.t } [\mathbf{B}\mathbf{B}^H]_{(n,n)} \leq \check{p}_n,\ \mathbf{b}_{ks}^H\mathbf{b}_{ks} \leq \check{p}_{ks},\ \forall n,k,s \\
\equiv: & \min_{\{\mathbf{B}_k, \mathbf{W}_k\}_{k=1}^{K}} \max \xi_{ks}^{DL} \\
& \text{s.t } [\mathbf{B}\mathbf{B}^H]_{(n,n)} \leq \check{p}_n,\ \mathbf{b}_{ks}^H\mathbf{b}_{ks} \leq \check{p}_{ks},\ \forall n,k,s
\end{aligned} \tag{4.71}$$

where $SINR_{ks}(R_{ks})$ is the SINR (rate) of the $k$th user $s$th symbol, and we use the fact that $R_{ks} = \log(1+SINR_{ks})$ and $\xi_{ks}^{DL} = (1+SINR_{ks})^{-1}$ [SSB08c].

## 4.10 Extension of the proposed duality based algorithms for robust transceiver design problems

In this section, the extension of the proposed duality-based algorithms for robust transceiver design problems will be discussed. The robustness against imperfect CSI is incorporated into our designs using stochastic approach [BCV11]. To this end, the channel can be modeled as (see Chapter 2.3)

$$\mathbf{H}_k^H = \widehat{\mathbf{H}}_k^H + \mathbf{R}_{mk}^{1/2}\mathbf{E}_{wk}^H\mathbf{R}_{bk}^{1/2} = \widehat{\mathbf{H}}_k^H + \mathbf{E}_k^H,\ \forall k \tag{4.72}$$

where $\mathbf{H}_k^H$ ($\widehat{\mathbf{H}}_k^H$) is the true (estimated) channel, $\mathbf{R}_{bk} \in C^{N\times N}$ ($\widetilde{\mathbf{R}}_{mk} \in C^{M_k\times M_k}$) antenna correlation matrix at the BS ($k$th MS), $\mathbf{R}_{mk} = (\mathbf{I}_{M_k} + \sigma_{ek}^2\widetilde{\mathbf{R}}_{mk}^{-1})^{-1}$, $\mathbf{E}_k^H$ is the estimation error and the entries of $\mathbf{E}_{wk}^H$ are i.i.d with $\mathcal{CN}(0, \sigma_{ek}^2)$.

Like in Chapter 2, we assume that each MS estimates its channel and feeds the estimated channel back to the BS without any error and delay. Thus, both the BS and MSs have the same channel imperfections. With these assumptions, the downlink AMSEs of the $k$th user $s$th symbol and $k$th user are given by

$$\bar{\xi}_{ks}^{DL} = \mathrm{E}_{\mathbf{E}_{wk}^{H}}\{\xi_{ks}^{DL}\} = \mathbf{w}_{ks}^{H}\Gamma_{k}^{DL}\mathbf{w}_{ks} - \mathbf{w}_{ks}^{H}\widehat{\mathbf{H}}_{k}^{H}\mathbf{b}_{ks} - \mathbf{b}_{ks}^{H}\widehat{\mathbf{H}}_{k}\mathbf{w}_{ks} + 1 \quad (4.73)$$

$$\bar{\boldsymbol{\xi}}_{k}^{DL} = \mathrm{E}_{\mathbf{E}_{wk}^{H}}\{\boldsymbol{\xi}_{k}^{DL}\} = \mathbf{I}_{S_k} + \mathbf{W}_{k}^{H}\Gamma_{k}^{DL}\mathbf{W}_{k} - \mathbf{W}_{k}^{H}\widehat{\mathbf{H}}_{k}^{H}\mathbf{B}_{k} - \mathbf{B}_{k}^{H}\widehat{\mathbf{H}}_{k}\mathbf{W}_{k} \quad (4.74)$$

where $\Gamma_{k}^{DL} = \widehat{\mathbf{H}}_{k}^{H}\mathbf{B}\mathbf{B}^{H}\widehat{\mathbf{H}}_{k} + \sigma_{ek}^{2}\mathrm{tr}\{\mathbf{R}_{bk}\mathbf{B}\mathbf{B}^{H}\}\mathbf{R}_{mk} + \mathbf{R}_{nk}$. Using these two equations, the symbol-wise and user-wise WSAMSEs can be expressed as

$$\begin{aligned}\bar{\xi}_{ws}^{DL} &= \sum_{k=1}^{K}\sum_{s=1}^{S_k}\eta_{ks}\bar{\xi}_{ks}^{DL} \\ &= \mathrm{tr}\{\boldsymbol{\eta}\} + \sum_{k=1}^{K}\sum_{s=1}^{S_k}\{\mathbf{w}_{ks}^{H}\Gamma_{k}^{DL}\mathbf{w}_{ks} - \mathbf{w}_{ks}^{H}\widehat{\mathbf{H}}_{k}^{H}\mathbf{b}_{ks} - \mathbf{b}_{ks}^{H}\widehat{\mathbf{H}}_{k}\mathbf{w}_{ks}\} \quad (4.75)\end{aligned}$$

$$\begin{aligned}\bar{\xi}_{wu}^{DL} &= \sum_{k=1}^{K}\tilde{\eta}_{k}\bar{\xi}_{k}^{DL} \\ &= \mathrm{tr}\{\tilde{\boldsymbol{\eta}}\} + \sum_{k=1}^{K}\tilde{\eta}_{k}\mathrm{tr}\{\mathbf{W}_{k}^{H}\Gamma_{k}^{DL}\mathbf{W}_{k} - \mathbf{W}_{k}^{H}\widehat{\mathbf{H}}_{k}^{H}\mathbf{B}_{k} - \mathbf{B}_{k}^{H}\widehat{\mathbf{H}}_{k}\mathbf{W}_{k}\} \quad (4.76)\end{aligned}$$

From equations (4.5) - (4.6) and (4.73) - (4.76), we can realize that the AMSE expressions are slight modification of MSE expressions. Moreover, by applying the approaches of Sections 4.4 - 4.8, one can examine the the robust MSE-based problems by applying the duality approach of this chapter (see for example the sum AMSE problems solved in [BV11b]). This shows that the robust MSE-based problems (for example $\mathcal{P}4.1$ - $\mathcal{P}4.5$) can be solved by slight modification of **Algorithm 4.I**. However, under stochastic robust design technique, the relationship between robust rate SINR and MSE-based problems are not known. Due to this reason, the duality approach of this chapter can not be extended straightforwardly to the robust rate and SINR-based problems (for example, the robust versions of $\mathcal{P}4.6$ - $\mathcal{P}4.8$).

We would like to mention here that the application of **Algorithm 4.II** is not limited to the problems of this chapter. Moreover, the distributive implementation of this algorithm can be analyzed like the algorithms of Chapter 3.

## 4.11 Simulation Results

In this section, we present simulation results for $\mathcal{P}4.1$ - $\mathcal{P}4.4$ (i.e., by assuming perfect CSI). All of our simulation results are averaged over 100 randomly chosen channel realizations. We set $K = 2$, $N = 4$ and $\{M_k = S_k = 2, \eta_{ks} = \rho_{ks} = \tilde{\eta}_k = \tilde{\rho}_k = 1, \forall s\}_{k=1}^{K}$. It is assumed that $\mathbf{R}_{n1} = \sigma_1^2 \mathbf{I}_{M_1}$, $\mathbf{R}_{n2} = \sigma_2^2 \mathbf{I}_{M_2}$ and $\sigma_2^2 = 2\sigma_1^2$. The maximum power of each BS antenna is set to $\{\breve{p}_n = 2.5mW\}_{n=1}^{N}$. And the maximum power allocated to each symbol and user are set to $\{\breve{p}_{ks} = 2.5mW, \forall s\}_{k=1}^{K}$ and $\{\hat{p}_k = 5mW\}_{k=1}^{K}$, respectively. For better exposition, we define the SNR as $P_{\max}/K\sigma_{av}^2$ and it is controlled by varying $\sigma_{av}^2$, where $P_{\max} = 10$mW is the total maximum BS power and $\sigma_{av}^2 = (\sigma_1^2 + \sigma_2^2)/2$. We also compare **Algorithm 4.II** and the algorithm in [SSVB08].

Note that the algorithm in [SSVB08] is designed for coordinated BS systems scenario. And, the iterative algorithm of [SSVB08] is based on the per BS power constraint. However, according to [SSVB08] and [TSC07] (see also Chapter 3 of this thesis), $B$ coordinated BS systems each with $Z$ antennas can be treated as one multiuser MIMO system with $BZ$ antennas. Thus, when $Z = 1$, the considered problem has exactly the same structure as that of [SSVB08]. To the best of our knowledge, there is no other general linear algorithm that can solve the problem types $\mathcal{P}1$ - $\mathcal{P}4.4$. On the other hand, in all problems, since there are more than one power constraints (i.e., per antenna and symbol (user) powers), all power constraints may not be active at the optimal solution. Due to these reasons, we compare **Algorithm 4.II** and the algorithm in [SSVB08] both in terms of the achieved MSE (i.e., minimized MSE) and total utilized BS power at the achieved MSE.

### 4.11.1 Simulation results for problems $\mathcal{P}4.1$ - $\mathcal{P}4.2$

In this subsection, we compare the performance of our proposed algorithm with that of [SSVB08]. As can be seen from Fig. 4.2, the proposed algorithm and the algorithm in [SSVB08] achieve the same symbol-wise and user-wise WSMSEs. Next, we plot the total utilized powers at the BS to achieve these WSMSEs which is shown in Fig. 4.3. From these two figures, one can see that to achieve the same WSMSE, the proposed duality based iterative algorithm requires less total BS power than that of [SSVB08]. This scenario fits to that of [BV11b] and [BV11c] where the sum MSE minimization constrained with a per BS antenna power problem has been examined by duality approach.

From the figures 4.2 and 4.3, one can notice that the total utilized power of $\mathcal{P}4.2$ is higher than that of $\mathcal{P}4.1$ at high SNR regions. And the achieved weighted sum MSE seems to be exactly the same for both $\mathcal{P}4.1$ and $\mathcal{P}4.2$. However, as the constraint of $\mathcal{P}4.2$ is more relaxed than that of $\mathcal{P}4.1$, one may expect that the weighted sum MSE of $\mathcal{P}4.2$ is strictly less than that of $\mathcal{P}4.1$. To verify this, we merge and zoom the achieved weighted sum MSEs of both of these problems on the same plot in the desired SNR values (see Fig. 4.4). As can be seen from this figure, the weighted sum MSE of $\mathcal{P}4.2$ is slightly less than that of $\mathcal{P}4.1$ which is expected.

### 4.11.2 Simulation results for problems $\mathcal{P}4.3$ - $\mathcal{P}4.4$

Like in the above subsection, here we compare the performance of our proposed algorithm with that of [SSVB08]. For these problems, we also observe from Fig. 4.5 that the proposed algorithm and the algorithm in [SSVB08] achieve the same maximum symbol and user MSEs. And from Fig. 4.6 the proposed duality based algorithms utilize less total BS power compared to that of [SSVB08]. For all of our problems, we observe that to achieve the same MSE, the proposed duality based iterative algorithm utilizes less total BS power compared to the algorithm of [SSVB08]. This scenario has also been observed for other MSE and rate-based problems in [BV11c, BV11b, BV12].

### 4.11.3 Convergence speed of Algorithm 4.II

As can be seen from Section 4.8, the overall computational complexity of **Algorithm 4.II** depends on the number of iterations to achieve convergence. In general, the number of iterations to achieve convergence may not be the same for all problems. On the other hand, for each problem, getting the exact number of iterations to achieve convergence analytically is very difficult. Due to these reasons, we provide numerical simulations to demonstrate the convergence speed of **Algorithm 4.II** for $\mathcal{P}4.1$. As can be seen from Fig. 4.7, the proposed algorithm converges within few iterations in low, medium and high SNR regions.

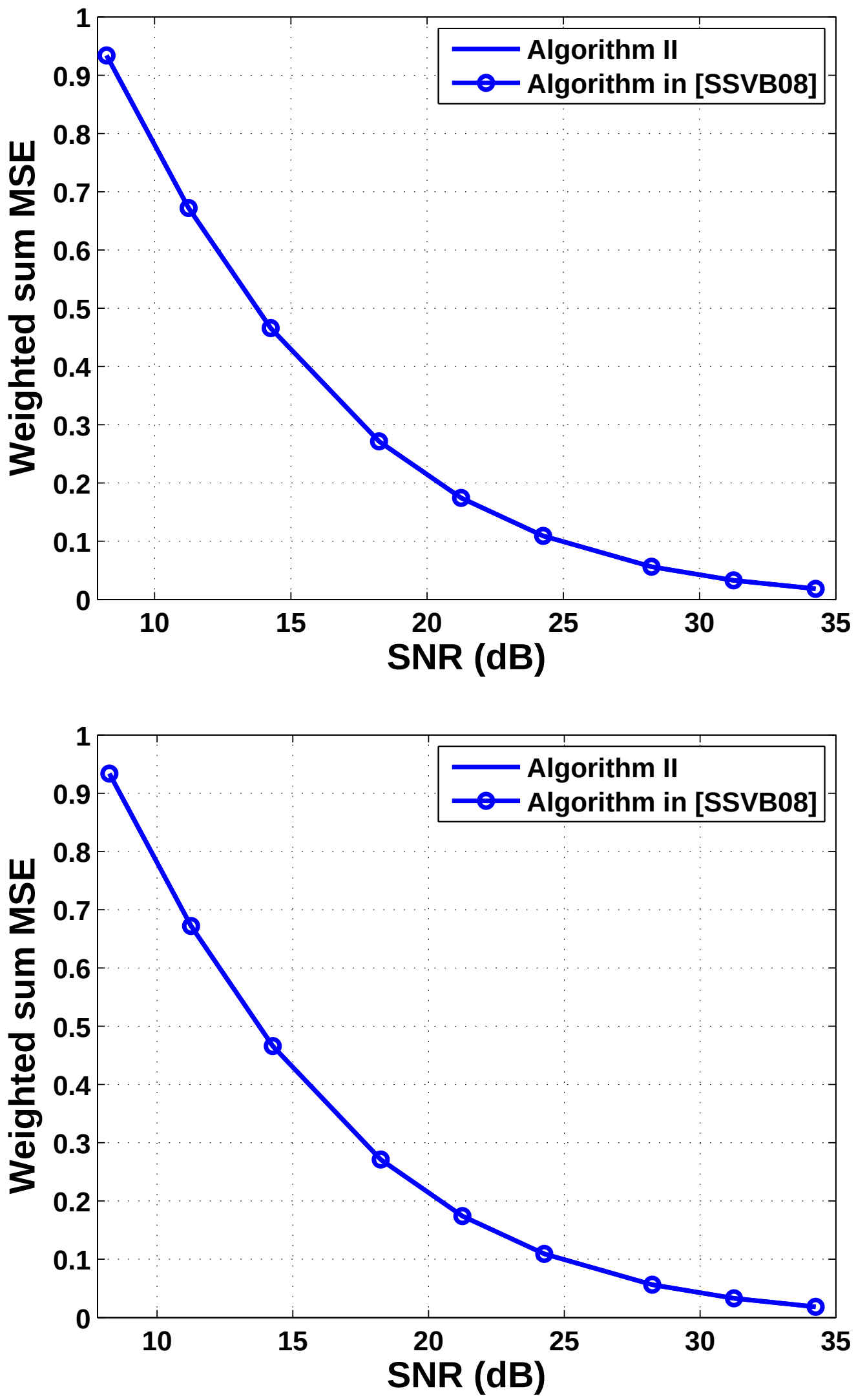

**Figure 4.2** Comparison of the proposed algorithm (**Algorithm 4.II**) and the algorithm of [SSVB08] in terms of WSMSE for: [upper] $\mathcal{P}4.1$, [lower] $\mathcal{P}4.2$.

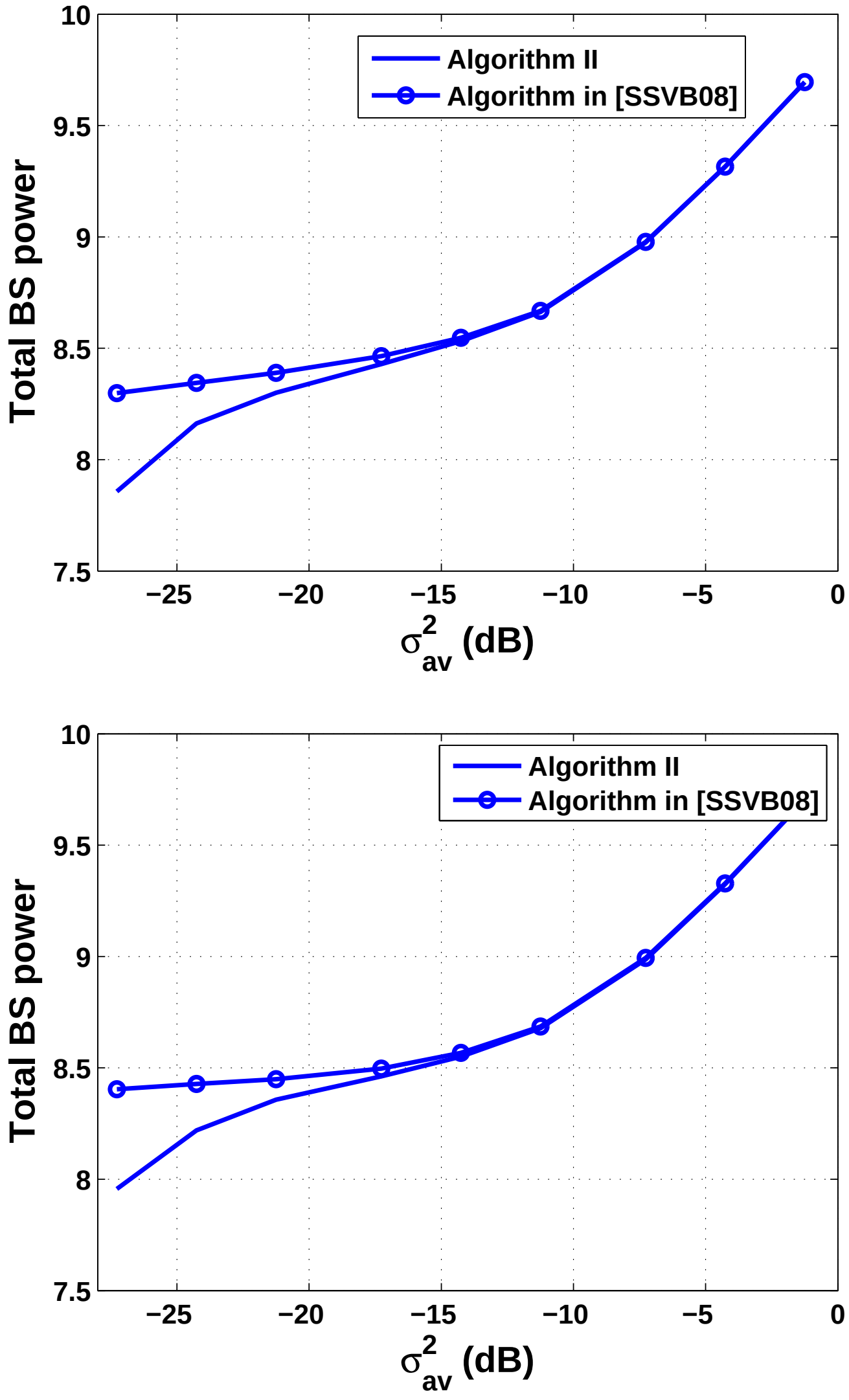

**Figure 4.3** Comparison of the proposed algorithm (**Algorithm 4.II**) and the algorithm of [SSVB08] in terms of total BS power for: [upper] $\mathcal{P}4.1$, [lower] $\mathcal{P}4.2$. For this figure we compute $\sigma_{av}^2(dB)$ as $\sigma_{av}^2$(dB)=$10\log\frac{\sigma_{av}^2}{1mW}$.

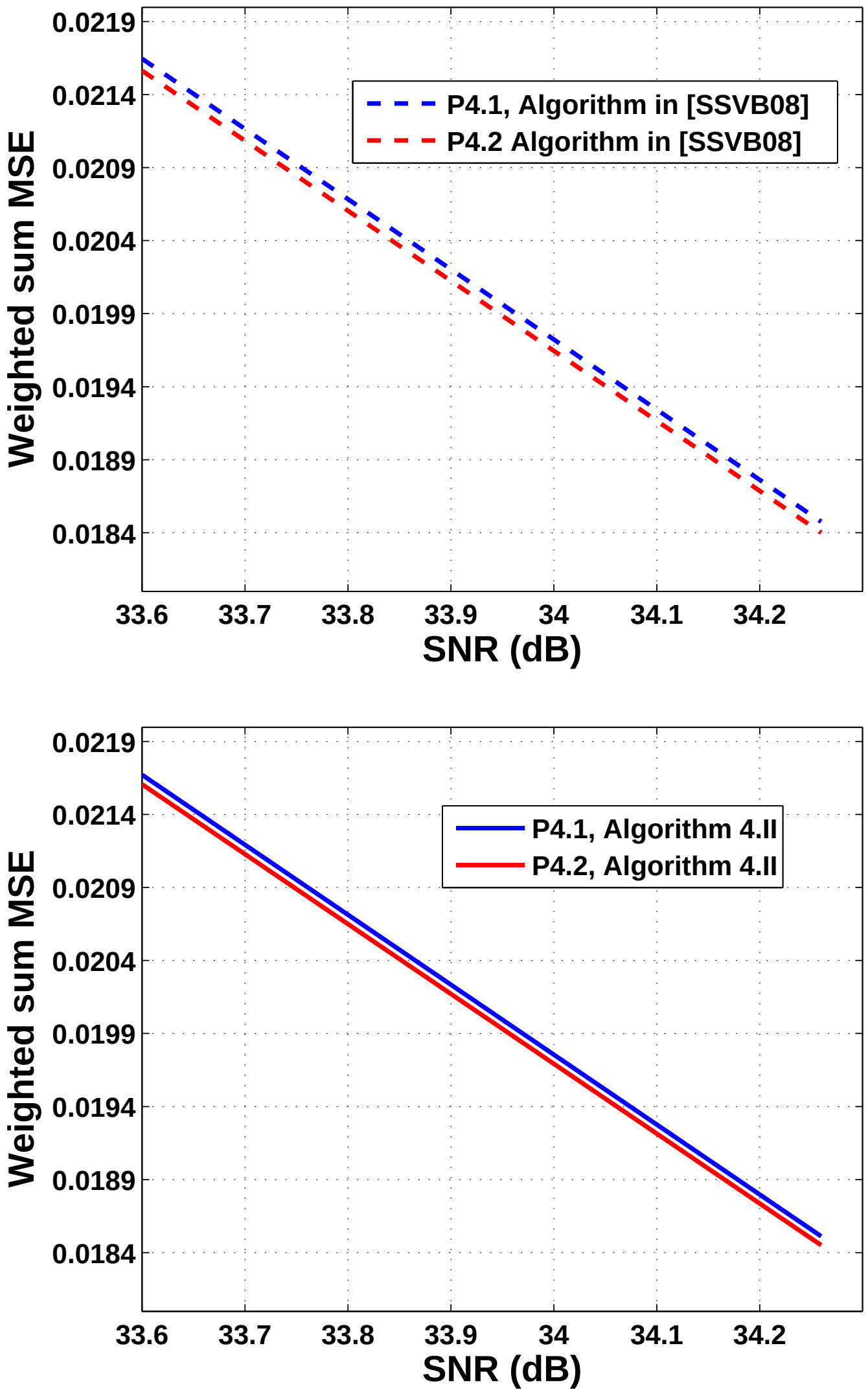

**Figure 4.4** Comparison of the weighted sum MSEs obtained from $\mathcal{P}4.1$ and $\mathcal{P}4.2$: [upper] The algorithm in [SSVB08], [lower] The proposed algorithm (**Algorithm 4.II**).

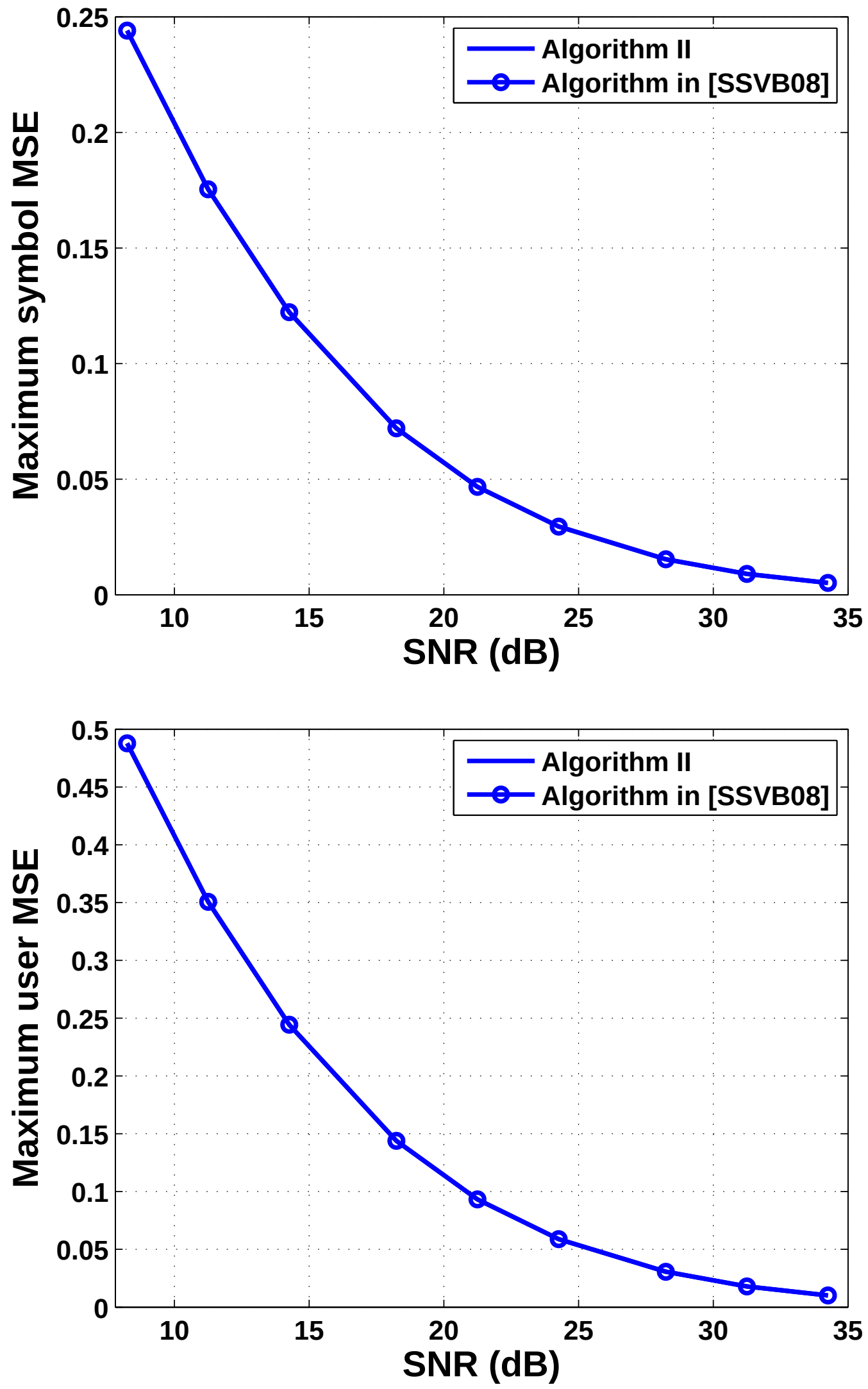

**Figure 4.5** Comparison of the proposed algorithm (**Algorithm 4.II**) and that of in [SSVB08] in terms of maximum achieved MSE for: [upper] $\mathcal{P}4.3$, [lower] $\mathcal{P}4.4$.

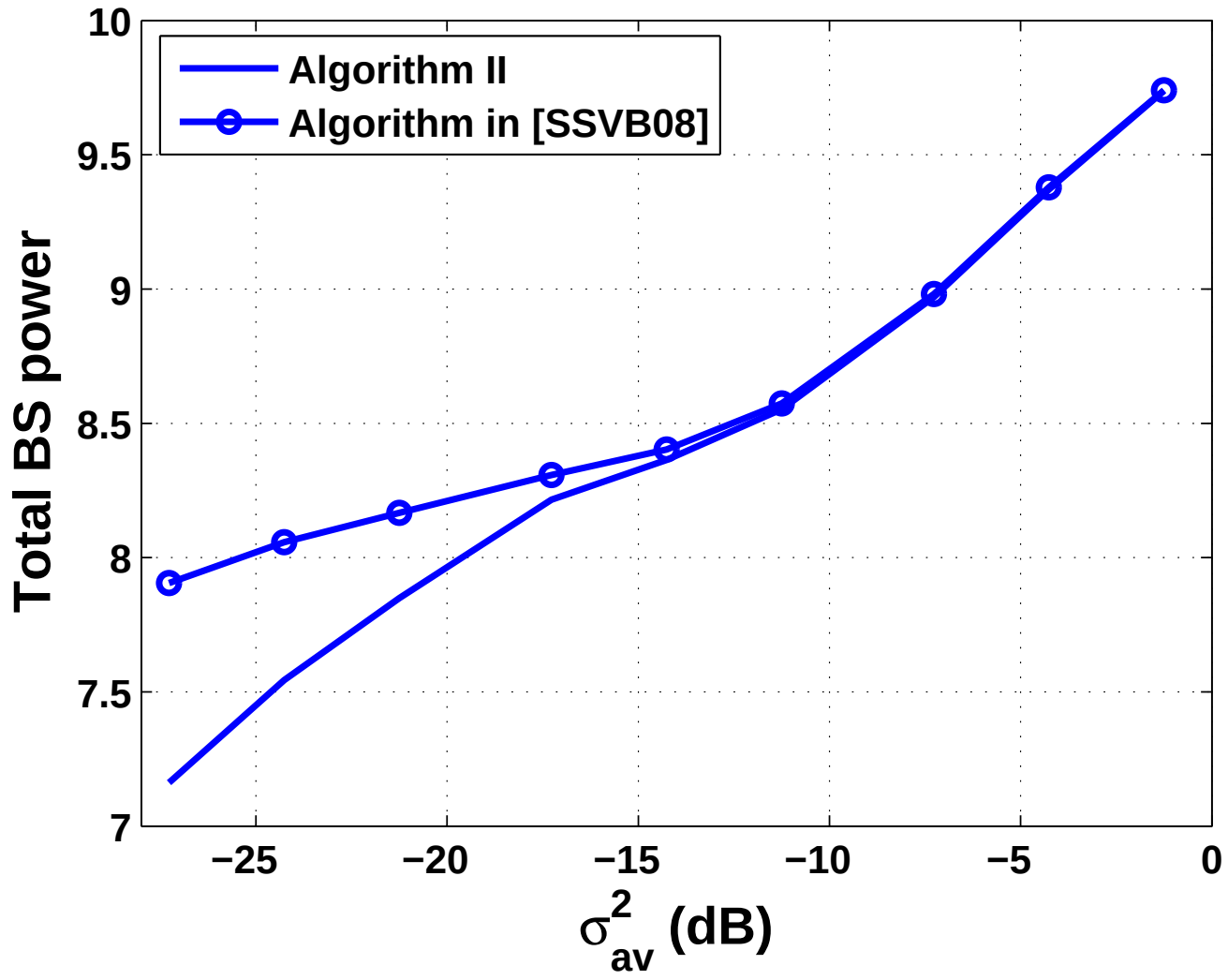

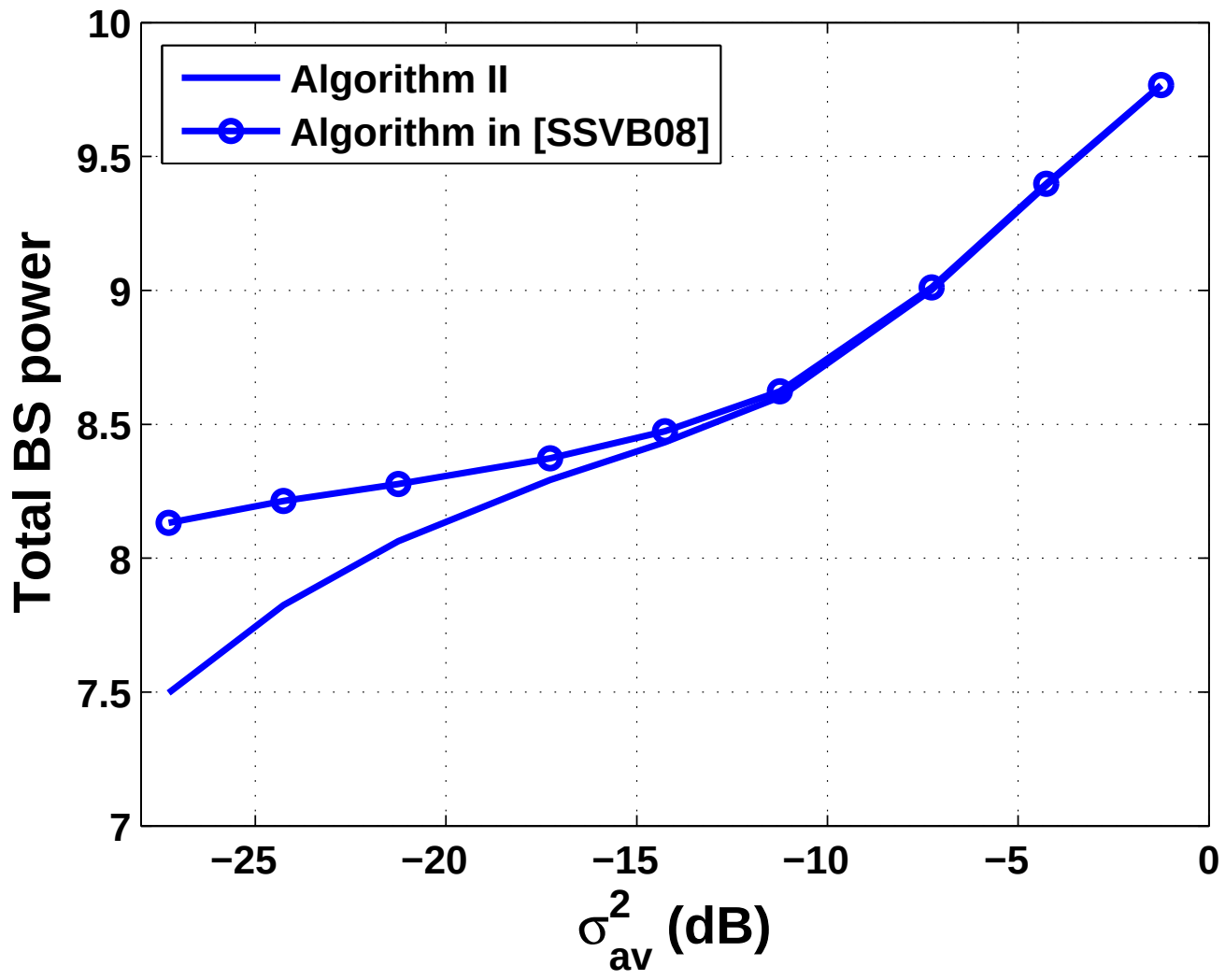

**Figure 4.6** Comparison of the proposed algorithm (**Algorithm 4.II**) and the algorithm of [SSVB08] in terms of total BS power for: [upper] $\mathcal{P}4.3$, [lower] $\mathcal{P}4.4$. In this figure we compute $\sigma_{av}^2(dB)$ as $\sigma_{av}^2$(dB)=10 log $\frac{\sigma_{av}^2}{1mW}$.

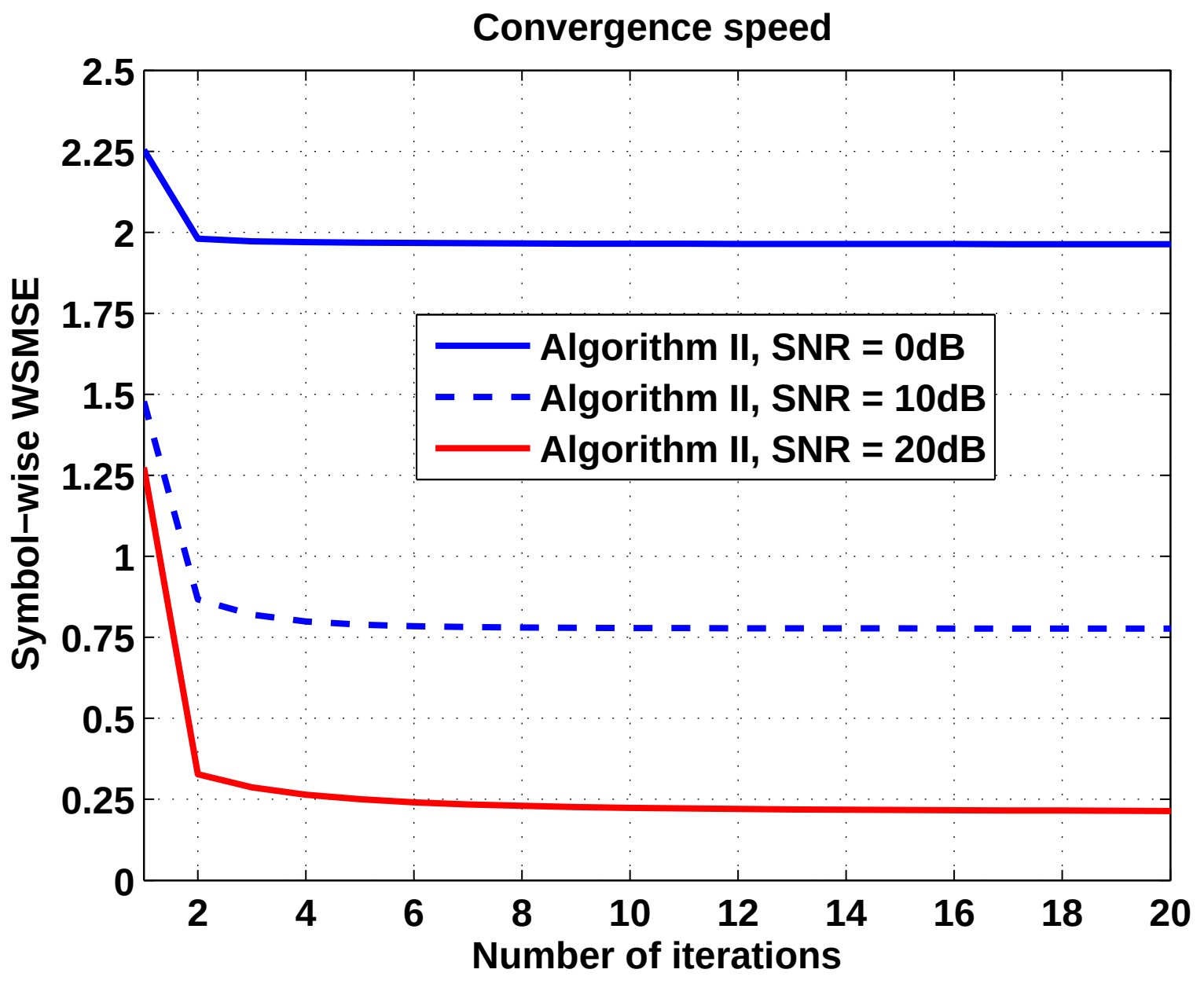

**Figure 4.7** Convergence speed of **Algorithm 4.II** for $\mathcal{P}1$.

## 4.12 Conclusions

In this chapter, we examine different transceiver design problems for multiuser MIMO systems under generalized linear power constraints. The problems are solved for the practically relevant scenario where the noise vector of each MS is a ZMCSCG random variable with arbitrary covariance matrix. For all of our problems, we propose new downlink-interference duality based iterative solutions. The current duality are established by formulating the noise covariance matrices of the dual interference channels as fixed point functions and marginally stable (convergent) discrete-time-switched systems. We show that the proposed duality based iterative algorithms can be extended straightforwardly to solve several practically relevant linear transceiver design prob-lems. We also show that our new MSE downlink-interference duality unify all existing MSE duality. Simulation results demonstrate that the proposed duality based algorithms utilize less total BS power than that of existing algorithms.

## 4.13 Appendix 4.A: Proof of Theorem 4.2

**Proof.** Define $\tilde{\mathbf{D}} \triangleq \text{diag}(\mathbf{A}_{1,1}, \cdots, \mathbf{A}_{n,n})$ and $\tilde{\mathbf{A}} \triangleq \tilde{\mathbf{D}} - \mathbf{A}$. It follows

$$\mathbf{A} = \tilde{\mathbf{D}} - \tilde{\mathbf{A}} \Rightarrow \mathbf{A}^{-1} = \tilde{\mathbf{D}}^{-1}(\mathbf{I} - \bar{\mathbf{A}})^{-1}$$

where $\bar{\mathbf{A}} = \tilde{\mathbf{A}}\tilde{\mathbf{D}}^{-1}$. Since $(\mathbf{I} - \bar{\mathbf{A}})$ is strictly diagonally dominant matrix, $(\mathbf{I} - \bar{\mathbf{A}})^{-1}$ exists [HJ85] (page 349 of [HJ85]). Furthermore, if $\rho(\bar{\mathbf{A}}) < 1$, $(\mathbf{I} - \bar{\mathbf{A}})^{-1}$ can be expressed as

$$(\mathbf{I} - \bar{\mathbf{A}})^{-1} = \sum_{k=0}^{\infty} \bar{\mathbf{A}}^k. \tag{4.77}$$

It follows

$$\bar{\mathbf{A}}^{-1} = \tilde{\mathbf{D}}^{-1}(\mathbf{I} - \bar{\mathbf{A}})^{-1} = \tilde{\mathbf{D}}^{-1} \sum_{k=0}^{\infty} \bar{\mathbf{A}}^k \geq 0 \tag{4.78}$$

From this equation we can see that if $\rho(\bar{\mathbf{A}}) < 1$, the nonnegativity of $\bar{\mathbf{A}}^{-1}$ can be ensured. Next we show that $\rho(\bar{\mathbf{A}})$ is indeed less than 1. For any $n \times n$ matrix $\mathbf{X}$, we have [HJ85] (pages 294 and 297 of [HJ85])

$$\rho(\mathbf{X}) \leq |||\mathbf{X}|||, \;\; |||\mathbf{X}||_1 \triangleq \max_{1 \leq j \leq n} \sum_{i=1}^{n} |x_{ij}| \tag{4.79}$$

where $|||.|||$ is any matrix norm and $|||.|||_1$ is a matrix one norm. By using (4.79), we get the following bound [HJ85]

$$\rho(\bar{\mathbf{A}}) \leq |||\bar{\mathbf{A}}|||_1 < 1. \tag{4.80}$$

Since $\mathbf{A}^{-1}$ has nonnegative elements, $\mathbf{A}$ is also an M-matrix [PB74]. By defining $\mathbf{S} \triangleq \mathbf{A}^{-1}$ and $\mathbf{e} \triangleq \mathbf{1}^{n \times 1}$, we get

$$\begin{aligned} \mathbf{e}^T\mathbf{A} = \mathbf{e}^T \;\Rightarrow \mathbf{e}^T = \mathbf{e}^T\mathbf{S} &= [\sum_{j=1}^{n} \mathbf{S}_{j,1}, \cdots, \sum_{j=1}^{n} \mathbf{S}_{j,n}] \\ \Rightarrow |||\mathbf{S}|||_1 &= 1 \end{aligned} \tag{4.81}$$

where the third equality follows from the fact that $\mathbf{S}$ is a nonnegative matrix. □

# Conclusions and Future Works

# 5

## 5.1 Conclusions

In this thesis, we accomplish the following key tasks:

1. We generalize the existing MSE uplink-downlink duality to the more practically relevant power constraint scenarios. The generalized duality are presented as MSE downlink-interference duality. For MSE-based problems, the duality are established by formulating the noise covariance matrices of the interference channels as marginally stable (convergent) discrete-time-switched systems. To express the noise covariance matrices as discrete time switched systems, this thesis employs the bound (4.55) and Theorem 2 of Chapter 4. Note that the second property of this Theorem is originally formulated and proved by us. For WSMSE-based problems, computationally less complex duality (compared to the latter duality) are established by formulating the noise covariance matrices of the interference channels as fixed point functions. The proposed duality can be applied to solve many classes of SINR, rate and MSE-based problems for multiuser MIMO uncoordinated and coordinated BS systems. The extensions of the MSE downlink-interference duality to imperfect CSI scenario has also been discussed.
2. We have shown that the weighted sum rate maximization problem can be equivalently formulated as weighted sum MSE minimization problem with additional optimization variables and constraints. This prob-lem reformulation employs two novel Lemmas (i.e., Lemma 1 and 2 of Chapter 3). We would like to mention here that Lemma 2 of Chapter 3 is originally formulated and proved by us.

3. We also develop distributed precoder/decoder design algorithms to solve weighted sum rate and MSE optimization problems. The distributed precoder/decoder design algorithms employ matrix fractional minimization and Lagrangian dual decomposition methods. The extension of the proposed distributed algorithms for solving robust weighted sum MSE optimization problem has also been discussed.

## 5.2 Future works

1. As mentioned throughout this thesis, the robust rate and SINR-based problems (i.e, in stochastic design approach) have not been examined. Thus, solving such problems is open research topic.
2. The duality results of this thesis are developed under the assumption of perfect CSI (and imperfect CSI under stochastic robust design approach). From this explanation, we can notice that the duality approach of this thesis to imperfect CSI under worst-case robust design approach is an open research topic. For transceiver design problems, the worst-case robust design approach can be found in [ECV10, VBS09]. Furthermore, it is also interesting to examine the duality of this thesis for several classes of channel uncertainty as discussed in [WBOP13]
3. In all of the proposed distributive algorithms, we assume that the global channel knowledge is available at the central controller (or at all BSs) prior to optimization. However, when the number of BSs are very large, getting the global channel knowledge at the central controller (each BS) appears to be difficult. Thus, the extension of the proposed distributive algorithms by assuming that each BS designs its precoder/decoder using the local CSI knowledge is interesting for future research.
4. All of the proposed algorithms of this thesis are linear. However, not all of these algorithms are able to guarantee global optimality for their respective problems. And to the best of our knowledge, for downlink multiuser MIMO systems (both uncoordinated and coordinated BSs), we are not aware of any optimal linear algorithm to solve general transceiver design problems. So it is interesting to conduct future research on the global optimal linear algorithm for solving general transceiver design problems. One potential approach could be to extend the majorization theory of [PLC04][1] for the downlink multiuser MIMO uncoordinated and coordinated BS systems.

[1]For single user MIMO systems, the optimal precoder/decoder structure can be obtained by employing Majorization theory of this paper.

# Convex Optimization Basics

A

## A.1 Convex Optimization Theory

In this appendix a brief overview of convex optimization theory will be discussed.

**Definition** (Convex set and convex function [BV04]): Let $\mathcal{S}$ be a vector space over some field. A set $\mathcal{C} \in \mathcal{S}$ is convex, if for any $\mathbf{x} \in \mathcal{C}, \mathbf{y} \in \mathcal{C}$ and $0 \leq \theta \leq 1$, we have

$$\theta\mathbf{x} + (1-\theta)\mathbf{y} \in \mathcal{C}.$$

A function $f(.)$ is convex if

$$f(\theta\mathbf{x} + (1-\theta)\mathbf{y}) \leq \theta f(\mathbf{x}) + (1-\theta)f(\mathbf{y})$$

and strictly convex if

$$f(\theta\mathbf{x} + (1-\theta)\mathbf{y}) < \theta f(\mathbf{x}) + (1-\theta)f(\mathbf{y}).$$

**Definition** (Affine function [BV04]): A function $f : \Re^n \to \Re^m$ is affine if it has of the form

$$f(\mathbf{x}) = \mathbf{A}\mathbf{x} + \mathbf{b}$$

where $\mathbf{A} \in \Re^{m\times n}$ and $\mathbf{b} \in \Re^m$.

Consider the following optimization problem:

$$\begin{aligned}\min_{\mathbf{x}} \ & f_0(\mathbf{x}) \\ \text{s.t } & g_i(\mathbf{x}) \leq 0, \ i = 1, 2, \cdots \\ & h_j(\mathbf{x}) = 0, \ j = 1, 2, \cdots \end{aligned} \tag{A.1}$$

where $f_0(\mathbf{x})$ is the objective function, $g_i(\mathbf{x}), \forall i$ and $h_j(\mathbf{x}), \forall j$ are the constraint functions.

This optimization problem is said to be convex if $f_0(\mathbf{x})$ and $g_i(\mathbf{x}), \forall i$ are convex functions and $h_j(\mathbf{x}), \forall j$ are affine functions [BV04].

The beauty of any convex optimization problem is that its global optimal solution can be obtained very efficiently either in closed form or by numerical methods. For most convex optimization problems, the closed form solution is obtained by examining the Lagrangian dual problem, and numerical solution is obtained using the well known interior point methods [BV04]. For this reason, the mathematics of convex optimization has received a lot interest in modern optimization theory. In the following, we summarize the most common convex optimization problems [BV04]:

## A.2 Examples of convex optimization problems

### A.2.1 Linear programming (LP)

A linear programming optimization problem is mathematically formulated as

$$\begin{aligned}\min_{\mathbf{x}} \ & \mathbf{c}^T\mathbf{x} + d \\ \text{s.t } & \mathbf{G}\mathbf{x} \preceq \mathbf{h} \\ & \mathbf{A}\mathbf{x} = \mathbf{b}\end{aligned}$$

where $\mathbf{c} \in \Re^{n\times 1}, \mathbf{x} \in \Re^{n\times 1}, d \in \Re, \mathbf{G} \in \Re^{k\times n}, \mathbf{h} \in \Re^{k\times 1}, \mathbf{A} \in \Re^{m\times n}$ and $\mathbf{b} \in \Re^{m\times 1}$. The objective and constraint functions of this problem is linear which is a special case of convex function. Thus, this problem is a convex optimization problem.

### A.2.2 Quadratic programming (QP)

A quadratic programming optimization problem is formulated as

$$\begin{aligned}\min_{\mathbf{x}}\ & 0.5\mathbf{x}^T\mathbf{P}\mathbf{x} + \mathbf{c}^T\mathbf{x} + d \\ \text{s.t}\ & \mathbf{G}\mathbf{x} \preceq \mathbf{h} \\ & \mathbf{A}\mathbf{x} = \mathbf{b}\end{aligned}$$

where $\mathbf{x} \in \Re^{n\times 1}, \mathbf{P}_{+} \in \Re^{n\times n}, \mathbf{c} \in \Re^{n\times 1}, d \in \Re, \mathbf{G} \in \Re^{k\times n}, \mathbf{h} \in \Re^{k\times 1}, \mathbf{A} \in \Re^{m\times n}$ and $\mathbf{b} \in \Re^{m\times 1}$, and $(.)_{+}$ represents positive semi-definite [BV04]. The objective function of this problem is quadratic and the constraint functions are linear which are special cases of convex functions. Thus, this problem is a convex optimization problem.

### A.2.3 Second-order cone programming (SOCP)

A second-order cone programming optimization problem is expressed as

$$\begin{aligned}\min_{\mathbf{x}}\ & \mathbf{f}^T\mathbf{x} \\ \text{s.t}\ & \|\mathbf{g}^T\mathbf{x} + d\|_2 \preceq \mathbf{c}^T\mathbf{x} + z \\ & \mathbf{A}\mathbf{x} = \mathbf{b}\end{aligned}$$

where $\mathbf{f} \in \Re^{n\times 1}, \mathbf{x} \in \Re^{n\times 1}, \mathbf{g} \in \Re^{n\times 1}, d \in \Re, \mathbf{c} \in \Re^{n\times 1}, z \in \Re, \mathbf{A} \in \Re^{m\times n}$ and $\mathbf{b} \in \Re^{m\times 1}$, and $\|.\|_2$ denotes the 2 norm of a vector.

### A.2.4 Semi-definite programming (SDP)

A semi-definite programming problem is mathematically formalized as

$$\begin{aligned}\min_{x_i}\ & \mathbf{c}^T\mathbf{x} \\ \text{s.t}\ & \sum_{i=1}^{n}\mathbf{F}_i x_i + \mathbf{G} \preceq \mathbf{0} \\ & \mathbf{A}\mathbf{x} = \mathbf{b}\end{aligned}$$

where $\mathbf{c} \in \Re^{n\times 1}, x_i \in \Re, \forall i, \mathbf{F}_i \in \Re^{k\times s}, \forall i, \mathbf{G} \in \Re^{k\times s}, \mathbf{A} \in \Re^{m\times n}$ and $\mathbf{b} \in \Re^{m\times 1}$.

## A.3 Geometric program

In this section we describe Geometric program problems that are not convex in their own form. These problems can, however, be transformed to convex optimization problems by a change of variables and a transformation of the objective and constraint functions.

### A.3.1 Monomial and posynomials

A Monomial function is a function of the form

$$f(\mathbf{x}) = cx_1^{a_1} x_2^{a_2} \cdots x_n^{a_n}$$

where $c \geq 0$ and $a_i \in \Re, \forall i$. A sum of monomials is called posynomials. Thus, a posynomial function is a function of the following form:

$$f(\mathbf{x}) = \sum_{k=1}^{N} c_k x_1^{a_{k1}} x_2^{a_{k2}} \cdots x_n^{a_{kn}}$$

where $c_k \geq 0, \forall k$ and $a_{ki} \in \Re, \forall k, i$.

### A.3.2 Geometric programming (GP)

A GP problem has the following form

$$\begin{aligned} \min_{\mathbf{x}>\mathbf{0}} \; & f(\mathbf{x}) \\ \text{s.t } & g_i(\mathbf{x}) \leq 1, \; i = 1, 2, \cdots \\ & h_j(\mathbf{x}) = 1, \; j = 1, 2, \cdots \end{aligned} \tag{A.2}$$

where $f(\mathbf{x})$ and $g_i(\mathbf{x}), \forall i$ are posynomial functions and $h_j(\mathbf{x}), \forall j$ are monomial functions.

## A.3.3 Geometric programming in convex form

In this subsection, we summarize how the above GP problem is transformed into a convex optimization problem. Towards this end, we express

$$\begin{aligned} f(\mathbf{x}) &= \sum_{k=1}^{N} c_k x_1^{a_{k1}} x_2^{a_{k2}} \cdots x_n^{a_{kn}} \\ g_i(\mathbf{x}) &= \sum_{k=1}^{M_i} d_{ki} x_1^{b_{ki1}} x_2^{b_{ki2}} \cdots x_n^{b_{kin}},\ \forall i \\ h_j(\mathbf{x}) &= z_j x_1^{t_{j1}} x_2^{t_{j2}} \cdots x_n^{t_{jn}},\ \forall j. \end{aligned}$$

Now by defining $\tilde{x}_i \triangleq \log(x_i)$, $\tilde{c}_k \triangleq \log(c_k)$, $\tilde{d}_k \triangleq \log(d_k)$, $\tilde{z}_k \triangleq \log(z_k)$, $\mathbf{a}_k = [a_{k1}, \cdots, a_{kn}]$, $\mathbf{b}_{ki} = [b_{ki1}, \cdots, b_{kin}]$ and $\mathbf{t}_j = [t_{j1}, \cdots, t_{jn}]$, we can reexpress $f(\mathbf{x})$, $g_i(\mathbf{x})$ and $h_j(\mathbf{x})$ as

$$\begin{aligned} f(\tilde{\mathbf{x}}) &= \sum_{k=1}^{N} \exp^{(\mathbf{a}_k^T \tilde{\mathbf{x}} + \tilde{c}_k)} \\ g_i(\tilde{\mathbf{x}}) &= \sum_{k=1}^{M_i} \exp^{(\mathbf{b}_{ki}^T \tilde{\mathbf{x}} + \tilde{d}_{ki})},\ \forall i \\ h_j(\tilde{\mathbf{x}}) &= \exp^{(\mathbf{t}_j^T \tilde{\mathbf{x}} + \tilde{z}_j)},\ \forall j. \end{aligned}$$

Using the transformed variable $\tilde{\mathbf{x}}$, and taking the logarithm of the objective and constraint functions, the GP problem (A.2) can be equivalently reformulated as[1]

$$\begin{aligned} \min_{\mathbf{x}>\mathbf{0}}\ & \log\Big(\sum_{k=1}^{N} \exp^{(\mathbf{a}_k^T \tilde{\mathbf{x}} + \tilde{c}_k)}\Big) \\ \text{s.t}\ & \log\Big(\sum_{k=1}^{M_i} \exp^{(\mathbf{b}_{ki}^T \tilde{\mathbf{x}} + \tilde{d}_{ki})}\Big) \le 0,\ i = 1, 2, \cdots \\ & \log\big(\exp^{(\mathbf{t}_j^T \tilde{\mathbf{x}} + \tilde{z}_j)}\big) = 0,\ j = 1, 2, \cdots. \end{aligned} \tag{A.3}$$

It can be shown that this problem is a convex optimization problem and (A.3) is referred as a GP problem in convex form [BV04].

[1] This is due to the fact that $\min x$ is equivalent to $\min \log(x)$.

# Duality Basics

# B

For better understanding of the duality, we consider multiuser MIMO system where a BS is serving $K$ decentralized MSs. We assume that the BS has N antennas and the $k$th MS has $M_k$ antennas. The overall transmitted symbols are given by $\mathbf{d} = [\mathbf{d}_1, \cdots, \mathbf{d}_K]$, where $\mathbf{d}_k \in \mathcal{C}^{S_k \times 1}$ is the transmitted signal for the $k$th user. Under these assumption, the multiuser transmission system model with precoding and decoding operations are shown in Fig. B.1.

For these system models, the dimensions of the precoder and decoder matrices of both channels are $\mathbf{B}_k \in \mathcal{C}^{N \times S_k}, \mathbf{T}_k \in \mathcal{C}^{N \times S_k}, \mathbf{W}_k \in \mathcal{C}^{M_k \times S_k}$ and $\mathbf{V}_k \in \mathcal{C}^{M_k \times S_k}$, and $\mathbf{B} = [\mathbf{B}_1, \cdots, \mathbf{B}_K]$ and $\mathbf{T} = [\mathbf{T}_1, \cdots, \mathbf{T}_K]$. The estimated symbol of the $k$th user in the downlink and uplink channels are given by

$$\widehat{\mathbf{d}}_k^{DL} = \mathbf{W}_k^H (\mathbf{H}_k^H \sum_{i=1}^{K} \mathbf{B}_k \mathbf{d}_k + \mathbf{n}_k^{DL}) \tag{B.1}$$

$$\widehat{\mathbf{d}}_k^{UL} = \mathbf{T}_k^H (\sum_{i=1}^{K} \mathbf{H}_i \mathbf{V}_i \mathbf{d}_i + \mathbf{n}^{UL}) \tag{B.2}$$

where $(.)^{DL}$ and $(.)^{UL}$ denote downlink and uplink, respectively.

If we assume $\mathbf{d}_k \sim \mathcal{CN}(\mathbf{0}, \mathbf{I})$, $\mathbf{n}_k^{DL} \sim \mathcal{CN}(\mathbf{0}, \sigma^2 \mathbf{I})$ and $\mathbf{n}^{UL} \sim \mathcal{CN}(\mathbf{0}, \tilde{\sigma}^2 \mathbf{I})$, the total sum MSE of all users is given by [BV11a]

$$\xi^{DL} = S + \mathrm{tr}\{\mathbf{W}^H \mathbf{H}^H \mathbf{B} \mathbf{B}^H \mathbf{H} \mathbf{W} + \sigma^2 \mathbf{W}^H \mathbf{W} - \mathbf{W}^H \mathbf{H}^H \mathbf{B} - \mathbf{B}^H \mathbf{H} \mathbf{W}\} \tag{B.3}$$

$$\xi^{UL} = \mathrm{tr}\{\mathbf{T}^H \mathbf{H} \mathbf{V} \mathbf{V}^H \mathbf{H}^H \mathbf{T} + \tilde{\sigma}^2 \mathbf{T}^H \mathbf{T} - \mathbf{T}^H \mathbf{H} \mathbf{V} - \mathbf{V}^H \mathbf{H}^H \mathbf{T}\} + S \tag{B.4}$$

where $S = \sum_{i=1}^{K} S_i$, $\mathbf{W} = \mathrm{blkdiag}(\mathbf{W}_1, \cdots, \mathbf{W}_K)$, $\mathbf{V} = \mathrm{blkdiag}(\mathbf{V}_1, \cdots, \mathbf{V}_K)$.

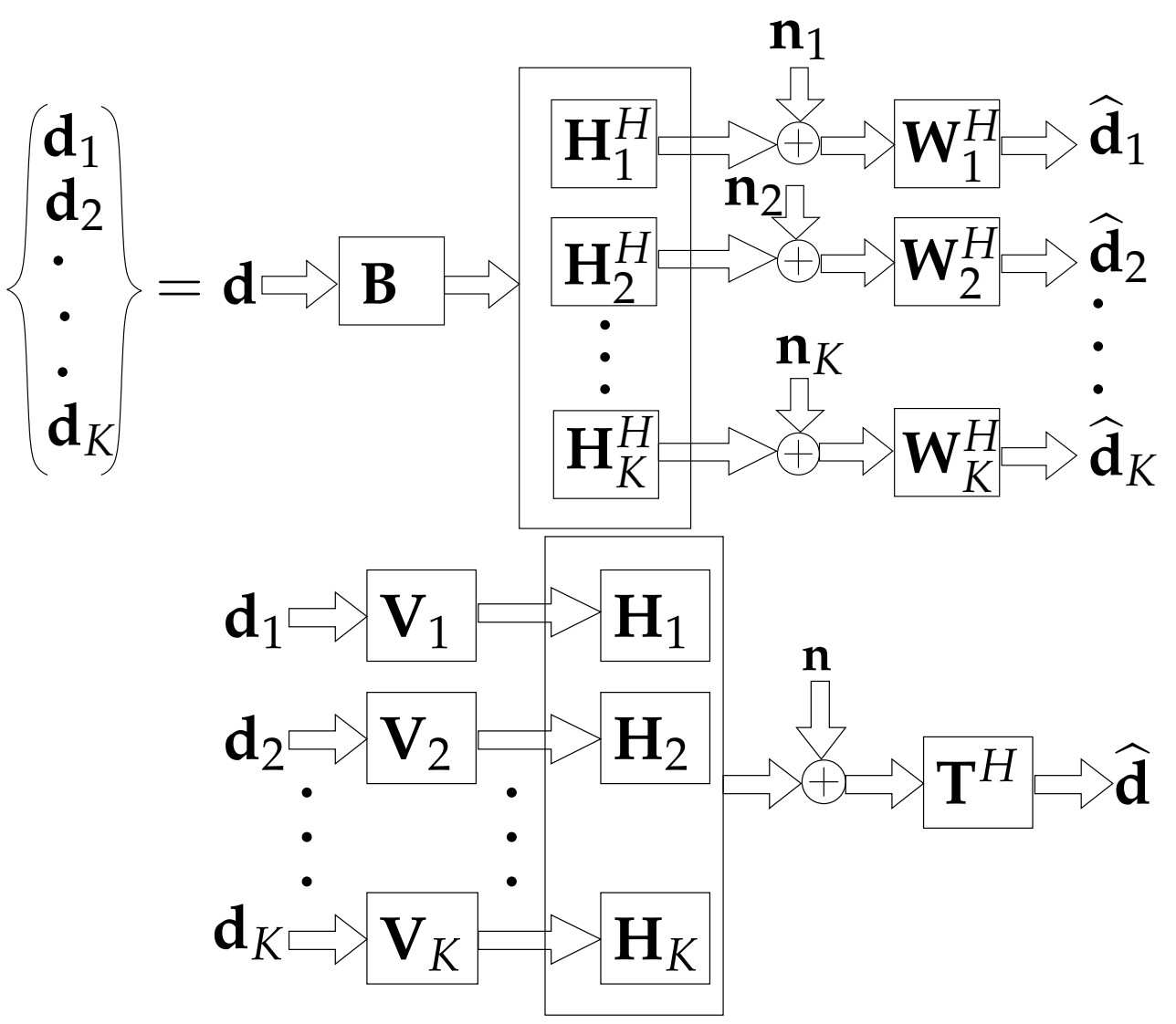

**Figure B.1** Multiuser MIMO system model: [upper] Downlink channel, [lower] Uplink channel.

Now let us consider the sum MSE minimization problem under total transmission power constraint: This problem can be expressed as [BV11a]

$$\mathcal{P}^{DL} : \min_{\mathbf{B},\mathbf{W}} \xi^{DL}, \quad \text{s.t}\ \text{tr}\{\mathbf{B}\mathbf{B}^{H}\} \leq P_{max}^{DL} \tag{B.5}$$

$$\mathcal{P}^{UL} : \min_{\mathbf{V},\mathbf{T}} \xi^{UL}, \quad \text{s.t}\ \text{tr}\{\mathbf{V}\mathbf{V}^{H}\} \leq P_{max}^{UL} \tag{B.6}$$

where $P_{max}^{DL}$ and $P_{max}^{UL}$ are the maximum transmission powers in the downlink and uplink channels.

Assume that we would like to solve the problem $\mathcal{P}^{DL}$ (i.e., our original problem).

From mathematical optimization point of view it appears that $\mathcal{P}^{UL}$ is a convex optimization problem (see Appendix A for a brief summary of convex optimization), where its global optimal solution can be obtained by existing convex optimization algorithms. However, unfortunately, $\mathcal{P}^{DL}$ is not convex problem and its global optimization can not be obtained by existing convex optimization algorithms[1].

Due to this reason, researchers propose to get the global optimal solution of $\mathcal{P}^{DL}$ from the the global optimal solution of $\mathcal{P}^{UL}$. Such approach of solving transceiver design problem is called Uplink-downlink duality approach[2]. From this discussion, we can understand that the uplink channel (i.e. in Fig. B.1) is created just to solve the original downlink channel problem. For this reason, the uplink channel is termed as a "virtual uplink channel" [FLT98, VM99, SB04, BV11a]. Recently, this notion is even more developed to a new duality which is termed as Downlink-interference duality in [BV13].

[1]This phenomena arises for many other classes of transceiver design problems for multiuser networks.

[2]We would like to mention here that for some design criteria, neither the uplink problem nor its downlink problem is convex. In such a case, the suboptimal solution of the original problem can be obtained by iteratively switching from the downlink to uplink channel problems and vice versa. This solution approach is also referred as Uplink-downlink duality approach.

# List of Publications

## Journal papers

1. Tadilo Endeshaw Bogale and Luc Vandendorpe, "Max-Min Signal Energy based Spectrum sensing Algorithms for Cognitive radio Networks under Noise Variance Uncertainty," in IEEE Trans. Wireless Communications, Sep. 2013.

2. Tadilo Endeshaw Bogale and Luc Vandendorpe, "Linear Transceiver Design for Downlink Multiuser MIMO Systems: Downlink-Interference Duality Approach," in IEEE Trans. Signal Processing, Jul. 2013.

3. Tadilo Endeshaw Bogale and Luc Vandendorpe, "Robust Sum MSE Optimization for Downlink Multiuser MIMO Systems with Arbitrary Power Constraint: Generalized Duality Approach," in IEEE Trans. Signal Processing, Dec. 21, 2011.

4. Tadilo Endeshaw Bogale and Luc Vandendorpe, "Weighted Sum Rate Optimization for MIMO Coordinated Base Station Systems: Centralized and Distributed Algorithms," in IEEE Trans. Signal Processing, Dec. 2011.

5. Tadilo Endeshaw Bogale, Batu Krishna Chalise and Luc Vandendorpe, "Robust Transceiver Optimization for Coordinated Base Station Systems: Distributed algorithm," IEEE Trans. Signal Processing, Oct. 2011.

6. Tadilo Endeshaw Bogale, Batu Krishna Chalise and Luc Vandendorpe, "Robust Transceiver Optimization for Downlink Multiuser MIMO Systems", in IEEE Trans. on Signal Processing, Sep. 2010.

## Conference papers